\documentclass{gtart_a}
\usepackage{pinlabel}
\pdfoutput=1

\title{Stabilization in the braid groups I: MTWS}

\author{Joan S Birman}
\givenname{Joan S}
\surname{Birman}
\address{Department of Mathematics\\
Barnard College\\
Columbia University\\\newline
2990 Broadway\\
New York NY 10027\\
USA}
\email{jb@math.columbia.edu}
\urladdr{}

\author{William W Menasco}
\givenname{William W}
\surname{Menasco}
\address{Department of Mathematics\\
University at Buffalo\\
Buffalo NY 14260\\
USA}
\email{menasco@math.buffalo.edu}
\urladdr{}

\volumenumber{10}
\issuenumber{}
\publicationyear{2006}
\papernumber{13}
\lognumber{0687}
\startpage{413}
\endpage{540}

\doi{10.2140/gt.2006.10.359}


\received{23 June 2005}
\revised{}
\accepted{25 January 2006}
\published{27 April 2006}
\publishedonline{27 April 2006}
\proposed{Robion Kirby}
\seconded{Benson Farb, David Gabai}

\numberwithin{equation}{section}

\def\reals{\mathbb{R}}
\def\ints{\mathbb{Z}}

\newtheorem{theorem}{Theorem}
\newtheorem{proposition}{Proposition}[subsection]
\newtheorem{lemma}{Lemma}[subsection]
\newtheorem{corollary}{Corollary}[subsection]

\theoremstyle{remark}
\newtheorem{remark}{Remark}[subsection]

\newtheorem{example}{Example}[subsection]
\def\pp{{\prime\prime}}
\def\be{\begin{enumerate}}
\def\ee{\end{enumerate}}
\def\bi{\begin{itemize}}
\def\ei{\end{itemize}}
\def\fib{\mathbf{H}}
\def\axis{\ensuremath{\mathbf{A}}}
\def\seif{\mathbf{F}}
\def\nhb{{N}}
\def\cA{{\cal A}}
\def\cB{{\cal B}}
\def\cC{{\cal C}}
\def\cD{{\cal D}}
\def\cE{{\cal E}}

\def\cH{{\cal H}}

\def\cM{{\cal M}}
\def\cN{{\cal N}}
\def\cP{{\cal P}}
\def\cR{{\cal R}}
\def\cS{{\cal S}}
\def\cShigh{{\cal S_+}}
\def\cSlow{{\cal S_-}}
\def\cT{{\cal T}}

\def\cX{{\cal X}}

\def\bB{{\bf B}}
\def\bD{{\bf D}}

\def\BS{{\cal BS}}
\def\bX{{\bf X}}
\def\bx{{\bf x}}
\def\ba{{\bf a}}
\def\bb{{\bf b}}
\def\bc{{\bf c}}
\def\bt{{\bf t}}

\def\ann{\cal A}
\def\e{\epsilon}
\def\ta{\ensuremath{\mathcal{CA}}}
\def\pa{\ensuremath{\mathcal{PA}}}
\def\glow{\gamma_-}
\def\ghigh{\gamma_+}
\def\ge{\gamma_\epsilon}
\def\d{\delta}
\def\D{\Delta^2}
\def\g{\gamma}
\def\connect{\longleftrightarrow}

\def\Dhigh{D_+}
\def\Dlow{D_-}
\def\Xhigh{X_+}
\def\Xlow{X_-}

\def\Xmed{X_0}
\def\Thigh{{\bf T}_+}
\def\Tlow{{\bf T}_-}
\def\Tboth{{\bf T}_\e}
\def\Nhigh{{\bf N}_+}
\def\Nlow{{\bf N}_-}
\def\Nboth{{\bf N}_\e}
\def\thintile{{\bf \wwbar{N}}}
\def\thintilehigh{{\bf \wwbar{N}}_+}
\def\thintilelow{{\bf \wwbar{N}}_-}
\def\thintileboth{{\bf \wwbar{N}}_\pm}
\def\annhigh{\cA_+}
\def\annlow{\cA_-}

\def\qhigh{{\bf q}_+}
\def\qlow{{\bf q}_-}

\def\pp{{\prime\prime}}
\def\seifhigh{{\bf F}_+}
\def\seifmed{{\bf F}_0}
\def\seiflow{{\bf F}_-}
\def\alow{a_-}
\def\ahigh{a_+}
\def\xhigh{{\bf x}}
\def\xlow{{\rm x}}
\def\rB{{\rm B}}
\def\rX{{\rm X}}
\def\rx{{\rm x}}
\def\ra{{\rm a}}
\def\rb{{\rm b}}
\def\rc{{\rm c}}
\def\rt{{\rm t}}
\def\rD{{\rm D}}

\newcommand{\intr}{\operatorname{int}}
\makeop{min}

\begin{document}

\begin{htmlabstract} 
Choose any oriented link type X and closed braid
representatives X<sub>+</sub>, X<sub>-</sub> of X, where X<sub>-</sub> has minimal braid
index among all closed braid representatives of X.  The main
result of this paper is a `Markov theorem without stabilization'. It
asserts that there is a complexity function and a finite set of
`templates' such that (possibly after initial complexity-reducing
modifications in the choice of X<sub>+</sub> and X<sub>-</sub> which replace them with
closed braids X<sub>+</sub>', X<sub>-</sub>') there is a sequence of closed
braid representatives X<sub>+</sub>' =
X<sup>1</sup>&rarr;X<sup>2</sup>&rarr;&hellip;&rarr;X<sup>r</sup> =
X<sub>-</sub>'
such that each passage X<sup>i</sup>&rarr;X<sup>i+1</sup> is strictly complexity  reducing
and non-increasing on braid index. The templates which define the passages
X<sup>i</sup>&rarr;X<sup>i+1</sup> include 3 familiar ones, the destabilization,  exchange
move and flype templates, and in addition, for each braid index m&ge; 4
a finite set T(m) of new ones.  The number of templates in
T(m) is a non-decreasing function of m.  We give examples of
members of T(m), m&ge; 4, but not a complete listing. There
are applications to contact geometry, which will be given in a separate
paper.
\end{htmlabstract}

\begin{webabstract} 
Choose any oriented link type $\mathcal{X}$ and closed braid
representatives $X_+, X_-$ of $\mathcal{X}$, where $X_-$ has minimal braid
index among all closed braid representatives of $\mathcal{X}$.  The main
result of this paper is a `Markov theorem without stabilization'. It
asserts that there is a complexity function and a finite set of
`templates' such that (possibly after initial complexity-reducing
modifications in the choice of $X_+$ and $X_-$ which replace them with
closed braids $X_+^\prime, X_-^\prime$) there is a sequence of closed
braid representatives $X_+' = X^1\to X^2\to \cdots \to X^r = X_-^\prime$
such that each passage $X^i\to X^{i+1}$ is strictly complexity  reducing
and non-increasing on braid index. The templates which define the passages
$X^i \to X^{i+1}$ include 3 familiar ones, the destabilization,  exchange
move and flype templates, and in addition, for each braid index $m\geq 4$
a finite set $\mathcal{T}(m)$ of new ones.  The number of templates in
$\mathcal{T}(m)$ is a non-decreasing function of $m$.  We give examples of
members of $\mathcal{T}(m), m\geq 4$, but not a complete listing. There
are applications to contact geometry, which will be given in a separate
paper.
\end{webabstract}

\begin{abstract} 
Choose any oriented link type $\mathcal{X}$ and closed braid
representatives $X_+, X_-$ of $\mathcal{X}$, where $X_-$ has minimal braid
index among all closed braid representatives of $\mathcal{X}$.  The main
result of this paper is a `Markov theorem without stabilization'. It
asserts that there is a complexity function and a finite set of
`templates' such that (possibly after initial complexity-reducing
modifications in the choice of $X_+$ and $X_-$ which replace them with
closed braids $X_+^\prime, X_-^\prime$) there is a sequence of closed
braid representatives $X_+' = X^1\to X^2\to \cdots \to X^r = X_-^\prime$
such that each passage $X^i\to X^{i+1}$ is strictly complexity  reducing
and non-increasing on braid index. The templates which define the passages
$X^i \to X^{i+1}$ include 3 familiar ones, the destabilization,  exchange
move and flype templates, and in addition, for each braid index $m\geq 4$
a finite set $\mathcal{T}(m)$ of new ones.  The number of templates in
$\mathcal{T}(m)$ is a non-decreasing function of $m$.  We give examples of
members of $\mathcal{T}(m), m\geq 4$, but not a complete listing. There
are applications to contact geometry, which will be given in a separate
paper \cite{B-M_transverse}.
\end{abstract}

\maketitle

\section{Introduction}
\label{section:introduction}
\subsection{The problem}
\label{subsection:the problem}
Let $\cX$ be an oriented link type in
the oriented 3--sphere $S^3$ or $\reals^3 = S^3 \setminus \{\infty\}.$  A
representative $X \in \cX$ is said to be a {\it closed braid}
if there is an unknotted curve $\axis\subset (S^3 \setminus X)$ (the {\it axis}) 
\index{braid axis $\axis$}
and a choice of
{\it  fibration}
$\fib$ of the open solid torus $S^3 \setminus \axis$  by meridian discs
$ \{ H_{\theta}; \ \theta \in [0,2\pi] \}$,  such that whenever
$X$ meets a fiber
$H_\theta$ the intersection is transverse. We call the pair ($\axis,\fib$) a
{\it  braid structure}
  The fact that 
$X$ is a closed braid with respect to $\fib$ implies that the number of points in
$X \cap H_\theta$ is independent of
$\theta$. We call this number the {\it braid index} of $X$, and denote
it by the symbol $b(X)$. The
{\it braid index} of $\cX$,
denoted $b(\cX)$, is the minimum value of $b(X)$ over all closed braid representatives $X \in \cX$.

\begin{figure}[htpb]
\labellist\small
\pinlabel {$P$} at 60 586
\pinlabel {$+$} [r] at 75 558
\pinlabel {\tiny $A$} [tr] at 100 582
\pinlabel {\tiny $w$} at 103 639
\pinlabel {$1$} at 101 625
\pinlabel {\tiny destabilize} [b] at 168 642
\pinlabel {\tiny stabilize} [b] at 168 492

\pinlabel {$P$} at 196 586
\pinlabel {\tiny $A$} [tr] at 235 578
\pinlabel {\tiny $w$} at 239 639
\pinlabel {$1$} at 240 625

\pinlabel {$P$} at 338 586
\pinlabel {\tiny $A$} [tr] at 380 582
\pinlabel {\tiny $w$} at 380 639
\pinlabel {$1$} at 379 625
\pinlabel {$-$} [r] at 355 558
\pinlabel {\tiny destabilize} [b] at 446 642
\pinlabel {\tiny stabilize} [b] at 446 492

\pinlabel {$P$} at 475 586
\pinlabel {\tiny $A$} [tr] at 513 578
\pinlabel {\tiny $w$} at 517 639
\pinlabel {$1$} at 518 625
\endlabellist
\centerline{\includegraphics[scale=.7]{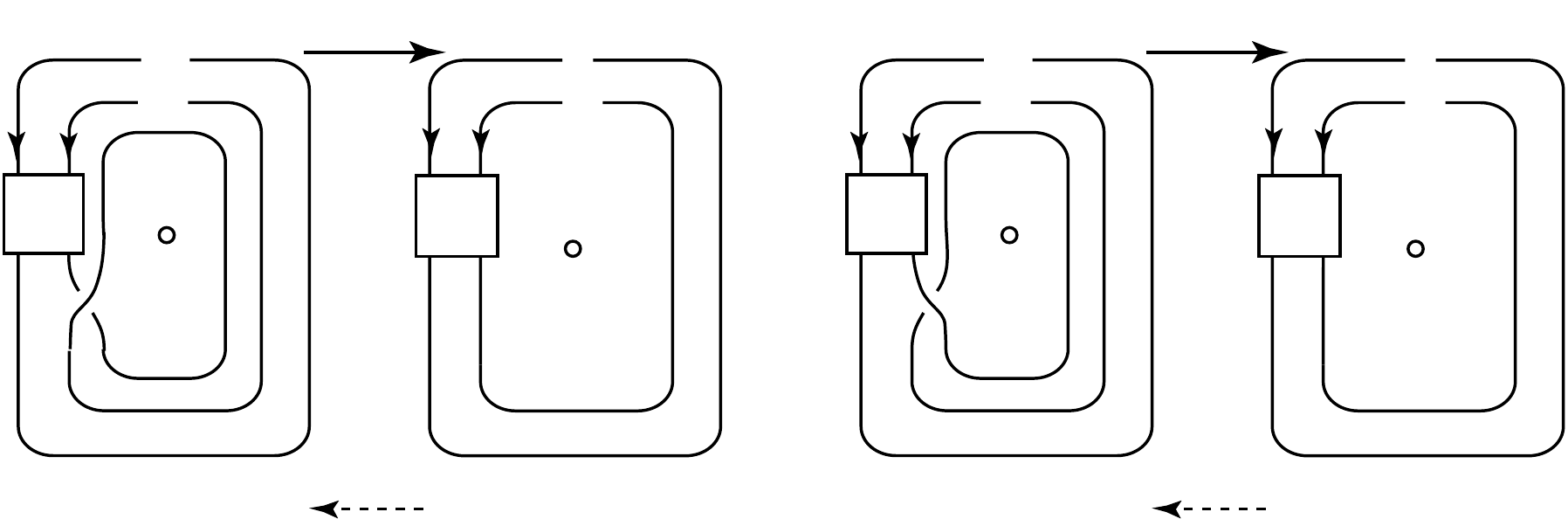}}
\caption{The two destabilization moves}
\label{figure:template1}
\end{figure}
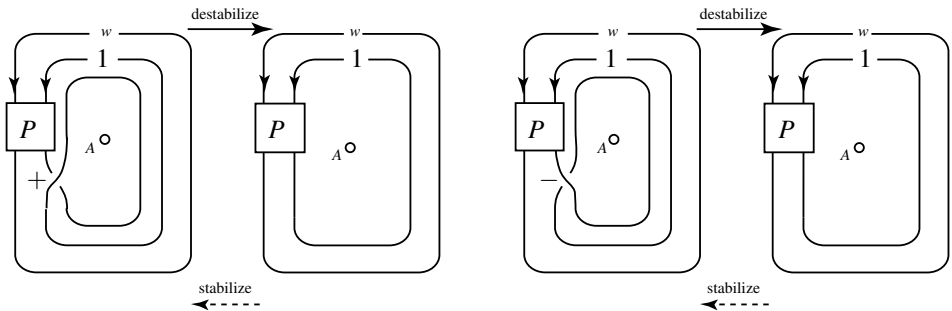

Closed braid representations of $\cX$ are not unique, and Markov's well-known 
theorem (see the book by Birman \cite{Birman1974}, and the papers by
Birman and Menasco \cite{B-M_Markov}, Lambropoulou and Rourke
\cite{Lambropoulou-Rourke}, Markov \cite{Markov}, Morton
\cite{Morton1986b} and Traczyk \cite{Traczyk})
asserts that any two are related by a finite
sequence of elementary moves. 
One of the moves is {\it  braid isotopy}
by which we mean an isotopy of the pair ($X, \reals^3\setminus\axis$) which preserves the condition
that $X$ is transverse to the fibers of $\fib$.
The
other two moves are mutually inverse, and are illustrated in \fullref{figure:template1}. 
Both take closed braids to closed braids. We call them
{\it  destabilization} and
{\it  stabilization}
where the former decreases braid
index by one and the latter increases it by one. The `weight' $w$ denotes $w$ parallel strands,
relative to the given projection. The braid inside the box which is labeled
$P$ is an arbitrary $(w+1)$--braid. Later, it will be necessary to distinguish between
positive and negative destabilizations, so we illustrate both now.  

\begin{theorem}[Markov's Theorem (MT) \cite{Markov}]
\label{theorem:MT}
Let $\Xhigh,\Xlow$ be closed
braid representatives of the same oriented link type $\cX$ in oriented
3--space, with the same braid axis $\axis$.  Then there exists a sequence
\begin{equation} \Xhigh = X_1 \to X_2 \to \cdots \to X_r = \Xlow 
\label{eq:markov}
\end{equation}
of closed braid representatives of $\cX$
such that, up to braid isotopy, each $X_{i+1}$ is obtained from $X_i$ by a single stabilization or
destabilization. 
\end{theorem}

It is easy to find examples of subsequences $X_j \to \cdots \to X_{j+k}$ of 
\eqref{eq:markov} in \mbox{\fullref{theorem:MT}} such that
$b(X_j) = b(X_{j+k})$, but  $X_j$ and $X_{j+k}$
 are not braid isotopic.
Call such a sequence a {\it  Markov tower}.
The stabilization and destabilization moves are very simple, but sequences of stabilizations, braid isotopies and
destabilizations can have unexpected consequences. In the braid groups these moves are `site
dependent', unlike the stabilization--destabilization move in the
Reidemeister--Singer Theorem.  (For an example
the reader should refer ahead to the specified site
of the stabilization in the sequence in \fullref{figure:exchange2}.)
Until now these moves have been predominately used to develop link invariants, but the Markov towers
themselves have been `black boxes'.  One of the main motivating ideas of this work is to open up
the black box and codify Markov towers.

Markov's Theorem is typical of an entire class of
theorems in topology where some form of stabilization and destabilization play a
central role. Other examples are:  
\begin{enumerate}
\item The Reidemeister--Singer Theorem \cite{Singer} relates any
two Heegaard
diagrams of the same 3--manifold, by a finite sequence of very simple
elementary changes on Heegaard diagrams.
The stabilization--destabilization move adds or deletes
a pair of simple closed curves $a,b$ in the defining
Heegaard diagram, where $a \cap b$ =
1 point and neither $a$ nor $b$ intersects any other curve $a_i, b_j$ in the
Heegaard diagram.  
\item The Kirby Calculus \cite{Kirby} gives a finite number of moves
which, when applied repeatedly, suffice to change any surgery
presentation of a given 3--manifold into any other, at the same time
keeping control of the topological type of a 4--manifold which the
given 3--manifold bounds. The stabilization--destabilization move is
the addition-deletion of an unknotted component with framing $\pm 1$
to the defining framed link.
\item Reidemeister's Theorem (see Burde--Zieschang \cite{B-Z}) relates any two diagrams of the same knot or link, by a
finite sequence of elementary moves which are known as RI, RII, RIII. The stabilization--destabilization
move is RI.  It is easy to see that Markov's Theorem implies Reidemeister's Theorem.   
\end{enumerate}

\begin{figure}[htpb]
\labellist\small
\pinlabel {\tiny $w$} [b] at 192 694
\pinlabel {\tiny $1$} [rb] at 181 677
\pinlabel {$P$} at 169 645
\pinlabel {\rotatebox{180}{\reflectbox{$Q$}}} at 252 645

\pinlabel {\tiny $w$} [b] at 372 694
\pinlabel {\tiny $1$} [rb] at 362 677
\pinlabel {$P$} at 350 645
\pinlabel {\rotatebox{180}{\reflectbox{$Q$}}} at 432 645
\endlabellist
\centerline{\includegraphics[scale=.8]{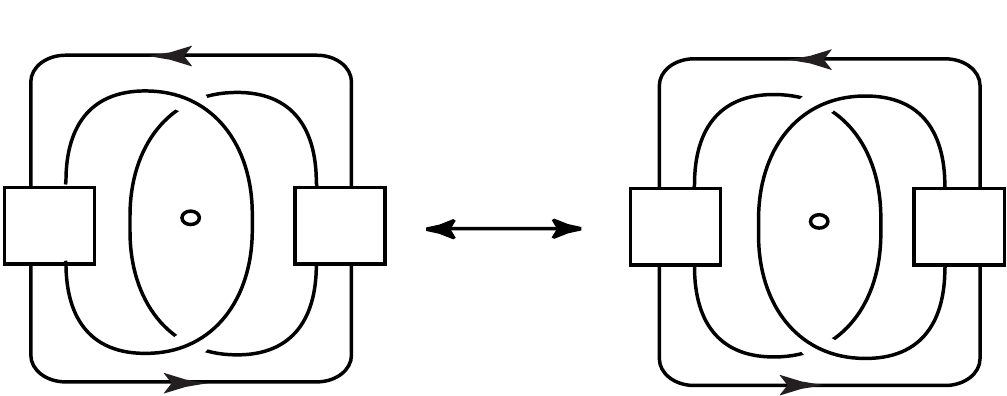}}
\caption{The exchange move}
\label{figure:exchange1}
\end{figure}
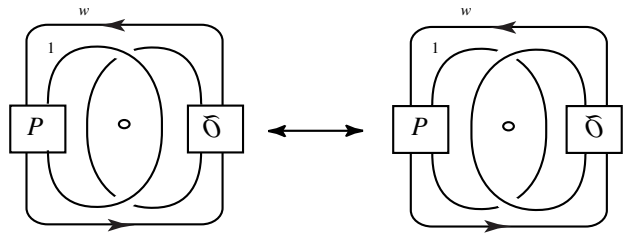

These theorems are all
like Markov's Theorem in the sense that while the
stabilization and destabilization  moves are very simple, a sequence
of these moves, combined with the appropriate
isotopy, can have very non-trivial consequences. Here are other examples in
which the stabilization move is not used, at the expense of restricting
attention to a special example: 
\begin{enumerate}
\item[(4)]  W Haken proved that any Heegaard
diagram for a non-prime 3--manifold is equivalent to a Heegaard diagram which is
the union of two separate Heegaard diagrams, one for each summand, supported on
disjoint subsets of the given Heegaard surface. See Scharlemann--Thompson \cite{Scharlemann-Thompson} for a
very pleasant proof.  
\item[(5)] Waldhausen \cite{Waldhausen} proved that any
two Heegaard diagrams of
arbitrary but fixed genus $g$ for the 3--sphere $S^3$ are
equivalent. 
\end{enumerate}
 In the course of an effort which we began in 1990 to discover the theorem which will be
the main result of this paper (see \fullref{theorem:MTWS} below) the authors made several
related contributions to the theory of closed braid representatives of knots and links:
\begin{enumerate}
\item [($4^{\prime}$)] A {\it  split} (resp.\  {\it  composite})
closed
$n$--braid is an $n$--braid which factorizes as a product $X Y$ where the sub-braid $X$ involves only
strands
$1,\dots,k$ and the sub-braid $Y$ involves only strands $k+1,\dots,n$ \
(resp.\  $k,\dots,n$).
In the manuscript \cite{B-M_IV} the authors proved that if
$X$ is a closed $n$--braid representative of a split or composite link, then up to
(braid-index preserving) isotopy and {\it  exchange} moves, 
\index{exchange move template}
as in \fullref{figure:exchange1}, 
$X$ may be assumed to be a split or composite closed braid.  
\end{enumerate}

\begin{enumerate}
\item [($5^\prime$)] In the manuscript \cite{B-M_V} the authors proved that
if $X$ is a closed braid representative of the $\mu$--component unlink $\cX$,
then a finite sequence of braid isotopies, exchange moves  and destabilization
can be found which change
$X$ to the closure of the identity braid in the braid group $B_\mu$.
\item [(6)] In the manuscript \cite{B-M_III} the authors discovered 
that there is another move, the {\it  3--braid flype} (see \fullref{figure:flype1})
 with the property that if $X$
is a closed 3--braid representative of a knot or link type $\cX$ which cannot be represented by a
$1$--braid or $2$--braid, then either $X$ has a unique conjugacy class or $X$ has
exactly two conjugacy classes, and these two classes are related by a 3--braid flype.
They also showed that the exchange move can be replaced by braid isotopy for prime links of braid index 3.
\end{enumerate}
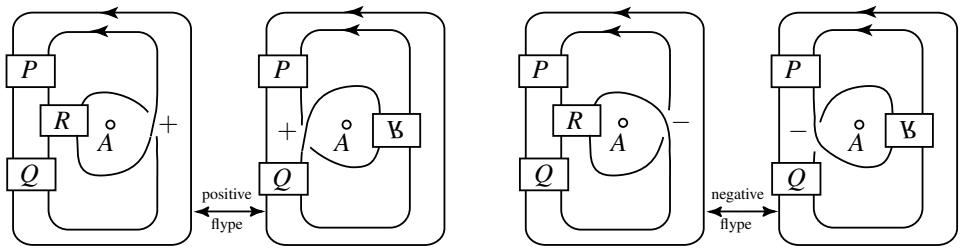
\begin{figure}[htpb]
\labellist\small
\pinlabel {$P$} at 60 663
\pinlabel {$Q$} at 60 605
\pinlabel {$R$} at 78 635
\pinlabel {$A$} [t] at 102 632
\pinlabel {$+$} [l] at 126 633

\pinlabel {$P$} at 198 663
\pinlabel {$Q$} at 198 603
\pinlabel {\reflectbox{\rotatebox{180}{$R$}}} at 260 630
\pinlabel {$A$} [t] at 231 632
\pinlabel {$+$} [r] at 211 629

\pinlabel {$P$} at 339 663
\pinlabel {$Q$} at 341 605
\pinlabel {$R$} at 358 635
\pinlabel {$A$} [t] at 382 632
\pinlabel {$-$} [l] at 406 634

\pinlabel {$P$} at 477 663
\pinlabel {$Q$} at 479 603
\pinlabel {\reflectbox{\rotatebox{180}{$R$}}} at 540 630
\pinlabel {$A$} [t] at 512 632
\pinlabel {$-$} [r] at 490 629

\pinlabel {\tiny positive} [b] at 168 586
\pinlabel {\tiny flype} [t] at 168 586

\pinlabel {\tiny negative} [b] at 447 586
\pinlabel {\tiny flype} [t] at 446 586

\endlabellist
\centerline{\includegraphics[scale=.7]{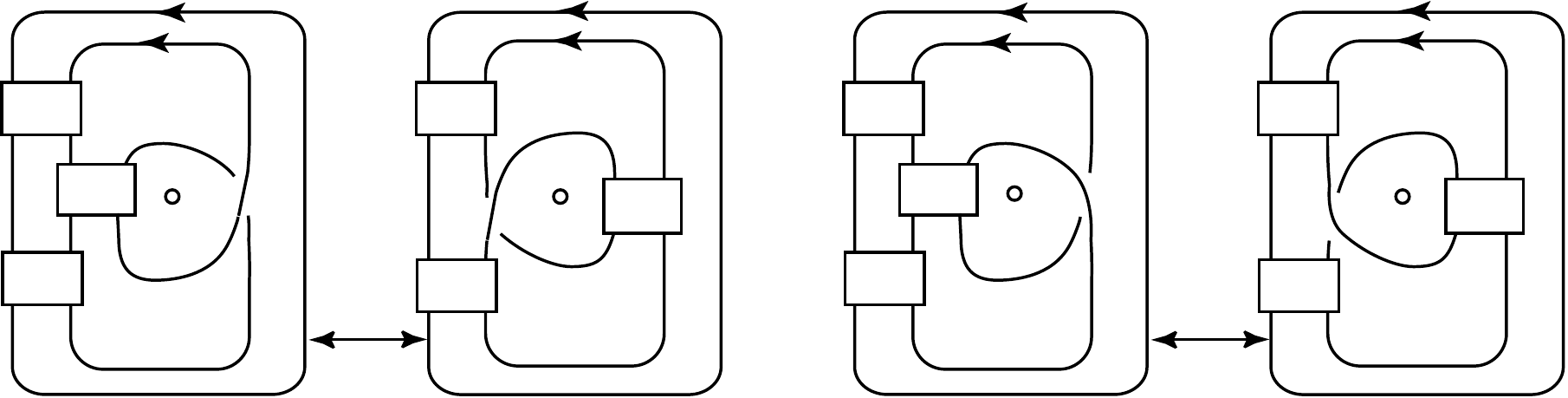}}
\caption{The two flype moves}
\label{figure:flype1}
\end{figure}
 The authors also established two fundamental facts which gave strong evidence that
a more general result might be true:
\begin{enumerate}
\item [(7)] In \cite{B-M_I} the authors introduced a complexity function on closed braid
representatives of $\cX$ and proved that, up to exchange moves, there are 
at most finitely many conjugacy classes of representatives of minimum complexity.

\item [(8)] In \cite{B-M_VI} the authors proved that if a link type $\cX$
has infinitely many conjugacy classes of closed braid representatives
of the same braid index, then, up to exchange moves,  they fall into
finitely many equivalence classes.
\end{enumerate}

The goal of this paper is to generalize examples (4$^\prime$),
(5$^\prime$) and (6), taking into account (7) and (8), to arbitrary
closed braid representatives of arbitrary oriented knots and links. We
call our main theorem {\em Markov's Theorem Without Stabilization}
(MTWS), because it is a direct modification of Markov's Theorem, but
with his stabilization move replaced by other moves  which allow one
to jump from one isotopy class in the complement of $\axis$ to another,
while keeping the braid index constant or decreasing it.

\subsection{Block-strand diagrams and templates}
\label{subsection:block strand diagrams and templates}
Before we can state our main result, we need to introduce new concepts. 
Our moves will be  described in terms of certain pairs of `block-strand
diagrams' which we call `templates'.  Examples are the block-strand
diagram pairs which make up the templates in Figures
\ref{figure:template1}, \ref{figure:exchange1} and \ref{figure:flype1}.   The reader may
wish to look ahead to the boxed pairs of block-strand diagrams in Figures
\ref{figure:template0}, \ref{figure:template2} and 
\ref{figure:exchange6} for  examples of more complicated templates.

A {\it  block} $B$ in $\reals^3-\axis$ is a 3--ball having the structure of
a $2$--disk $\Delta$ crossed with an interval $[0,1]$ such that (i) for any fiber
$H_\theta \in \fib$ the intersection $H_\theta\cap B$ is either $\emptyset$ or $\Delta\times\{\theta\}$
for some $\theta\in [\theta_1,\theta_2]$, and (ii) there exists some $\theta\in )0,2\pi)$ such that $H_\theta\cap B
= \emptyset$.  The disc $ \bt = B \cap H_{\theta_1} $ is the {\it  top} of $B$ and the disc
$ \bb = B \cap H_{\theta_2} $ is the {\it  bottom} of $B$.
A {\it  strand} $l$ is homeomorphic to an interval $[0,1]$ or a circle $S^1$. It is oriented 
and transverse to each fiber of $\fib$ such that its orientation agrees with the
forward direction of $\fib$.  When $l$ is homeomorphic to an interval, $\partial l = l_0 \cup l_1 $,
where $l_0$ is the {\it  beginning endpoint} of $l$ and $l_1$ is the {\it  ending endpoint} of $l$.
A {\it  block-strand diagram}
$\cD$ is a collection of pairwise disjoint blocks 
$\{B^1, \ldots, B^{\rm k} \}$ and pairwise disjoint strands $\{l^1,\ldots,l^{\rm l}\}$
which together have the following structure:
\bi
\item[(1)] If $l^i \cap B^j \not= \emptyset$
then $l^i \cap B^j = (l_1^i \cap \bt^j ) \cup (l_0^i \cap \bb^j) $ where $\bt^j$ is the top of $B^j$
and $\bb^j$ is the bottom of $B^j$.  (We allow for the
possibility that either $(l_1^i \cap \bt^j )$ or $(l_0^i \cap \bb^j)$ is empty.)
\item[(2)] For each $l_0^i$ (resp.\  $l_1^i$) there is some $\bb^j
\subset B^j$ (resp.\  $\bt^j \subset B^j$)
such that $l_0^i \subset \bb^j$ (resp.\  $l_1^i \subset \bt^j$).
\item[(3)] For each block $B^j$ we have
$| \bt^j \cap (\cup_{1\leq i \leq {\rm l}} l_1^i ) | = | \bb^j \cap (\cup_{1\leq i \leq {\rm l}} l_0^i ) | \geq 2 $.
\ei
The fact that for each $j=1,\dots,k$ there is a fiber which misses $B^j$ shows that, by rescaling, we may
find a distinguished fiber $H_{\theta_0}$ which does not intersect any block. We define the {\it  braid index}
$b(\cD)$ of the block-strand diagram
$\cD$
to be the number of times  the strands of
$\cD$ intersect the distinguished fiber $H_{\theta_0}$.  Condition (3) above makes $b(\cD)$ well defined.
For a specified block $B^j \subset \cD$ we define its {\it  braid index} $b(B^j) =
| \bt^j \cap (\cup_{1\leq i \leq {\rm l}} l_1^i ) |$.  Continuing the definition of a block-strand diagram, we
assume:
\bi
\item [(4)] If $B^j \subset \cD$ then the braid index $b(B^j)$ is strictly
less than the braid index $b(\cD)$, with one exception. The exception
occurs when  $\cD$ is the block strand diagram that results after a
destabilization has been performed, as in \fullref{figure:template1}.  
\ei

A {\it  template} $\cT$ 
\index{template $\cT$} 
is a pair of block-strand diagrams
$(\cD_+,\cD_-)$, both with blocks $B^1\!,\dots,B^{\rm k}$ and an isotopy which
takes $\cD_+$  to $\cD_-$, in such a way that 
for every fixed choice of braiding assignments to the blocks
$B^1,\dots,B^{\rm k}$ the resulting closed braids $X,X'$ represent the same oriented link
type $\cX$.  The diagrams $\cD_+$ and $\cD_-$ are the {\it  initial} and {\it   final} \ block-strand
diagrams in the pair.  The {\it  fixed blocks} and {\it  fixed strands}  in
$\cT = (\cD_+,\cD_-)$  are the blocks  and strands where the isotopy is pointwise the identity. All other
blocks and strands are {\it  moving}. For example, in \fullref{figure:flype1} the blocks $P$ and $Q$ are
fixed blocks, whereas
$R$ is a moving block.
A {\it  braiding assignment} to a block-strand diagram $\cD$ 
is a choice of a braid on $m_j$ strands for each
$B^j\in\cD$. That is, we replace $B^j$ with the chosen braid, so that $B^j$ with
this braiding assignment becomes a braided tangle with $m_j$ in-strands and $m_j$
out-strands.  In this way a block strand diagram gives us a closed braid
representative of a link $\cX$.  

Let $X$ be a closed $m$--braid. We say that $X$ is   {\it  carried} by $\cD$
 if there exists a
braiding assignment for the blocks in $\cD$ such that the resulting closed braid is braid-isotopic to
$X$. 

When we first began to understand that templates were
the appropriate settings for our work on the MTWS we wondered whether our definition
was so broad (because the diagrams in question support so many knot and link types) as
to be content-free! In this regard, the following fact is fundamental: 

\begin{proposition} 
\label{proposition:block-strand diagrams don't carry all links}
Let $\cD$ be a block-strand diagram of braid index $n$.  Then there exist
$n$--braids that $\cD$ does not carry.  
\end{proposition}

\begin{proof} 
Up to conjugation, a block-strand diagram may
be described by a word $V_1 W_1 V_2 W_2 .... V_k W_k$ in the standard elementary braid generators
$\sigma_1,\dots,\sigma_{n-1}$ of the
$n$--strand braid group, where each $V_j$ represents a word which describes the braid carried by the $j^{th}$
block (after making a braiding assignment to the block) and each  $W_i$ is a braid word on $n$ strands which
describes the strands that connect the blocks.   By hypothesis no block has more than $n-1$ strands entering or
leaving it, so by modifying the
$W_i$'s we may assume without loss of generality that each $V_j$ is a braid on the first
$q_j$--strands, where $q_j<n$.  After this modification,  the only places where the
elementary braid generator $\sigma_{n-1}$ appears is in the braid words that describe the strands that join the
blocks, ie the words $W_1,\dots,W_k$.

Now let $|W_i|$ be the number of times $\sigma_{n-1}$ occurs in $W_i$.  
The $|W_i|'s$  are fixed numbers since we were handed a block-strand diagram.
For an arbitrary conjugacy class $\{C\}$ of $n$--braids, let $|C|$ be the minimum number
of times the generator $\sigma_{n-1}$ is used, in all possible words which represent
$\{C\}$. Our block diagram can only carry closed n-braids $C$ such that (up to
conjugacy) $|C| \leq |W_1| + |W_2| + ... + |W_k|$.  But there are closed n-braids 
$X$ such that  $|Y|$ is arbitrarily large for all $Y\in \{X\}$. An example
is $X = (\sigma_1\sigma_2\dots\sigma_{n-1})^{nN}$ for a large positive integer $N$. For, the braid word 
$(\sigma_1\sigma_2\dots\sigma_{n-1})^{n}$ is a full twist of the braid strands. It generates the center of
the n-string braid group.  It cannot be represented by any braid word that does not use all the elementary
braid generators. Therefore it's $N^{th}$ power uses the generator $\sigma_{n-1}$ at least $nN$ times. 
\end{proof} 

 Having \fullref{proposition:block-strand
diagrams don't carry all links} on hand, we proceed to define the templates that we will use in the statement of
the MTWS. Our main theorem begins with an arbitrary closed
$n$--braid representative
$\Xhigh$ of an arbitrary oriented knot or link type $\cX$ in 3--space. Let $\Xlow$ be a second such
representative, where $b(\Xlow) = b(\cX)$.  Our goal in the subsections which follow is to describe some of
the templates that we need, and at the same time to describe the building blocks of all of
them.
Note that we regard braid isotopy as a trivial move, sometimes even
forgetting to mention it. By a result of Morton
\cite[Theorem~1]{Morton1978} braids $\beta,\beta^\prime$ in the $n$--string
braid group $B_n$ are conjugate if and only if the associated closed
braids are isotopic in the complement of the braid axis.  In keeping
with our motivating idea of codifying Markov towers, the names that we
give some of our templates $\cT = (\Dhigh, \Dlow)$ correspond to the
name of the isotopy that is used to move $\Dhigh$ to $\Dlow$.

\subsubsection {The two destabilization templates} 
\label{subsubsection:destabilization}
Our two {\it  destabilization} templates were defined in
\fullref{figure:template1}.  We distinguish the cases of positive and negative
destabilization because the strands which join the fixed blocks are different, and so
the templates are different.  The
destabilization templates do not have any moving blocks. They occur at every braid index.

\subsubsection{The admissible flype templates}
\label{subsubsection:flypes}
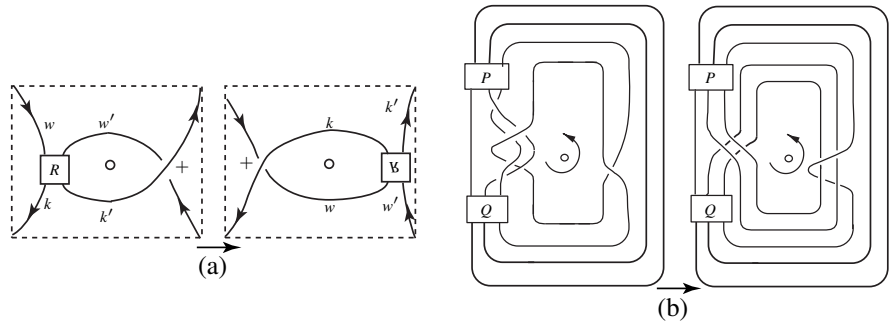
\begin{figure}[htpb!]
\labellist\tiny
\pinlabel {$R$}      at 72 484
\pinlabel {$w$}  [l] at 61 510
\pinlabel {$w'$} [b] at 103 504
\pinlabel {$k$}  [l] at 61 464
\pinlabel {$k'$} [t] at 103 466
\pinlabel {$+$}  [l] at 138 486

\pinlabel {\reflectbox{\rotatebox{180}{$R$}}} at 272 486
\pinlabel {$+$}  [r] at 194 488
\pinlabel {$k$}  [b] at 234 504
\pinlabel {$w$}  [t] at 234 468
\pinlabel {$k'$} [r] at 280 522
\pinlabel {$w'$} [r] at 280 464

\pinlabel {\small(a)}  [t] at 165 438

\pinlabel {$P$}      at 326 538
\pinlabel {$Q$}      at 326 460

\pinlabel {$P$}      at 458 538
\pinlabel {$Q$}      at 458 460

\pinlabel {\small(b)}  [t] at 436 413

\endlabellist
\centerline{\includegraphics[scale=.65]{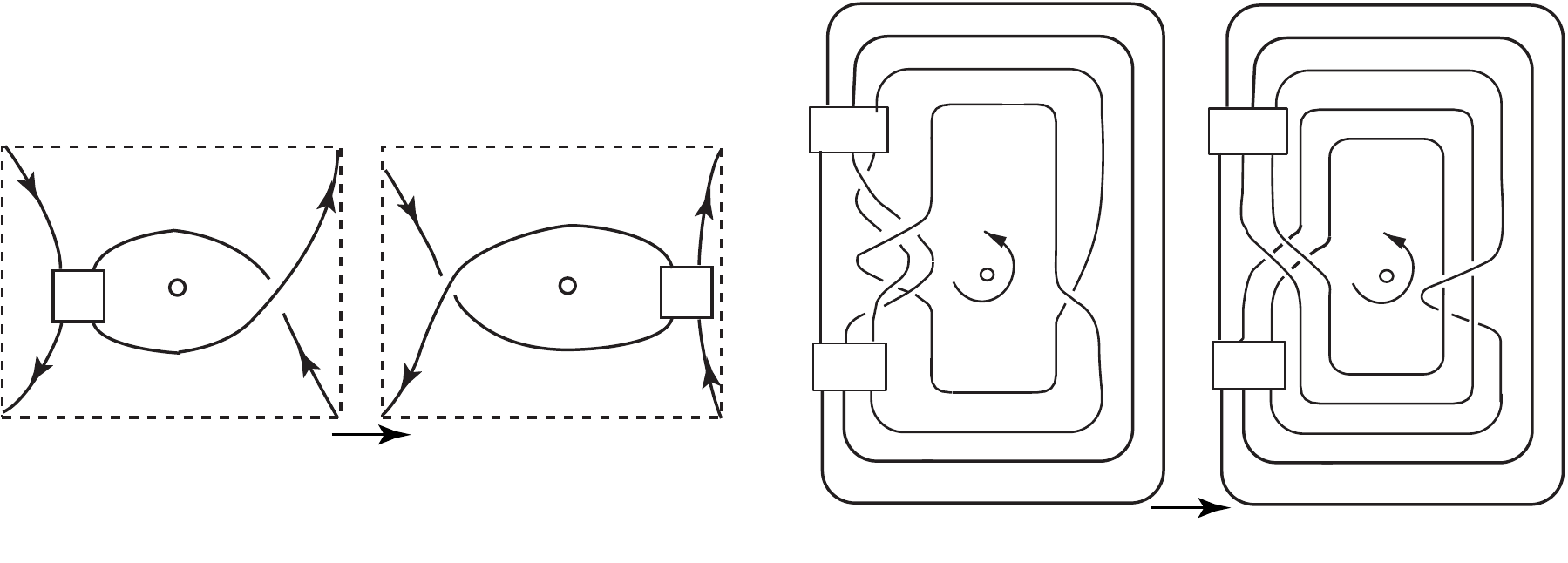}}
\caption{(a) The support of a positive flype with weighted strands. (b) Example of an inadmissible flype. Notice the
extra twists introduced because of the weighted strands.}  
\label{figure:flype2}
\end{figure}

Flypes first enter the picture when the braid index is 3, and we already illustrated the two
3--braid flype templates in \fullref{figure:flype1}. There is an obvious way to generalize it to any braid index $n$, namely declare
the strands to be weighted strands.  See  \fullref{figure:flype2}(a), which shows
the support of the flype with weighted strands.  From now on, the term {\it  flype}
\index{flype, defined as a move on a closed braid diagram}
will always have this meaning.  The  {\it  sign of a flype}
\index{sign of a flype} is the sign of
the single crossing (possibly weighted) which is not in the braid block.  Both positive and negative
flypes are illustrated in   \fullref{figure:flype1}. They have distinct templates. 

There is a subtle point: 
 Let $\Xhigh$ and $\Xlow$ be the closed braid before and after a flype, which
we shall consider (for the purpose of describing our moves) as acting left to right.  The flype motion 
is supported in a 3--ball $B^3$. In \fullref{figure:flype2}(a) observe that the fiber
$H_\theta$ at $\theta=\pi/2$ intersects $\Xhigh\cap B^3$ in $w'$ points, but intersects $\Xlow$ in $k$
points. Observe that $w+w' = k+k'$.  We have shown that
$$b(\Xhigh) - b(\Xlow) = w^\prime - k  = k' - w.$$
Thus flypes with weighted strands are non-increasing on braid index if and only if
$w'-k = k'-w \geq 0$. We will refer to a flype which is non-increasing on braid
index as an {\it  admissible flype}. An example of an inadmissible flype is given in \fullref{figure:flype2}(b).
While we are obviously interested in the admissible flypes, it will turn out that the {\it  inadmissible
flypes} 
\index{inadmissible flype} 
are important too, as they lead to additional, more complicated templates. 

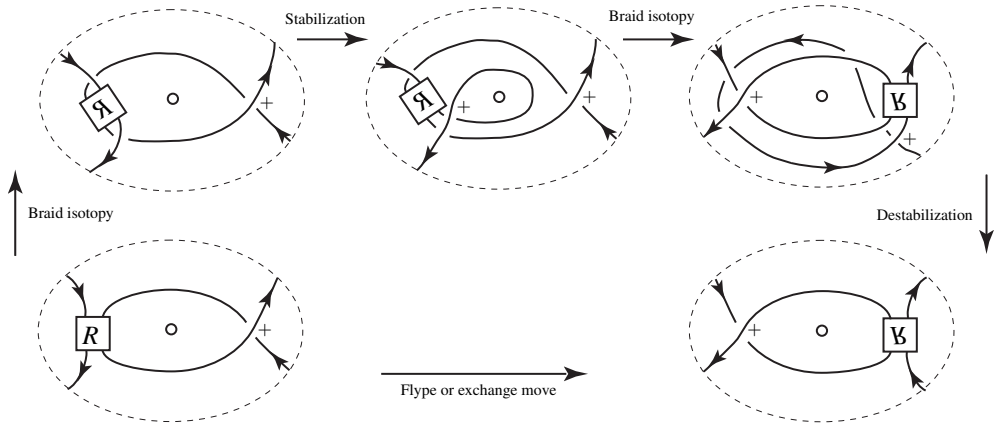
\begin{figure}[htpb!]
\labellist\tiny
\hair=5pt
\pinlabel {{Braid isotopy}}          [l] at 61 551
\pinlabel {{Stabilization}}          [b] at 217 640
\pinlabel {{Braid isotopy}}          [b] at 381 640
\pinlabel {{Destabilization}}        [r] at 549 551
\pinlabel {{Flype or exchange move}} [t] at 294 472

\hair=3pt
\pinlabel {\rotatebox{30}{\reflectbox{\small $R$}}} at 105 604
\pinlabel {\small $R$} at 99 490
\pinlabel {\rotatebox{30}{\reflectbox{\small $R$}}} at 266 608
\pinlabel {\rotatebox{180}{\reflectbox{\small $R$}}} at 505 608
\pinlabel {\rotatebox{180}{\reflectbox{\small $R$}}} at 505 490

\pinlabel {$+$} [l] at 180 607
\pinlabel {$+$} [bl] at 280 599
\pinlabel {$+$} [l] at 342 610
\pinlabel {$+$} [l] at 426 610
\pinlabel {$+$} [l] at 503 589
\pinlabel {$+$} [l] at 179 493
\pinlabel {$+$} [l] at 425 493
\endlabellist
\centerline{\includegraphics[scale=.75]{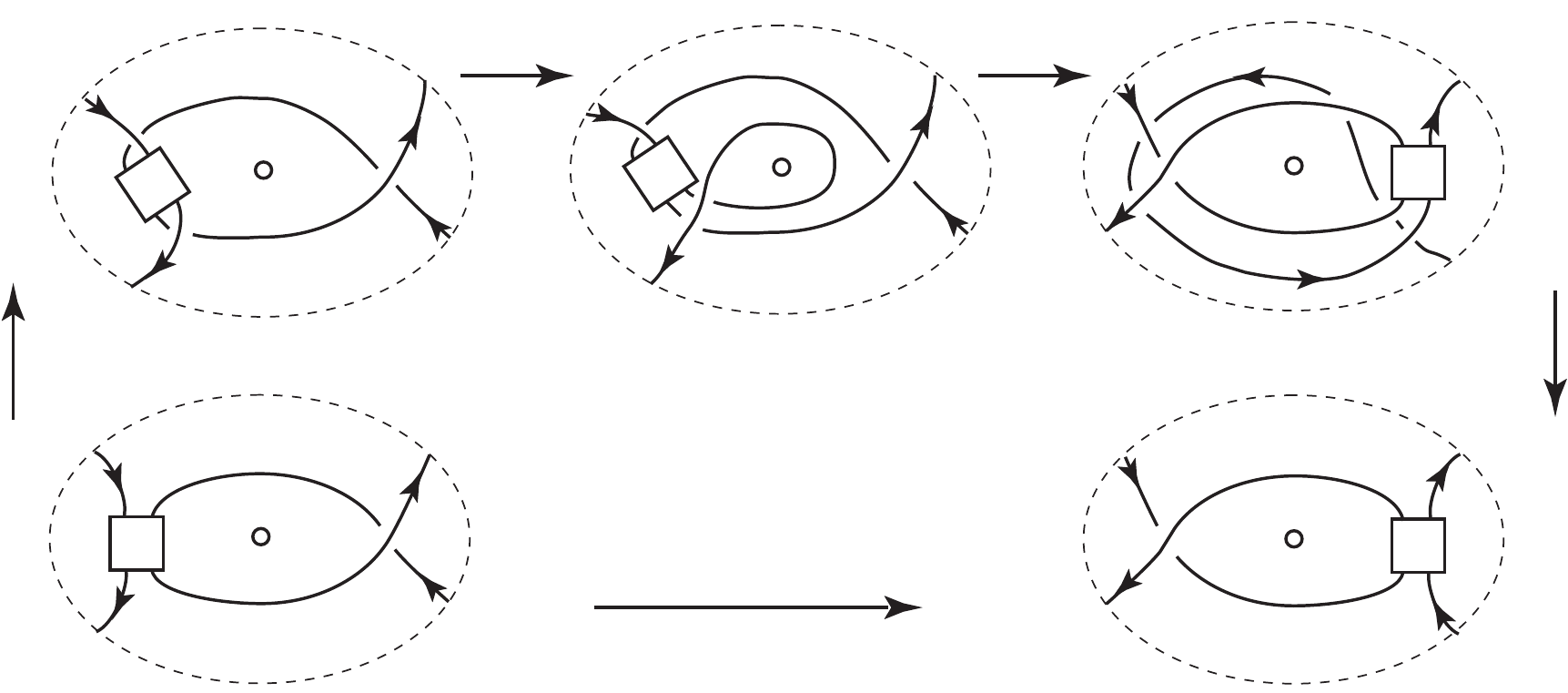}}
\caption{A very simple Markov tower}
\label{figure:exchange2}
\end{figure}

By Markov's Theorem, the 
left and right braids in every admissible flype template 
must  be related by a Markov tower. \fullref{figure:exchange2} shows such
a 2--step tower, in the case when the braid index
is 3.   The moves used in the sequence
are (up to braid isotopy) a single stabilization and a single
destabilization. Thus flypes arise in a very natural way in the
study of stabilization in the braid groups: they replace a sequence 
stabilization, destabilization  by a single braid-index preserving (or possibly
reducing) move. 
Notice that when
$R$ is a negative half-twist, the tower can be replaced by an exchange move.
  Observe that since flypes  are
replacements for Markov towers, we are now free to use them to construct more general Markov towers.

\subsubsection{The exchange move template and sequences of exchange moves}
\label{subsubsection:exchange moves and sequences of exchange moves}
The exchange move template was defined by the block-strand diagram
in \fullref{figure:exchange1}.  It was proved in our earlier paper
\cite{B-M_III}  that for $n\leq 3$ it is equivalent to braid isotopy, and
by Fiedler \cite{Fiedler} that for $n\geq 4$, and generic choices of the
braids $P$ and $Q$, the exchange move cannot be replaced by braid isotopy.

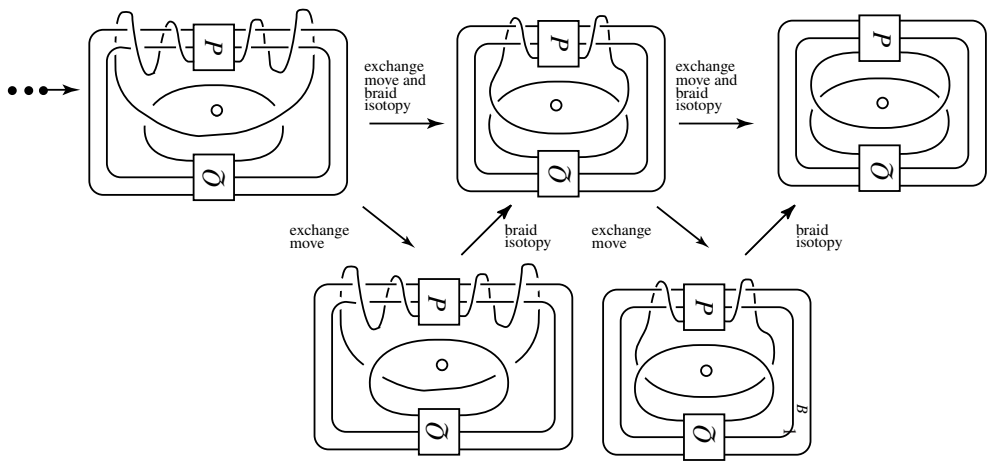
\begin{figure}[htpb!]
\labellist\tiny
\pinlabel {\parbox{30pt}{\smash{exchange} move and\\braid\\isotopy}}
  [b] at 260 517
\pinlabel {\parbox{30pt}{\smash{exchange} move and\\braid\\isotopy}}
  [b] at 430 517
\pinlabel {\parbox{30pt}{exchange\\move}} [tr] at 244 464
\pinlabel {\parbox{30pt}{braid\\isotopy}} [tl] at 312 464
\pinlabel {\parbox{30pt}{exchange\\move}} [tr] at 409 464
\pinlabel {\parbox{30pt}{braid\\isotopy}} [tl] at 471 464
\pinlabel {\rotatebox{270}{\small $P$}} at 156 558
\pinlabel {\rotatebox{90}{\small $Q$}} at 156 486
\pinlabel {\rotatebox{270}{\small $P$}} at 345 560
\pinlabel {\rotatebox{90}{\small $Q$}} at 345 488
\pinlabel {\rotatebox{270}{\small $P$}} at 521 562
\pinlabel {\rotatebox{90}{\small $Q$}} at 521 490
\pinlabel {\rotatebox{270}{\small $P$}} at 278 418
\pinlabel {\rotatebox{90}{\small $Q$}} at 278 348
\pinlabel {\rotatebox{270}{\small $P$}} at 424 416
\pinlabel {\rotatebox{90}{\small $Q$}} at 424 346
\pinlabel {\rotatebox{270}{\tiny $B$}} at 477.5 360
\endlabellist
\centerline{\includegraphics[scale=.7]{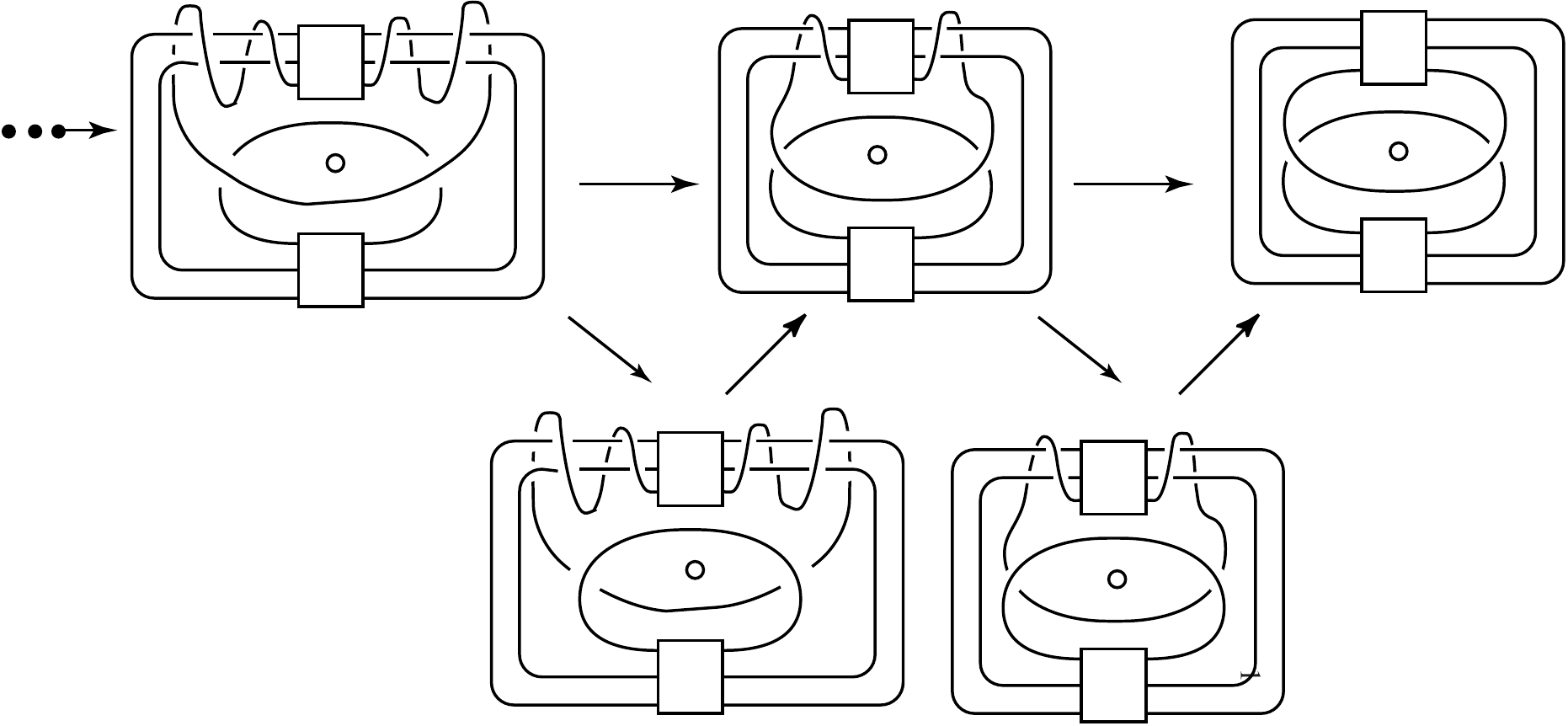}}
\caption{The exchange move can lead
to arbitrarily many distinct braid isotopy classes of closed
$n$--braid representatives of a single knot type, $n\geq 4$.}
\label{figure:exchange3}
\end{figure}
\fullref{figure:exchange3} shows how exchange moves, together with
braid isotopy,  can lead to infinitely many conjugacy classes of closed braid
representatives of the same knot or link (see our earlier paper
\cite{B-M_V}).   Indeed, in
\cite{B-M_VI} the authors proved that if a link has infinitely many conjugacy
classes of closed $m$--braid representatives for any fixed value
of $m$ then all but finitely many of them are related by exchange moves. 
This fact will shape the form of our main theorem. More precisely, our main theorem
shows exactly how far one may go, using only exchange moves and destabilizations, and then identifies the
finitely many moves which are needed {\it in addition to} exchange moves and destabilizations, to take one
closed
$n$--braid representative of a knot to another of the same braid index. 

\begin{figure}[htpb]
\labellist\tiny
\pinlabel {\parbox{30pt}{unfold tree}} [t] at 302 452
\pinlabel {\parbox{35pt}{\smash{sequence} of exchange moves}} [l] at 437 446
\pinlabel {\parbox{30pt}{collapse tree}} [t] at 316 321
\pinlabel {\small $B_1$} at 345 418
\pinlabel {\small $B_2$} at 395 416
\pinlabel {\rotatebox{15}{$B_1''$}} at 342 502.5
\pinlabel {\rotatebox{345}{$B_2''$}} at 396 502.5
\pinlabel {$B'$} at 371 466
\endlabellist
\centerline{\includegraphics[scale=.8]{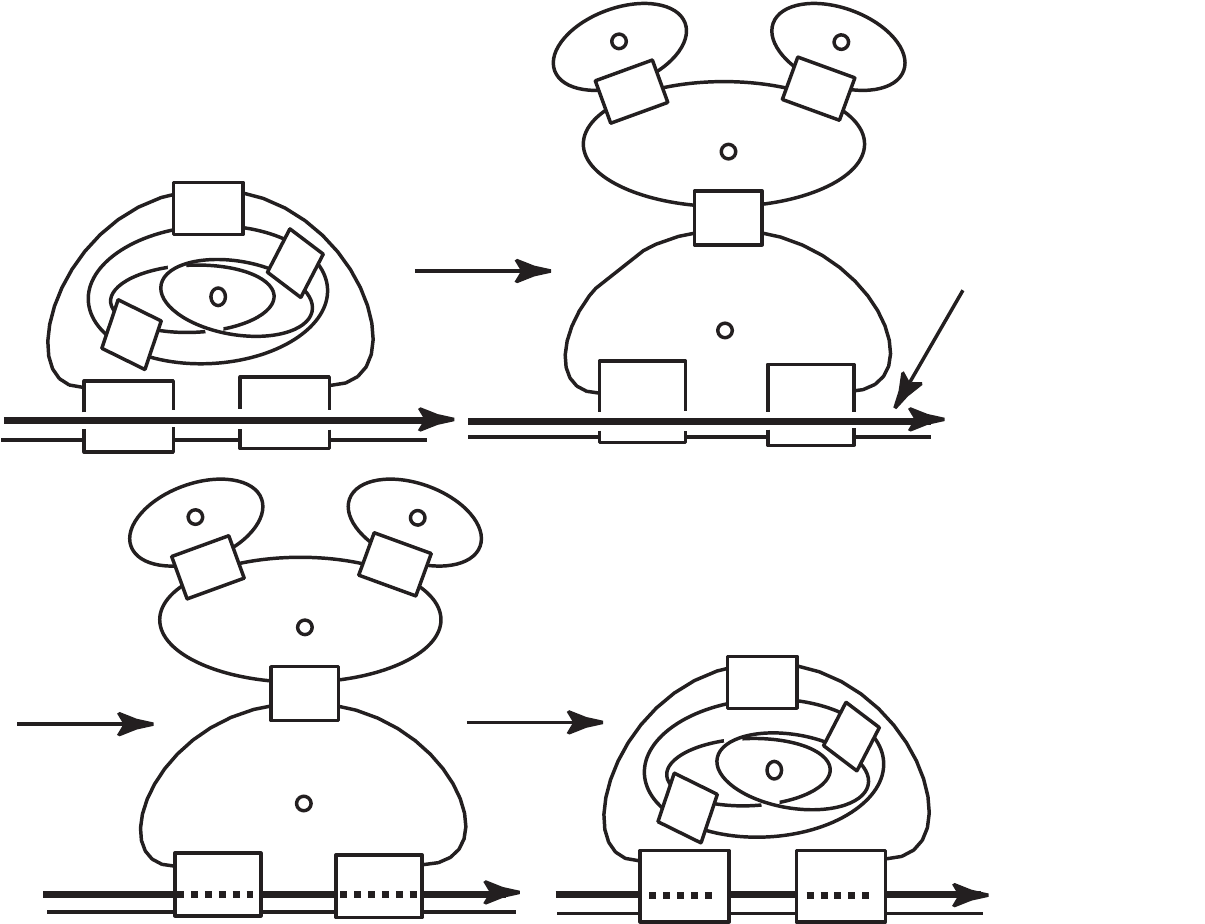}}
\caption{A sequence of exchange moves passes a distinguished strand
(the thick black one) over a block-strand tree}
\label{figure:handle0}
\end{figure}
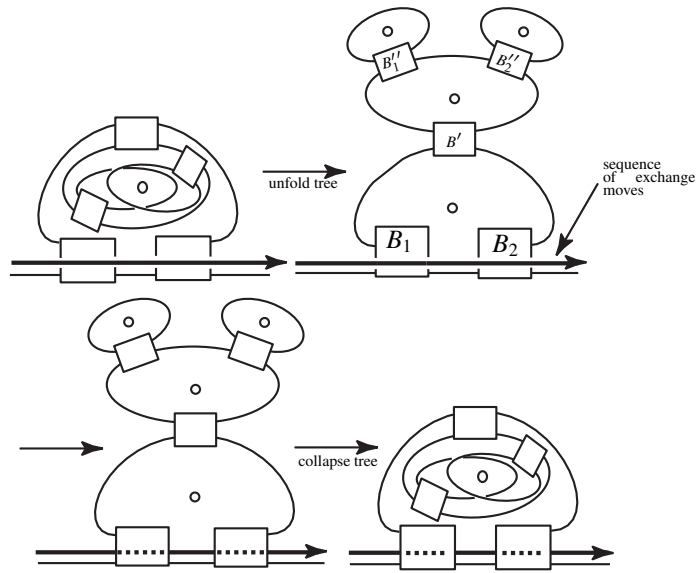
The sequences of exchange moves which we next define are very useful and important.  \fullref{figure:handle0} shows how we unfold a piece of a closed
braid to reveal that it has the structure of a `rooted block and strand
tree', and then `loop' a distinguished subarc of the braid (always of weight 1)
over the tree.  (These concepts will be defined
in \fullref{subsection:pushing across regions with a G-exchange foliation}.   See \fullref{figure:handle1} and the nearby text. 
We hope the reader will be patient.  Our initial goal is to state our main result.)  
Notice that, while the unlooping process does not preserve closed braids, 
we have retained the
closed braid structure by the device of cutting the
braid axis into 4 little `axis pieces'. Of course the fibers of $\fib$ are
arranged radially around these little axis pieces, in a sufficiently small
neighborhood, so that when we open up the tree we can retain a local
picture of the braid structure. During the looping
motion the distinguished strand cuts each axis piece twice.  

\subsubsection{The cyclic templates}
\label{subsubsection:the cyclic templates}
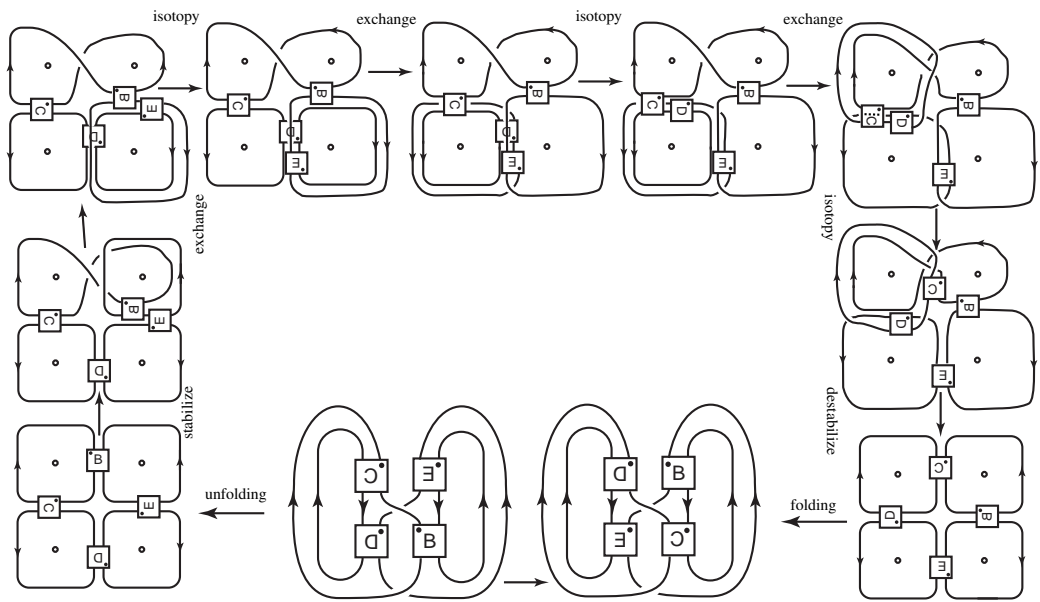
\begin{figure}[htpb]
\labellist\tiny
\pinlabel {isotopy} at 131 575
\pinlabel {exchange} at 230 574
\pinlabel {isotopy} at 330 575
\pinlabel {exchange} at 431 574
\pinlabel {\rotatebox{270}{isotopy}} at 438 481
\pinlabel {\rotatebox{270}{destabilize}} at 440 385
\pinlabel {folding} [b] at 431 338
\pinlabel {unfolding} [b] at 159 343
\pinlabel {\rotatebox{90}{stabilize}} at 136 389
\pinlabel {\rotatebox{90}{exchange}} at 142 478
\endlabellist
\centerline{\includegraphics[scale=.80]{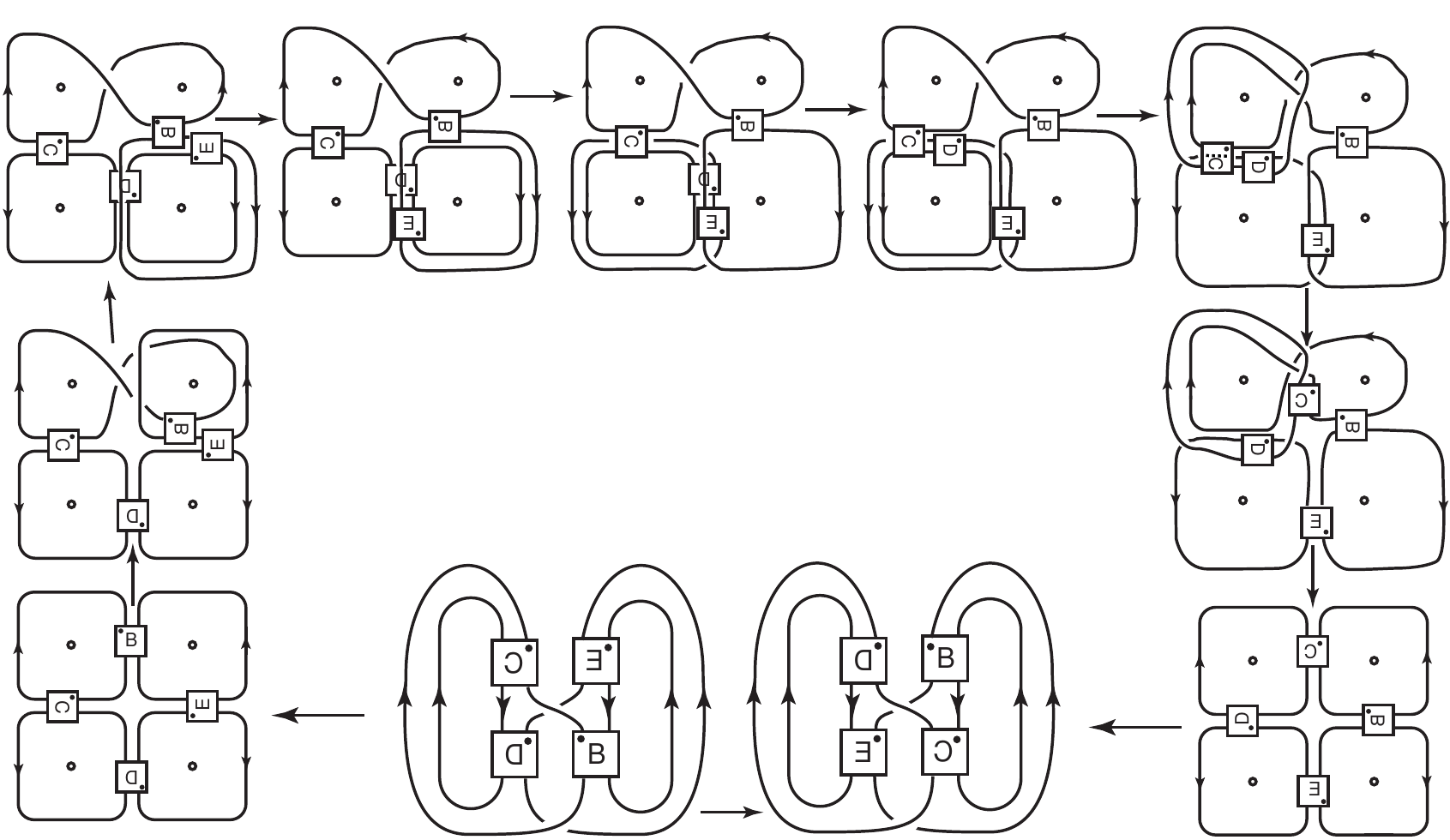}}
\caption{Example of a 4--braid cyclic template, and the Markov tower
(stabilize, exchange, exchange, exchange, destabilize) that it replaces}
\label{figure:template0}
\end{figure}

\fullref{figure:template0} gives an example of the cyclic template.  
We have singled it out because it shows an interesting way in which stabilization introduces flexibility
into the manipulation of closed braids, by allowing us to permute the blocks in a rather special and highly
symmetric block-strand diagram. The associated closed braid diagrams have been unfolded to make it easier to follow
the sequence of moves.  A more general Markov tower for a cyclic template uses weighted strands, the entire tower being equivalent to permuting the blocks and weighted strands in a cycle.  The resulting move on closed braids will be referred to as the
{\it  cyclic move}.

\subsubsection{The G--flype and G--exchange templates} 
\label{subsection:the G-flype and G-exchange templates}
The moves that we next describe are gathered together into the set of templates that we call $\cT(m),
m=4,5,\dots$, where $m$ is the braid index $b(\Xhigh)$ of $\Xhigh$.  

\begin{figure}[htpb]
\labellist\tiny
\pinlabel {\parbox{30pt}{exchange move}} [r] at 223 567
\pinlabel {\parbox{28pt}{exchange move and isotopy}} [bl] at 443 567
\pinlabel {\parbox{25pt}{braid isotopy}} [r] at 142 352
\pinlabel {\parbox{25pt}{braid isotopy}} [l] at 515 364
\pinlabel {\parbox{27pt}{negative flype to a 7--braid}} [r] at 140 220
\pinlabel {\parbox{27pt}{positive flype to a 6--braid}} [l] at 515 220 
\pinlabel {G--flype} [b] at 325 118
\pinlabel {$n=7$} [t] at 331 562
\pinlabel {$n=7$} [t] at 230 382
\pinlabel {$n=7$} [t] at 428 386
\pinlabel {$n=7$} [t] at 224 226
\pinlabel {$n=7$} [t] at 430 226
\pinlabel {$n=6$} [t] at 220 58
\pinlabel {$n=6$} [t] at 432 58

\pinlabel {$W$} at 262 651
\pinlabel {$X$} at 331 581
\pinlabel {$Y$} at 375 698
\pinlabel {$Z$} at 300 685

\pinlabel {$W$} at 160 470
\pinlabel {$X$} at 230 400
\pinlabel {$Y$} at 273 518
\pinlabel {$Z$} at 198 505

\pinlabel {$W$} at 358 475
\pinlabel {$X$} at 427 405
\pinlabel {$Y$} at 470 522
\pinlabel {$Z$} at 395 510

\pinlabel {$W$} at 155 302
\pinlabel {$X$} at 224 232
\pinlabel {$Y$} at 268 349
\pinlabel {$Z$} at 193 338

\pinlabel {$W$} at 362 302
\pinlabel {$X$} at 430 232
\pinlabel {$Y$} at 474 349
\pinlabel {$Z$} at 399 337

\pinlabel {$W$} at 151 140
\pinlabel {$X$} at 220 70
\pinlabel {$Y$} at 264 188
\pinlabel {$Z$} at 189 176

\pinlabel {$W$} at 363 139
\pinlabel {$X$} at 432 69
\pinlabel {$Y$} at 476 187
\pinlabel {$Z$} at 401 173
\endlabellist
\centerline{\includegraphics[scale=.75]{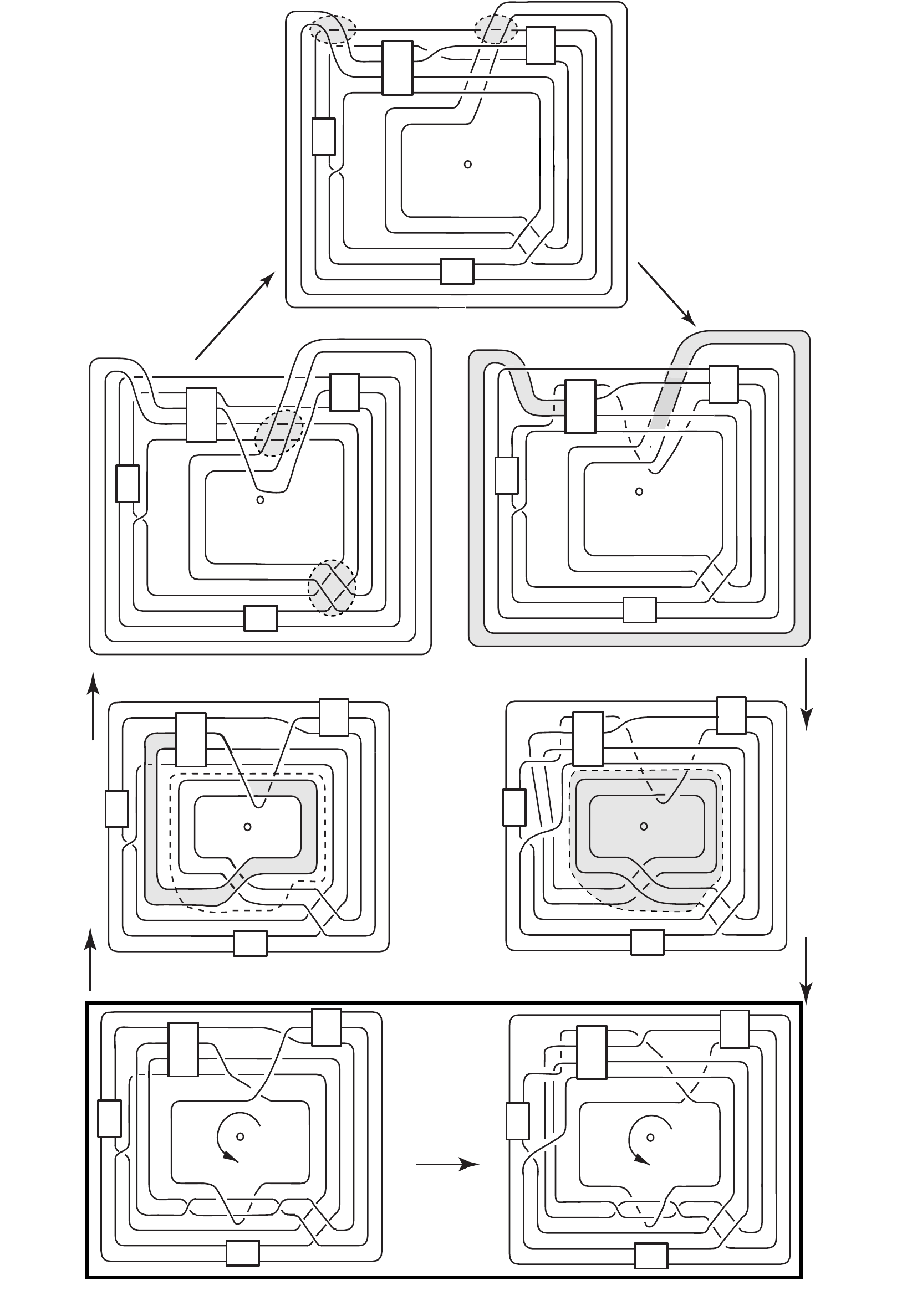}}
\caption{Example of a G--flype on a 6--braid}
\label{figure:template2}
\end{figure}
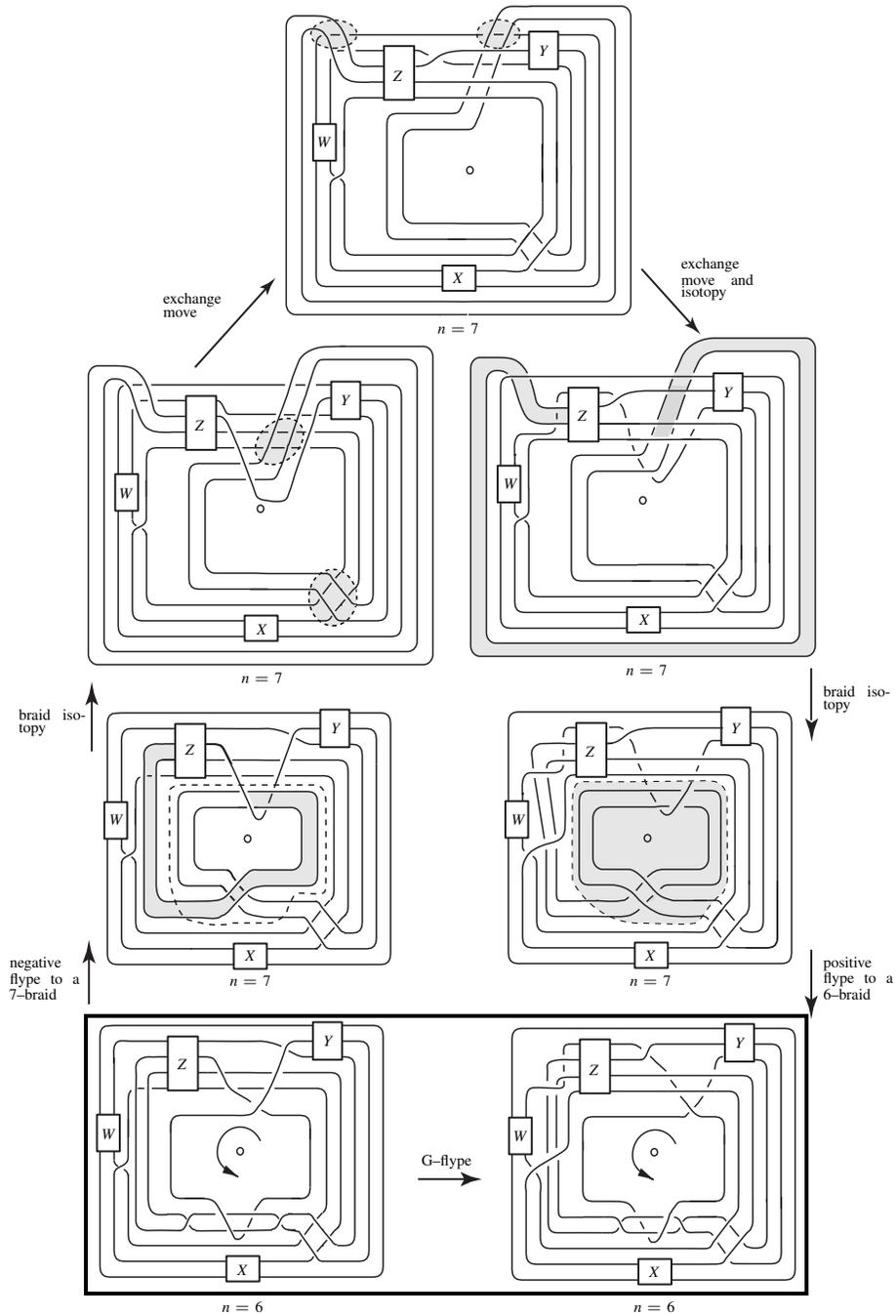

A  {\it  G--flype} template (for `generalized' flype template) is a block-strand diagram, the support
of which is the replacement for a Markov tower in which (i) the first move $X_1\to X_2$ is an inadmissible
flype; (ii) more generally adjacent terms in the sequence differ by destabilizations, exchange moves and not
necessarily admissible flypes; (iii)
$b(X_j) > b(X_1)$ for every $j=2,\dots,k-1$; and (iv) $b(X_k)\leq b(X_1)$.     
 An example  is the boxed pair of 6--braid block-strand diagrams at the
bottom of \fullref{figure:template2}. 
It can be understood by running around the diagram clockwise.  As can be seen,
the first step in that sequence is an inadmissible flype which increases the braid
index by 1.  The intermediate steps are exchange moves and the final step is an admissible flype that
reduces braid index. In more general examples the final step could also be a destabilization. 
See \fullref{subsection:G-flype moves} for a definition which shows
precisely how G--flypes arise and
describes the fixed and moving blocks that always occur.

A {\it  G--exchange move} is the template that results from a sequence of $k$
 interrelated exchange moves, each of which moves a distinguished subarc $\{\alpha_1,\dots,\alpha_k\}$ of
the closed braid across $\axis$ and has the property that the exchange move on $\alpha_i$ cannot be
completed before at least part of the exchange move on some other $\alpha_j$ is started, for $1\leq i\leq
k$. See \fullref{subsection:pushing across regions with a G-exchange foliation} for a definition which shows
precisely how G--exchange moves arise and shows the fixed blocks that always occur.

\begin{figure}[htpb!]
\labellist\tiny
\pinlabel {G--exchange} [b] at 325 100

\pinlabel {$A$} at 131 639
\pinlabel {$B$} at 131 606
\pinlabel {$C$} at 239 621
\pinlabel {$D$} at 287 621
\pinlabel {$E$} at 194 685
\pinlabel {$F$} at 195 555
\pinlabel {$a$} at 167 681
\pinlabel {$b$} at 158 690
\pinlabel {$c$} at 133 666
\pinlabel {$c$} at 130 576
\pinlabel {$d$} at 249 691

\pinlabel {$A$} at 391 633
\pinlabel {$B$} at 391 601
\pinlabel {$C$} at 501 615
\pinlabel {$D$} [l] at 550 615
\pinlabel {$E$} at 453 679
\pinlabel {$F$} at 455 549
\pinlabel {$a$} at 428 675
\pinlabel {$b$} at 417 685
\pinlabel {$c$} at 479 675
\pinlabel {$c$} at 483 555
\pinlabel {$d$} at 527 661
\pinlabel {$d$} at 527 570

\pinlabel {$A$} at 131 421
\pinlabel {$B$} at 131 388
\pinlabel {$C$} at 237 405
\pinlabel {$D$} at 286 405
\pinlabel {$E$} at 194 468
\pinlabel {$F$} at 194 337
\pinlabel {$a$} at 167 463
\pinlabel {$b$} at 158 474
\pinlabel {$c$} at 231 463
\pinlabel {$d$} at 247 474

\pinlabel {$A$} at 391 420
\pinlabel {$B$} at 391 388
\pinlabel {$C$} at 507 403
\pinlabel {$D$} at 552 403
\pinlabel {$E$} at 458 467
\pinlabel {$F$} at 458 336
\pinlabel {$a$} at 433 463
\pinlabel {$b$} at 422 472
\pinlabel {$c$} at 498 469
\pinlabel {$d$} [r] at 366 465
\pinlabel {$a$} at 422 343
\pinlabel {$b$} at 428 332
\pinlabel {$c$} at 494 340
\pinlabel {$d$} [r] at 366 332

\pinlabel {$A$} at 112 197
\pinlabel {$B$} at 100 164
\pinlabel {$C$} at 228 181
\pinlabel {$D$} at 277 181
\pinlabel {$E$} at 185 244
\pinlabel {$F$} at 185 115
\pinlabel {$a$} at 158 241
\pinlabel {$b$} at 139 250
\pinlabel {$a$} at 157 121
\pinlabel {$b$} at 140 110
\pinlabel {$c$} at 247 215
\pinlabel {$d$} at 257 232

\pinlabel {$A$} at 392 197
\pinlabel {$B$} at 381 164
\pinlabel {$C$} at 510 181
\pinlabel {$D$} at 558 181
\pinlabel {$E$} at 464 245
\pinlabel {$F$} at 465 116
\pinlabel {$a$} at 440 241
\pinlabel {$b$} at 420 251
\pinlabel {$c$} at 539 233
\pinlabel {$d$} at 528 218
\pinlabel {$a$} at 440 121
\pinlabel {$b$} at 423 110
\pinlabel {$c$} at 537 127
\pinlabel {$d$} at 528 145

\endlabellist
\centerline{\includegraphics[scale=.65]{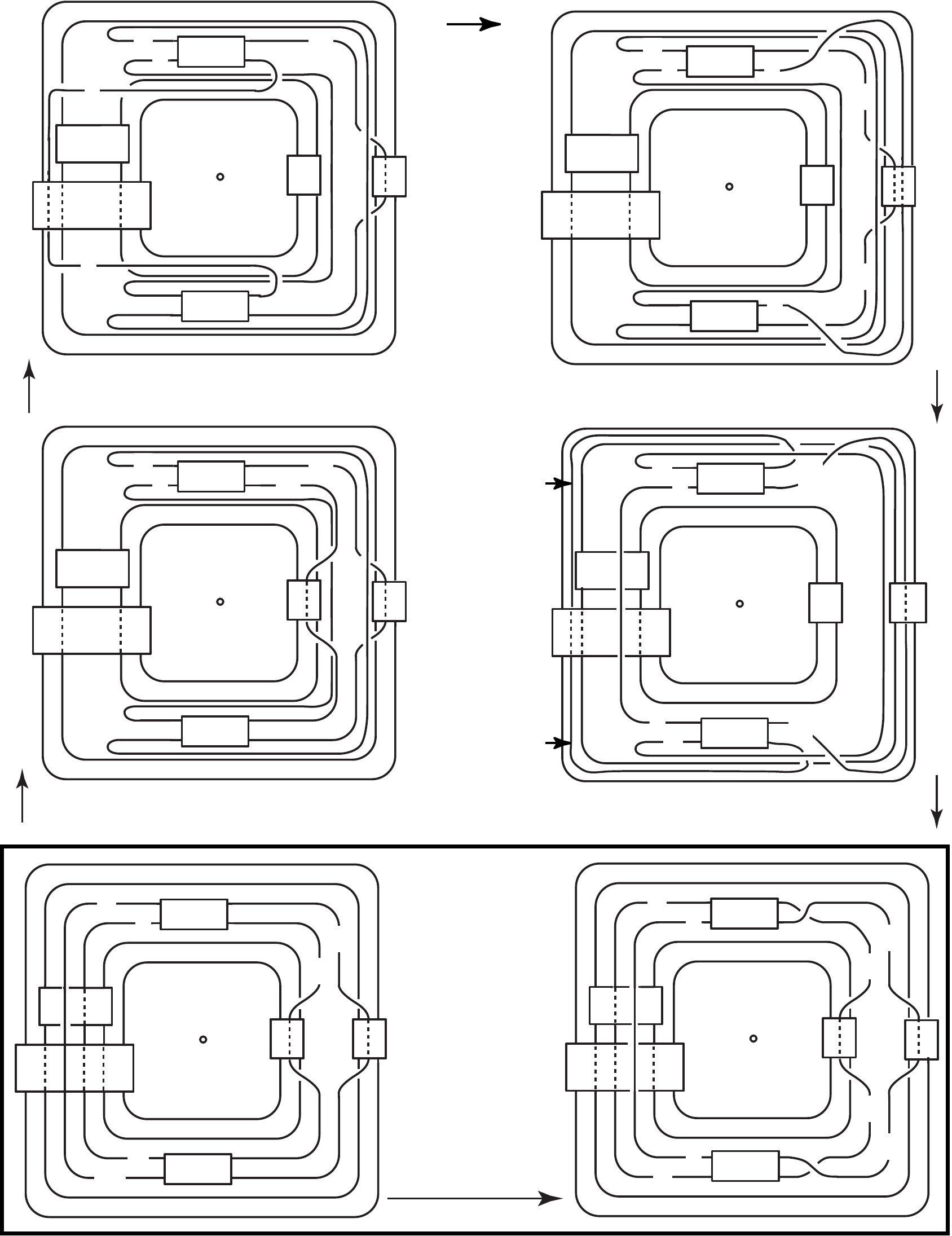}}
\caption{Example of a G--exchange move on a 6--braid}
\label{figure:exchange6}
\end{figure} 
The boxed pair of block-strand diagrams in the bottom row of
\fullref{figure:exchange6} are an example of a template for a G--exchange move on a
6--braid.  There are 6 braid blocks: $A,B,C,D,E,F$. Running around the  figure clockwise we
show how a coordinated sequence of partial exchange moves, each of which can be completed
as soon as enough of the other arcs are moved out of the way, achieves the same goal. 
In the first passage we have pushed strand $a$ under the braid blocks $A$ and
$B$ and across \axis, to a position just to the right of braid block $C$.  We
have also lifted strand $b$ above the braid blocks $A$ and
$B$ and pulled it across \axis \ to a position just to the left of braid block
$D$.  Then we begin our G--exchange move on arc $c$. In the fourth sketch we
complete it. In the fifth sketch we begin the G--exchange move on arc $d$ and complete the
G--exchange move on arc $a$. In the final sketch  we complete the G--exchange moves on arcs $b$
and $d$.  This example was discovered in the course of our proof. It illustrates the ideas developed
in \fullref{subsection:pushing across regions with a G-exchange foliation}.   Observe that, since exchange
moves preserve link type and braid index, it follows that G--exchange moves do too.

\subsection{Statement of results}
\label{subsection:statement of results}
We are finally ready to state our main result, the Markov theorem without stabilization.
Let $\cB$ be the collection of all braid isotopy classes of closed braid representatives of oriented
knot and link types in oriented 3--space. Among these, consider the subcollection
$\cB(\cX)$ of representatives of a fixed link type $\cX$.  Among these,
let $\cB_{\min}(\cX)$ be the
subcollection of representatives whose braid index is equal to the braid index of $\cX$.  

\begin{theorem}[Markov's Theorem Without Stabilization (MTWS)]
\label{theorem:MTWS}
Choose any
$\Xhigh\in \cB(\cX)$  and any $\Xlow\in \cB_{\min}(\cX)$.  Then
 \bi
\item there is  a
complexity function  with values in $\ints_+\times\ints_+\times\ints_+$ which is associated
to $\Xhigh,\Xlow$, and 
\item for each braid index $m$ there is a finite set $\cT(m)$ of templates, each
determining a move which is non-increasing on braid index, \ei
such that  the following hold:

First, there are initial modifications in the choices of $X_+$ and $X_-$, which replace them by $X_+'$ and $X_-'$:
\begin{eqnarray}
\label{sequence:unwinding Xlow}  
\Xlow & = & \Xlow^1 \to \cdots\to \Xlow^p = \Xlow' \\
\label{sequence:exchange-reducibility on Xhigh}
\Xhigh & = & \Xhigh^1\to\dots \to \Xhigh^q = \Xhigh'  
\end{eqnarray}
In \eqref{sequence:unwinding Xlow} and
\eqref{sequence:exchange-reducibility on Xhigh} each passage
$\Xlow^j\to\Xlow^{j+1}$  and $\Xhigh^j\to\Xhigh^{j+1}$ is strictly
complexity-reducing.  In \eqref{sequence:unwinding Xlow} each
passage  is realized by an exchange move, so that $b(\Xlow^{j+1}) =
b(\Xlow^j)$.  In \eqref{sequence:exchange-reducibility on Xhigh}
each passage  is realized by either an exchange move or a destabilization,
so that $b(\Xhigh^j) \geq b(\Xhigh^{j+1})$.

After these initial modifications, there is another  sequence  taking $\Xhigh'$ to $\Xlow'$:
\begin{equation}
\label{sequence:flype sequence on Xhigh'}
\Xhigh'=X^q\to\cdots \to X^r = \Xlow'
\end{equation}
In {\rm \eqref{sequence:flype sequence on Xhigh'}} each passage $X^j\to X^{j+1}$   is 
strictly complexity-reducing, and is realized by either an exchange move,
a destabilization, an admissible flype or a move defined by one of the
templates $\cT\in\cT(m)$, where $m=b(X^j)$.  The inequality $b(X^j)
\geq b(X^{j+1})$ holds for each $j = q,\dots,r-1$.
\end{theorem}

\begin{remark}
We explain why we refer to \fullref{theorem:MTWS} as
{\it Markov's Theorem without Stabilization}.  Setting $X^j = \Xhigh^j$
when $j=1,\dots,q-1$, we have the combination of sequences
\eqref{sequence:exchange-reducibility on Xhigh}
and \eqref{sequence:flype sequence on Xhigh'}:
$$ \Xhigh = X^1\to\dots \to X^q \to\cdots \to X^r = \Xlow'$$
In this combined sequence $b(X^j) \geq b(X^{j+1})$ for each
$j=1,\dots,q,\dots r.$  On the other hand, $b(X_-) = b(X_-^1) = b(X_-^2)
= \cdots = b(X_-^p) = b(X_-^\prime)$, so that the braid index of the
second representative  remains fixed during all of our modifications.
Thus \fullref{theorem:MTWS} is, indeed, a version of \fullref{theorem:MT}
which avoids stabilization. \hfill$\diamondsuit$
\end{remark}

\begin{remark}
The sequence  \eqref{sequence:unwinding Xlow} deals with
the phenomenon which was exhibited in \fullref{figure:exchange3}. It
must be treated separately because if we only allowed modifications to
$\Xhigh$ then the complexity would be forced to increase as $\Xhigh$
approached $\Xlow$, if $\Xlow$ happened to be wound up as on the
left in \fullref{figure:exchange3}. Since
\eqref{sequence:unwinding Xlow} and
\eqref{sequence:exchange-reducibility on Xhigh} are
inter-related, we treat \eqref{sequence:exchange-reducibility on Xhigh}
(which uses a limited subset of the moves in
\eqref{sequence:flype sequence on Xhigh'}) and
\eqref{sequence:flype sequence on Xhigh'} separately.  \hfill$\diamondsuit$
\end{remark}

\begin{remark}
When a passage $X^j\to X^{j+1}$ is realized by a template $\cT$, there
are braiding assignments to the blocks in $\cT$ such that the initial
and final diagrams of $\cT$ carry the pair $(X^j,X^{j+1})$. However,
the template $\cT$ also carries infinitely many other knots and links,
for other braiding assignments to the blocks. \hfill $\diamondsuit$
\end{remark}

\begin{remark} {\rm   The templates in the sets
$\cT(m), m\in\ints^+$, are precisely the additional moves which were not needed for the work in our
earlier papers \cite{B-M_I}--\cite{B-M_VI}, but are needed for our particular
proof of the MTWS.  We discuss them briefly,
starting with $m=2$.} 
\be

\item [{\rm (1)}] The
2--string braid group is an infinite cyclic group. Let $\sigma_1$ denote its generator. An arbitrary
element is then $\sigma_1^k, k\in\ints$. It is easy to see that links
which are closed 2--braids are
either (i) the 2--component unlink ($k=0$),  or (ii) the unknot ($k=\pm 1$)
 \  or (iii) the type
$(2,k)$ torus knots and links ($ |k|\geq 2$). It is clear that the
2--component unlink and the type
$(2,k)$ torus knots and links have unique closed 2--braid representatives. The unknot has exactly 2
closed 2--braid representatives, with $\sigma_1$ (resp.\  $\sigma_1^{-1}$) admitting a positive
(resp.\  negative) destabilization. Since the set $\cT(m)$ does not include the two destabilization
templates, it follows that
$\cT(2) =
\emptyset$. In our paper
\cite{B-M_transverse}, which contains applications of \fullref{theorem:MTWS} to transverse knots, we will prove that, as a consequence of the main theorem in
our paper \cite{B-M_III}, $\cT(3) = \emptyset$.

\item [{\rm (2)}] It was proved by Fiedler
\cite{Fiedler} that  closed 4--braids include infinitely many inequivalent
4--braid
representatives of the unknot. His basic one is the example discovered
by  Morton \cite{Morton1983}, with the others obtained from it by the
winding process which we illustrated in
\fullref{figure:exchange3}.  Fifteen other families of 4--braid unknot
examples were uncovered in the paper by Birman, Boldi, Rampichini and
Vigna \cite{B-B-R-V}, in the course of a computer implementation of the unknot recognition  algorithm of
Birman and Hirsch \cite{B-H}. All of them can be simplified to braids
which admit a destabilization with the help of
exchange moves.  
We do not know whether more general 4--braids are too complicated to be simplified with the use of the
2 destabilization templates, the exchange move template, the cyclic templates and the admissible
flype templates. 

\item [{\rm (3)}] Note that, given
any template $\cT$ of braid index
$m$, other templates for braid index $>m$ may be obtained from it by declaring the strands to be
weighted, and also by replacing some  of the blocks by other templates.  From this it follows that
the cardinality $|\cT(m)|$ of $\cT(m)$ is an increasing function of $m$.  However, we do not have a
precise description of
$\cT(m)$ for any $m>3$, although we do not expect any fundamental difficulty in doing the actual
enumeration for, say, $m=4,5$ and perhaps $6$.  For the special case
$m=6$ two examples were given in the boxed pairs of block-strand diagrams at the
bottom of Figures \ref{figure:template2} and \ref{figure:exchange6}.   The general
picture seems to be quite complicated. \hfill $\diamondsuit$
\ee
\end{remark}

\subsection{Plan of the paper}
\label{subsection:plan of the paper}

In \fullref{section:getting started} we
set up the topological construction which will be the basis for our work. 
We will show that there is a very special isotopy that takes us from $\Xhigh$ to
$\Xlow$.  The trace of the isotopy will be seen to be an immersed annulus $\ta$ whose double
point set is the union of finitely many pairwise disjoint clasp arcs. We call it a `clasp annulus'. 
 
The principle tool in our proof of \fullref{theorem:MTWS} is the study of certain `braid
foliations' of the immersed annulus $\ta$ and its preimage $\pa$. Braid foliations were used by the
authors in earlier work \cite{B-M_I}--\cite{B-M_VI}, but always in the setting of embedded surfaces. In
\fullref{section:preliminaries} we review the ideas that we need from the
literature on braid foliations.  Readers who are familiar with the literature will probably
want to pass quickly over \fullref{section:preliminaries}, referring to it
instead, as needed, later in the paper.
In \fullref{section:braid foliations of the immersed annulus} we study braid foliations of our
immersed annulus.  We will need to do hard technical work to arrange things so that the 
clasp arcs are close to or contained in a union of
leaves and have nice neighborhoods  (we call them normal neighborhoods) on the preimage
annulus $\pa$.   

In \fullref{section:pushing across ta} we learn how to
translate data in the braid foliation of $\ta$ and the induced foliation of $\pa$ into
data about the passage from the closed braid $\Xhigh$ to the closed braid $\Xlow$.  The tools
that are needed become increasingly complicated as we proceed. First, we ask how far we can get with exchange
moves and destabilizations.  Flypes enter the picture next, but in the form of very rudimentary examples which we
call `microflypes'.  A rather surprising use of stabilization becomes apparent in \fullref{subsection:pushing across a thin annulus}. Briefly, we learn that stabilization is the tool for
creating flypes with weighted strands and complicated braiding assignments in the moving blocks out of
microflypes. G--exchange come into play next.  There are hints in this part of the work about the need for
G--flypes, however the reasons for needing G--flype templates will not become clear until we are part way through the
proof of the MTWS. 

The proof of \fullref{theorem:MTWS} is given in \fullref{section:the
proof}. We will see how G--flype
templates arise.  The most difficult part of the argument will be the proof that for each fixed braid
index $m = b(\Xhigh)$ the cardinality $|\cT(m)|$ is bounded. The finiteness can, perhaps, be understood by
appreciating that the `infinite parts' are pushed into the blocks in the block-strand diagrams of the templates
in $\cT(m)$. This is, perhaps, the key point about block-strand diagrams and templates: they are at the same time
both very flexible and very inflexible. A given template supports a huge family of knots and links, because there
are no restrictions on the braiding assignments in the blocks, but on the other hand a template always supports at
most a special family of links.

The paper ends, in \fullref{section:open problems} with a discussion of open problems suggested by the MTWS and by its proof.  In a
separate paper \cite{B-M_transverse} several applications will be given to the study of transversal knot types in
the standard contact structure on $\reals^3$.   

\subsubsection*{Conventions}
Results which will be used explicitly in the proof of the MTWS are highlighted by calling them `propositions' rather than `lemmas'.    There are many remarks, scattered throughout the paper. Most of them are guides to the reader.  The end of a remark is indicated by the symbol $\diamondsuit$.

\subsubsection*{Acknowledgments}
The authors single out for special thanks John Etnyre, Keiko Kawamuro and
the unknown G\&T referee. All three studied our manuscript with care and
attention to details in the proof,  and asked the kind of questions that
can only be asked by an expert,  leading us to clarify many tricky points.

This work has been in progress for so many years that it's difficult for us to even remember everyone else who helped. With apologies for oversights we also thank Tara Brendle, Oliver Dasbach, Ivan Dynnikov, Brian Mangum, Walter Neumann,  Nancy Wrinkle,  Joel Zablow, 
Xingru Zhang,  Joseph Birman and Melissa Menasco.   

The first author acknowledges partial support from the US National
Science Foundation, under grants DMS-9402988, 9705019, 9973232 and
0405586. She also wishes to thank the Mathematics Department at the
Technion (Israel Institute of Technology) for hospitality during several
visits when this paper was the main focus of her work.  The second
author acknowledges partial support from the US National Science
Foundations, under grants DMS-9200881, 9626884 and 0306062; and the
Mathematical Sciences Research Institute, where he was a Visiting Member
during winter/spring of 1997.

\section{Getting started} 
\label{section:getting started}

In this section we develop the basic construction which will allow us to prove
\fullref{theorem:MTWS}.  For ease in
presentation, we give our construction first for the special case when $\cX$ is a knot. After
that it will be easy to modify it in the more general case when we begin with a link. 
The section will end with a key example. 

\subsection{The basic construction for knots}
\label{subsection:basic construction for knots}

\begin{figure}[htpb]
\labellist\small
\pinlabel {positive clasp} at 180 564
\pinlabel {negative clasp} at 421 564
\pinlabel {\parbox{45pt}{clasp\\intersection}} [b] at 300 651
\endlabellist
\centerline{\includegraphics[scale=.7]{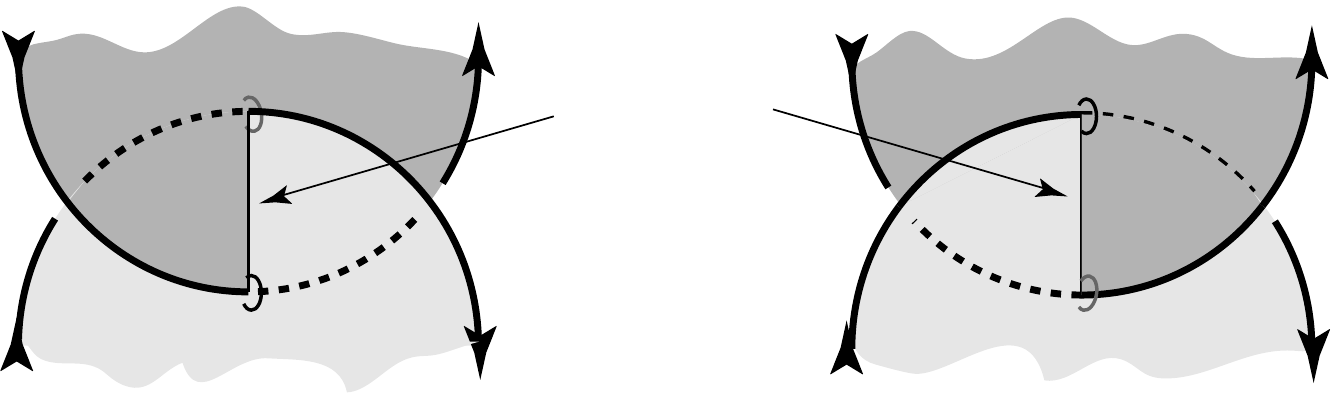}}
\caption{Clasp intersections}
\label{figure:getstarted1}
\end{figure}

\begin{lemma}
\label{lemma:construction of the two annuli} 
Let
$\Xhigh, \Xlow$ be arbitrary disjoint closed braid representatives of
the same knot type $\cX$.   Then
there is an intermediate representative $\Xmed$ of
$\cX$ such that the following hold:
\begin{enumerate}
\item[{\rm (1)}] $\Xmed$ is the braid connected sum of $\Xhigh$ and $k$ closed braid
representatives of the  unknot, for some $k\geq 0$. These $k$ representatives of the unknot bound pairwise
disjoint discs.
\item[{\rm (2)}] $\Xmed-\Xlow$ (resp.\  $\Xhigh-\Xmed$) is the
boundary of an embedded annulus $\annlow$ (resp.\  $\annhigh$).  It will be seen from the construction that $X_0$ is the boundary of a Seifert surface $\seif_0$, and $\annlow$ is a collar neighborhood of $\Xmed$ on $\seif_0$.  
\item[{\rm (3)}] The intersections $\annhigh\cap \annlow$ are precisely $k$ clasp arcs. See \fullref{figure:getstarted1}.

\end{enumerate}
\end{lemma}

\begin{remark}
A word is in order on the basic construction. There
are very simple constructions which yield everything in
\fullref{lemma:construction of the two annuli} except the fact that
$\annlow \subset \seifmed$.  Astute readers will notice that $\seifmed$
plays an almost invisible role in the pages that lie ahead, and ask why we
needed it?   Our reasoning was that in our earlier papers we had developed
extensive machinery regarding braid foliations of Seifert surfaces for
knots and links, including the case of discs bounded by unknots. In this
paper we need related braid foliations of the immersed annulus of
\fullref{lemma:construction of the two annuli}.  The easiest way to
develop the tools that we need seemed to be to make one of our annuli
(ie $\annlow$) a subset of a Seifert surface, and arrange that the other
(ie $\annhigh$) is a finite family of foliated discs banded together
with narrow  strips to form a second annulus that intersects the first in
a controlled fashion.  That will allow us to make full use of the existing
braid foliation machinery for Seifert surfaces and discs, rather than to
develop everything anew.  This issue will become clear in
\fullref{section:preliminaries} and \fullref{section:braid foliations of
the immersed annulus} below.  \hfill $\diamondsuit$ 
\end{remark}

\begin{proof}   Without loss of generality we may assume that $\Xhigh\subset\reals^3_+$ and $\Xlow\subset\reals^3_-$, with $\Xhigh$ far above $\reals^2=\reals^3_-\cap\reals^3_+$, and $\Xlow$ a little bit below $\reals^2$.  We may further assume that $\Xhigh,\Xlow$ are closed braids with respect to the same braid structure $(\axis, \fib)$.  Our first task is to construct a series of knots
$\Xlow^\prime,
\Xlow^{\prime\prime}, \Xmed$, all representing $\cX$, with $\Xmed$ the braid-connected
sum of $\Xhigh$ and $k$ pairwise disjoint and pairwise
unlinked copies $U_1,\dots,U_k$ of the unknot.  
 
Choose a Seifert surface $\seiflow$  for $\Xlow$.  Let $\Xlow^\prime
\subset \seiflow$ be a preferred longitude for $\Xlow$, chosen close enough to $\Xlow$, so that the annulus that
they cobound in $\seiflow$ does not intersect $\axis$, and so that $\Xlow'$ is also a closed
braid.  The knots
$\Xlow$ and $\Xlow'$ will have algebraic linking number 0, but $\Xlow\cup \Xlow'$ will not be a split link unless
$\cX$ is the unknot.  Therefore, if we try to push $\Xlow'$ out of $\reals^3_-$ and into $\reals^3_+$ it will get
stuck, ie there will be a finite number of `undercrossing hooks'  where $\Xlow'$  is forced to dip back
into $\reals^3_-$ to pass under $\Xlow$, as in \fullref{figure:getstarted2}(a).  

Our first change is to modify
$\Xlow'$ (holding $\Xlow$ fixed) to a new closed braid $\Xlow^\pp$ which has the same knot
type as
$\Xlow'$ and is entirely in $\reals^3_+$. This can be accomplished by pushing $\Xlow'$
across  pairwise disjoint discs $D_1\cup \dots \cup D_k$, as in \fullref{figure:getstarted2}(b) to $\Xlow^\pp$. By placing $\Xlow$ very close to $\reals^2$ and choosing the discs to
be very `thin' we may assume that each subarc
$\beta_i$ is in braid position and in $\reals^3_+$.   By construction $\Xlow^\pp$ is a  closed braid (because $\Xlow'$ is a closed braid and the isotopy $\Xlow' \to \Xlow^\pp$ is arbitrarily small), also it represents $\cX$, and it is entirely in $\reals^3_+$.

\begin{figure}[htpb]
\labellist\small
\pinlabel {(a)} at 255 523
\pinlabel {(b)} at 382 523
\pinlabel {(c)} at 538 523
\pinlabel {(d)} at 212 365
\pinlabel {(e)} at 370 365
\pinlabel {(f)} at 537 365

\pinlabel {$\mathbb{R}^2$} [br] at 106 581
\pinlabel {$\mathbb{R}^2$} [l] at 371 576
\pinlabel {$\mathbb{R}^2$} [br] at 500 546
\pinlabel {$\mathbb{R}^2$} [br] at 178 391
\pinlabel {$R_i$} at 152 448

\tiny
\pinlabel {$X_-'$} [l] at 231 613
\pinlabel {$X_-$} [l] at 245 580
\pinlabel {$\alpha_i$} [l] at 167 533
\pinlabel {$\beta_i$} [b] at 154 614
\pinlabel {$D_{i,+}$} at 154 596
\pinlabel {$D_{\!i\!,\!-}$} at 154 558

\pinlabel {$X_-''$} [l] at 376 612
\pinlabel {$\beta_i$} at 313 619
\pinlabel {$D_{i,+}$} at 314 594
\pinlabel {$D_{i,-}$} at 314 561

\pinlabel {$\beta_i'$} [b] at 462 649
\pinlabel {$X_+$} [l] at 513 651
\pinlabel {\rotatebox{60}{$g(D_{i,+})$}} at 476 606

\pinlabel {$\beta_i'$} [b] at 147 491
\pinlabel {$X_+$} [l] at 194 494
\pinlabel {$D_{i,-}$} [tl] at 149 372

\pinlabel {$X_+$} [l] at 364 494
\pinlabel {$X_0$} [l] at 348 485
\pinlabel {$X_-$} [l] at 353 407

\pinlabel {$X_0$} [l] at 519 488
\pinlabel {$X_-$} [l] at 522 409

\endlabellist
\centerline{\includegraphics[scale=.65]{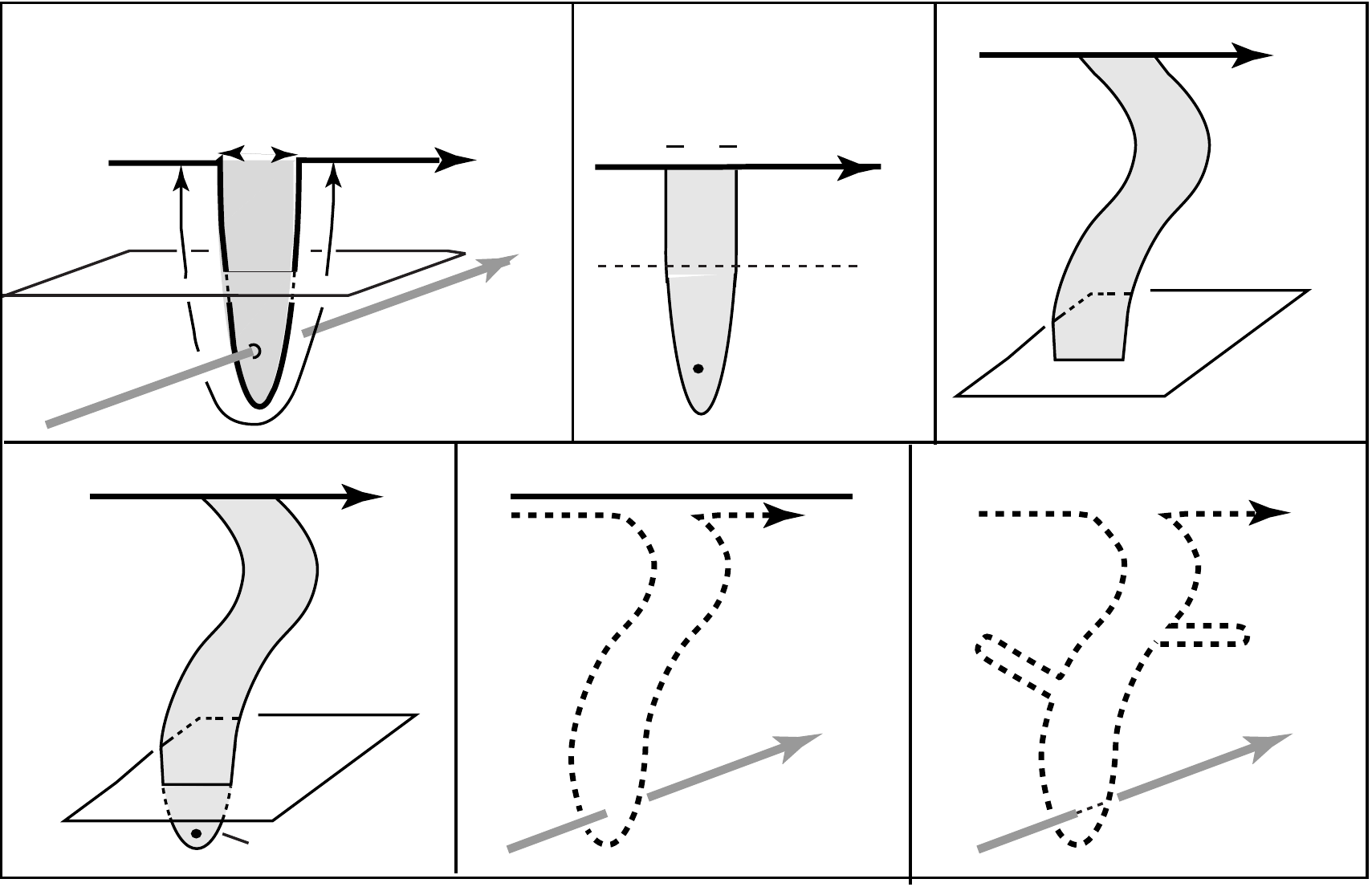}}
\caption{Constructing $\Xmed$}
\label{figure:getstarted2}
\end{figure}

We are now ready to bring $\Xhigh$ into the picture.
The fact that $\Xlow^\pp$ and $\Xhigh$ both represent $\cX$ and are both
in the interior of $\reals^3_+$ shows that we may find a homeomorphism
$g\co \reals^3_+\to\reals^3_+$ which is the identity on $\reals^2$
with $g(\Xlow^\pp)
= \Xhigh$. Extend $g$ by the identity on $\reals^3_-$ to a homeomorphism
$G\co\reals^3\to\reals^3$.  Let $R_i =  G(D_{i,+}\cup D_{i,-}) = g(D_{i,+})\cup D_{i,-}.$ The facts that
(i)
$G$ is a homeomorphism which is the identity in $\reals^3_-$  and (ii) if
$i\not=j$ then $ D_i\cap D_j = \emptyset$ tell us that the
$R_i's$ are pairwise disjoint embedded discs, and also that $\Xlow$ pierces each $R_i$ exactly once. 
The fact that $\Xhigh$ was well above $\reals^2$ shows that we may assume that each $R_i$ intersects
$\Xhigh$ in a single arc $\beta_i' = g(\beta_i)$.

Let $\Xmed'$ be the
knot which is obtained from $\Xhigh$ by replacing each $\beta_i'\subset \Xhigh$ by
$\partial R_i \setminus \beta_i'$.  Then $\Xmed'$ is the connected sum of $\Xhigh$ and $k$ copies of the
unknot,  the $i^{\rm th}$ copy being $\partial R_i$.  By construction the $k$ discs $R_1,\dots,R_k$
are pairwise disjoint, so that our $k$ unknots represent the $k$--component unlink. 
It may happen that $\Xmed'$ is not a closed braid.   The only subarcs which might not be in braid position are the
$2k$ `vertical' arcs in each $\partial g(D_{i,+})$. To overcome this
problem, recall that J\,W Alexander
\cite{Alexander} 
introduced a very simple way to change an arc $\delta$
which is not in braid position to one that is.  See \fullref{figure:getstarted3}(a). Alexander showed that whenever a knot
$X$ is not transverse to the fibers of $\fib$ it can be divided into small segments which can
then be pushed across the braid axis, one at a time, avoiding unwanted intersections with the rest of $X$, to
change $X$ to a closed braid.  
Using this construction,  we 
change all the wrongly oriented subarcs of the $2k$ vertical arcs in the boundaries of
the `ribbons' $R_1\cup\dots\cup R_k$ to subarcs that are everywhere transverse to the fibers
of $\fib$.   We have proved
(1). 
\begin{figure}[htpb]
\labellist\small
\pinlabel {$\delta$} [l] at 177 559
\pinlabel {$\delta'$} [l] <-1.15pt,3pt> at 274 559
\pinlabel {$X_0$} [l] at 141 495
\pinlabel {(a)} at 174 460
\pinlabel {clasps} [b] at 406 584
\pinlabel {$X_0$} [b] at 474 553
\pinlabel {$X_-$} [l] at 552 606
\pinlabel {$X_+$} [l] at 558 508
\pinlabel {(b)} at 442 460
\pinlabel {$\axis$} [b] at 256 651
\endlabellist
\centerline{\includegraphics[scale=.75]{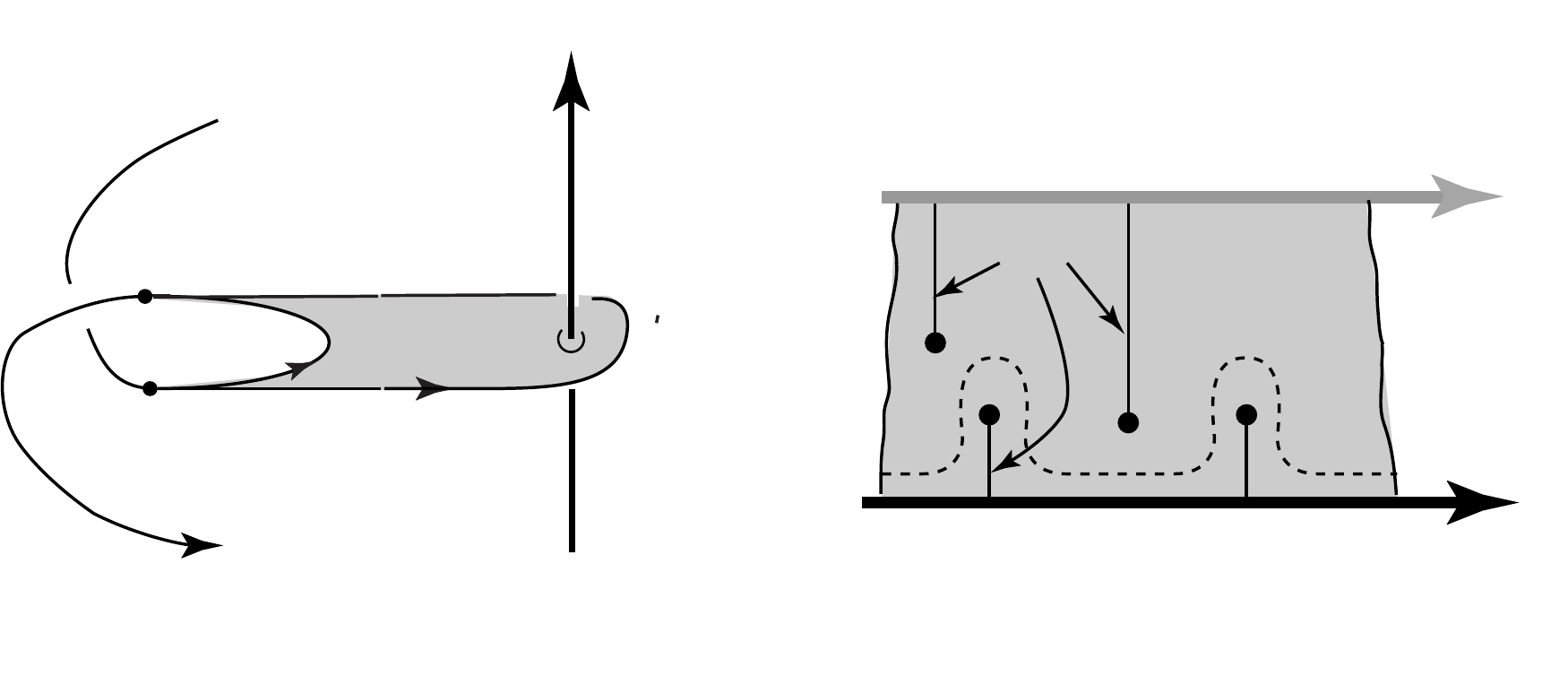}}
\caption{(a) Alexander's trick\quad (b) A fragment of the preimage $\pa$ of $\annhigh\cup\annlow$}
\label{figure:getstarted3}
\end{figure}

There is an important aspect of our construction, which will give us the part of (2) that relates to $\annlow$: 
\bi
\item $\Xlow\cup \Xmed$ has the same link type as $\Xlow\cup
\Xlow'$.  For, by construction, the homeomorphism $G^{-1}\co\reals^3\to\reals^3$, being the
identity on $\reals^3_-$, sends
$\Xlow\cup \Xmed$ to $\Xlow\cup \Xlow'$. 
\ei
This simple fact gives us the annulus $\annlow$, in the
following way: Since $\Xlow'$ is a preferred longitude for
$\Xlow$, and since $\Xlow\cup \Xmed$ has the same link type as $\Xlow\cup \Xlow'$, it 
follows that  $\Xmed$ is also a preferred longitude for $\Xlow$.  From this it follows that 
$\Xlow$ is also a preferred longitude for $\Xmed$.
Choose a Seifert surface $\seifmed$ for $\Xmed$.  Holding $\Xmed$ and $\Xlow$
fixed, isotope the interior of $\seifmed$ until $\Xlow$  lies on $\seifmed$ as a preferred longitude.
Let
$\annlow\subset\seifmed$ be the annulus in $\seifmed$ which $\Xmed$ and $\Xlow$ bound. 
This annulus is embedded because $\seifmed$ is embedded.  Thus we have proved the part of (2) that relates to $\annlow$.

In fact, a small modification in $\Xmed$ also gives us $\annhigh$. The discs $R_1
\cup \cdots \cup R_k$ have a natural order which is determined by the order of the
subarcs $\beta_1'\cup\cdots\cup\beta_k'$ along $\Xhigh$. Using this order, and the
framing provided by $\seifmed$, join the discs $R_1\cup \cdots \cup R_k$ by $k$ very narrow bands in
$\seifmed$, each having one edge on $\Xhigh$. Modify $\Xmed$ by pushing it a little bit
into $\seifmed$ along the bands and along the $R_i$'s. (By an abuse of notation, we use the same
names for the modified $R_i'$s and the modified $\Xmed$). The union of the new $R_i'$s and the
bands is our annulus $\annhigh$. The annulus $\annhigh$ is embedded because the $R_i'$s are
disjoint and embedded, and the bands are too. We have proved (2).  Since $\Xlow\subset\partial\annlow$ and since 
$\Xlow$ pierces each $R_i\subset\annhigh$ once, it follows that $\annlow\cap\annhigh$ always
contains $k$ clasp arc intersections.  Thus we have also proved (3), and the proof of \fullref{lemma:construction of the two annuli} is complete. \end{proof}

We establish conventions that will be used throughout this paper. The symbol $\ta$ denotes the immersed annulus $\annhigh\cup\annlow$. We shall refer to it as a {\it clasp annulus}. We will also be interested in its preimage $\pa$ under the immersion.  \fullref{figure:getstarted3}(b) is a schematic that illustrates our  basic construction,  as described in \fullref{lemma:construction of the two annuli}.  The knot  $\Xmed$ will only be needed in the beginning of our argument, therefore we show it as a dotted curve. The closed braid
$\Xhigh$ will be the primary focus of our attention, and so we show it as a thick black line. We will
eventually modify it to $\Xlow$, which we illustrate as a thick grey line (to suggest that it is
a distant goal). To avoid clutter in our figures we will, whenever the meaning is unambiguous,
suppress the labels $\Xhigh$ and $\Xlow$. Most of the time the black-grey convention will
enable us to recognize them without labels.

\begin{remark}
The reader may wonder why the
isotopy which we constructed from
$\Xhigh$ to $\Xlow$ in \fullref{lemma:construction of the two annuli}  required us to increase the
braid index of $\Xhigh$ by connect-summing it with some number of closed braid representatives of the
unknot, because one expects these to have braid index greater than 1. This does indeed sound wrong, in view of the
fact that our goal is to have the braid index go down and not up, but in fact it is to be expected
because we are attempting to prove the MTWS, and for that we first need to have in hand a proof of the
MT.  Indeed, in
\cite{B-M_Markov} we showed that a variation on the construction given here in
\fullref{lemma:construction of the two annuli} and \fullref{proposition:get started for links} (but without taking any note of the clasp arcs) can be used to
give a new proof of Markov's classical result.  When we began the work in this paper we started with
that new proof of the Markov theorem and proceeded to modify the Markov tower that it gave to us.
Later, we realized that it was not necessary to literally prove the Markov theorem, all we needed was
the isotopy encoded by $\ta$. We hope that explains the logic.
\hfill$\diamondsuit$
\end{remark}

\subsection{The general case}
\label{subsection:the basic construction, general case}

\begin{proposition}
\label{proposition:get started for links}
Choose any $\mu$--component oriented link type $\cX$ in oriented $\reals^3$. Let
$\Xhigh,\Xlow$ be closed braid representatives sharing a braid structure $(\axis, \fib)$.
Then an intermediate closed braid representative
$\Xmed\in\cX$ exists such that the following hold:
\begin{enumerate}

\item[(1)] Each component $\Xmed^i$ of $\Xmed$ is the braid	connected sum of $\Xhigh^i$ and $k^i$ 
closed braid representatives of the unknot.  These $k = k^1+\cdots + k^\mu$ representatives of the
unknot bound $k$ pairwise disjoint discs.

\item[(2)] There is a pairing of the components of $\Xhigh$ and $\Xlow$ such that for each pair
$(\Xhigh^j,\Xlow^j)$ there exists a closed braid $\Xmed^j$ which represents the same component of $\cX$. Moreover,
 for each $j=1,\dots,\mu$ the following holds: $\Xhigh^j-\Xmed^j$
(resp.\  $\Xmed^j-\Xlow^j$) is the boundary an
embedded annulus $\annhigh^j$ (resp.\  $\annlow^j$).  The union of these $\mu$ embedded annuli forms an
embedded surface $\annhigh$ (resp.\  $\annlow$). 

\item[(3)] The intersections $\annhigh\cap\annlow$ are a finite collection
of pairwise disjoint clasp arcs. These intersections can be clasp intersections between the annuli
$\annhigh^i$ and $\annlow^i$ associated to a single component $\cX^i$ of $\cX$ and also clasp
intersections between the annuli $\annhigh^i$ and $\annlow^j$ associated to distinct components
$\cX^i, \cX^j$ of $\cX$. 

\item[(4)] Moreover, without loss of generality we may assume that $\cX$ is a non-split link.  If, instead, $\cX$ is a
split link, then the basic construction may be applied to each non-split component.

\end{enumerate}
\end{proposition}

\begin{proof}
We leave it to the reader to check that, except for (4),
the proof is a minor adaptation of the argument that we gave in the proof
of \fullref{lemma:construction of the two annuli}. As for (4), we note
that in the manuscript \cite{B-M_IV} the authors used tools that are
closely related to the tools that will be used in this paper to prove
that if $X$ is an arbitrary closed $n$--braid representative of a split
(resp.\  composite) link, then  there is a complexity function $C(X)$ with
values in $\ints_+$ which is associated to $X$ such that after a strictly
complexity-reducing sequence of exchange moves, all of which preserve
braid index, $X$ can be changed to a split (resp.\  composite) closed
$n$--braid representative of $\cX$.  Thus, for split (resp.\  composite)
links, the proof of the MTWS can be preceded by applying the results in
\cite{B-M_IV} to find the non-split (resp.\  prime) summands.  However,
we note that non-prime knots and links do not require any special
consideration in our work. Therefore there would not be any point in
making the assumption that $\cX$ is prime. On the other hand, there are
two places where the assumption that $\cX$ is not split simplifies our
work a little bit, so in what follows we make the assumption (4) that
$\cX$ is non-split.
\end{proof}

\subsection{A key example of the basic construction}
\label{subsection:a key example of the basic construction}

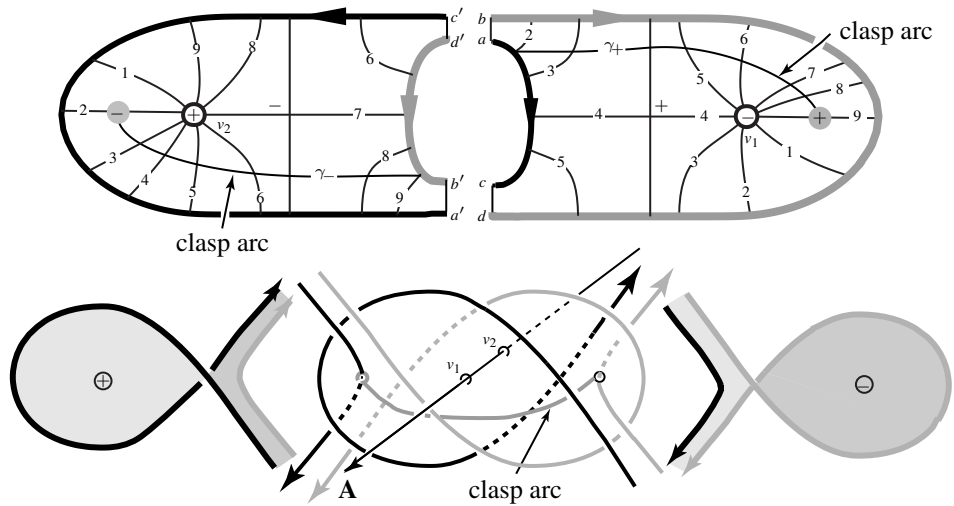
\begin{figure}[htpb]
\labellist\small
\pinlabel {clasp arc} [l] at 484 507
\pinlabel {clasp arc} [t] at 178 410
\pinlabel {clasp arc} [t] at 325 286
\pinlabel {$\axis$} [t] at 241 286
\tiny
\pinlabel {$\gamma_{\!-}$} at 230 435
\pinlabel {$\gamma_{\!+}$} at 375 499
\pinlabel {$-$} [br] at 211 464
\pinlabel {$+$} [bl] at 392 464
\pinlabel {$+$} at 163.5 465.5
\pinlabel {$-$} at 442 464
\pinlabel {$+$} at 478 464
\pinlabel {$-$} at 125 466
\pinlabel {$+$} at 118 332
\pinlabel {$-$} at 500 329

\pinlabel {$1$} at 129 487
\pinlabel {$2$} at 108 468
\pinlabel {$3$} at 123 444
\pinlabel {$4$} at 139 433
\pinlabel {$5$} at 163 427
\pinlabel {$6$} at 197 424
\pinlabel {$8$} at 193 500
\pinlabel {$9$} at 165 499

\pinlabel {$2$} at 332 506
\pinlabel {$3$} at 343 488
\pinlabel {$4$} at 367 466
\pinlabel {$5$} at 348 442
\pinlabel {$6$} at 252 494
\pinlabel {$7$} at 246 466
\pinlabel {$8$} at 260 445
\pinlabel {$9$} at 267 424

\pinlabel {$1$} at 463 447
\pinlabel {$2$} at 441 427
\pinlabel {$3$} at 415 442
\pinlabel {$4$} at 420 465
\pinlabel {$5$} at 418 484
\pinlabel {$6$} at 440 504
\pinlabel {$7$} at 474 487
\pinlabel {$8$} at 488 479
\pinlabel {$9$} at 497 465

\pinlabel {$v_1$} [t] at 444 459
\pinlabel {$v_2$} [t] at 179 466

\hair=2pt
\pinlabel {$a$} [r] at 313 503
\pinlabel {$b$} [r] at 313 514
\pinlabel {$c$} [r] at 313 431
\pinlabel {$d$} [r] at 313 415

\pinlabel {$a'$} [l] at 290 416
\pinlabel {$b'$} [l] at 290 432
\pinlabel {$c'$} [l] at 290 516
\pinlabel {$d'$} [l] at 290 503

\pinlabel {$v_1$} [br] at 300 334
\pinlabel {$v_2$} [br] at 318 347
\endlabellist
\centerline{\includegraphics[scale=.75]{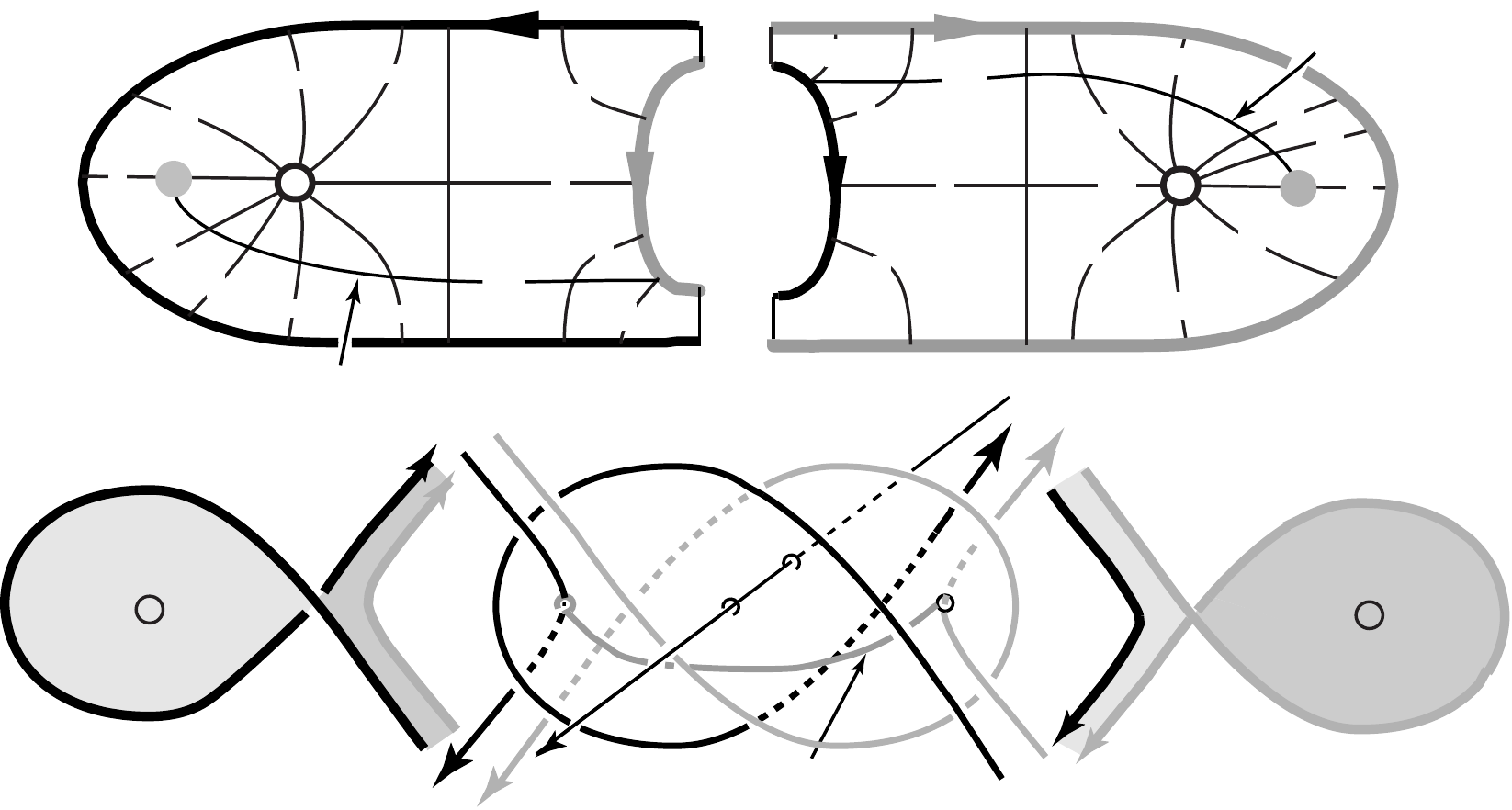}}
\caption{A key example}
\label{figure:getstarted4}
\end{figure}

The reader is referred to \fullref{figure:getstarted4}.
It is a key example, and we will study it in full detail during the course of this manuscript.
We explain those features which can be understood at this time. 

In our example $\cX$ is a knot, and the isotopy from $\Xhigh$ to $\Xlow$ is realized by a flype. The top sketch
shows
$\pa$ as the union of two discs which are identified along bands which join the arcs $ab$ and $a'b'$, and
also $cd$ and $c'd'$.  \fullref{figure:getstarted4} depicts the disc neighborhood of one of its clasp arcs.
The black (resp.\  grey) boundary is
$\Xhigh$ (resp.\ 
$\Xlow$). To avoid clutter we do not show $\Xmed$. There is one clasp arc whose two preimages
$\glow$ and $\ghigh$ are shown. The clasp arc $\ge,$ where $\epsilon =
\pm$, has one endpoint on $X_\epsilon$ and the other in the interior of the annulus, at the
point where $X_{-\epsilon}$ pierces the annulus. 
The bottom sketch illustrates the clasp annulus $\ta$. To visualize $\ta$
in 3--space, first give each of the discs 
in the top sketch a half-twist, as in the leftmost and rightmost bottom sketches. Then
identify the two half-twisted discs along the clasp arc, as in the middle bottom
sketch. The passage from $\Xhigh$ (the black boundary) to $\Xlow$ (the grey boundary) is realized
by a push  of $\Xhigh$ across the immersed annulus $\ta$ to $\Xlow$.  While the annulus is immersed,
there are no self-intersections of the boundary braid as it is pushed across the annulus. 

In fact, our example illustrates a (negative) 3--braid flype, as in the passage
from the 3rd to the 4th sketch in \fullref{figure:flype1}, in the
special case when   
$$\Xhigh = \sigma_1^p\sigma_2^{-2}\sigma_1^q\sigma_2^{-1}, \ \  
\Xlow = \sigma_1^p\sigma_2^{-1}\sigma_1^q\sigma_2^{-2}$$
with $p$ and $q$ integers, also $p+q$ odd (so that we get a knot) and
absolute value at least 3 (so that, by the work in our paper \cite{B-M_III}, we know that the flype
cannot be replaced by braid isotopy).  Notice that in the passage from the left to the
right sketch in \fullref{figure:flype1} the isotopy is supported in a 3--ball $B^3$ whose
boundary intersects $\Xhigh$ in 4 points. These are the endpoints of the 2 subarcs of
$\Xhigh$ in the bottom middle sketch in \fullref{figure:getstarted4}. The
3--ball
$B^3$ contains in its interior the braid box $R$ (which is our case is a single negative full
twist), a little subarc of the braid axis and the single crossing to its right (which in our
case is also negative). The signs of the singular points were chosen to correspond to the fact that the exponents of $\sigma_2$ in
the braids which represent $\Xhigh$ and $\Xlow$ are both negative.   The braiding in those blocks depends, of course, on the choices of the exponents $p$ and
$q$.   They are examples of the `fixed blocks' which are formed from the
`braiding of bands of $s$--arcs'.  (See the proof of \fullref{proposition:sequence of flypes}). 

It should be clear to the reader that (except in very special cases) the braiding of the long thin
bands will lead to geometric linking between $\Xhigh$ and $\Xlow$, and this implies that there will be
additional `short clasp arcs' in the bands. That matter will be discussed in \fullref{subsection:short clasp arcs}. 

We remark that the puncture endpoints of the two clasp arcs in \fullref{figure:getstarted4} are illustrated as being signed, but their signs have not been defined or discussed. We will need them later, so we define them now. The puncture  point on
a clasp arc is said to be a {\it positive} puncture point (resp.\ 
negative puncture point) according as the orientation of $X_\epsilon$ at
the puncture agrees (resp.\  disagrees) with
the orientation of the outward drawn normal, when the surface has the orientation induced by the orientation on $\Xhigh$.  We will see later that the sign on one of the punctures determines the sign on the other.   In \fullref{figure:getstarted4} we illustrated the
case when the left puncture is negative and the right puncture is positive, but the opposite
choices are also possible.

\section{Introducing braid foliations}
\label{section:preliminaries}

\fullref{theorem:MTWS} is about the relationship
between two closed braid diagrams which represent the same link.  However, the method
we use to prove it will not be in the setting of link diagrams.
Rather, we will be working with the immersed annulus $\ta$, and with
certain `braid foliations' of $\ta$.  Foliated surfaces have been
used before, in our earlier papers
\cite{B-M_I}--\cite{B-M_VI}.  In this section we will review and describe the
machinery which we use from these other papers.  The reader who has seen
these foliations before will be able to omit this section and go directly
to \fullref{section:braid foliations of the immersed annulus}, possibly pausing
to refer back to this review to refresh his/her memory of details.  A
more detailed review may be found in the review article by Birman and
Finkelstein \cite{B-F}. 

\subsection{Braid foliations of Seifert surfaces}
\label{subsection:braid foliations of Seifert surfaces}
We are given $X\in\cB(\cX)$ with $b(X)=n$ and with braid structure $(\axis, \fib)$ in $\reals^3$.
To make this review as simple as possible,
we assume that $X$ is a knot.  Choose 
a Seifert surface $\seif$ of minimal genus for $X$.  After modifying ${\seif}$ we will show that it supports
a special type of singular foliation which was studied and used by the authors in
\cite{B-M_I}--\cite{B-M_VI}.  We call it  a {\it braid foliation}.

There are choices of orientation which determine the sign conventions in braid foliations.
First, we assign the standard orientation to $\reals^3$ and choose the braid axis $\axis$
to be the positively oriented $z$--axis.  Using cylindrical coordinates, this determines a sense of increase
of the polar angle coordinate $\theta$.  Next, the closed braid
$X$ is oriented so that it points in the direction of increasing
$\theta$ at each point of $X \cap H_\theta$.  The orientation on $X$ induces an orientation on
$\seif$, and so determines a positive normal at each interior point on $\seif$.

\begin{figure}[htpb]
\labellist\small
\pinlabel {$\axis$} [b] at 283 436
\pinlabel {$+$ vertex} [t] at 140 358
\pinlabel {$-$ vertex} [t] at 200 358
\pinlabel {viewed on the surface} at 170 318
\pinlabel {\parbox{40pt}{$+$ side\\of surface}} [l] at 307 402
\pinlabel {\parbox{40pt}{$-$ side\\of surface}} [l] at 307 360
\pinlabel {viewed in 3--space} at 300 318
\pinlabel {near the boundary} [t] at 471 330
\pinlabel {(a)} at 240 300
\pinlabel {(b)} at 470 300
\endlabellist
\centerline{\includegraphics[scale=.8]{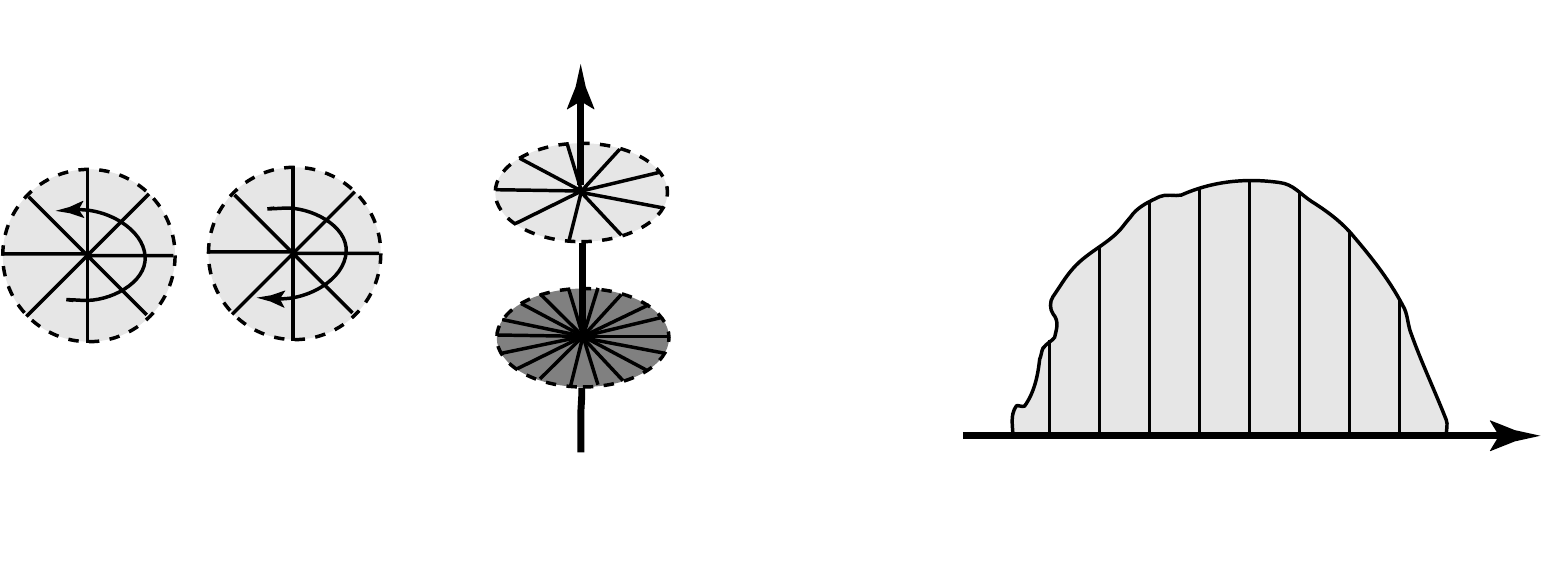}}
\caption{(a) Foliation of surface near a vertex\quad
(b) Foliation of surface near the boundary}
\label{figure:foliations0}
\end{figure}
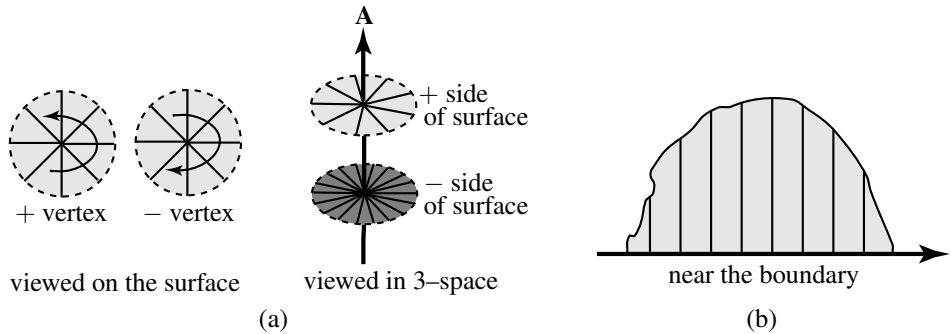

 Since $X = \partial\seif$ is a closed braid, $\axis\cap\seif$ is non-empty.  The intersections of $\fib$ with
$\seif$ may be assumed to be (a) radial in a neighborhood of each point of
$\axis \cap {\seif}$ and
(b) transverse to the boundary in a neighborhood of $\partial {\seif}$.  See \fullref{figure:foliations0}.  

{\it Vertices} in the foliation are points where $\axis$ pierces $\seif$.  
We call a vertex positive or negative, according as $\axis$ intersects $\seif$ from the 
negative or positive side of $\seif$ respectively.  The sketches in \fullref{figure:foliations0}
illustrate the positive side of $\seif$, so that the flow induced by the
fibration is anticlockwise (resp.\ 
clockwise) about a positive (resp.\  negative) vertex.

\emph{Singularities} 
in the foliation occur at points where $\seif$ is tangent to one of the fibers of
$\fib$.  The singularities may be assumed to be finite in number and to occur on distinct
fibers of $\fib$.  By Morse's classification theorem every singularity may be assumed
to result from a local maximum or minimum or a saddle point tangency between ${\seif}$ and a fiber of
$\fib$.  Let $s$ be a singular point of the foliation of ${\seif}$, and
let $H_\theta$ be the disc fiber which contains $s$.  We say that the singularity $s$ is 
positive if the outward-drawn oriented normal to the
oriented surface ${\seif}$ coincides in direction with the normal to
$H_\theta$ in the direction of increasing $\theta$.  Otherwise $s$ is 
negative. 

{\it Leaves} in the foliation 
are components of 
intersection of $H_{\theta}$ with the surface  ${\seif}$.  A {\it singular leaf} is a leaf
which contains a singularity of the foliation.  Every other leaf is a {\it non-singular
leaf}.

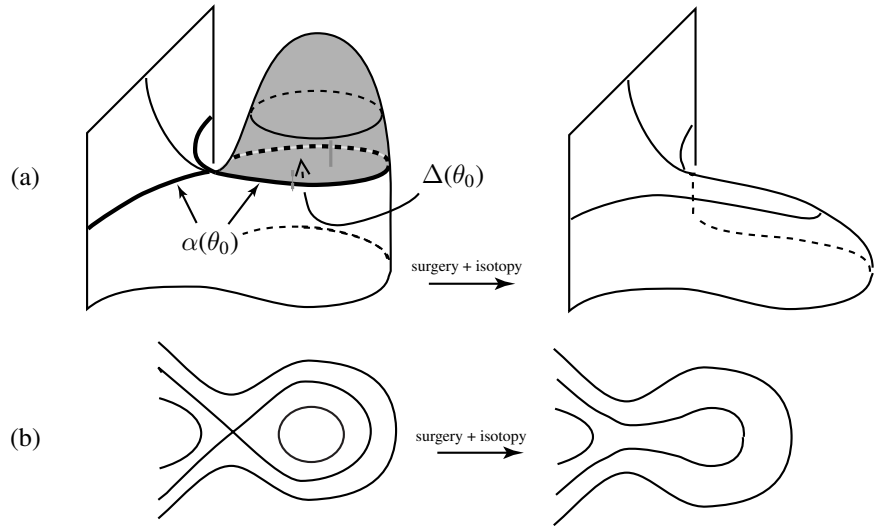
\begin{figure}[htpb!]
\labellist\small
\pinlabel {(a)} at 93 437
\pinlabel {(b)} at 93 314
\pinlabel {\tiny surgery + isotopy} [b] at 301 388
\pinlabel {\tiny surgery + isotopy} [b] at 303 308
\pinlabel {$\Delta(\theta_0)$} [bl] at 277 430
\pinlabel {$\alpha(\theta_0)$} [t] at 180 416
\endlabellist
\centerline{\includegraphics[scale=.8]{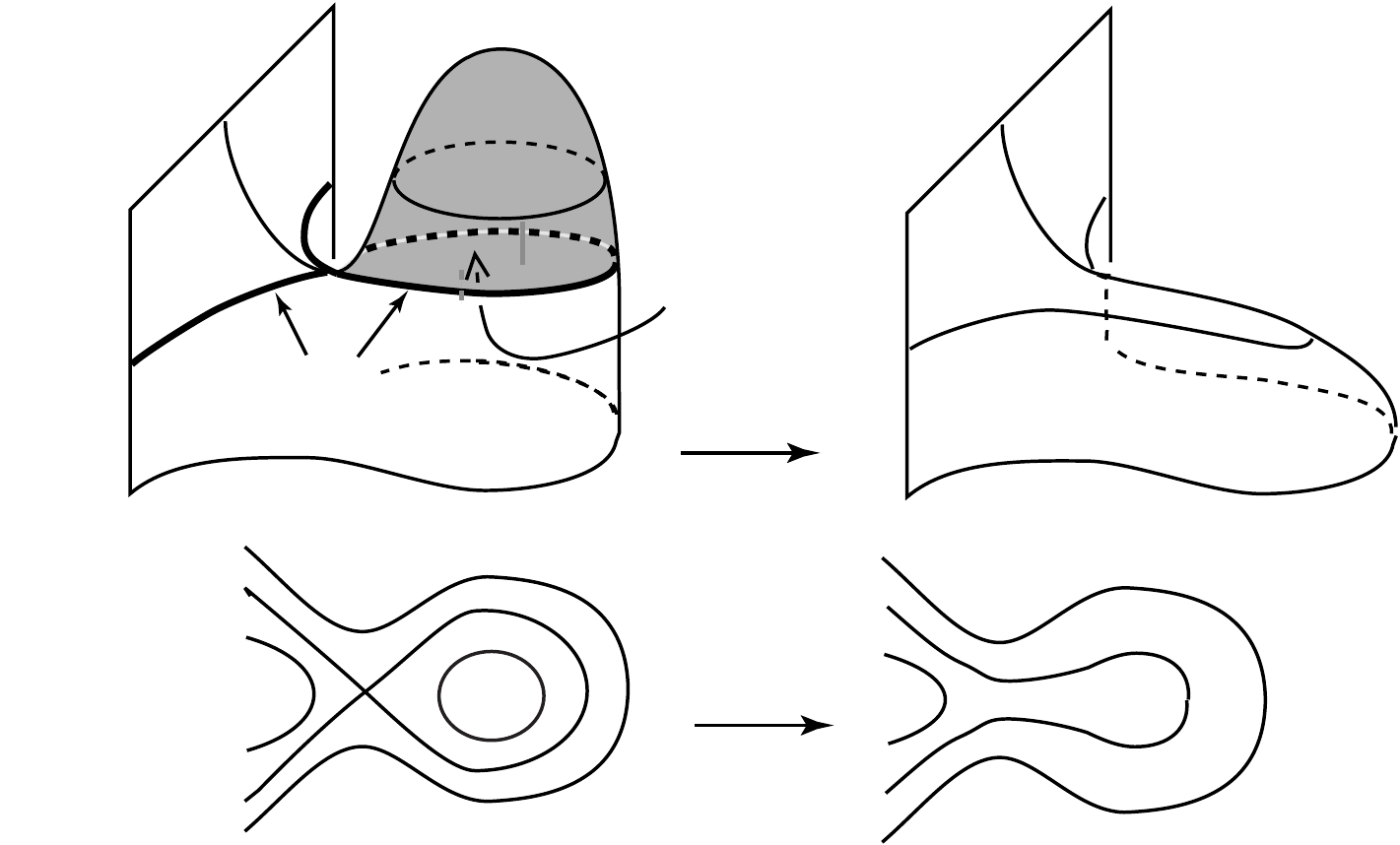}}
\caption{In this sketch the fibers of $\fib$ are to be thought of as 
horizontal planes.  (a) Surgering
$\seif$ along the disc $\Delta$ removes the singularity.  (b) The change in foliation
on $\seif$.}
\label{figure:foliations1}
\end{figure}

A very basic property of our braid foliations of Seifert surfaces is that non-singular leaves may all
be assumed to be arcs.  We review the reasons (which goes back to
Bennequin \cite{Bennequin}).
Suppose that there is a simple closed curve $\alpha(\theta)$ in $\seif\cap H_\theta$ for some
non-singular polar angle $\theta$.  The fact that $\seif$ is pierced non-trivially by
the braid axis, and that the foliation is transverse to the boundary, shows that
$\seif$ cannot be foliated entirely by simple closed curves, so if we follow the
sequence of arcs $\alpha(\theta)$ as
$\theta$ increases or decreases we must arrive at a singularity.  Let $H_{\theta_0}$ be
the singular fiber.  The singularity may be assumed to be either
be a center or a saddle point, but if it is a center, then by following
$\alpha(\theta)$ in the opposite direction we will arrive at another singularity,
and that one cannot also be a center because $\seif$ is not a 2--sphere, so it must
be a saddle point, and the singularity must be a homoclinic point (a singularity which is formed when
a generic leaf has a saddle point singularity with itself), as illustrated in \fullref{figure:foliations1}(a).  Note that the singular leaf
$\alpha(\theta_0)$ lies in both $\seif$ and $H_{\theta_0}$, and necessarily bounds
a disc $\Delta$ in the latter.  Assuming that $\Delta$ is innermost, we surger $\seif$ along $\Delta$ as in \fullref{figure:foliations1}.

The surgered surface has two components.  By
assumption $\seif$ has maximal Euler characteristic, which implies that one of the two
components is a 2--sphere.  Discarding it, and smoothing the new $\seif$, we can eliminate
the singularity.

It remains to consider the case when the interior of the disc $\Delta$ intersects
$\seif$.  Since $s_0$ is the only singularity in $H_\theta$, there are no
singularities in the interior of $\Delta$.  But then each component of $\seif\cap
\intr(\Delta)\subset H_\theta$ must be a simple closed curve.  Let $c$ be an innermost such
simple closed curve.  Then we can surger $\seif$ along $c$, and then smooth the surgered
surface by an isotopy.  
This will introduce center singularities but no additional saddles.  Any ${S^2}'s$ that are formed we
discard.
(We leave it to the reader to draw appropriate pictures).  After a
finite number of such surgeries we obtain a new surface
$\seif^\prime$ which has the same homoclinic point in its foliation, but does not
meet the disc $\Delta$.  We then do the surgery which is illustrated in \fullref{figure:foliations1}, reducing the number of saddle point singularities.  In this way all leaves which are
simple closed curves can be eliminated.  

Since each non-singular leaf is an arc, one of its endpoints could either be at a vertex of the foliation or a point
on the boundary, however we now claim that non-singular leaves which have both endpoints on the
boundary do not occur.  For, suppose there is a leaf $l$ which has both endpoints on the
boundary.  Let $N(l)$ be a neighborhood of the leaf on $\seif$.  Thinking of
$N(l)$ as a rectangle, it will have a pair of opposite edges on the boundary, and
the orientation on these edges will be consistent with the orientation on $N(l)$.
However, $l$ also lies in a fiber $H_\theta$ of $\fib$, and the boundary of $\seif$
is a closed braid.  But then, the orientation on one of the edges of $N(l)$ will
agree with that of the normal to
$H_\theta$ and (since $l$ is by hypothesis non-singular) the orientation of the
other must disagree with the that of the normal, but that is impossible because the
boundary curve is a closed braid.

\begin{figure}[htpb]
\labellist\small
\pinlabel {$\axis$} [b] at 260 672
\pinlabel {$\axis$} [r] at 290 654
\pinlabel {$\axis$} [l] at 344 654
\pinlabel {$X$} [l] at 386 630
\pinlabel {$a$} at 261 646
\pinlabel {$b$} at 315 654
\tiny
\pinlabel {$+$} at 261 667
\pinlabel {$+$} at 294 653
\pinlabel {$-$} at 341 652
\endlabellist
\centerline{\includegraphics[scale=.9]{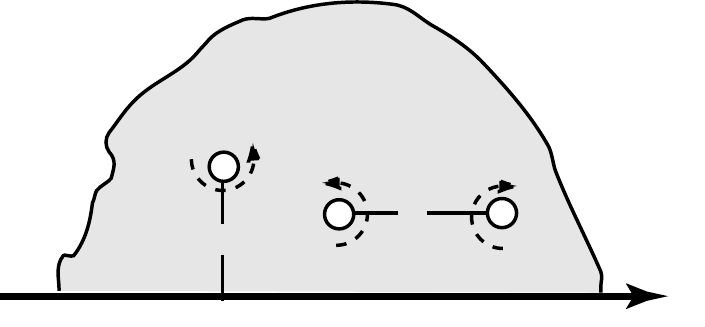}}
\caption{Non-singular leaves in the foliation of $\seif$}
\label{figure:foliations2}
\end{figure}

It follows from this that the non-singular leaves
have two types: those which have one endpoint on the boundary and the
other at a vertex (we call them {\it $a$--arcs}) and those which have both endpoints being vertices of
the foliation (we call them {\it $b$--arcs}).
See \fullref{figure:foliations2}.

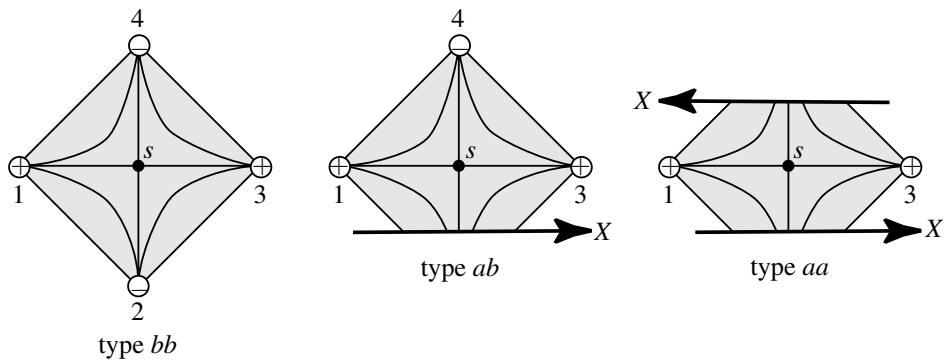
\begin{figure}[htpb]
\labellist\small
\pinlabel {$1$} [t] at 52 313
\pinlabel {$2$} [t] at 117 249
\pinlabel {$3$} [t] at 184 313
\pinlabel {$4$} [b] at 117 391
\pinlabel {$-$} at 118 383
\pinlabel {$-$} at 118 251
\pinlabel {$+$} at 53 319
\pinlabel {$+$} at 185 319
\pinlabel {$s$} [bl] at 117 320
\pinlabel {type $bb$} at 117 220

\pinlabel {$1$} [t] at 227 313
\pinlabel {$3$} [t] at 359 313
\pinlabel {$4$} [b] at 292 391
\pinlabel {$s$} [bl] at 293 320
\pinlabel {$+$} at 227 319
\pinlabel {$+$} at 359 319
\pinlabel {$-$} at 293 383
\pinlabel {type $ab$} at 293 262
\pinlabel {$X$} [l] at 362 283

\pinlabel {$1$} [t] at 407 313
\pinlabel {$3$} [t] at 539 313
\pinlabel {$s$} [bl] at 472 320
\pinlabel {$+$} at 408 319
\pinlabel {$+$} at 540 319
\pinlabel {$X$} [r] at 402 355
\pinlabel {$X$} [l] at 543 285
\pinlabel {type $aa$} at 473 262
\endlabellist
\centerline{\includegraphics[scale=.7]{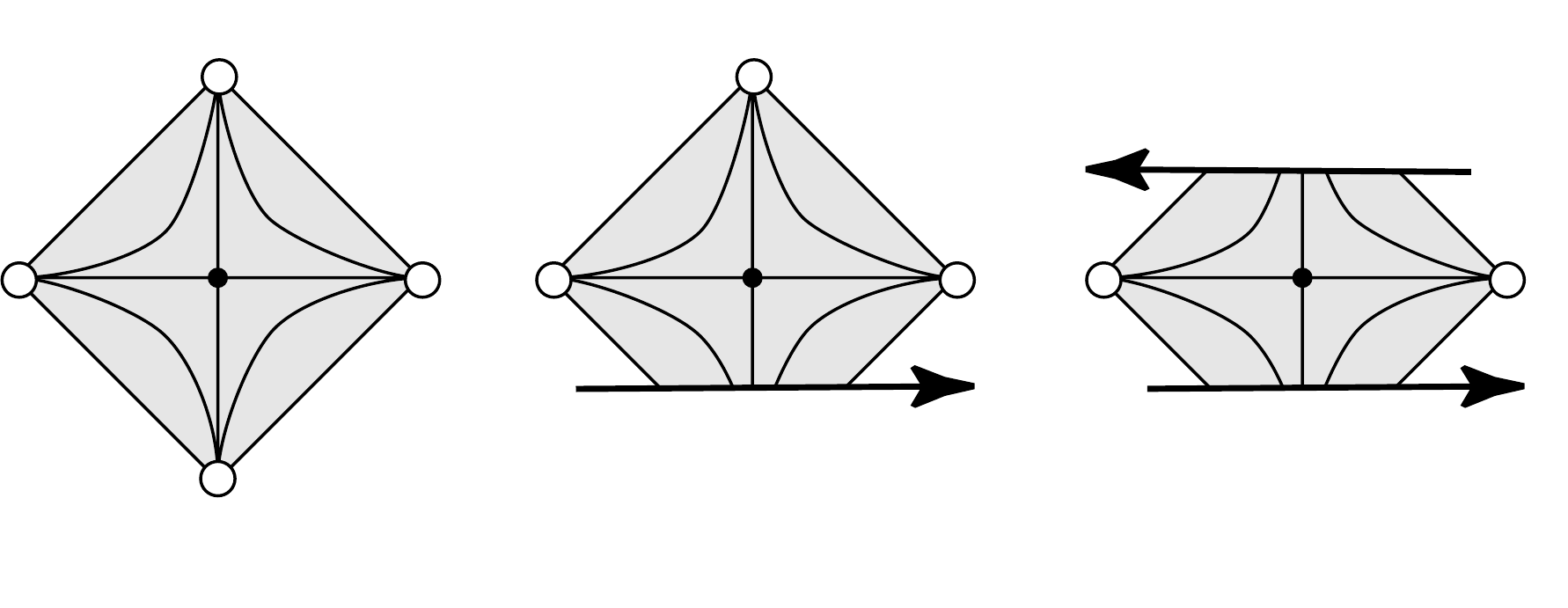}}
\caption{Tiles of type $bb,ab,aa$}
\label{figure:foliations3}
\end{figure}

Singularities  fall into three types, which we call types 
$aa, ab$ and $bb$, the notation indicating that just before an $aa$--singularity
(resp.\  $ab$, $bb$--singularity) the non-singular leaves were both type $a$
(resp.\  types $a$ and $b$, resp.\  both type $b$).  We shall refer to the 2--cells which are foliated
neighborhoods of the singular leaves as `tiles'.  See \fullref{figure:foliations3}.

The foliation may be used to decompose the surface ${\seif}$ into
a union of foliated $2$--cells, each of which contains exactly one singularity of
 the foliation.  Each 2--cell is a regular neighborhood on
$\seif$ of a singular leaf.  These foliated
$2$--cells are our {\it tiles} and the resulting decomposition of ${\seif}$ is a  {\it tiling}.
See the three sketches in \fullref{figure:foliations3}.  The 
{\it tile vertices} are
the points where the braid axis $\axis$ intersects the  surface ${\seif}$.
(There are also other vertices on $\partial\seif$, but we prefer
to exclude them when we refer to tile vertices.) The {\it tile edges} 
are arbitrary choices of $a$--arcs or $b$--arcs.  (There are also other tile edges
which are  subarcs of $X$, but it will be convenient to ignore those, just as we ignored the vertices 
which are on $X$.)

Assume that the vertices and singularities in the braid foliation of $\seif$ have been labeled.
The {\it combinatorial data} which is associated to the foliation is a
listing of the signed vertices, in the natural cyclic order in which they occur on $\axis$, and a
listing of the signed singularities, in the natural cyclic order in which they occur in the fibration.
The following proposition is very natural, because leaves in the foliation of $\seif$ are level sets for the
embedding of $\seif$: 

\begin{proposition}[{{Birman--Menasco \cite{B-M_VI},
Birman--Finkelstein \cite[Theorem~4.1]{B-F}}}]
\label{proposition:embedding of X} 
Let $T$ be any tile in the braid foliation of $\seif$.  Label the vertices and singularities of $T$ by their signs
and their cyclic orders as above.  Then this decorated braid foliation determines the embedding of $T$.  More
globally, the embedding of all of the $aa$ and $ab$ tiles in the foliation of $\seif$  determines the embedding
of $X$ in $\reals^3 \setminus \axis$.
\end{proposition}

\begin{figure}[htpb]
\labellist\small
\pinlabel {$\axis$} [b] at 79 718
\pinlabel {$\axis$} [b] at 245 718
\pinlabel {$\axis$} [b] at 412 718
\pinlabel {$X$} [b] at 349 672
\pinlabel {$X$} [b] at 539 695
\pinlabel {$X$} [b] at 545 639
\pinlabel {$s$} [l] at 142 632
\pinlabel {$s$} [tl] at 304 618
\pinlabel {$s$} [l] at 505 617
\pinlabel {$1$} [r] at 71 590
\pinlabel {$1$} [r] at 239 602
\pinlabel {$1$} [r] at 403 599
\pinlabel {$2$} [r] at 71 629
\pinlabel {$2$} [r] at 239 639
\pinlabel {$3$} [r] at 71 663
\pinlabel {$3$} [r] at 239 693
\pinlabel {$3$} [r] at 403 690
\pinlabel {$4$} [r] at 71 699
\pinlabel {type $bb$} at 110 555
\pinlabel {type $ab$} at 279 555
\pinlabel {type $aa$} at 509 555
\endlabellist
\centerline{\includegraphics[scale=.75]{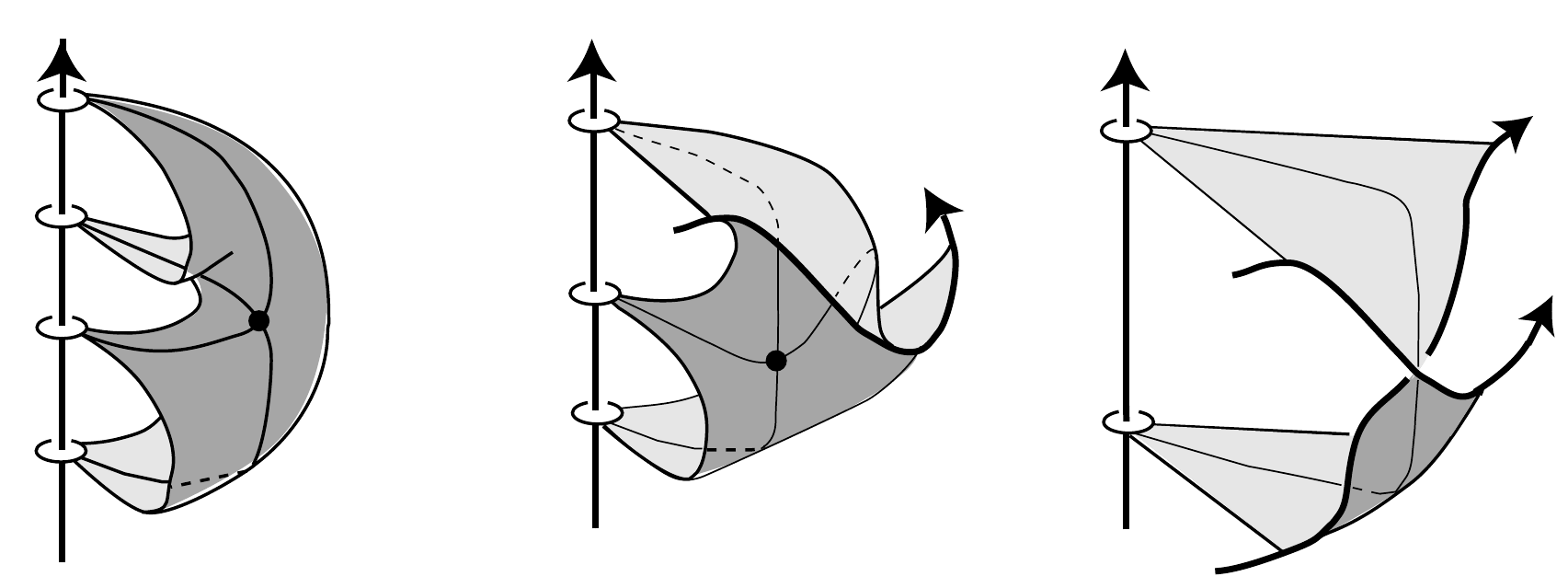}}
\caption{Embeddings of tiles of type $bb,ab,aa$ in 3--space}
\label{figure:foliations3b}
\end{figure}
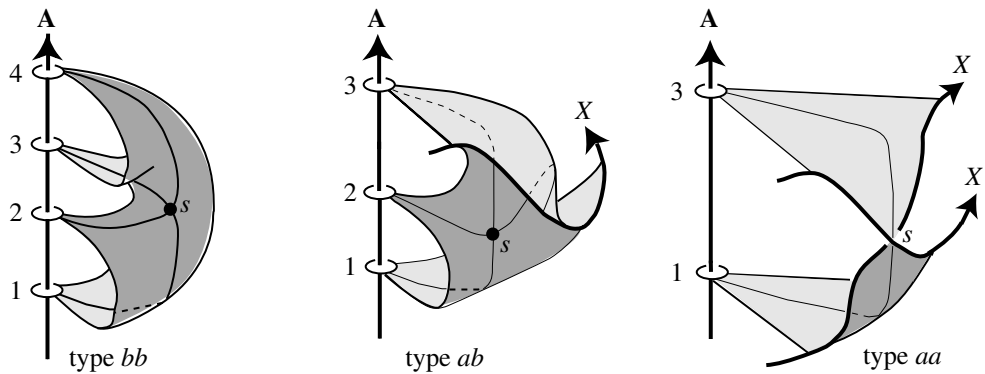

We illustrate in \fullref{figure:foliations3b} the embeddings of the three tiles in
\fullref{figure:foliations3}, for one of the finitely many choices of the combinatorial data on
these tiles.

\subsection{Control over the foliations}
\label{subsection:control over the tiling}
Braid foliations are not unique, and in this section we describe
some of the ways we have discovered to modify them.  This is an essential part
of the argument in the proof of the MTWS, because the vertices of valence 1 and 2 that we use
to recognize destabilizations and exchange moves may not be present initially, but after a
change in foliation they may be present.  The existence of a vertex of valence 1 or 2 will
be one of our ways to learn when the complexity can be reduced.
 
The foliation of ${\seif}$ depends
upon the choice of half-planes $H_\theta$ in the fibration of $\reals^3 \setminus \axis$.  A
{\it change in fibration} 
is the choice of a new set of half-planes $\fib =\{H_\theta : 0 \leq \theta \leq 2\pi\}$.
Equivalently, one could fix the fibers of $\fib$ and the link, and move the surface.  This induces a 
{\it change in foliation}.

The changes which we shall make are always very controlled and very local in terms of changing the tiling.  In
particular they do not change the braid, and are supported in a neighborhood in $\reals^3$ of subarcs of singular
leaves.  The question of when such changes are possible has been studied.
They were used in earlier papers by the authors, and we use those results as needed here.
We describe two changes in foliation.  For details, see
Birman--Finkelstein \cite{B-F}. 
\begin{figure}[htpb]
\centerline{\includegraphics[scale=.7]{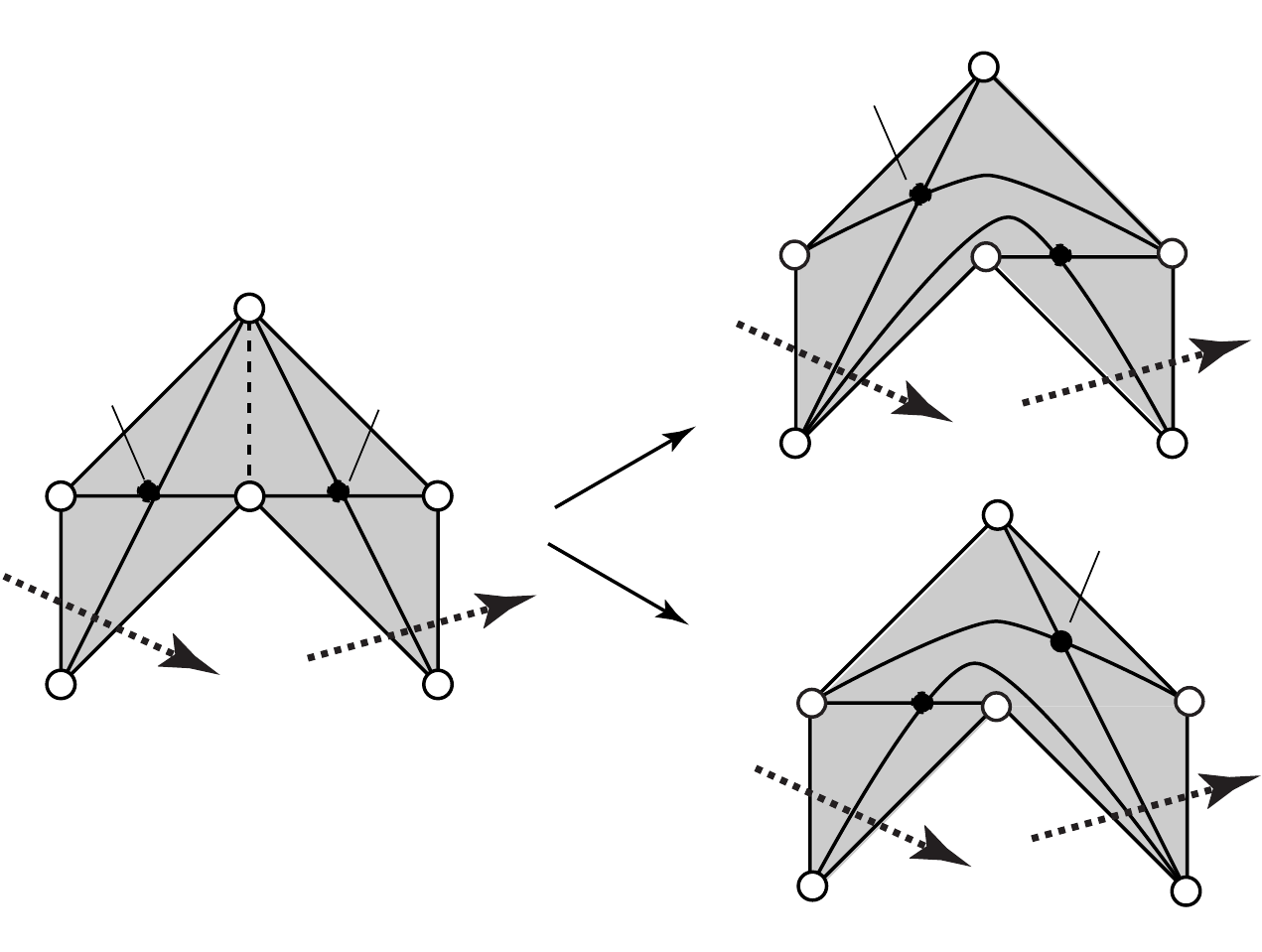}}
\caption{The first change in foliation, in the case when both tiles have type $bb$.  There are two other cases,
obtained from the one which is illustrated by deleting part of the $bb$--tile and adding one or both dotted arcs
as boundary, to convert to an $ab$--tile.  The 3 cases are: case (1), two
$bb$--tiles; case (2), one $ab$ and one
$bb$--tile; case (3), two $ab$--tiles.}
\label{figure:first change in foliation}
\end{figure}
\begin{lemma}[First change in foliation]
\label{lemma:first change in foliation}
Let $s_1$ and
$s_2$ be singularities of the same sign $\delta = \pm$ in tiles $D_1$ and
$D_2$, where $D_1$ and $D_2$ intersect in a common leaf $vw$ of type $b$.  For example, see the left
sketch in \fullref{figure:first change in foliation}, which relates to the case when both tiles
have type $bb$.  Then after a change in foliation which is supported on a neighborhood  in
3--space of an arc $\gamma$ which joins the two singular points, the foliation of $D_1\cup D_2$
changes in one of the two ways which are illustrated in \fullref{figure:first change in foliation}.  In particular, the
valence of the vertices $v$ and $w$ decreases as a result of the change.
\end{lemma} 
\begin{proof}  See \cite[Theorem~2.1]{B-F} for a very detailed proof of the Lemma.  We note the following feature of the
proof.  There are 6 vertices in $D_1\cup D_2$, labeled $v,a,y,w,z,b$.  The proof in \cite{B-F} shows that they remain
fixed during the change in foliation.  There are also singular leaves $vy$
(resp.\  $vz$) with one endpoint at the
vertex
$v$ in the region of interest and the other at $y$ (resp.\  $z$).  In the passage from the left to the middle
(resp right) sketch in \fullref{figure:first change in foliation} there is, at every point in the isotopy, a
singular leaf $vy$ (resp.\  $vz$) which contains one of the two singularities.
This justifies our labeling that
singularity  $s_1'$ (resp.\  $s_2'$) in the middle (resp.\ 
right) sketch, because it evolves directly from $s_1$ (resp.\ 
$s_2$) during the isotopy which realizes the change.
We call the other singularity $s_2'$ (resp.\ 
$s_1'$).  In fact (see \cite{B-F})  if
$s_1$ and $s_2$ occur in an angular interval $[\theta_1,\theta_2]\subset [0,2\pi]$ in the fibration $\fib$, with
$s_1 < s_2$ in this interval, then after the change in foliation the new singularities will still be in the same
angular interval, only now, in both cases, we will have $s_1' > s_2'$.\end{proof}  

We will need one more change in foliation.  It is similar to that of \fullref{lemma:first change in foliation}, except that it holds without restriction as
to the signs of the two singularities which are involved.  We call it the {\it second}  change
in foliation.
The reader may find \fullref{figure:second change in foliation} 
helpful in understanding what it says.
\begin{figure}[htpb!]
\labellist
\tiny
\pinlabel {$1$} [b] at 189 657
\pinlabel {$2$} [b] at 168 651
\pinlabel {$3$} [t] at 154 553
\pinlabel {$4$} [t] at 182 542
\pinlabel {$5$} [tl] at 299 566
\pinlabel {$6$} [bl] at 310 621

\pinlabel {$1$} [b] at 406 661
\pinlabel {$2$} [b] at 376 651
\pinlabel {$3$} [t] at 367 559
\pinlabel {$4$} [t] at 403 546
\pinlabel {$5$} [tl] at 517 570
\pinlabel {$6$} [bl] at 528 625

\hair=1pt
\pinlabel {$1$} [br] at 352 434
\pinlabel {$2$} [bl] at 385 520
\pinlabel {$3$} [br] at 209 520
\pinlabel {$4$} [br] at 227 482
\pinlabel {$5$} [bl] at 494 473
\pinlabel {$6$} [bl] at 475 431

\hair=2pt
\pinlabel {$s_1$} [l] at 177 604
\pinlabel {$s_2$} [r] at 270 599
\pinlabel {$\gamma$} at 146 613
\pinlabel {$p$} [l] at 253 562

\pinlabel {$s_1'$} [bl] at 426 625
\pinlabel {$s_2'$} [bl] at 402 604
\pinlabel {$\gamma$} at 356 620
\pinlabel {$p$} [tl] at 462 571

\pinlabel {$\alpha$} [r] at 268 437
\pinlabel {$\gamma$} [bl] at 312 452
\pinlabel* {$p$} [bl] at 310 431
\pinlabel {$s_1$} [r] <2pt,0pt> at 288 465
\pinlabel {$s_2$} [b] at 402 408
\pinlabel {$\Delta$} at 292 442

\small
\pinlabel {(a)} at 93 599
\pinlabel {(b)} at 93 457

\endlabellist
\centerline{\includegraphics[scale=.60]{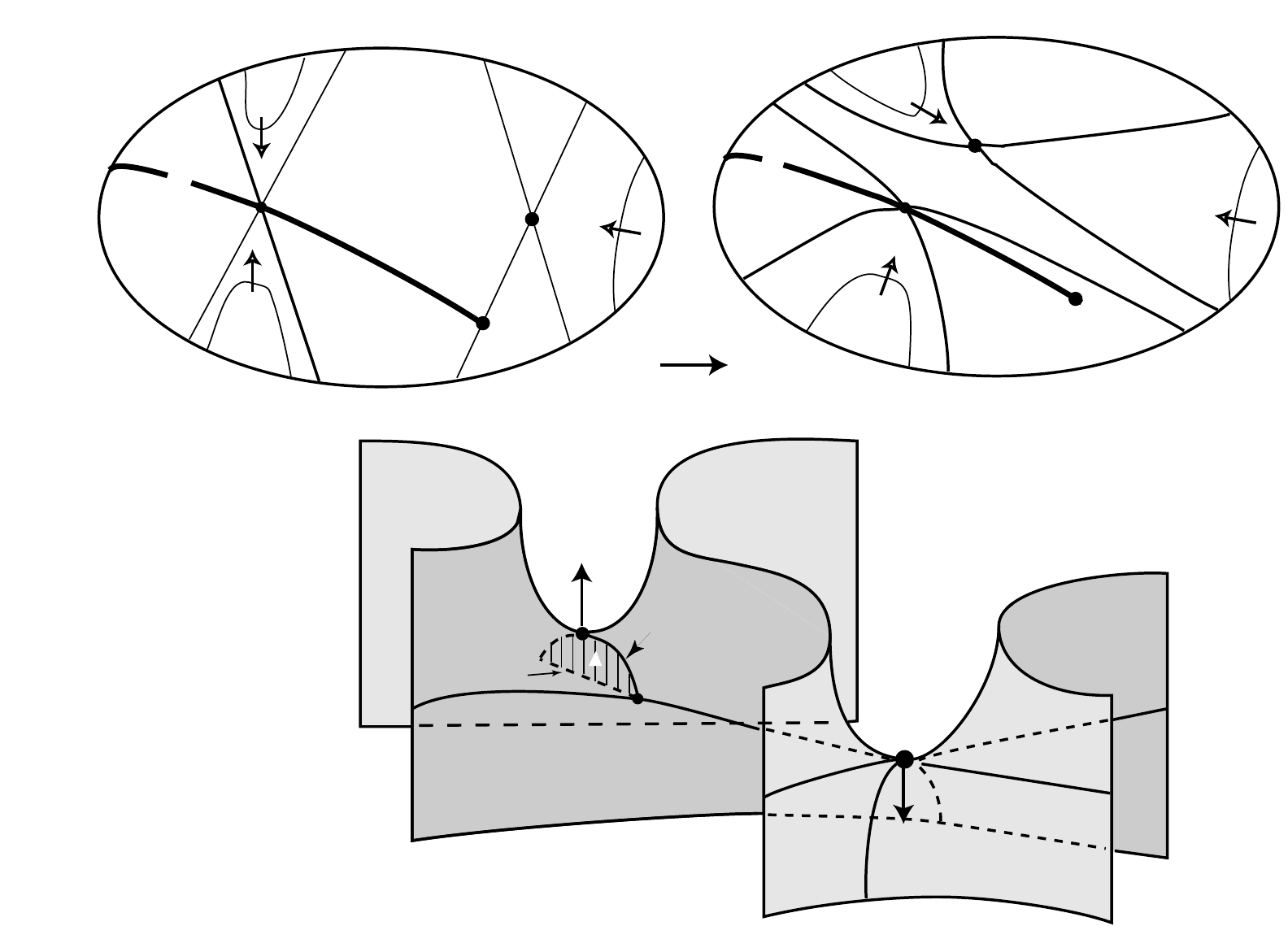}}
\caption{Sketch (a) shows the second change in foliation.  Sketch (b)
shows $N$, embedded in 3--space, before the change, illustrating the position of
$\Delta$.  The move is a push downward of $s_1$ along $\Delta$.}
\label{figure:second change in foliation}
\end{figure}
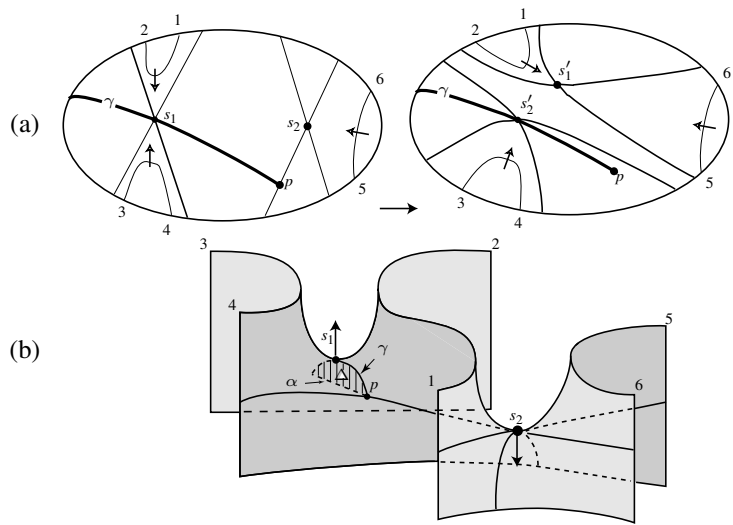
\begin{lemma}[Second change in foliation]
\label{lemma:second change in foliation} 
Let $\gamma$ be an arc which is located in a
foliated disc $\nhb$ which is a subset of a foliated surface.  Assume that the
foliation of $\nhb$  contains exactly two singular points $s_1$ and $s_2$.  Let $l_i$ be the singular leaf through
$s_i, i=1,2$.  Let $H_{\theta_i}$ be the fiber of $\fib$ which contains $s_i, i=1,2$.  Suppose that there is a disc
$\Delta$ in 3--space, such that:
\begin{enumerate}
\item[{\rm(i)}] The interior of $\Delta$ has empty intersection with $\seif$.

\item[{\rm(ii)}] $\partial \Delta = \alpha \cup \gamma$, with  $\Delta \cap \nhb = \gamma$  and  $\Delta
\cap H_{\theta_2} = \alpha $.

\item[{\rm(iii)}] $\Delta$ is trivially foliated, ie there are no vertices or singularities in $\Delta$.

\item[{\rm(iv)}] $\gamma \cap l_1 = s_1$ and $\gamma \cap l_2 =\{p\}$, where $p$ is a point.
\end{enumerate}
 Then after a change in foliation which is induced by pushing $N$ along
$\Delta$, changing the order of $s_1$ and $s_2$, we may assume that
$\gamma \cap l_1 = s_2^\prime$ and $\gamma \cap l_2 = \emptyset.$  
Moreover, the change in the foliation of $\seif$ may be assumed to be supported on
 an arbitrarily
small neighborhood on $\nhb$ of the subarcs $[s_1,p] \subset \gamma$ and
$[p,s_2] \subset l_2.$   There could also be several singularities 
$s_1,s_2,\dots s_k$, with associated neighborhoods, and if the conditions are
met for each of them in turn then the disc $\Delta$ may be used to push $s_1$ past
many singularities, one at a time. 
\end{lemma}

\begin{proof}  Figures \ref{figure:second change in foliation}(a)
shows the foliated disc $\nhb$ before and after the
change which we propose to make.  The arrows which are attached to the leaves
indicate the direction of increasing $\theta$.  Using the foliation of $\nhb$, and
knowing the signs of the
 singularities, one
may construct an embedding of $\nhb$ in 3--space, and we have done so in \fullref{figure:second change in foliation}(b) in the case
when the signs of the singularities at $s_1$ and $s_2$ are different.  (The other
case is similar).
In \fullref{figure:second change in foliation}(b)
fibers of $\fib$ are to be thought of as horizontal planes.  The auxiliary disc $\Delta$, is also 
illustrated.
The move which we make to realize the change in foliation in
\fullref{figure:second change in foliation}(a)
is to push $\nhb$ down
along the disc $\Delta$.
To understand how this changes the foliation of $\nhb$, we
have labeled
certain endpoints on $\partial \nhb$ with numbers 1,2,3,4,5,6.  There are non-singular
leaves which
we call 12, 34, 56, each with arrows directed inward (to illustrate the direction
of increasing $\theta$) and joining 1 to 2, 3 to 4 and 5 to 6 respectively.  In the left
picture the first
singularity occurs when leaf 12 approaches leaf 34, but in the right picture the
 first singularity
occurs when leaf 12 approaches leaf 56. \end{proof}

\subsection{Using braid foliations to detect destabilizations and exchange moves}
\label{subsection:recognizing destabilizations and exchange moves from the foliation}
In this section we will show how to recognize, from the foliation of $\seif$, when a closed braid admits
a destabilization or an exchange move.

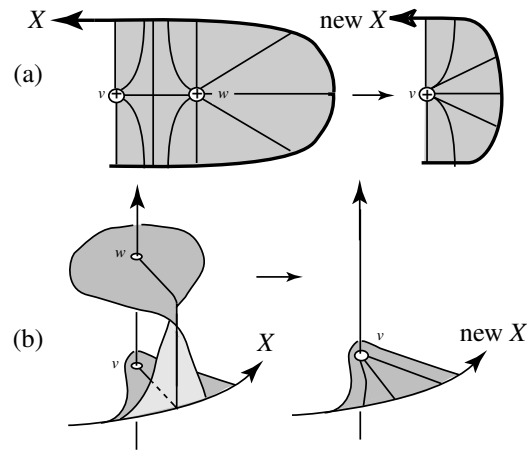
\begin{figure}[htpb!]
\labellist\small
\pinlabel {$X$} [r] at 215 689
\pinlabel {$X$} [b] at 340 489
\pinlabel {(a)} at 200 656
\pinlabel {(b)} at 200 500
\pinlabel {new $X$} [r] at 416 689
\pinlabel {new $X$} [b] at 474 495
\tiny
\pinlabel {$v$} [r] at 433 645
\pinlabel {$v$} [r] at 250 645
\pinlabel {$w$} at 315 645
\pinlabel {$w$} [r] at 262 551
\pinlabel {$v$} [r] at 258 486
\pinlabel {$v$} [bl] at 402 497
\endlabellist
\centerline{\includegraphics[scale=.65]{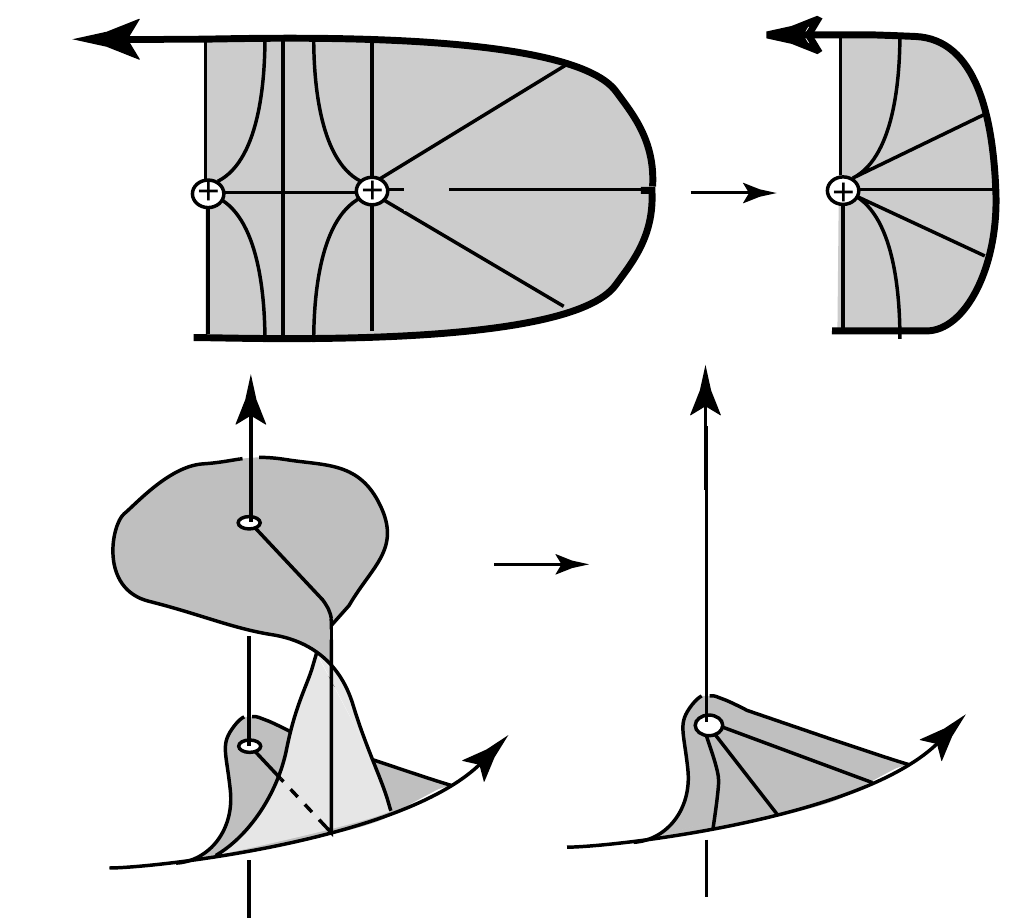}}
\caption{Destabilization along a valence 1 vertex, viewed (a) on the
foliated surface and (b) in 3--space}
\label{figure:foliations4}
\end{figure}

Destabilizations are easy.  It is shown in \cite{B-F} that $\Xhigh$ admits a
destabilization if the foliation has a vertex of valence 1, as in \fullref{figure:foliations4}(a).  The embedding of a tile which contains the vertex of
valence 1, for one of the two possible  choices of the sign of the singularity, is illustrated in
sketch (b).

\begin{lemma} 
\label{lemma:recognizing destabilizations}
Destabilizations which are predicated on the existence of a valence 1 vertex reduce the number of singularities
(resp.\  vertices) in the foliation of $\seif$ by 1 (resp.\  1).
\end{lemma}

\begin{proof}  Clear.  See \fullref{figure:foliations4}(a).\end{proof} 

Before we can describe our exchange moves, a new concept is needed.  We observe that, intuitively,
$b$--arcs in the foliation of $\seif$ arise when there are `pockets' in the surface, and we  are now interested in
the case when a `pocket is empty' and so can be removed.  We now make this precise.  A
$b$--arc $\beta$ in the foliation is an
{\it inessential $b$--arc} if it joins a pair of vertices $v,w$ which are consecutive vertices in the natural 
cyclic ordering of vertices along $\axis$.  We use the term `inessential' because, if we think of a
fiber $H_\theta$ of $\fib$ as a disc with the braid axis $\axis$
as its boundary, then an inessential $b$--arc  will cut off a disc $\Delta\subset H_\theta$ which
has the property $\Delta\cap\seif = \beta$.  The disc $\Delta$ cannot be
pieced by $X$, for if it were there would be an $a$--arc $\alpha$ which
either intersected our inessential $b$--arc, which is impossible because $\seif$  is embedded, or $\alpha$ has a vertex endpoint between $v$ and $w$, which is impossible because $v$ and $w$ are consecutive vertices on $\axis$. The disc $\Delta$ can be used to push
$\beta$ (and nearby leaves in nearby disc fibers) across $\axis$, reducing the number of vertices
in the foliation of
$\seif$.  Peek ahead to the right sketch in \fullref{figure:exchange4} for
examples of $b$--arcs which are
inessential and also essential.  An {\it essential} $b$--arc is one which is not inessential.

Exchange moves were defined in \fullref{figure:exchange1} as a move  on a block-strand diagram.
Our task now is to
detect them in the foliation.  The complexity function mentioned in \fullref{theorem:MTWS}
will include the number of singularities in
the foliation of a clasp annulus which is bounded by the given braids.  As will be seen shortly,
exchange moves always reduce this complexity function.
However, the exchange moves that are used in this paper come in two
flavors, and it's necessary to check both. 

The \emph{$ab$--exchange move}:
\begin{figure}[htpb!]
\labellist\small
\pinlabel {(a)} at 64 641
\pinlabel {(b)} at 64 496

\pinlabel {$X$} [l] at 260 662
\pinlabel {$w$} [br] at 157 643
\pinlabel {$v$} [tl] at 189 612
\pinlabel {$q$} [tr] at 157 613
\pinlabel {$s$} [bl] at 188 643
\pinlabel {$p$} [t] at 157 577
\pinlabel {$t$} [l] at 228 643
\pinlabel {$\alpha$} at 229 575

\pinlabel {new $X$} [l] at 459 664
\pinlabel {$\alpha'$} at 354 653

\pinlabel {$\axis$} [r] at 105 457
\pinlabel {$\alpha$} [t] at 165 406
\pinlabel {$X$} [l] at 202 394
\pinlabel {$X$} [l] at 252 416
\pinlabel {$v$} [t] at 154 497
\pinlabel {$w$} [b] at 177 519

\pinlabel {$\axis$} [r] at 357 463
\pinlabel {$\alpha'$} [b] at 451 423
\pinlabel {$X$} [l] at 448 400
\pinlabel {new $X$} [l] at 505 421
\endlabellist
\centerline{\includegraphics[scale=.65]{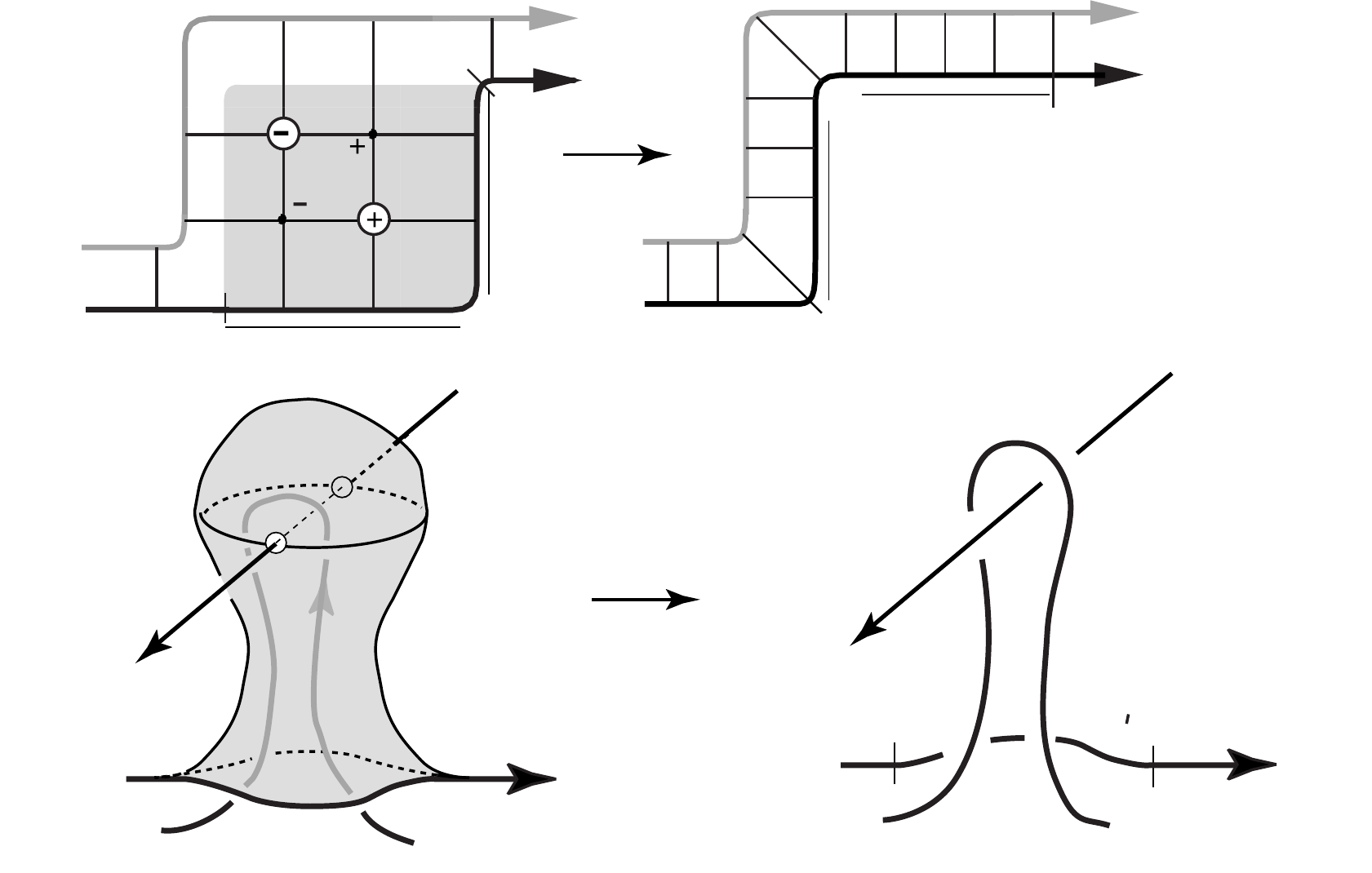}}
\caption{The type $ab$--exchange move}
\label{figure:exchange5}
\end{figure}
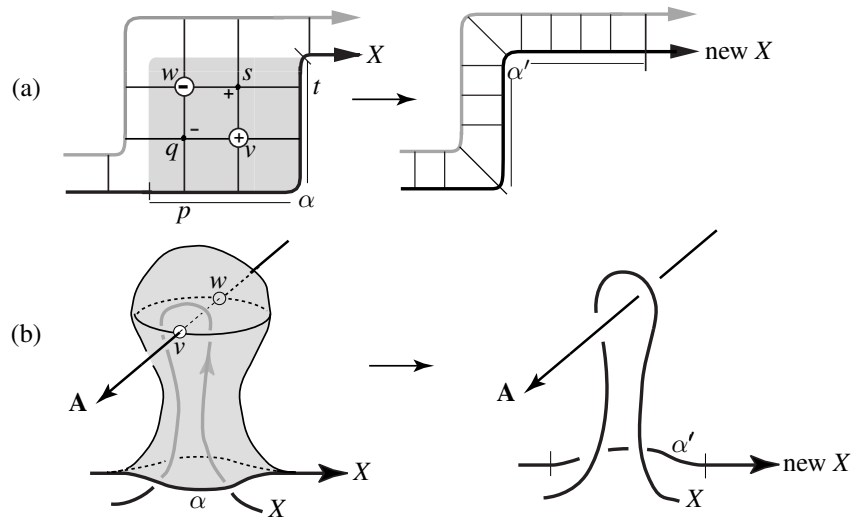
In \fullref{figure:exchange5}(a) we have illustrated $v$, a {\it type $ab$} valence 2 vertex.
The left and right sketches in (a) (resp.\  (b)) show the configurations on the foliated
surface $\seif$ (resp.\  in 3--space), the left sketches being before and the right sketches being after the exchange
move.  In both (a) and (b) the vertex $v$ is adjacent to two type $ab$ singularities,
and these singularities have opposite signs.
The $ab$--exchange move occurs in the passage from the left to the right sketch.  One pushes the subarc 
$\alpha\subset X$ across the shaded disc  to  $\alpha'$.
If we choose $\alpha^\prime$ so that it is everywhere transverse to leaves of the foliation (as is clearly
possible) the move takes closed braids to closed braids.
\fullref{figure:exchange5}(b) illustrates the embedding for one particular choice of signs and orderings.

Note that if the $b$--arcs which end at the vertex $v$ are all essential, then $X$ must
encircle the subarc $vw$ of the braid axis which is inside the pocket.  The
braid changes we make can then be understood by looking at the projection onto a plane orthogonal to the
braid axis.  After the exchange move the shaded disc will have vanished.  Peering down the braid axis
(as we did in \fullref{figure:exchange1}) we see that the projection has changed in the
predicted manner.  For full details, consult \cite{B-F}.

The \emph{$bb$--exchange move}:
In Figures \ref{figure:exchange4} and \ref{figure:foliations7} we have illustrated $v$, a
{\it type $bb$} valence 2 vertex.  The non-singular
leaves which have an endpoint at $v$ are all type $b$.  The left
and right sketches show the changes in 3--space.  The changes in the braid
projection can be understood by looking down the axis onto a plane orthogonal to
the braid axis.
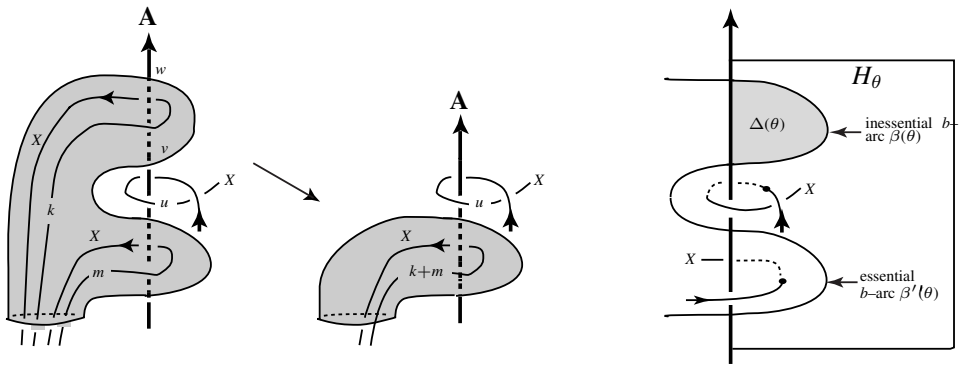
\begin{figure}[htpb]
\labellist\tiny
\hair=2pt
\pinlabel {\small $\axis$} [b] at 140 596
\pinlabel {$X$} [r] at 85 536
\pinlabel {$X$} [b] at 111 478
\pinlabel {$X$} [l] at 179 514
\pinlabel {$k$} at 87 498
\pinlabel {$v$} [b] at 149 526
\pinlabel {$w$} [b] at 147 569
\pinlabel {$m$} at 113 463
\pinlabel {$u$} at 149 501

\pinlabel {\small $\axis$} [b] at 310 550
\pinlabel {$X$} [b] at 281 478
\pinlabel {$X$} [l] at 348 514
\pinlabel {$u$} at 321 501
\pinlabel {$k{+}m$} at 292 464

\pinlabel {$\Delta(\theta)$} at 478 544
\pinlabel {$X$} [l] at 496 509
\pinlabel {$X$} [r] at 442 471
\pinlabel {\small $H_\theta$} at 533 570
\pinlabel {\parbox{35pt}{inessential $b$--arc $\beta(\theta)$}} [l] at 529 541
\pinlabel {\parbox{35pt}{essential\\$b$--arc $\beta'(\theta)$}} [l] <0pt,-2pt> at 526 459
\endlabellist
\centerline{\includegraphics[scale=.7]{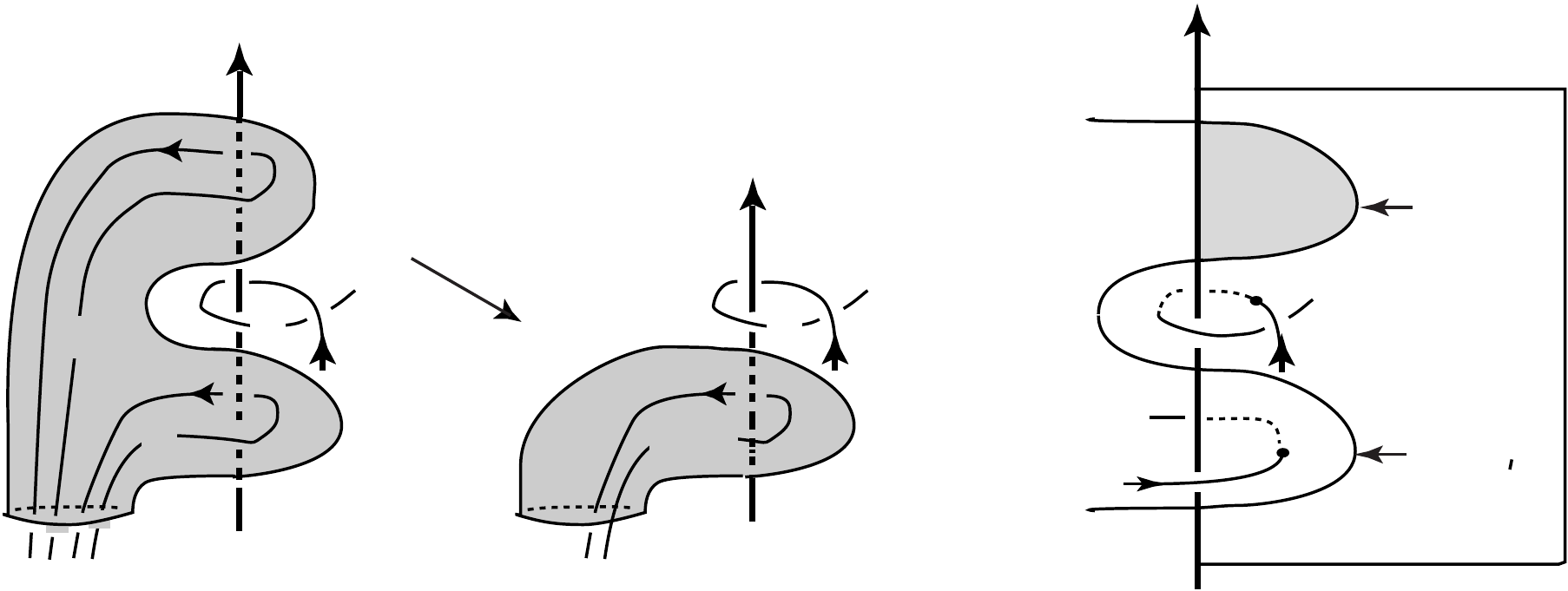}}
\caption{The type $bb$--exchange move, followed by an isotopy of $\seif$} 
\label{figure:exchange4}
\end{figure}
The foliation does not change at all after a $bb$--exchange move, however there
is a change in the order of the vertices along the braid axis.  The `pocket has
been emptied' and after the exchange move, the empty pocket can be collapsed by
pushing every $b$--arc $\beta(\theta)$ across  its disc $\Delta(\theta)$, as in \fullref{figure:foliations7}.  Remark: The pocket could of course have much more complicated
braiding inside it.  The proof that it can always be emptied in this way is
non-trivial; details may be found in our earlier paper \cite{B-M_IV} or
the article by Birman and Finkelstein \cite{B-F}.  The passage from the left to right sketches in \fullref{figure:foliations7} shows
the change in the foliation, after the removal of all inessential $b$--arcs.
\begin{figure}[htpb]
\labellist\small
\hair=2pt
\pinlabel {$p'$} [bl] at 91 410
\pinlabel {$p$} [bl] at 246 372
\pinlabel {$q$} [br] at 164 369
\pinlabel {$s$} [tr] at 164 421
\pinlabel {$v$} [l] at 172 396
\pinlabel {$w_1$} [r] at 218 395
\pinlabel {$w_2$} [l] at 109 395
\pinlabel {$p'$} [bl] at 379 405
\pinlabel {$p$} [tl] at 431 382
\pinlabel {$w_2$} [t] at 391 383
\endlabellist
\centerline{\includegraphics[scale=.7]{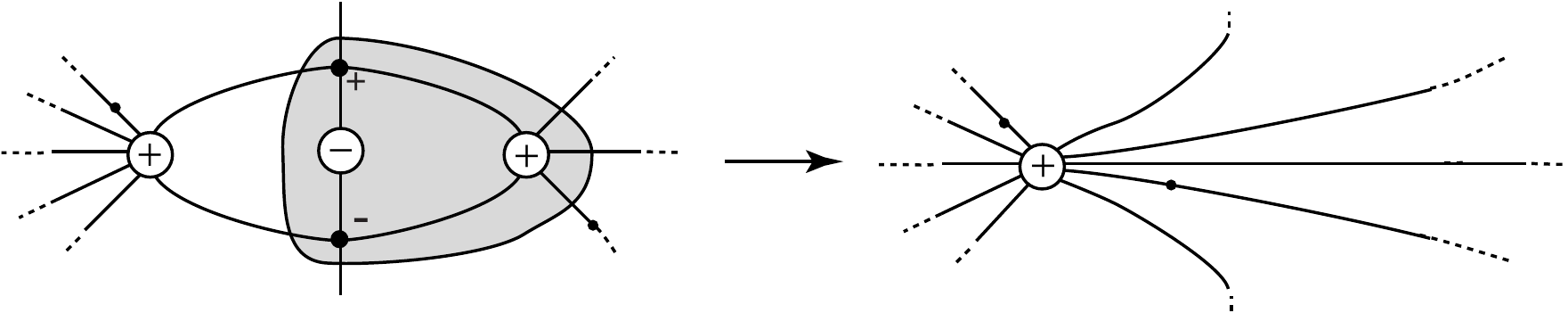}}
\caption{Changes in foliation after the removal of inessential $b$--arcs} 
\label{figure:foliations7}
\end{figure}
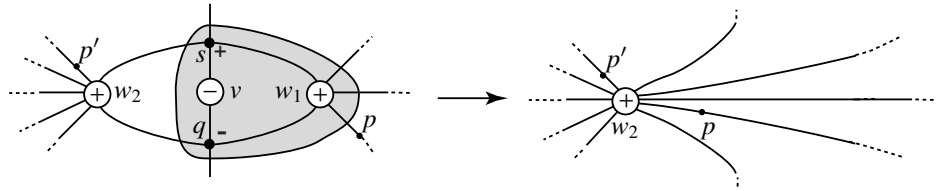 

\begin{lemma}
\label{lemma:ab and bb exchange}
Suppose that the foliation of $\seif$ has a vertex
$v$ of  valence 2 and type $ab$ or $bb$.  Assume that the adjacent singularities have
opposite signs.  Then the closed braid  $X=\partial\seif$ admits an exchange move.
After the move, there is a surface ${\seif}^{'}$, isotopic to
${\seif}$, with a decomposition containing two fewer vertices and two fewer
singularities than the decomposition of ${\seif}$.
\end{lemma} 
\begin{proof}  See Figures \ref{figure:foliations4}(a), \ref{figure:exchange4} and 
\ref{figure:foliations7}.\end{proof}  

\subsection{Using braid foliations to detect stabilizations}
\label{subsection:using braid foliation to detect stabilizations}
As elaborated on in the  introduction to this paper, one of the questions that motivated our work was a
desire to understand why stabilization played such an important role in the classical Markov Theorem.  Having braid foliations of Seifert surfaces bounded by knots and links in hand, one answer to that question became clear:  they  allow one to simplify a Seifert surface by eliminating negative vertices and their associated singularities.
\begin{figure}[htpb!]
\labellist\small
\pinlabel {(a)} at 213 387
\pinlabel {(b)} at 423 387
\pinlabel {stabilize} [b] at 211 484
\pinlabel {stabilize} [b] at 422 484
\pinlabel {\parbox{50pt}{negative\\singularity}} [b] at 413 433
\pinlabel {\parbox{50pt}{positive\\stabilization}} [b] at 512 436

\pinlabel {$v_i$} [br] at 129 454
\pinlabel {$v_j$} [bl] at 199 454
\pinlabel {$\alpha'$} [r] at 154 441
\pinlabel {$\alpha$} [t] at 163 410
\pinlabel {$p$} [t] at 143 409
\pinlabel {$q$} [t] at 184 409
\pinlabel {$K$} [l] at 205 408

\pinlabel {$p$} [t] at 243 410
\pinlabel {$q$} [t] at 283 410
\pinlabel {$v_i$} [br] at 225 454
\pinlabel {$v_j$} [bl] at 295 454
\pinlabel {$K$} [l] at 303 408

\pinlabel {$K$} [l] at 408 408
\pinlabel {$K$} [l] at 508 408
\endlabellist
\centerline{\includegraphics[scale=.90]{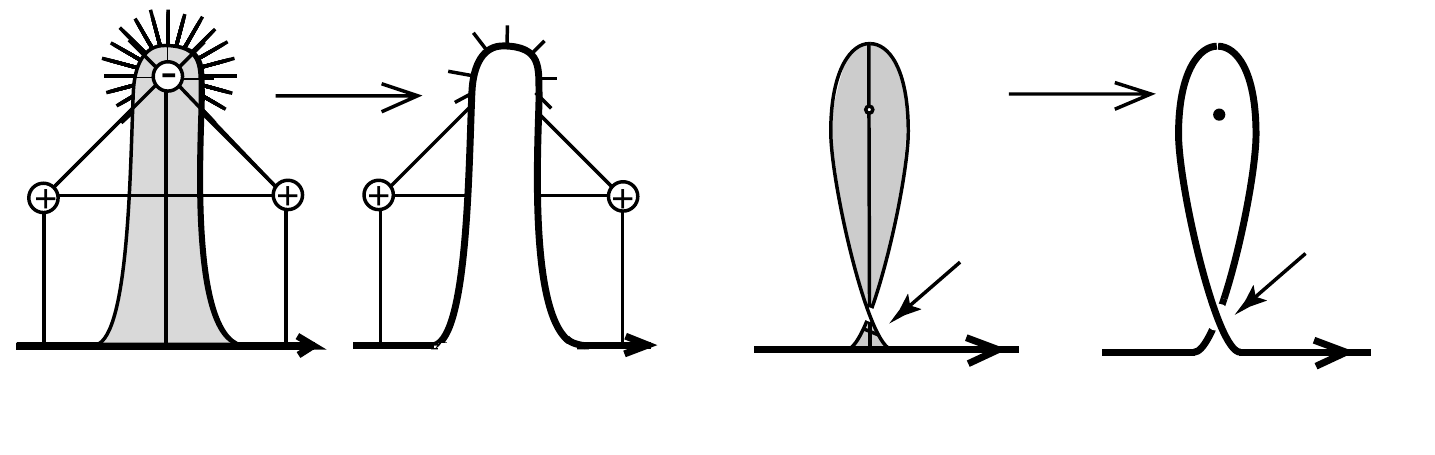}}
\caption{Stabilization along an $ab$--tile, viewed (a) as a move on the foliated surface, and (b) concentrating on
how it alters the boundary} 
\label{figure:foliation5}
\end{figure}
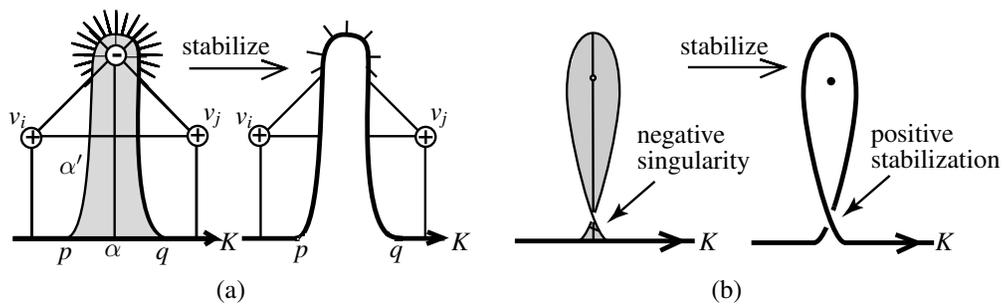
The first sketch in \fullref{figure:foliation5} shows how: When
there is an $ab$--singularity in the foliation, we may push $X$ across the
associated negative vertex and its singularity, in a neighborhood of
the singular leaf, to a new position which is again
everywhere transverse to the foliation.  It follows that after we do
this move we will have a new closed braid representative.
This moves simplifies $\seif$ because it eliminates a vertex and a singularity.  The right sketch
shows why 
the move is actually a stabilization. \fullref{figure:foliations6} shows
our stabilization move on the embedded surface in 3--space.  If one looks carefully one can see the
half-twist which has been introduced in the course of the push.  We
note that the pictures of $ab$--tiles in \fullref{figure:foliations6} are deformations of the
picture in \fullref{figure:foliations3}: we stretched out the top
sheet to make visible a neighborhood of the singular leaf.
\begin{figure}[htpb!]

\centerline{\includegraphics[scale=.90]{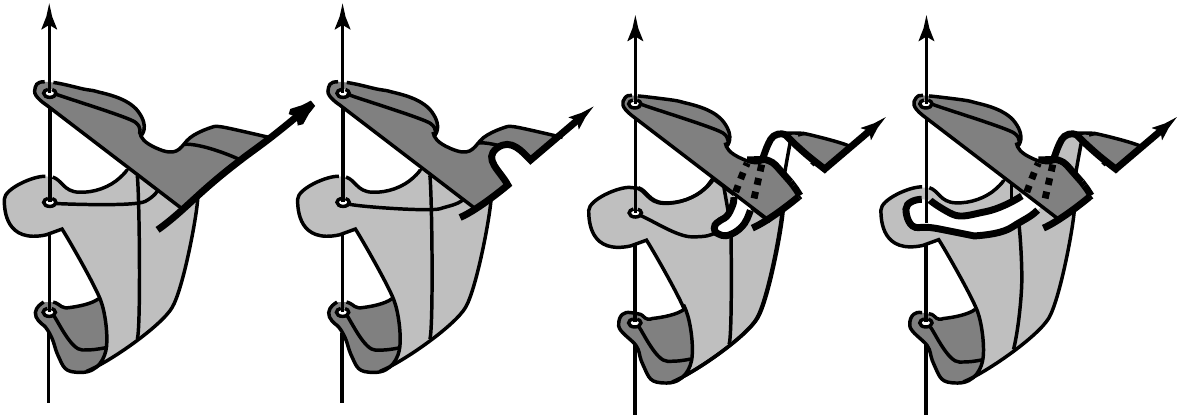}}
\caption{Stabilization along an $ab$--tile, viewed in 3--space} 
\label{figure:foliations6}
\end{figure}

\section{Braid foliations of the immersed annulus}
\label{section:braid foliations of the immersed annulus} 

We are now ready to investigate braid foliations of the clasp annulus $\ta$ which we constructed
in \fullref{section:getting started}.  We will need to confront two new aspects of the geometry:
\be
\item There is a (basically trivial) new aspect to our geometry,  which unfortunately will lead to
new bookkeeping: The clasp annulus $\ta$ has two boundary components, $\Xhigh$ and $\Xlow$.  By the basic
construction, the
$\Xlow$ boundary is a curve in the interior of the chosen Seifert surface $\seifmed$ for the
closed braid
$\Xmed$.  As was shown in the previous section, $\seifmed$ admits a braid foliation.  Clearly $\Xlow$
can cut through the tiles in that foliation in any way as long as 
it  is always transverse to the leaves of the foliation.  Therefore  we will allow for new
tile types, to account for the partial tiles at the $\Xlow$ boundary of $\ta$. 
\item The second new aspect of the geometry is central to the work in this paper.  Our annulus $\ta$
is not embedded.
\ee

A preliminary modification in the clasp arcs will be helpful in what follows:

\begin{lemma}
\label{lemma:clasp arcs transverse to foliation}
We may assume that the $k$ clasp arcs are transverse to fibers of $\fib$, and so also to the leaves in the braid
foliation of $\ta$ and $\pa$.
\end{lemma}

\begin{proof}   
Let $\Gamma =
\gamma_1\cup\dots\cup\gamma_k$ be the union of the clasp arcs.  Then  the graph
$\Xhigh\cup \Xmed\cup \Xlow\cup\Gamma$ is embedded in $\ta$, which is a subset of
$\reals^3$.  We focus now on that graph.  By our earlier construction, its subsets $\Xhigh, \Xmed$
and
$\Xlow$ are in braid position, but in general $\Gamma$ is not in braid
position, ie the interior
of some clasp arc may not be transverse to the fibers of $\fib$..
From the proof of \fullref{lemma:construction of the two annuli}
we employ Alexander's braid trick to every wrongly oriented subarc of $\Gamma$,
doing it so as to avoid intersections with  $\Xhigh, \Xmed$ and $\Xlow$.
The construction allows us to find an orientation-preserving PL homeomorphism 
$f\co S^3\to S^3$ which changes $\Gamma$ to braid position.  A classical result of Gugenheim
(see Hempel \cite[Theorem~1.5]{Hempel}) then tells us that we may assume
that $f$ is isotopic to the identity.  Replacing $\ta$ by $f(\ta)$, and
using the fact that $f$ leaves $\Xhigh, \Xmed$ and $\Xlow$ invariant,
it follows that we may assume that every subarc of the graph $\Xhigh\cup
\Xmed\cup \Xlow\cup\Gamma \subset \ta$ is in braid position. \end{proof}

\subsection{Tile types in $\pa$}
\label{subsection:tile types in pa} 

Our work begins with the two closed braid  representatives $\Xhigh$ and $\Xlow$ of $\cX$.
Our basic construction in \fullref{section:getting started} gave us the immersed annulus $\ta =
\annhigh\cup\annlow$ and the clasp arcs.  We begin by studying the braid
foliations of the two annuli.  The key point which will allow us to apply the machinery of
\fullref{section:preliminaries} is that each is embedded.  The really new
feature is the presence of the clasp arcs.  We first make the clasp arc as nice as possible (we
change them to `tabs').  After that we will remove `short clasp arcs', ie ones which do not pass
through any singular leaves.  Then we will create  `normal neighborhoods' about the doubly modified
clasp arcs, which will isolate them from the rest of
$\ta$.  Our device for isolating them is to introduce lots of inessential
$b$--arcs, 
which will give us the freedom we need to modify the clasp arcs when we need to do so.    

In this section we are interested in the foliation of
$\annhigh$ and $\annlow$ which are induced by intersections of these annuli with the half-planes
of our braid structure.
The closed braids $\Xhigh, \Xmed, \Xlow$ are all oriented so that they
point in the direction of increasing $\theta$ at each point of intersection with an $H_\theta$.
We choose an orientation on the annulus $\annhigh$ in such a way that it induces the given
orientation on its boundary component $\Xhigh$.  Notice that this means that the 
orientation on
$\Xmed$ does not agree with that induced by the chosen orientation on
$\annhigh$.  Similarly, we choose an orientation on the annulus $\annlow$ in such a way 
that it induces the given orientation on its boundary component $\Xmed$, which implies
that the orientation on $\Xlow$ does not agree with that induced by the chosen
orientation on $\annlow$. 

As in the situation of Seifert surfaces, the
foliation may be assumed to be radial in a neighborhood of each point of $\axis \cap \cA_\pm$
(see \fullref{figure:foliations0}(a)) and
transverse to the boundary in a neighborhood of $\partial \cA_\pm$ (\fullref{figure:foliations0}(b)).  The braid axis
$\axis$ pierces $\cA_\pm$ from either the negative or the positive side at each pierce
point, and we have indicated this by attaching positive or negative signs to the
pierce points on $\cA_\pm$.  As before, leaves in the foliation are  singular if they
contain a singularity of the foliation, otherwise they are  non-singular.
The singularities may be assumed to be finite in number and to occur on distinct fibers of $\fib$. 

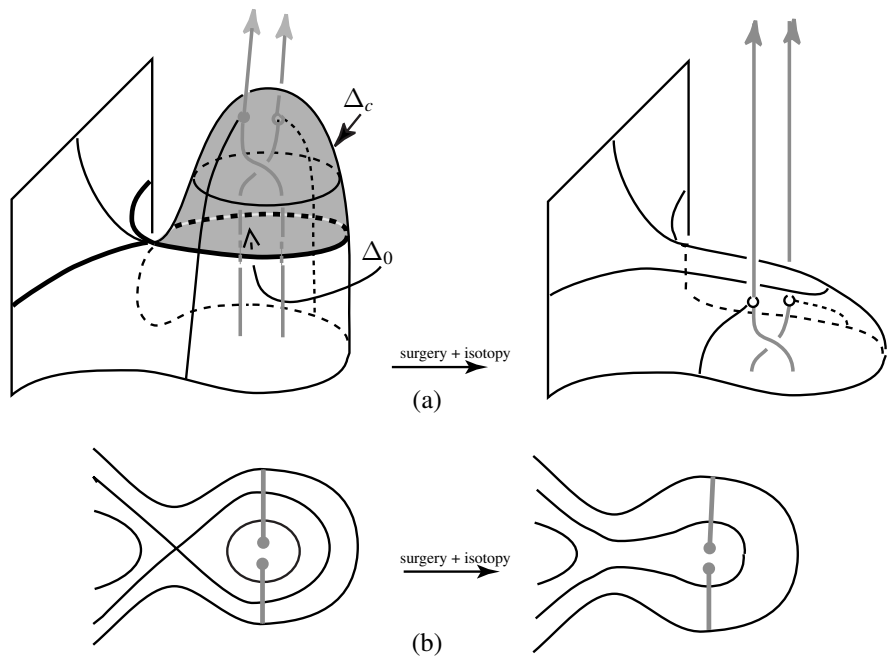
\begin{figure}[htpb!]
\labellist\small
\hair=1pt
\pinlabel {$\Delta_0$} [b] at 277 430
\pinlabel {$\Delta_c$} [b] at 269 495
\pinlabel {\tiny surgery + isotopy} [b] at 310 389
\pinlabel {\tiny surgery + isotopy} [b] at 310 302
\pinlabel {(a)} at 298 373
\pinlabel {(b)} at 298 270
\endlabellist
\centerline{\includegraphics[scale=.9]{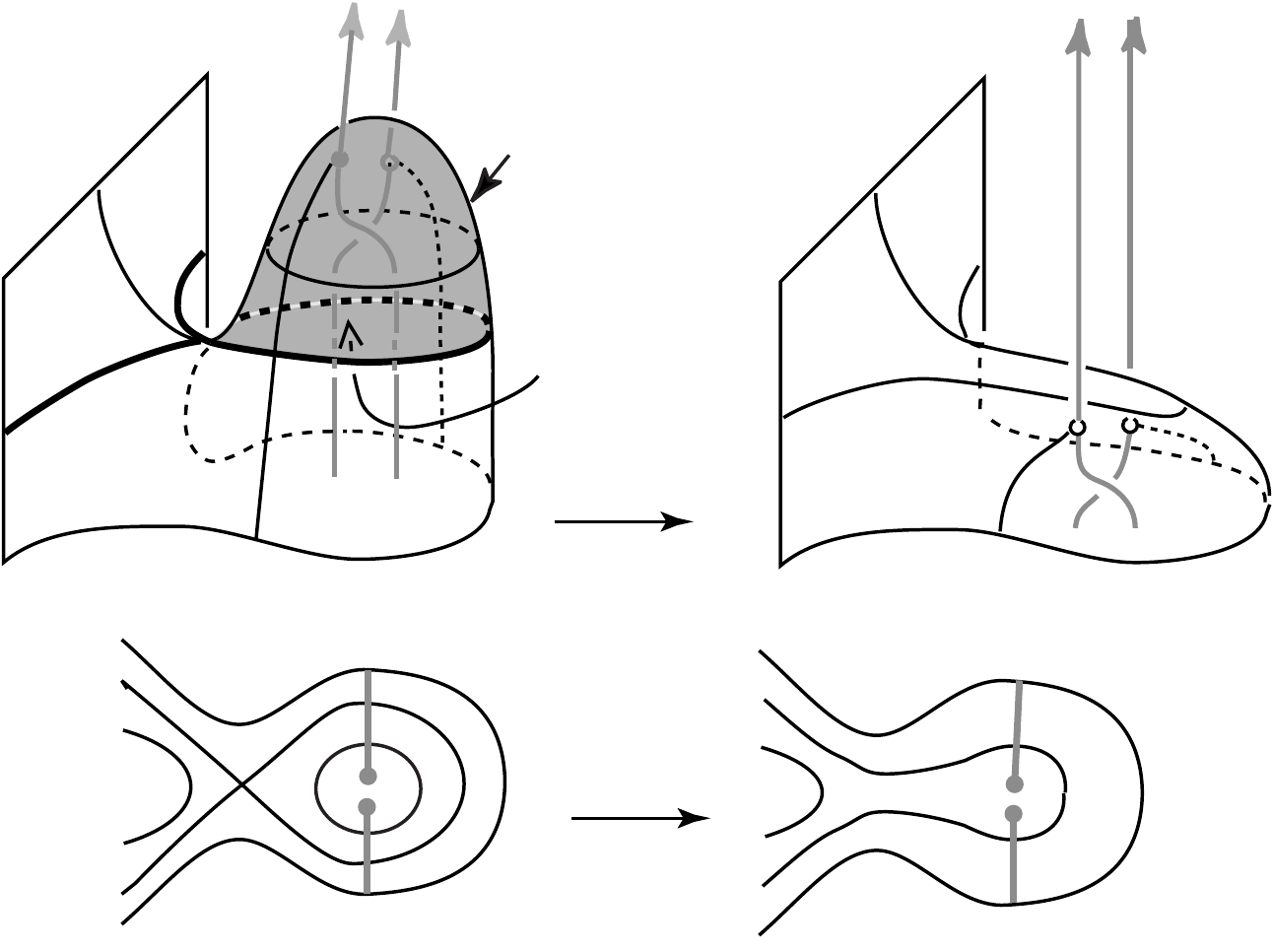}}
\caption{The elimination of a homoclinic singularity and associated simple closed curves, 
in the presence of clasp arcs. Sketch (a) is the geometric realization. Sketch (b) shows the preimage in $\pa$.}
\label{figure:annulus0}
\end{figure}

We now show that, as in the situation of Seifert surfaces, we may assume that there are
no leaves in the foliation of $\ta = \annhigh\cup\annlow$ which are simple closed curves.
See \fullref{figure:annulus0}.
There are new issues to settle.  The first question we ask is whether
a simple closed curve $c$ could intersect both
$\annhigh$ and $\annlow$?  Assume that the homoclinic loop determined by $c$ is innermost
in $H_\theta$, and also innermost in the foliation of $\ta$.  Then $c$  bounds a disc on
$\annhigh\cup\annlow$ which is foliated by simple closed curves, so it would intersect $\Xmed$ twice if it
intersects it at all.  This would force a tangency between $\Xmed$ and a fiber of
$\fib$, but $\Xmed$ is in braid position, so this cannot happen.
Therefore any leaf in the foliation which is a simple closed curve must lie in the interior of $\cA_\e$,
$\e=+$ or $-$.

Next, for simplicity assume that the foliation of $\annhigh$ contains the homoclinic loop $c$ and refer
to \fullref{figure:annulus0}.
Our leaf $c$ bounds  a disc $\Delta_\theta$ in the fiber $H_\theta$ which
contains $c$.  See 
\fullref{figure:annulus0}(a).  (If $c$ is homologically nontrivial in $\annhigh$
then, since it bounds a disc in $H_\theta$, $\Xhigh$ would have a component that was the unknot.)
If $c$ is not intersected by any clasp arcs then we can
apply the argument of \fullref{section:preliminaries} (see \fullref{figure:foliations1}).  Since both
$\Delta_\theta$ and $\Delta_c$ are embedded their union $\Delta_\theta
\cup_c \Delta_c$ is an embedded $2$--sphere. 

Now consider a clasp arc that intersects $c$.  Since clasp arcs are transverse to the foliation of
$\ta$, any clasp arc that intersects $c$ intersects it exactly once.  Thus, any clasp arc that intersects
$c$ must have a puncture endpoint in $\Delta_c$, where $\Xlow$ intersects $\annhigh$.  (Referring
to \fullref{figure:annulus0}(a), the gray braid strands represent $\Xlow$.)  Since any braid that
enters the $2$--sphere $\Delta_\theta \cup_c \Delta_c$ must puncture $\Delta_\theta$, and can only
exit by puncturing $\Delta_c$, we conclude that only $\Xlow$ can intersect $\Delta_\theta$. 

Next, we consider the intersection of $\Delta_\theta \cap (\annhigh \cup \annlow) \subset \Delta_\theta$.
In principle, this intersection set can contain three types of arcs/curves: arcs that have an endpoint
on both $\Xlow$ and $c$; arcs that have both endpoints on $c$; and simple closed curves.  Notice that there
can be no arcs 
that have both endpoints on $\Xlow$ because this would violate the orientation of $\Xlow$.
Since we are assuming that $c$ is innermost for the moment we ignore the issue of simple closed curves.
Referring to \fullref{figure:annulus0}(b) notice that for every puncture
point in $\Delta_c$ there is a point on
$c$ that is an intersection with a clasp arc.  But, also for every puncture point of $\Xlow$ with
$\Delta_c$ there must be a puncture point of $\Xlow$ with $\Delta_\theta$.  So the second type of
arc of intersection---having both endpoints on $c$---cannot occur.  If we then perform the surgery
illustrated in \fullref{figure:annulus0}(a), we see that this corresponds to a truncation of
the clasp arcs.
That is, this surgery on a simple closed
curve in the homoclinic leaf replaces the immersed annulus with a new immersed surface 
which is union of an annulus and 
a $2$--sphere, where the latter can be discarded.  It replaces a clasp intersection with a 
clasp intersection and ribbon intersection on the discarded $2$--sphere.  Afterwards, clasp 
intersections are still in braid position.
This surgery eliminates at least one saddle singularity (and, possibly some points of
$\ta \cap \axis$).  Reiterating this procedure we will arrive at a point where the foliation of
$\ta$ has no leaves that are simple closed curves, except in one situation. 

\begin{figure}[htpb!]
\labellist\tiny
\pinlabel {type $a_-$} at 236 642
\pinlabel {type $a_+$} at 235 568
\pinlabel {type $b$} at 319 602
\pinlabel {type $s$} at 354 641

\pinlabel {$X_-$} [l] at 395 665
\pinlabel {$X_+$} [l] at 395 548
\endlabellist
\centerline{\includegraphics[scale=.6]{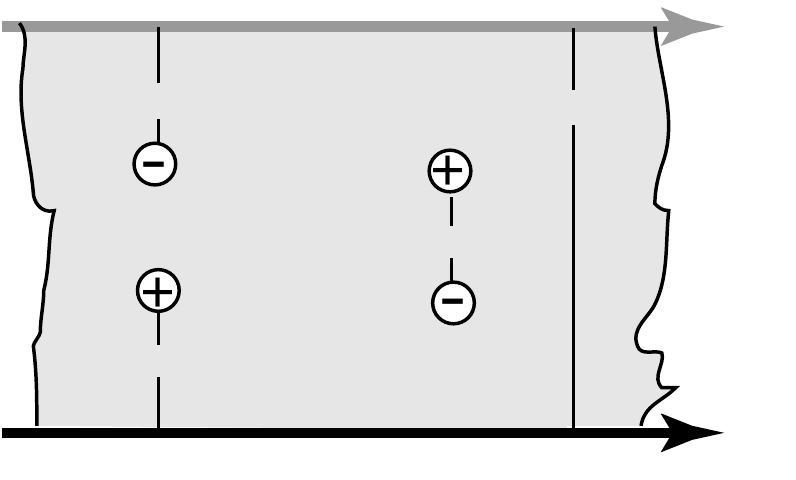}}
\caption{Non-singular leaves in the foliation of $\annhigh\cup\annlow$, as
viewed in the preimage of $\ta$}
\label{figure:annulus1} 
\end{figure}

The situation where we will not be able to perform the surgery illustrated in
\fullref{figure:annulus0} is when the homoclinic point involves a singularity between an
$s$--arc (see \fullref{figure:annulus1})
and a simple closed curve where the disc $\Delta_\theta$ contains the endpoints
of the resulting singular leaf which
are in $\Xhigh$ and $\Xlow$.  However, we can alter the foliation of
$\ta$ in the following manner to eliminate the leaf $c$.
Let $\ge$ be any clasp arc that has an endpoint $ e_{-\e} \subset \Delta_c$.
Let $\alpha \subset \ta $ be an arc that starts at $e_{-\e}$ and ends at
$X_\e$ such that $\intr(\alpha)$ does not
intersect any clasp arcs.  It is convenient to take the support of $\alpha$ to be in a union of
$\Delta_c$ and a regular neighborhood of the homoclinic singular leaf; and it can be assumed that
$\alpha$ is transverse to the foliation of $\ta$ except at a single point in $\Delta_c$ near $e_{-\e}$.
We now perform Alexander's braid trick on $\alpha$ to make $\alpha \cup \ge$ transverse to $\fib$.  In the
resulting new foliation of $\ta$ the path $\alpha \cup \ge$ cannot intersect any leaf that is a
simple closed curve for
reasons of orientation.  It is easily seen that once the homoclinic point is eliminated for $\ge$ it will be
eliminated for all clasp arcs that intersect $\Delta_c$. No new homoclinic points are introduced. 

Finally, if $c$ is not innermost and $\Delta_\theta \cap (\annhigh \cup \annlow) \subset \Delta_\theta$
does contain other simple closed curves, we can achieve the assumption that $c$ is innermost by first
performing a surgery which is similar to that illustrated in \fullref{figure:annulus0} on all of the simple
closed curves in $\Delta_\theta$, starting with the innermost.  Such surgeries may or may not eliminate
any saddle singularities, but they will create center points in the foliation.  But, then we can perform
the surgery to eliminate the homoclinic saddle singularity on $c$. 

We have learned that we may assume that every leaf is an
arc.  The situation is a little bit more complicated than it was in the case of Seifert
surfaces.  Consult \fullref{figure:annulus1}.  There are now 5 possible types of
non-singular arcs in the foliation of $\annhigh$ and
$\annlow$:  (i) arcs which have both endpoints at vertices (type $b$); (ii) arcs which have
exactly one endpoint at a positive vertex  (type $a_+$); (iii) arcs which have exactly one endpoint at a
negative vertex (type $a_-$);  (iv) arcs which have one endpoint on one boundary
component and the other on the opposite boundary component (type $s$); and (v) arcs
which have both endpoints on the same boundary component.  As for type (v), the exact
argument that we used in the case of Seifert surfaces applies, because
$\Xhigh,\Xmed,\Xlow$ are all closed braids, so type (v) does not occur. 

\begin{remark}
With regard to \fullref{figure:annulus1} we may need to reinterpret the boundary components, temporarily,
as
$(\Xhigh,\Xmed)$ or $(\Xmed,\Xlow)$ instead of $(\Xhigh,\Xlow)$.  This ambiguity will be removed shortly.
See \fullref{remark:pm notation} below.\hfill $\diamondsuit$
\end{remark}

The annulus $\annhigh \cup \annlow$ is said to be {\it trivially foliated} or {\it trivial} if it is
foliated without clasp arcs and the leaves in its foliation are all $s$--arcs.

 An exhaustive list of the
singularities which could, in principle, occur in the foliation of $\annhigh$ and $\annlow$ are types 
$\ahigh\ahigh, \ahigh b, \
\ahigh s,
\ \ahigh\alow,\ \alow\alow, \ \alow b,  \ \alow s, \ sb, \ ss,$   and $bb$, where the
notation is consistent with that used in \fullref{section:preliminaries}.  However,
we have:

\begin{lemma}
\label{lemma:tile types on the annulus}
Singularities of type $ss$ do not occur in either $\annhigh$ or $\annlow$.
\end{lemma}
\begin{proof}   There is an oriented flow on $\annhigh \cup \annlow$ that is transverse to the foliation.  In particular,  we have arranged that the foliation is transverse to the boundary in small neighborhoods of $\Xhigh$ and $\Xlow$, so that the flow will be parallel to the boundary and oriented to agree with the orientations on $\Xhigh$ and $\Xlow$, both of which are oriented closed braids.  On the other hand, the fact that every singularity in the foliation is a saddle point implies that every singularity in the associated flow is hyperbolic. Suppose that there was a singularity of type $ss$ at polar angle $0$. Then, for an arbitrarily small negative polar angle $-\epsilon$ there must be a pair of type $s$ non-singular leaves $s_1, s_2$ moving toward one another on a collision course.   However, a type $s$ leaf separates the annulus and as we just noted the flow near the two endpoints of both $s_1$ and $s_2$ is in the direction of the orientation of $\Xhigh$ and $\Xlow$.  It is then impossible for the images of $s_1$ and $s_2$ in the flow to approach one-another, as the flow evolves. \end{proof}

\begin{lemma}
\label{lemma:not type alow alow or alow s}
If a singularity of type $\alow\alow$ or $\alow s$ occurs, then the associated singular
leaf is always intersected by a clasp arc.
\end{lemma}
\begin{proof}  Suppose that a singularity
of either type  $\alow\alow$ or $\alow s$ occurs, and that no clasp arc intersects the
associated singular leaf.  See \fullref{figure:annulus1a}.  
\begin{figure}[htpb!]
\labellist\small
\pinlabel {$\alpha$} [r] at 216 649
\pinlabel {$\beta$} at 369 640
\endlabellist
\centerline{\includegraphics[scale=.6]{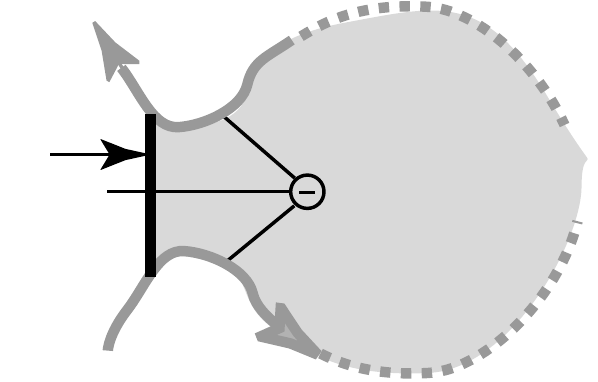}}
\caption{Type $a_-a_-$ and $a_-s$ singularities do not occur}
\label{figure:annulus1a}
\end{figure}
Notice that there is an arc,
$\alpha$, contained in the singular leaf which has both of its endpoints on  $\Xlow$ and which,
together with one of the two subarcs $\beta$ of
$\Xlow \setminus (\Xlow \cap \partial \alpha)$ forms a simple closed curve $C$ which 
bounds a disc $D \subset \annlow$.  Since $C$ bounds a disc it represents the unknot.
Now observe that if we perturb the  endpoints  of $\alpha$ slightly along 
$\Xlow$ we may change $\alpha$ to an arc  $\alpha^\prime$ which is transverse to fibers
of $\fib$.  Thus, after an   arbitrarily small perturbation, we may change $C$ to
$C^\prime = \alpha^\prime \cup \beta$, where $C^\prime$ is a closed braid and where  
$\beta \subset \Xlow$ is a subarc of the closed braid $\Xlow$.  Let $b$ be the 
braid index of $C^\prime$.
The fact that $C^\prime$ is  everywhere transverse to fibers of $\fib$ 
guarantees that $b \geq 1$.
But then, we may reduce the braid index of $\Xlow$ by at least one 
without changing its knot type by replacing $\beta$ by  $\alpha^\prime$.  However that
is  impossible because by hypothesis the braid index of $\Xlow$ is minimal.\end{proof}

\begin{figure}[htpb!]
\labellist\small
\pinlabel {$a_+a_+$} at 127 615
\pinlabel {$a_+b$} at 222 615
\pinlabel {$a_+s$} at 325 615
\pinlabel {$a_-a_+$} at 460 615
\pinlabel {$a_-a_-$} at 84 475
\pinlabel {$a_-b$} at 180 475
\pinlabel {$a_-s$} at 290 475
\pinlabel {$sb$} at 415 475
\pinlabel {$bb$} at 530 475
\endlabellist
\centerline{\includegraphics[scale=.7]{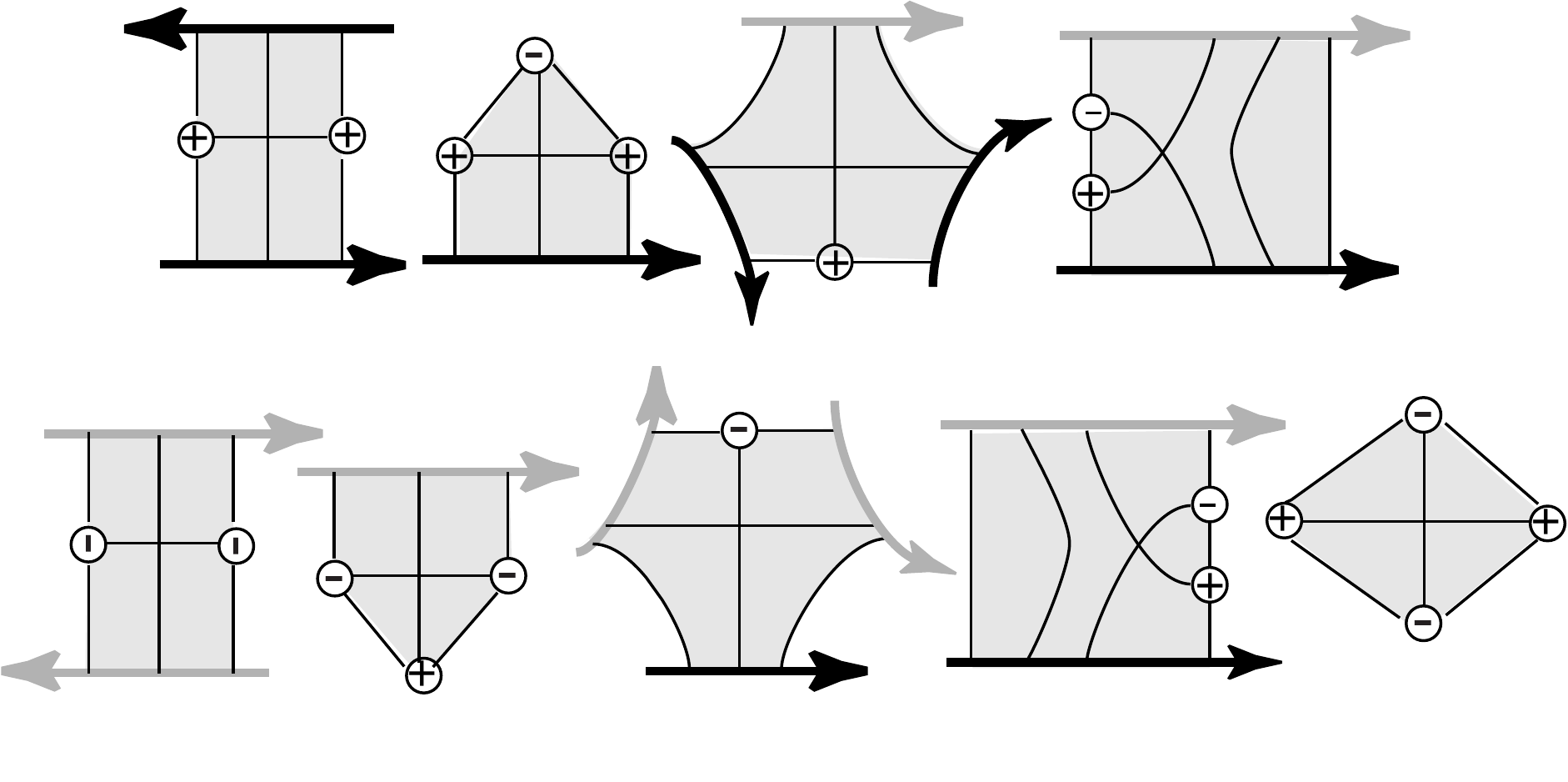}}
\caption{Possible tile types in the foliation of $\cA_+\cup\cA_-$.  The
black (resp.\  grey) arcs always
represent subarcs of $X_+$ (resp.\  $X_-$), with one exception: when we construct the tabs we
will be working in $\cA_-$ and $X_0$ will enter the picture as a (dotted) grey arc.}
\label{figure:annulus2}
\end{figure}
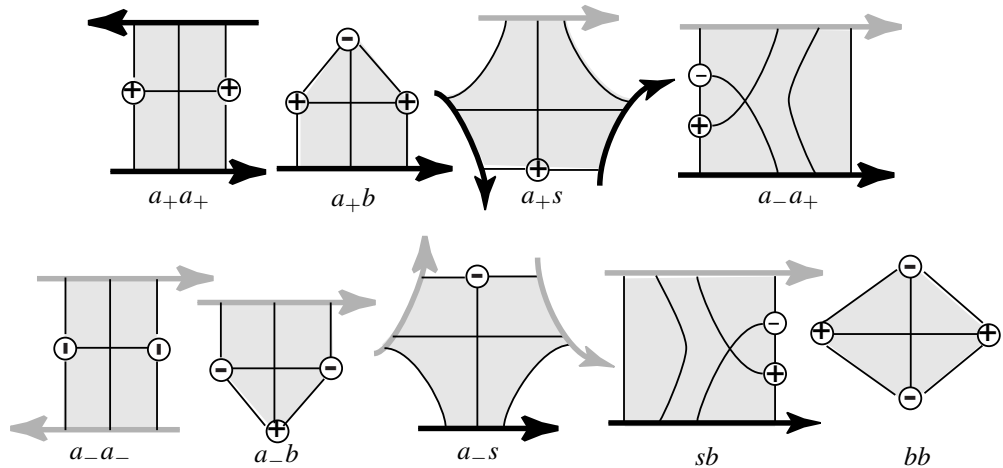

\subsection{Preliminary modifications in the clasp arcs}
\label{subsection:short clasp arcs}

We assume from now on that our clasp annulus $\ta$
supports a braid foliation, so that its preimage $\pa$ supports the lifted foliation.  We continue to
use the symbols $\ta$ and $\pa$, but from now on $\ta$ means the {\it foliated} clasp annulus and $\pa$
means its {\it foliated} preimage.  Each clasp arc $\gamma$ in $\ta$ will have two preimages
$\ghigh$ and $\glow$ in $\pa$, where $\ge$ indicates the preimage of $\gamma$ which begins on 
$X_\e$ in the boundary of $\pa$ and ends at its {\it puncture endpoint} on $X_{-\e}$ in the
interior of $\pa$.

 We say that a clasp arc is {\it short} if $\ghigh$ or $\glow$ does not
cross any singular leaves.  Since our ultimate goal is to push
$\Xhigh$ across
$\ta$ to
$\Xlow$, every clasp arc will ultimately become short.

\begin{lemma}
\label{lemma:eliminating short clasp arcs}
After a braid isotopy, we may assume that there are no short clasp arcs.  Even more, suppose that $\ge$ is a
short clasp arc in $\pa$.  Suppose further that the only non-singular leaves that $\ge$ intersects are
$s$--arcs and $a_\e$--arcs (without dividing the vertices of any $a_\e
a_\e$--singular leaves).  Then we may eliminate $\ge$. \end{lemma}

\begin{proof}  Since a short clasp arc does not cross any
singular leaves and has one endpoint on $\Xhigh$ it must be in the boundary of a region in $\pa$ which is near
$\Xhigh$ and foliated entirely by $a$--arcs and $s$--arcs (call it an
$a/s$--region).  See \fullref{figure:clasp1}.
\begin{figure}[htpb]
\labellist\small
\pinlabel {$e_+$} [t] at 101 356
\pinlabel {$p$} [t] at 254 356
\pinlabel {$e_-$} [b] at 254 405
\pinlabel {$e_-$} [b] at 404 401
\pinlabel {$e_+$} [t] at 404 369
\endlabellist
\centerline{\includegraphics[scale=.8]{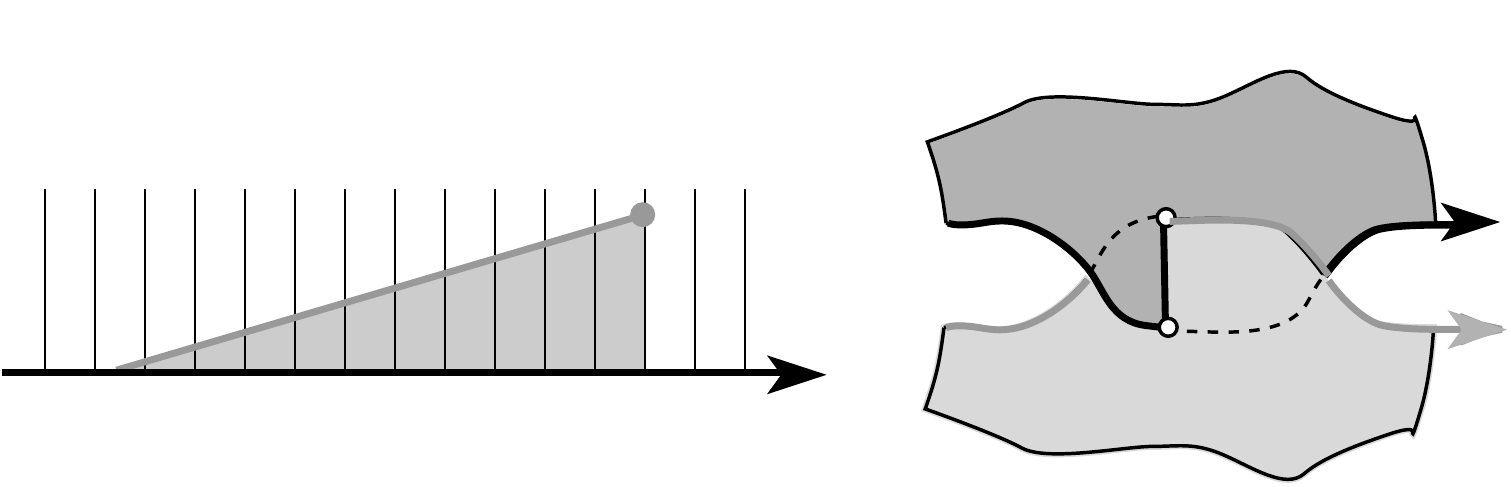}}
\caption {Eliminating short clasp arcs}
\label{figure:clasp1}
\end{figure}
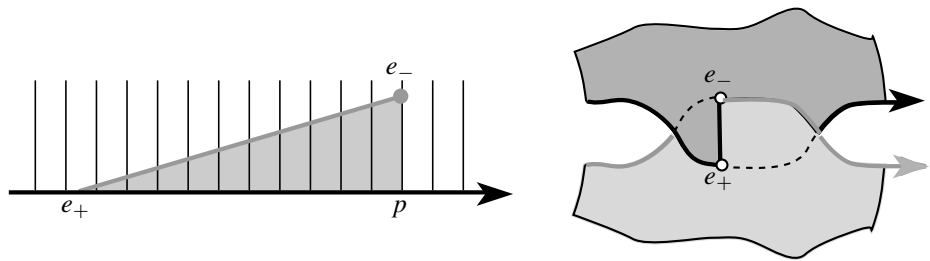
Let $e_+$ and $e_-$ be the endpoints of the clasp arc $\ghigh$, where $e_+\in\Xhigh$.  Then $e_-$ is on
an $a/s$--arc, and running along that arc to $\Xhigh$ we obtain a point $p\in\Xhigh$.  Modifying $e_-p$
slightly, we make it transverse to the fibers of $\fib$.  The `triangular'
region $e_+e_-p$, shaded in \fullref{figure:clasp1}, is  foliated entirely by subarcs of
$a/s$--arcs.  Rescale $\Xhigh$ in such a way as to shrink the shaded region, pulling $e_+$ to $p$
and $\ghigh$ into the arc $pe_-$. 

The clasp arc is still in an $a/s$--region,  as in the right sketch of \fullref{figure:clasp1}, which shows the
local picture in 3--space.  The black (resp.\  grey) boundary arc is
$\Xhigh$ (resp.\  $\Xlow$).  We now push
$\Xhigh$ along
$e_+e_-$ to eliminate the short clasp arc, changing the geometric linking
between $\Xhigh$ and $\Xlow$ in  the process.  Note that while it looks as if we are `unlinking' the two
curves, what we are really doing is to correct the geometric linking
between $\Xhigh$ (which started out by being geometrically unlinked from $\Xlow$) and $\Xlow$ by putting it into the
position of $\Xmed$.  \end{proof}  

\begin{remark}
By our basic construction, the only way that $\ta$ will fail to have clasp
arcs is if $\Xlow'$, the pushoff of $\Xlow$ onto a Seifert surface which was constructed during the proof
of \fullref{lemma:construction of the two annuli}, is geometrically unlinked from
$\Xlow$.  The only way that can happen is if $\cX$ is the unlink.  So, if $\cX$ is not the unlink,
then clasp arcs occur.  Short clasp arcs arise in the basic construction in situations where we could have
arranged for $\Xhigh$ to have the correct linking with $\Xlow$ locally by a braid isotopy.  If
it happened that every clasp arc was short, then we would know that, after a braid isotopy of $\Xhigh$,
it could be assumed that $\Xlow$ is a preferred longitude for $\Xhigh$,  ready to be pushed across
$\annlow$, an embedded annulus. \hfill$\diamondsuit$
\end{remark}

\subsection{Construction of the tabs}
\label{subsection:construction of the tabs}
 In this section we modify the discs
$R_1,\dots,R_k$ which we constructed in \fullref{section:getting started} to special foliated discs $\Thigh^1,\dots,\Thigh^k$ containing the clasp arcs. 

We say that
$\Thigh^i\subset \pa$ is a {\it tab} associated to the clasp arc $\ghigh^i$ if
the following hold:

\begin{itemize}
\item $\ghigh^i\subset\Thigh^i$.
\item $\partial \Thigh^i = \alpha \cup \beta$
where $\alpha \subset \Xhigh$ and $\beta$ is an arc that is transverse
to the foliation of $\pa$.
\item There is a simple path $l_+$ which is contained in a union of singular leaves in $\Thigh^i$.  The path
$l_+$ starts on $\Xhigh$, ends at a negative vertex, and contains all the negative vertices in $\Thigh^i$.

\item The arc $ \ghigh^i$
is the only clasp arc which intersects $\Thigh^i$.
It is everywhere transverse to the leaves in the foliation of $\pa$.  The arc
$\ghigh^i$ intersects each of the singular leaves in the induced
foliation of $\Thigh^i$ exactly once.
\end{itemize}
\index{tab} 

\begin{figure}[htpb!]
\labellist\small
\pinlabel {$e_-$} [b] at 135 475
\pinlabel {$X_0$} [b] at 187 479
\pinlabel {$\gamma_+^i$} [b] at 288 464
\pinlabel {$\delta_4$} [bl] at 195 436
\pinlabel {$\delta_2$} [bl] at 352 436
\pinlabel {$\delta_1$} [bl] at 409 436
\pinlabel {$e_+$} [b] at 453 439
\endlabellist
\centerline{\includegraphics[scale=.9]{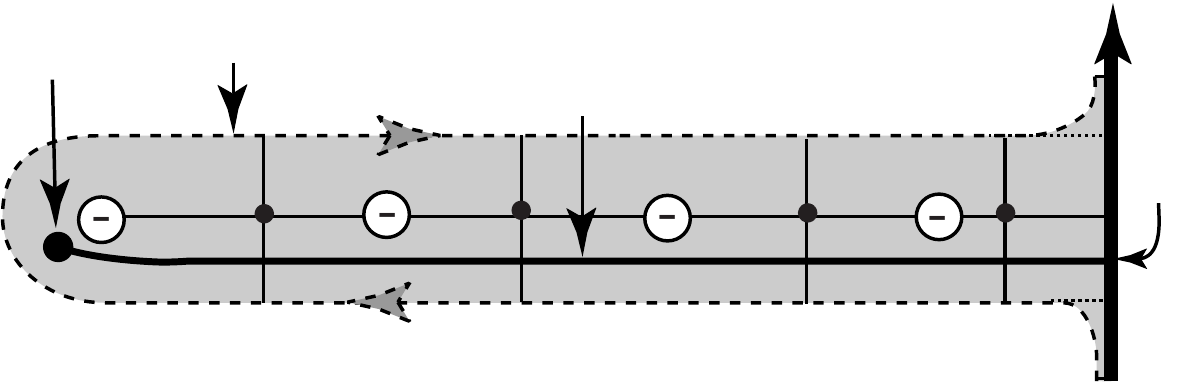}}
\caption{An example of a tab}
\label{figure:clasp2}
\end{figure}
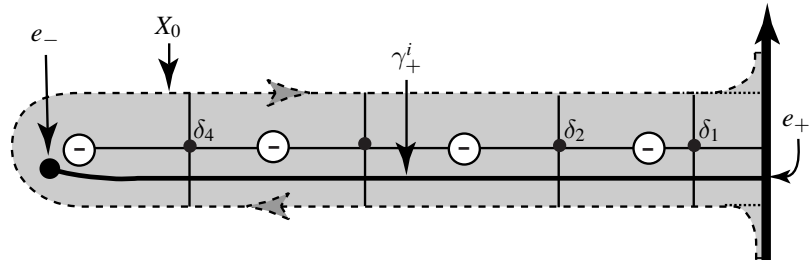

\begin{lemma}
\label{lemma:straightening out the tabs}
We may assume that each $R_i$ is a tab $\Thigh^i$.
\end{lemma}

\begin{proof}  We focus on a single $R_i\subset\pa$, where $R_i$ is one of the discs that we constructed in 
\fullref{section:getting started}.  We construct the tab $\Thigh^i$ as a subset of the
disc $R_i$ in the basic construction.
Going back to \fullref{section:preliminaries}, recall that the boundary of $R_i$ is a union of 4 arcs: 
\bi
\item a subarc of
the modified $\Xmed$, 
\item a subarc of
$\Xhigh$, and 
\item two modified $s$--arcs which join them.
\ei
It contains 
$\ghigh^i$ in its interior.  Using the methods described in \fullref{section:preliminaries}, we
may assume that $R_i$ supports a braid foliation.  Initially, the foliated disc $R_i$ will not look
anything like a tab.  What we wish to do is to rechoose the discs $R_1,\dots,R_k$  so that they are as simple
as possible.  The new discs that we choose will be subsets of the old ones, with possibly modified foliations.
In the course of the construction we will, in effect, be making a new choice of the separating curve
$\Xmed$ which divides $\pa$ into the embedded annuli
$\annlow$ and $\annhigh$.  We note that right after the completion of this proof we will discard $\Xmed$, as
we will have no further use for it.  

The arc $\ghigh^i$ has one of its endpoints at a point $e_+^i$ on $\Xhigh$ and its other
endpoint $e_-^i\in\Xlow$ in the interior of $R_i$.  In between, it winds around the
vertices of the foliated disc $R_i$, constrained to remain transverse to the leaves of the 
foliation.  We may assume without loss of generality that the point $e_-^i$ where $\Xlow$
pierces  the interior of $R_i$ is not on a singular leaf.  Therefore $e_-^i$ belongs to a
non-singular  leaf  of type $s,a_+,a_-$ or $b$. 

Let $\alpha \subset R_i \setminus \ghigh^i$ be a simple path that starts at 
$e_-^i$ and ends at a point $p \in \Xhigh$.  Then
$\ghigh^i \cup \alpha \subset R_i$
is a simple path that joins $e_+^i$ to $p$ in the interior of $R_i$.  Without
loss of generality we may assume that this simple path is in braid position.  For, if it is not, apply
Alexander's braiding trick to wrongly oriented subarcs.  This will change the interior of $R_i$ by an isotopy
in 3--space.  The changes will modify the foliation of
$R_i$ by introducing new vertices and singularities which allow $\alpha$ to avoid the points of
non-transversality.  After these modifications, let
$R_i' \subset R_i$ be the foliated subdisc that 
$\ghigh^i \cup \alpha$ splits off in $R_i$.  Reapply the argument for eliminating leaves that are circles in
the foliation of $R_i'$, so that $R_i'$ supports a braid foliation. Therefore we may assume that $R_i'$ supports a
braid foliation.  

We are now in position to re-choose $\ghigh^i,\alpha$ and $\Xmed$  so that every vertex in $R_i'$
is negative.  See \fullref{figure:clasp3}.
\begin{figure}
\labellist\small
\pinlabel {(a)} at 163 523
\pinlabel {(b)} at 125 410
\pinlabel {(c)} at 346 410
\pinlabel {$v$} [r] at 209 482
\pinlabel {$\gamma_+^i$} [b] at 376 442
\pinlabel {$e_+$} [b] at 467 365
\endlabellist
\centerline{\includegraphics[scale=.75]{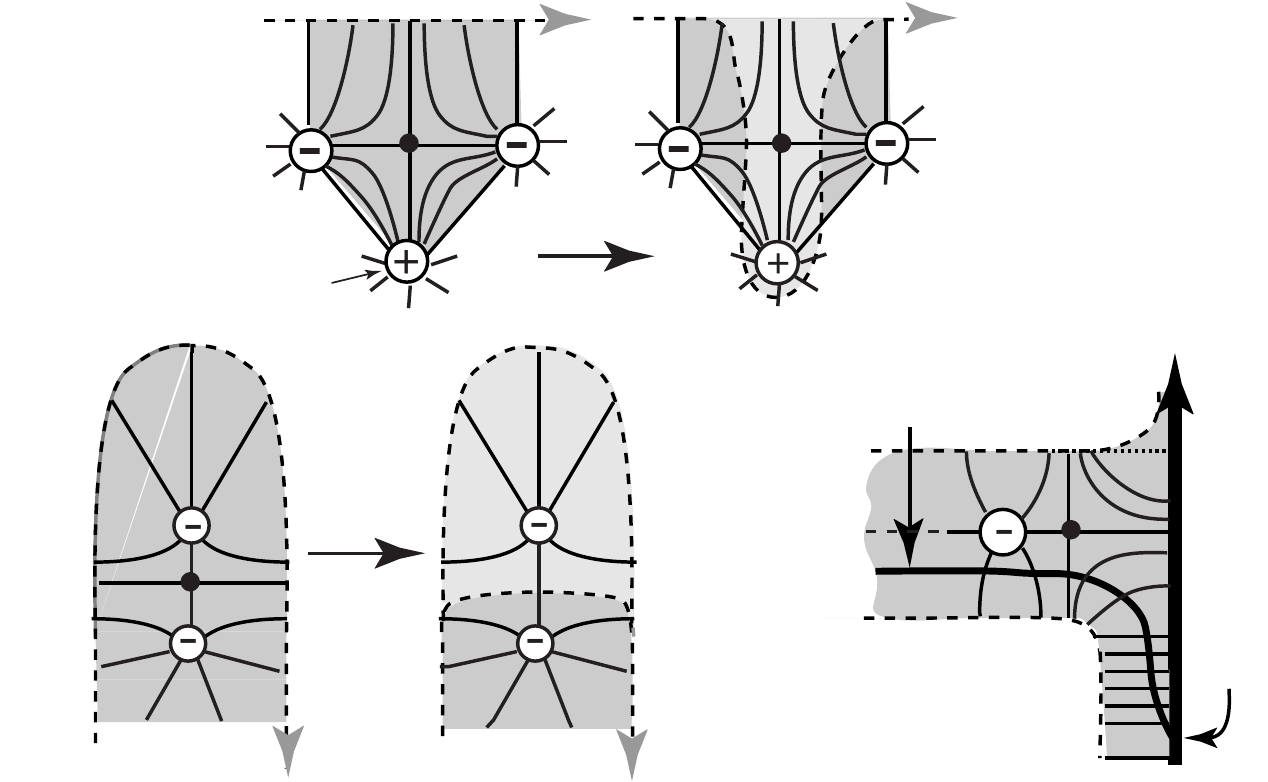}}
\caption{Rechoosing the disc $R_i$. In sketches (a) and (b) the dark area is in $R_i^\prime$. The passage from left to right shows how we modify $R_i^\prime$ to eliminate positive vertices. }
\label{figure:clasp3}
\end{figure}
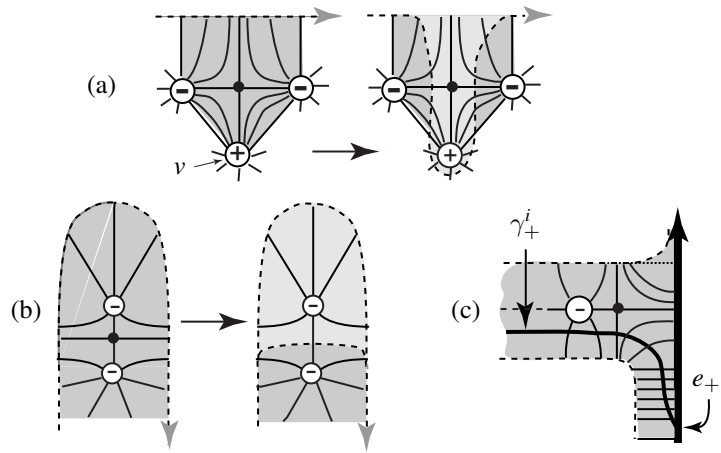 
If the induced foliation of $R_i'$ contains a positive vertex $v$ then this
sub-foliation  must also contain a singular
leaf that has its endpoint on $\ghigh^i \cup \alpha$ and is adjacent to a (possibly different) positive
vertex $v'$.   We can then push
$\ghigh^i \cup \alpha$ along the singular leaf and across $v$, staying transverse to the foliation, as in
\fullref{figure:clasp3}(a), moving $v$ out of $R_i^\prime$.  Inducting on the
number of positive vertices in the foliation of $R_i'$, we have arranged that the foliation of $R_i'$ contains
only negative vertices. 

Now the graph of singular leaves in $R_i'$ is either a linear tree, or it is a tree with
branches.  In the latter case, since
$e_-^i$ can be adjacent to only one negative vertex, either  $\ghigh^i$ or $\alpha$ intersects
a singular leaf in the foliation of $R_i$ twice.  We can then find
a  valence one vertex, as in \fullref{figure:clasp3}(b), with either 
$\ghigh^i$ or $\alpha$ in its boundary.
After a `destabilization'  we can eliminate this negative vertex from  $R_i$.  Iterating this
procedure, we alter  $\ghigh^i \cup \alpha$, until $R_i'$ is changed to a tab. 

There is only one more problem.  It may happen that near the $\Xhigh$ boundary of $R_i'$, either $\ghigh^i$ or
$\alpha$ or both run along a band which is foliated without singularities
by $a$ and/or $s$--arcs, as
depicted in \fullref{figure:clasp3}(c) in the case of $\ghigh^i$.  If so, we simply `rescale' $\Xhigh$,
as we did in the proof of \fullref{lemma:eliminating short clasp arcs} to
pull the long arcs back into the single $\alow s$--tile which
forms the base of the tab.  \end{proof}

\begin{remark} 
\label{remark:pm notation}
We make an important remark about simplified notation and simplified tile types.   
Having standardized the tabs, we will not have further use for
$\Xmed$.  This eliminates the dual meanings of the boundary arcs in \fullref{figure:annulus2} and 
allows the following simplification. If $\e = \pm$, an $a_\e$--arc has its interior endpoint
on a vertex of sign $\epsilon$ and has its boundary endpoint on $X_{\epsilon}$.  At the same time, we will no
longer need to depict $\Xmed$ in our figures, and will be free to use dotted arcs in other settings as we
proceed through the proof of the MTWS. \hfill $\diamondsuit$
\end{remark}

\begin{remark}
Notice that the work in this section uses the braid foliation of $R_i$.  If we had
braid foliations on hand as a tool, during the basic construction, we could have arranged at that time for
each
$R_i$ to be a tab.\hfill $\diamondsuit$
\end{remark}

\begin{remark}
In \fullref{figure:getstarted4} we gave an example of a tab with the clasp arc on one side of the singular
leaves.  Soon we will develop the ability to move it so that it is transverse to the foliation and to one side or the
other of a sequence of singular leaves, or alternatively to position it in a union of leaves (singular and/or
non-singular), at the expense of introducing many new $b$--arcs into the foliation.  See
\fullref{subsection:normal} below. \hfill $\diamondsuit$
\end{remark}

\subsection{The two finger moves}
\label{subsection:the finger moves}  
We need tools that will allow us to modify neighborhoods of the clasp arcs in $\ta$ in
a controlled manner, keeping track of the foliation on the two branches and making sure that
no new self-intersections are introduced.  The `finger moves' will help us to do that. 
See Figures \ref{figure:fingermoves1} (and eventually \ref{figure:fingermoves2}).
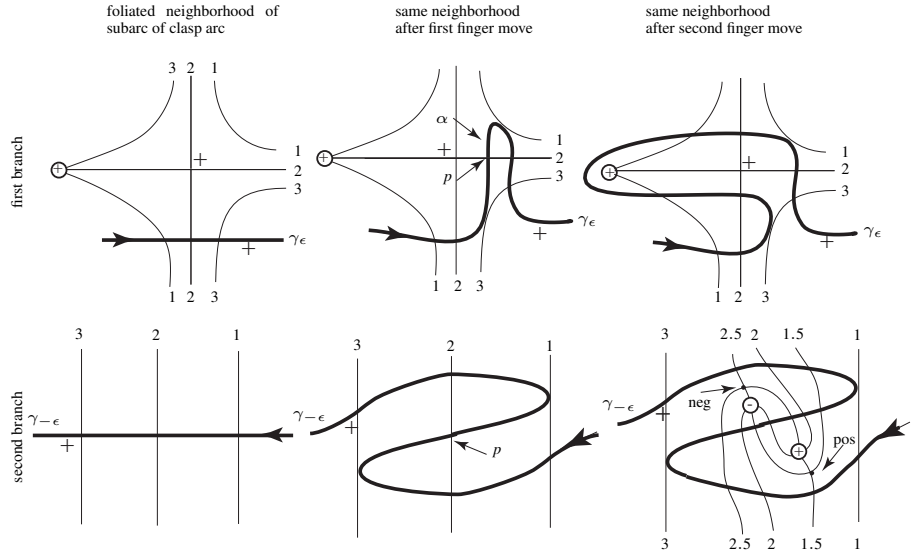
\begin{figure}[htpb]
\labellist\tiny
\pinlabel {\parbox{65pt}{foliated neighborhood of subarc of clasp arc}}
  [b] at 180 600
\pinlabel {\parbox{60pt}{same neighborhood\\after first finger move}}
  [b] at 334 600
\pinlabel {\parbox{70pt}{same neighborhood\\after second finger move}}
  [b] at 478 600
\pinlabel {\rotatebox{90}{first branch}} at 84 530
\pinlabel {\rotatebox{90}{second branch}} at 84 385

\pinlabel {$1$} [b] at 192 580
\pinlabel {$2$} [b] at 179 580
\pinlabel {$3$} [b] at 169 580
\pinlabel {$1$} [t] at 169 469
\pinlabel {$2$} [t] at 179 469
\pinlabel {$3$} [t] at 192 469
\pinlabel {$1$} [l] at 231 542
\pinlabel {$2$} [l] at 231 531
\pinlabel {$3$} [l] at 231 521
\pinlabel {$\gamma_\epsilon$} [l] at 228 493

\pinlabel {$1$} [t] at 313 475
\pinlabel {$2$} [t] at 325 475
\pinlabel {$3$} [t] at 337 475
\pinlabel {$1$} [l] at 374 548
\pinlabel {$2$} [l] at 374 538
\pinlabel {$3$} [l] at 374 528
\pinlabel {$\gamma_\epsilon$} [l] at 386 503
\pinlabel {$p$} [r] at 326 526
\pinlabel {$\alpha$} [r] at 324 559

\pinlabel {$1$} [t] at 467 469
\pinlabel {$2$} [t] at 478 469
\pinlabel {$3$} [t] at 493 469
\pinlabel {$1$} [l] at 530 541
\pinlabel {$2$} [l] at 530 532
\pinlabel {$3$} [l] at 530 520
\pinlabel {$\gamma_\epsilon$} [l] at 542 497

\pinlabel {$3$} [b] at 118 434
\pinlabel {$2$} [b] at 160 434
\pinlabel {$1$} [b] at 204 434
\pinlabel {$\gamma_{-\epsilon}$} [b] at 99 386

\pinlabel {$3$} [b] at 270 428
\pinlabel {$2$} [b] at 321 428
\pinlabel {$1$} [b] at 375 428
\pinlabel {$\gamma_{-\epsilon}$} [b] at 244 388
\pinlabel {$p$} [l] at 340 377

\pinlabel {$1$} [b] at 543 434
\pinlabel {$1.5$} [b] at 508 434
\pinlabel {$2$} [b] at 487 434
\pinlabel {$2.5$} [b] at 474 434
\pinlabel {$3$} [b] at 438 434
\pinlabel {$1$} [t] at 543 334
\pinlabel {$1.5$} [t] at 518 334
\pinlabel {$2$} [t] at 496 334
\pinlabel {$2.5$} [t] at 479 334
\pinlabel {$3$} [t] at 438 334
\pinlabel {neg} [t] at 457 410
\pinlabel {pos} [b] at 536 375
\pinlabel {$\gamma_{-\epsilon}$} [b] at 413 394

\endlabellist
\centerline{\includegraphics[scale=.70]{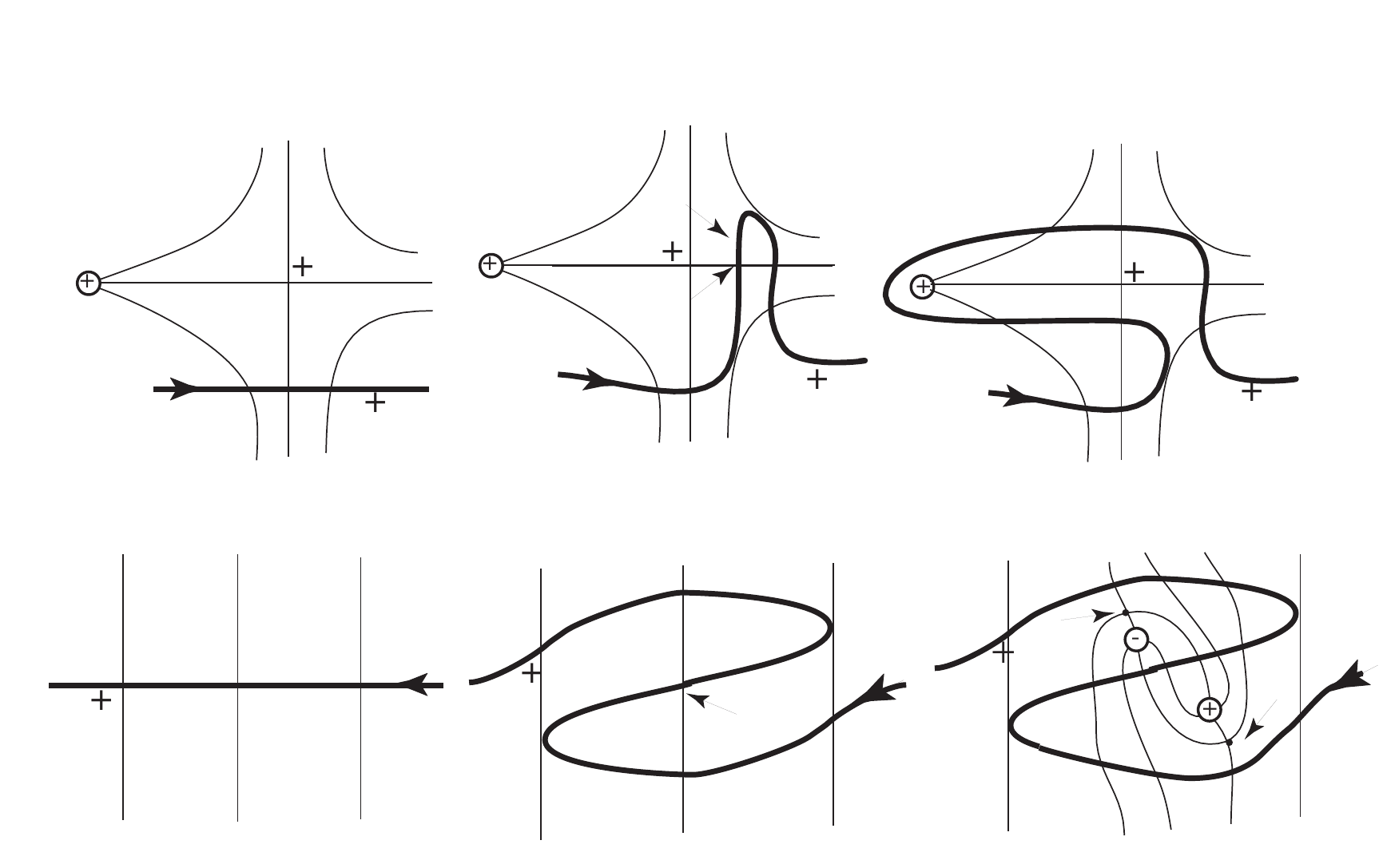}}
\caption{Controlled local changes in the clasp arc in a neighborhood of a
singularity after the first and second finger moves}
\label{figure:fingermoves1}
\end{figure} 
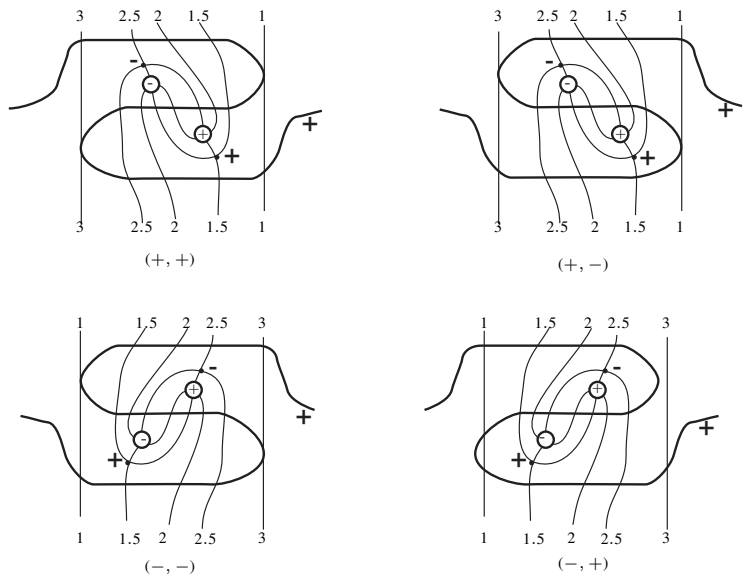
\begin{figure}[htpb!]
\labellist\tiny
\pinlabel{$(+,+)$} at 210 568
\pinlabel{$(+,-)$} at 385 566
\pinlabel{$(-,-)$} at 210 438
\pinlabel{$(-,+)$} at 385 439

\pinlabel {$1$} [b] at 249 666
\pinlabel {$1.5$} [b] at 222 666
\pinlabel {$2$} [b] at 205 666
\pinlabel {$2.5$} [b] at 193 666
\pinlabel {$3$} [b] at 172 666
\pinlabel {$1$} [t] at 249 588
\pinlabel {$1.5$} [t] at 230 588
\pinlabel {$2$} [t] at 212 588
\pinlabel {$2.5$} [t] at 198 588
\pinlabel {$3$} [t] at 172 588

\pinlabel {$1$} [b] at 426 666
\pinlabel {$1.5$} [b] at 399 666
\pinlabel {$2$} [b] at 381 666
\pinlabel {$2.5$} [b] at 370 666
\pinlabel {$3$} [b] at 349 666
\pinlabel {$1$} [t] at 426 588
\pinlabel {$1.5$} [t] at 407 588
\pinlabel {$2$} [t] at 390 588
\pinlabel {$2.5$} [t] at 374 588
\pinlabel {$3$} [t] at 349 588

\pinlabel {$3$} [b] at 249 536
\pinlabel {$2.5$} [b] at 230 536
\pinlabel {$2$} [b] at 217 536
\pinlabel {$1.5$} [b] at 200 536
\pinlabel {$1$} [b] at 172 536
\pinlabel {$3$} [t] at 249 456
\pinlabel {$2.5$} [t] at 225 456
\pinlabel {$2$} [t] at 207 456
\pinlabel {$1.5$} [t] at 193 456
\pinlabel {$1$} [t] at 172 456

\pinlabel {$3$} [b] at 420 536
\pinlabel {$2.5$} [b] at 400 536
\pinlabel {$2$} [b] at 387 536
\pinlabel {$1.5$} [b] at 369 536
\pinlabel {$1$} [b] at 343 536
\pinlabel {$3$} [t] at 420 456
\pinlabel {$2.5$} [t] at 397 456
\pinlabel {$2$} [t] at 380 456
\pinlabel {$1.5$} [t] at 364 456
\pinlabel {$1$} [t] at 343 456

\endlabellist
\centerline{\includegraphics[scale=.90]{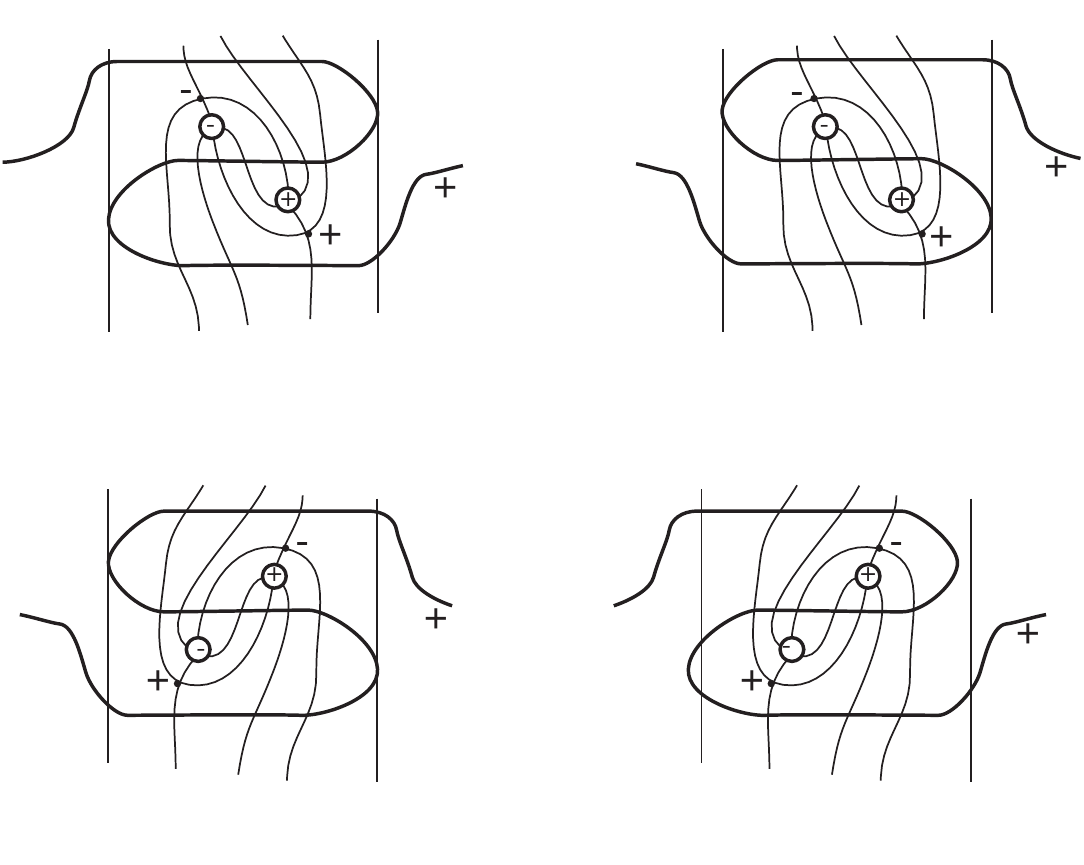}}
\caption{The 4 possible
choices of the signs of the pair (vertex,singularity)}
\label{figure:fingermoves2}
\end{figure} 
To begin, we explain the pictures in  \fullref{figure:fingermoves1} The first column 
shows foliated neighborhoods  $N$ of a subarc of one of the clasp arcs.  We give 
separate pictures  of the two foliated branches, $N_1$ (the {\it first branch}) and 
$N_2$ (the {\it  second branch}).  They intersect transversally.
By hypothesis the clasp arc (and so also its preimages $\ghigh$ and $\glow$) is
transverse to the leaves of $\fib$.  We have oriented the clasp arc (arbitrarily).  

We are interested in modifying the position of one of the clasp arcs in a
neighborhood of a singular leaf on one of two branches $N_1$ or $N_2$.  Since
singular leaves correspond to places where one of the two branches is tangent to a
fiber of $\fib$, and since the two branches intersect transversally along a clasp
arc, we may assume that there is a neighborhood of the singular point in which the
other branch is foliated without singularities.  The neighborhood has been chosen so
that $N_2$ is foliated  without singularities, but $N_1$ contains a singularity,
together with the vertex endpoint of one of the branches of the singular leaf.  There
are two sign choices: the sign of the vertex and the sign of the singularity.  We have
chosen these to be $(+,+)$, but in a moment we will consider the 4 possible sign
choices.  There are also little `+' signs next to the two components of the clasp
arcs.  They indicate which side of $N_2$ (resp.\  $N_1$) is the positive side, at
$\ge$ (resp.\  $\gamma_{-\e})$. 

Selected fibers $H_{\theta_i}$ have labels
$i=1,2,3$.  The labels on the fibers and the way that $\gamma$
intersects them and the little + signs next to $\ge$ and $\gamma_{-\e}$ completely
determine the position of $N$ in 3--space, relative to the 
coordinate system provided by the fibers of $\fib$.
We are ready to describe our two controlled changes in
the clasp arc, and the corresponding changes in the foliation of $N_1$ and $N_2$
and in the way that $\ge$ and $\gamma_{-\e}$ intersect the leaves in the
foliation.

We now define  the moves precisely:
\begin{enumerate}
\item Our {\it first finger move}, 
illustrated in the middle column of
\fullref{figure:fingermoves1}, pushes $\gamma$ across the horizontal branch of
the singular leaf in the first branch, creating two points where it is not 
transverse to
the fibers of $\fib$, one on fiber 1 and the other on fiber 3.
This move is always possible, because it occurs within an arbitrarily small 
neighborhood of
the singularity.  The corresponding change in the second branch can be understood by
noticing that before the change $\ge$ intersected fibers 1,2,3
transversally in that order.  After the change there are two points of 
tangency with fibers of
$\fib$, the first with fiber 3 and the second with fiber 1.  This explains the 
doubling back of $\gamma_{-\e}$ after the 
first finger move. 

\item Let
$\alpha$ denote the subarc of the clasp arc which is between the two points of
non-transversality, and let
$p\in\alpha\cap{\rm fiber \ 2}$ be the point which is closest to the singularity on
$\ge$.  The {\it second finger move}, illustrated in the right column of \fullref{figure:fingermoves1}, pushes a neighborhood
$N(p)$ of $p$ on $\alpha$ across the
 singularity and across the vertex, staying within a neighborhood of the
singular leaf.  The foliation on the first branch is unchanged.  The foliation of
the second branch changes in an arbitrarily small neighborhood of $p$ on the
second branch.  Two new vertices of opposite parity and two new singularities of
opposite parity are created, as illustrated in the bottom row sketch.
\end{enumerate}
See \fullref{figure:fingermoves2} for the local changes with
the four possible sign choices for the pair (vertex,singularity).

\begin{remark}
The two finger moves are always
possible because of their local nature and because of our control over the
geometry.  We have illustrated the case when the singularity and vertex are
both positive.  The other three cases differ from this one by local
symmetries.  Observe that we have given a great deal of
detailed data in Figures \ref{figure:fingermoves1} and \ref{figure:fingermoves2}, including the cyclic order of
leaves in the fibration, the signs of the singularities and the signs of the vertices.  Fibers of $\fib$ are level
sets for positioning the two branches of the immersed surface in 3--space.  The test for whether the finger moves
are realizable in 3--space is to examine them on a sequence of fibers of $\fib$, and the data in Figures
\ref{figure:fingermoves1} and  \ref{figure:fingermoves2} suffices for that purpose. 

If  braid foliations had been available as a tool during the basic construction, we could have
used the two finger moves then.  We stress this because later we will use an inductive argument and we
need to know that, after many changes, we are still in the situation of
the basic construction. \hfill $\diamondsuit$
\end{remark} 

\subsection{Creating symmetric normal neighborhoods of the clasp arcs}
\label{subsection:normal}

In \fullref{subsection:construction of the tabs} we modified the discs $R_1,\dots,R_k$ which had been
constructed as part of the basic construction to very special foliated discs
which we called tabs.  Thus we now know that each clasp arc $\ghigh^i$ is contained in a tab.  In this section
we study $\ghigh^i$ and $\glow^i$ and arrange that they are contained in much larger
neighborhoods which support a canonical foliation.  We will call them normal neighborhoods.  The modifications
will be made with the help of the finger moves of \fullref{subsection:the finger moves}, at
the expense of adding new $b$--arcs (see \fullref{subsection:recognizing destabilizations and
exchange moves from the foliation}).  Intuitively, the new $b$--arcs add pouches to the surface, and so give additional `room to move around'.  The normal
neighborhoods which we will construct will give us choices, which can be made one way or another as it is
convenient.  For example we will be able to regard a clasp arcs as being on either side of its associated chain of
singular leaves and in braid position, or as being contained  in a union of leaves and so lying in a union of
fibers of $\fib$. 
\begin{figure}[htpb]
\labellist\small
\pinlabel {$\varphi_1$} [t] at 137 310
\pinlabel {$\varphi_2$} [l] at 219 409
\pinlabel {$\varphi_3$} [b] at 174 487
\pinlabel {$\varphi_4$} [r] at 50 403
\pinlabel {$\gamma_+$} at 136 450
\pinlabel {$\gamma_-$} at 328 483
\pinlabel {$\gamma_-$} at 504 402

\tiny
\hair=2pt
\pinlabel {$\delta_k$} [r] at 74 436
\pinlabel {$-\delta_k$} [tr] at 114 464
\pinlabel {$\delta_k$} [tl] at 158 436
\pinlabel {$-\delta_k$} [tr] at 199 461
\pinlabel {$\delta_1$} [r] at 74 365
\pinlabel {$-\delta_1$} [tl] at 113 392
\pinlabel {$\delta_1$} [tl] at 158 365
\pinlabel {$-\delta_1$} [tr] at 199 392

\pinlabel {$-\delta_k$} [r] at 265 497
\pinlabel {$\delta_k$} [tr] at 307 524
\pinlabel {$\delta_{k-1}$} [br] at 307 466
\pinlabel {$-\delta_k$} [tl] at 351 497
\pinlabel {$\delta_k$} [tr] at 390 524
\pinlabel {$\delta_{k-1}$} [br] at 390 466
\pinlabel {$-\delta_1$} [r] at 265 364
\pinlabel {$\delta_1$} [tr] at 308 392
\pinlabel {$-\delta_1$} [tr] at 351 364
\pinlabel {$\delta_1$} [tr] at 389 392

\pinlabel {$-\delta_1$} [r] at 443 362
\pinlabel {$\delta_1$} [tr] at 483 389
\pinlabel {$-\delta_1$} [tr] at 528 362
\pinlabel {$\delta_1$} [tr] at 567 389
\pinlabel {$-\delta_2$} [r] at 443 415
\pinlabel {$-\delta_2$} [tr] at 528 415

\endlabellist
\centerline{\includegraphics[scale=.7]{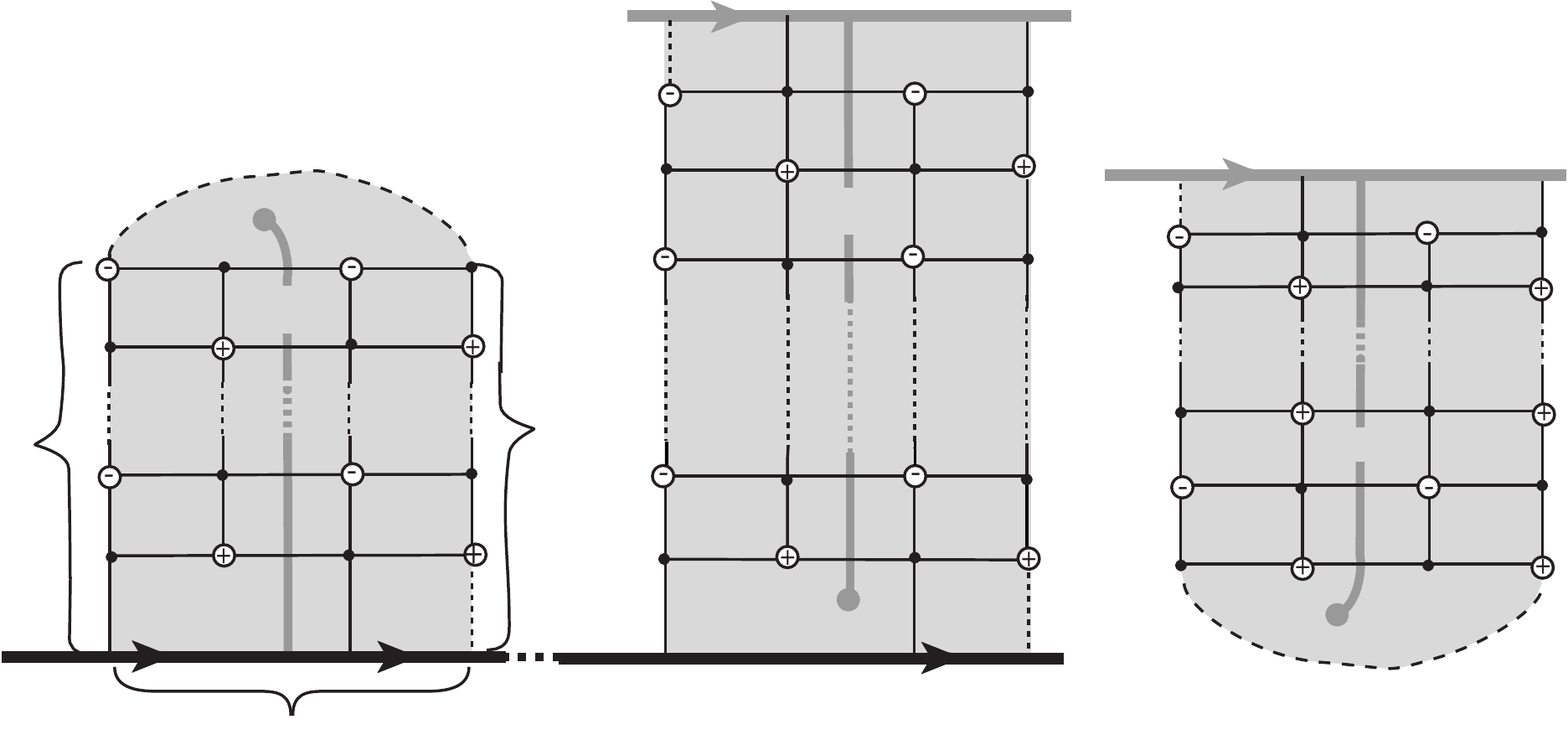}}
\caption{Normal neighborhoods of $\glow$ and $\ghigh$ in $\pa$.  The left
(resp.\  middle and right) sketch
shows
$\Nhigh$ (resp.\  $\Nlow,\Nlow$).  The arc $\glow$ ends near $\Xhigh$ in the
middle sketch and in the interior of $\pa$ in  the right sketch.  In all 3 sketches the foliation is defined up
to reflection about a vertical  axis which joins $\Xlow$ and $\Xhigh$.}
\label{figure:normal2}
\end{figure}
Let $\ge \subset
 \pa$ be a clasp arc with the property that $\ghigh$ crosses $k$ singular leaves.  A
neighborhood $\Nboth$ of $\ge$ is a 
{\it normal neighborhood} (\fullref{figure:normal2}, which gives an example, should be helpful) 
if the following conditions hold:

\be
\item $\Nboth\supset\ge$.  Also $\Nboth$  intersects no other clasp arcs.

\item  $\Nboth$ contains $2k$ disjoint `horizontal' paths $\varrho_1^\e, \dots , \varrho_{2k}^\e$, each
contained in the singular leaves of $\Nboth$ and each containing two vertices and two singularities.
Traversing
$\ge$, starting at the $X_\e$ endpoint, each $\varrho_i^\e$ is crossed once.  The vertices on
$\varrho_{2i-1}^\e$ (resp.\  $\varrho_{2i}^\e$) have sign $\e$
(resp.\  $-\e$).
For $1\leq i \leq k$, the singularities on $\varrho_{2i-1}^\e$
have parity $\e \delta_i$ and the singularities on $\varrho_{2i}^\e$ have parity $-\e \delta_i$, where
$\delta_i=\pm$.

\item  $\partial \Nboth = \varphi_1 \cup \varphi_2 \cup \varphi_3
\cup \varphi_4$ where:
\begin{enumerate}
\item $\varphi_1 \subset X_\epsilon$.
\item $\varphi_2$ is a path contained in one arc of type $a_{\epsilon}$ and
$k$ singular leaves.  It contains $k$ vertices of sign $\epsilon$.
\item $\varphi_3$ is transverse to the foliation of $\pa$.
\item $\varphi_4$ is a path contained in $k$ singular leaves and one arc of type
$b$ (or type $a_{\epsilon}$, in the special case when the puncture
point on $\ge$ is near $X_{-\epsilon}$).  It contains $k$ vertices of sign $-\epsilon$.
\item Traversing $\varphi_1 \subset X_\epsilon \cap \partial \Nboth$,
 beginning at the point $\varphi_1 \cap \varphi_4$, we pass 
through the endpoint of a type $a_\epsilon b$ singular leaf, and end at the
$X_\e$ endpoint of an $a_\epsilon$ leaf.
\end{enumerate} 

\ee

We say that $(\Nhigh,\Nlow)$ is a {\it normal neighborhood pair}
if each $\Nboth$ is a normal neighborhood of $\ge$.  Our main result in this section will be very important in the
detailed work we will need to do to push $\Xhigh$ across $\ta$ to $\Xlow$, and to prove the MTWS.
We stress this by
calling in a `proposition'.

\begin{proposition}
\label{proposition:normal}
For each clasp arc pair $(\ghigh,\glow) =
(\ghigh^i , \glow^i)$ we may assume that there is a normal neighborhood pair
$(\Nhigh,\Nlow)$.  Moreover, within the normal neighborhoods, we may assume that instead of being
transverse to the foliation each clasp arc is positioned in a finite union of leaves in its normal neighborhood
pair.
\end{proposition}

\begin{proof} The proof of \fullref{proposition:normal}  will occupy the rest of
this section.  Our work begins with the tab neighborhoods of $\ghigh^1,\dots,\ghigh^k$ which we 
constructed earlier,
when we simplified the discs $R_1,\dots,R_k$ of the basic construction.  The tab neighborhood 
construction told
us nothing about the other sheet of $\ta$, and our first goal is to modify it so that we have related tab
neighborhoods of $\glow^1,\dots,\glow^k$.  See \fullref{figure:normal1} for an example.  After that we will
iterate the construction to produce normal neighborhoods.

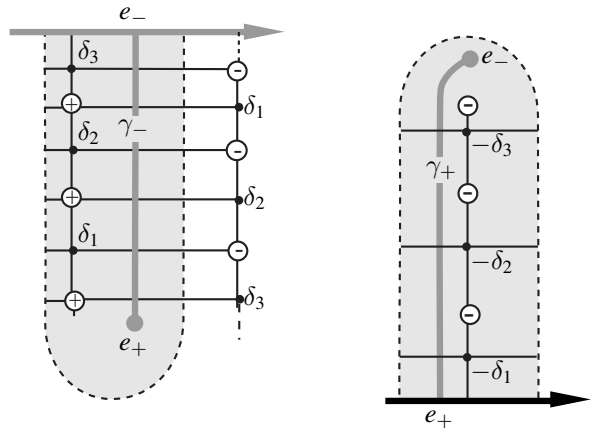
\begin{figure}[htpb]
\labellist\small
\pinlabel {$\delta_3$} [bl] at 188 464
\pinlabel {$\delta_2$} [bl] at 188 426
\pinlabel {$\delta_1$} [bl] at 188 378
\pinlabel {$e_-$} [b] at 217 481
\pinlabel {$e_+$} [t] at 217 341
\pinlabel {$\delta_1$} [l] at 265 446
\pinlabel {$\delta_2$} [l] at 265 403
\pinlabel {$\delta_3$} [l] at 265 356
\pinlabel {$-\delta_1$} [tl] at 373 329
\pinlabel {$-\delta_2$} [tl] at 373 381
\pinlabel {$-\delta_3$} [tl] at 373 436
\pinlabel {$e_-$} [l] at 377 468
\pinlabel {$e_+$} [t] at 361 307
\pinlabel {$\gamma_-$} at 217 436
\pinlabel {$\gamma_+$} at 362 416
\endlabellist
\centerline{\includegraphics[scale=.8]{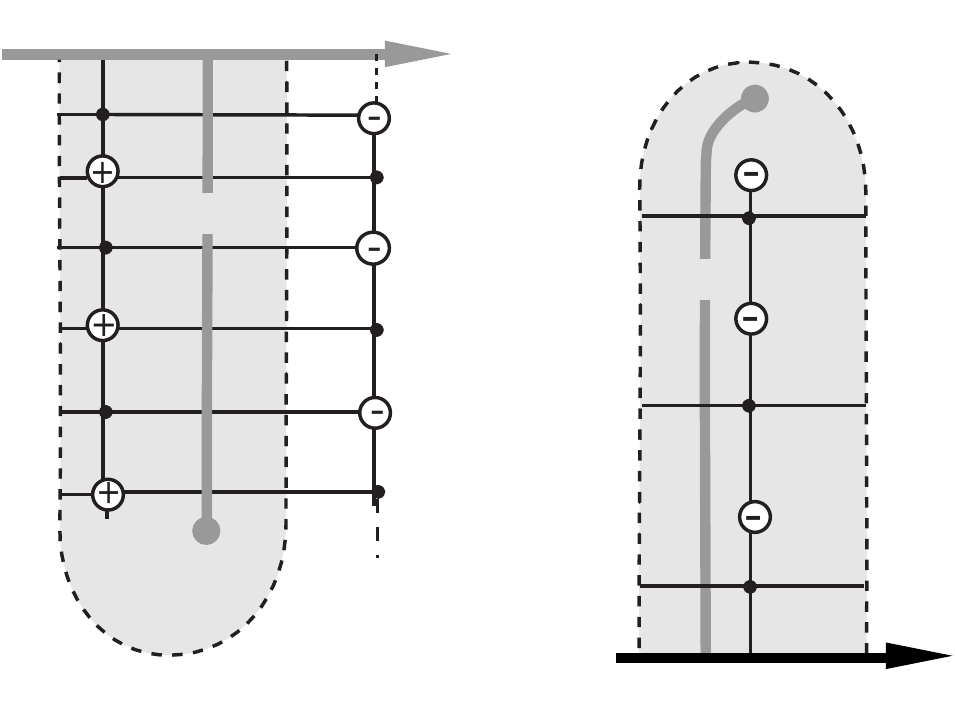}}
\caption{Tab neighborhood pair in $\pa$}
\label{figure:normal1}
\end{figure}
Let $\gamma$ be
a clasp arc and let $\ghigh,\glow$, be its preimages in $\pa$.  We say that 
$(\Thigh,\Tlow)$ is a
 {\it tab neighborhood} of $\gamma $ in $\pa$ if the following hold: 

\be
\item[(i)] $\partial \Tboth = \alpha \cup \beta$
where $\alpha \subset X_\epsilon$ and $\beta$ is an arc that is transverse
to the foliation of $\pa$.

\item[(ii)]  There is a simple path $l_\epsilon \subset \Tboth$, contained  in singular leaves.  It
starts on $X_\epsilon$, contains all $k$ vertices and all $k$ singularities.  If the
singularities on $l_+$ have sign $\delta_1, \delta_2, \dots,\delta_k$ on $\Thigh$, then the
singularities on $l_-$ have sign $-\delta_1, -\delta_2, \dots,-\delta_k$.

\item[(iii)] The arc $ \ge \subset \Tboth$
is the only clasp arc which intersects $\Tboth$.  It intersects each of the $k$ singular leaves in the
induced foliation of $\Tboth$ exactly once.

\item[(iv)] Let $e_+ e^\prime , e^\prime e^{\prime \prime}, \dots , e^{\prime
\dots
\prime} e_-
\subset \gamma_+ \subset \Thigh$ be a subdivision of $\ghigh$
into $k$ subarcs such that each subarc crosses one singular leaf.  
Then the corresponding induced subdivision of $\glow$ in $\Tlow$ given by the immersion $\ta$  also has the
property that each subarc crosses one singular leaf.
\ee

\begin{lemma}
\label{lemma:creating tab neighborhood pair}
We may assume that each pair of clasp arcs $\g$ has a tab neighborhood
$(\Thigh,\Tlow)$ for its associated pair $(\ghigh,\glow)\subset\pa$.
\end{lemma}

\begin{proof} 
To prove the lemma we will make repeated use of
the first and second finger moves.  Consult \fullref{figure:clasp4}, which shows the changes we will
make as they appear on $\Thigh$, and \fullref{figure:clasp5}, which depicts
the corresponding changes on $\pa$ in a neighborhood of $\glow$.  After the sequence of changes $\ghigh$
will have moved to the other side of the chain of singular leaves in $\Thigh$, but the key features
will not have changed. 

\begin{figure}[htpb]
\labellist\small
\pinlabel {(a)} at 66 450
\pinlabel {(b)} at 144 450
\pinlabel {(c)} at 223 450
\pinlabel {(d)} at 300 450
\pinlabel {(e)} at 380 450
\pinlabel {(f)} at 460 450
\pinlabel {(g)} at 542 450

\tiny
\pinlabel {$e_-$} [l] at 66 631
\pinlabel {$p$} [l] at 66 620
\pinlabel {$\gamma_+$} at 54 587
\pinlabel {$p'$} [l] at 66 527
\pinlabel {$e'$} [r] at 54 511

\pinlabel {$p$} [l] at 144 620
\pinlabel {$p'$} [l] at 144 527

\pinlabel {$p'$} [l] at 223 523
\pinlabel {$\Delta$} at 215 598

\pinlabel {$e'$} [tl] at 310 582
\pinlabel {$p'$} [l] at 300 556
\pinlabel {$p''$} [l] at 300 465

\pinlabel {$e'$} [tl] at 390 582
\pinlabel {$p'$} [l] at 380 556
\pinlabel {$p''$} [l] at 380 465

\pinlabel {$e'$} [tl] at 472 582
\pinlabel {$p'$} [tl] at 460 556
\pinlabel {$p''$} [l] at 461 465
\pinlabel {$\Delta\!'$} at 448 547

\pinlabel {$e'$} [l] at 557 588
\pinlabel {$p'$} [l] at 542 555
\pinlabel {$e''$} [tl] at 551 498
\pinlabel {$p''$} [l] at 542 464

\endlabellist
\centerline{\includegraphics[scale=.7]{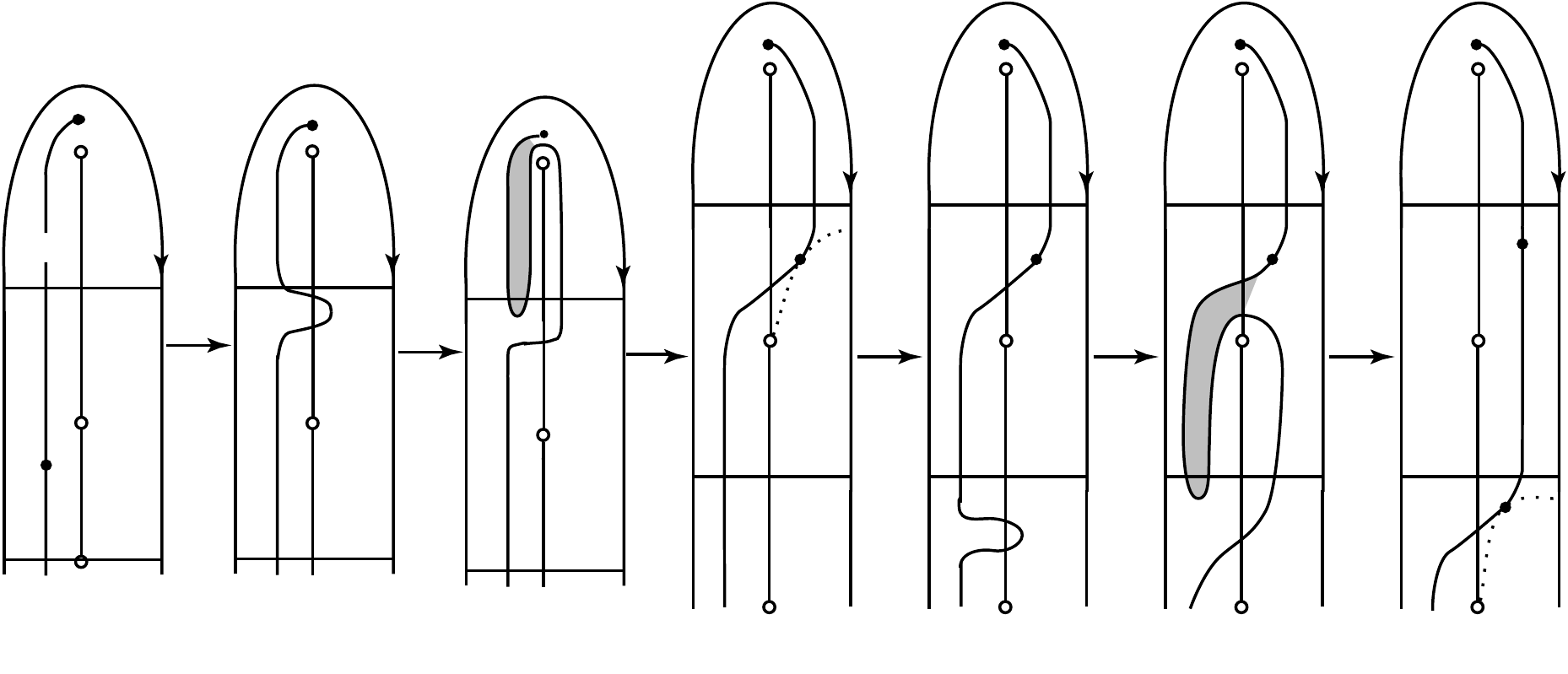}}
\caption{Local changes, as they appear on $\Thigh$}
\label{figure:clasp4}
\end{figure} 

By \fullref{lemma:straightening out the tabs} we may assume that
$\ghigh$ is already contained in a tab neighborhood $\Thigh$.  Since we have not changed the fact that
$\gamma$ is transverse to fibers of $\fib$, it follows that in a sufficiently small foliated neighborhood
of $\glow$ in $\pa$ the arc $\glow$ will also be transverse to the leaves of the foliation of $\pa$.  Some
of the leaves intersecting $\glow$ may be singular.  Label them
$l_1,\dots,l_t$, as  in \fullref{figure:clasp5}.
The first change we introduce is to perform the first finger move on a
small subarc of $\gamma$ which is
just below the singular leaf in the end-tile of $\Thigh$, pushing $\gamma$  across the singular leaf, as
in the passage from sketch (a) to (b) in  \fullref{figure:clasp4}(a) $\to$ (b).  This will induce a
corresponding change in $\glow$, as illustrated in \fullref{figure:clasp5}(a) $\to$ (b).

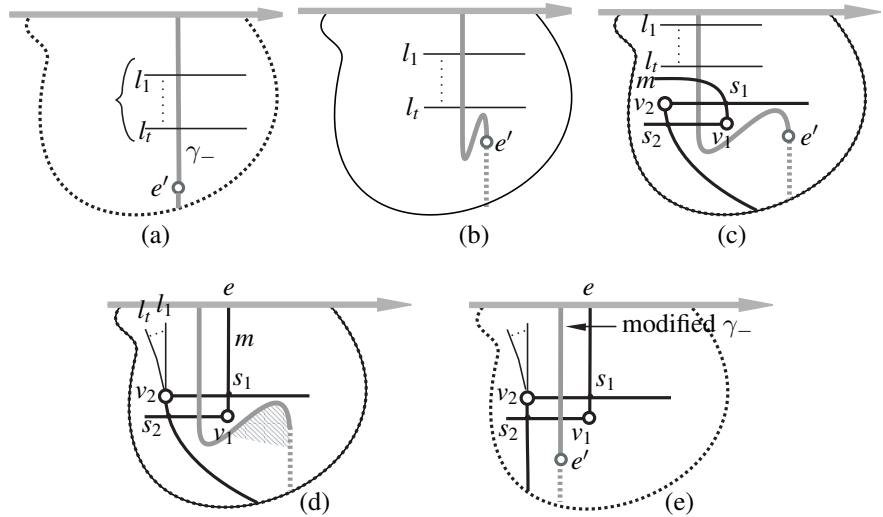
\begin{figure}[htpb]
\labellist\small
\pinlabel {(a)} at 123 420
\pinlabel {(b)} at 323 420
\pinlabel {(c)} at 485 420
\pinlabel {(d)} at 223 252
\pinlabel {(e)} at 452 252

\pinlabel {$l_1$} [r] at 126 518
\pinlabel {$l_t$} [r] at 126 486
\pinlabel {$\gamma_-$} [l] at 137 471
\pinlabel {$e'$} [r] at 137 452

\pinlabel {$l_1$} [r] at 294 535
\pinlabel {$l_t$} [r] at 294 501
\pinlabel {$e'$} [l] at 331 480

\pinlabel {$l_1$} [r] at 444 553
\pinlabel {$l_t$} [r] at 444 530
\pinlabel {$m$} [r] at 442 517
\pinlabel {$v_2$} [r] at 442 502
\pinlabel {$s_2$} [tr] at 445 490
\pinlabel {$v_1$} [t] at 479 491
\pinlabel {$s_1$} [bl] at 480 504
\pinlabel {$e'$} [l] at 523 484

\pinlabel {$l_t$} [b] at 116 361
\pinlabel {$l_1$} [b] at 129 365
\pinlabel {$e$} [b] at 169 378
\pinlabel {$m$} [l] at 168 355
\pinlabel {$s_1$} [bl] at 168 321
\pinlabel {$v_1$} [t] at 166 305
\pinlabel {$v_2$} [r] at 129 321
\pinlabel {$s_2$} [tr] at 131 308

\pinlabel {$v_2$} [r] at 356 320
\pinlabel {$v_1$} [t] at 391 305
\pinlabel {$e'$} [l] at 379 280
\pinlabel {$s_1$} [bl] at 396 320
\pinlabel {$s_2$} [tr] at 355 306
\pinlabel {$e$} [b] at 395 378
\pinlabel {modified $\gamma_-$} [l] at 411 364
\endlabellist
\centerline{\includegraphics[scale=.6]{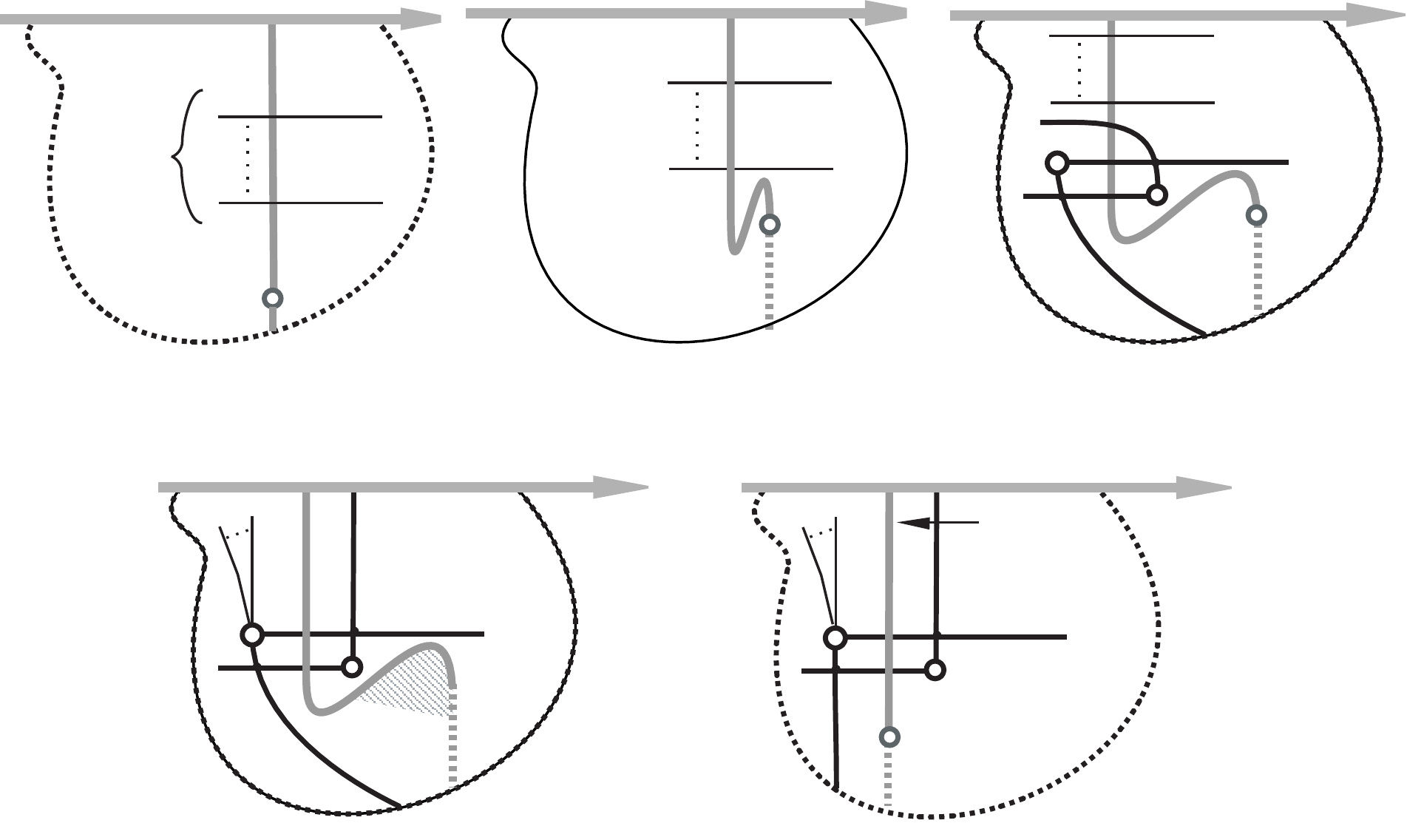}}
\caption{Local changes, depicted in a neighborhood of $\glow$ on $\pa$}
\label{figure:clasp5}
\end{figure}

We next perform the second finger move as shown in the passage 
 \fullref{figure:clasp4}(b)
$\to$\ref{figure:clasp4}(c) and the corresponding 
alteration
\ref{figure:clasp5}(b)$
\to$\ref{figure:clasp5}(c).  It will be helpful
to label the two new vertices and singularities introduced by this second finger move
as $v_1$ and $v_2$ (the vertices) and $s_1$ and $s_2$ (the singularities).  This second
finger move creates a disc $\Delta$ in the tab $\Thigh$ whose interior is necessarily
embedded because it is on a tab-neighborhood of
$\ghigh$ and there is only one clasp arc on each such tab-neighborhood.  We use it to
do the second change  in foliation (see \fullref{lemma:second change in foliation}
and \fullref{figure:second change in foliation}).  The notation has been chosen so that
$s_1$ in \fullref{figure:clasp5}(c) and $\Delta$ in
\fullref{figure:clasp4}(c) correspond to $s_1$ and $\Delta$ in
\fullref{figure:second change in foliation}.  The singularity $s_2$ of \fullref{figure:second change in foliation} is not shown in \fullref{figure:clasp5}(c), also the singularity
$s_2$ of \fullref{figure:clasp5}(c)
is similarly not part of the geometry of \fullref{figure:second change in foliation}.
After a series of such changes in foliation the singular leaf which is labeled $m$ in
\fullref{figure:clasp5}(c) will have 
exchanged order with the singularities in the leaves $l_1,\dots,l_t$ as illustrated in
\fullref{figure:clasp5}(d).  (Note that this implies that each of the leaves $l_1,\dots,l_t$ has one of its endpoints at $v_2$).   After the change in foliation the clasp arc can be tightened,
resulting in the picture we see in \fullref{figure:clasp5}(e). 

From the tab neighborhood definition, condition (ii) on the signs of the vertices follows from the fact that the
vertex endpoint of a singular leaf of the foliation of $\pa$ which begins at $X_\e$
necessarily has sign $-\e$.  But then, the tiles which make up
$\Tboth$ are all type $a_{-\e}a_{-\e}$ and never type
$a_\e a_\e$. 

The point $e^\prime$ in \fullref{figure:clasp4}(d) 
is defined to be the point where the deformed clasp arc is tangent to a fiber of
$\fib$.  Figures \ref{figure:clasp4}(d)
$\to$\ref{figure:clasp4}(e)
$\to$ \ref{figure:clasp4}(f)
$\to$\ref{figure:clasp4}(g) show how the argument can be
iterated.  If there are $k$ tiles on $\Thigh$,
then after $k$ iterations---finger move 1 followed by finger move 2
followed by the second change in foliation, or as we shall call it, the $FFF$ move,  we will have
created a tab neighborhood on $\pa$ for $\glow$.
The subdivision of the tabs described in condition (iv) is achieved automatically via
the iteration of the $FFF$ sequence.
Since we have not moved $\ghigh$ outside
of $\Thigh$ and since after the final iteration of our
second change in foliation $\ghigh$ will again be transverse to the foliation, we
will have created a tab neighborhood pair $(\Thigh,\Tlow)$ for $\g$ as its associated pair $(\ghigh, \glow)$. 

We refer the reader back to the discussion in
\fullref{subsection:the finger moves}, where we introduced the two
finger moves, to verify that the parity information of condition (ii) is satisfied.  The proof of \fullref{lemma:creating tab neighborhood pair} is complete. 
\end{proof}

We are ready to complete the proof of \fullref{proposition:normal}
By repeating the procedure of \fullref{lemma:creating tab neighborhood pair} we can replicate another $\Thigh$
neighborhood inside the  $\Thigh$ neighborhoods of \fullref{figure:normal1}.  This places $\ghigh$ inside a
normal neighborhood $\Nhigh$.  Since $\Nhigh$ contains a tab neighborhood 
$\Thigh$ which belongs to a tab neighborhood pair, we can
interchange the roles of $\ghigh$ and $\glow$ to produce a corresponding
normal neighborhood for $\glow$.  The main point is that the application of our
$FFF$ procedure does not move $\ghigh$
outside the tab neighborhood that is nested inside $\Thigh$.
The key properties of the normal neighborhoods follow. 
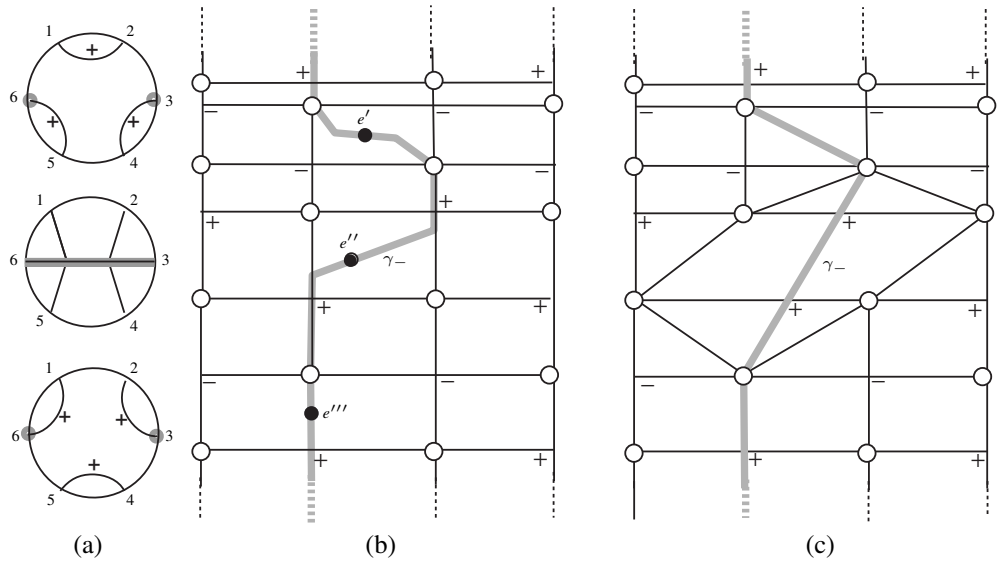
\begin{figure}[htpb!]
\labellist\tiny
\pinlabel {$1$} [br] at 64 521
\pinlabel {$2$} [bl] at 98 521
\pinlabel {$3$} [l] at 118 489
\pinlabel {$4$} [tl] at 100 461
\pinlabel {$5$} [tr] at 65 461
\pinlabel {$6$} [r] at 46 491

\pinlabel {$1$} [br] at 59 429
\pinlabel {$2$} [bl] at 100 429
\pinlabel {$3$} [l] at 116 402
\pinlabel {$4$} [tl] at 100 372
\pinlabel {$5$} [tr] at 59 375
\pinlabel {$6$} [r] at 46 402

\pinlabel {$1$} [br] at 64 337
\pinlabel {$2$} [bl] at 100 337
\pinlabel {$3$} [l] at 118 306
\pinlabel {$4$} [tl] at 98 277
\pinlabel {$5$} [tr] at 65 277
\pinlabel {$6$} [r] at 47 306

\pinlabel {$e'$} [b] at 231 472
\pinlabel {$e''$} [b] at 223 404
\pinlabel {$e'''$} [l] at 204 319
\pinlabel {$\gamma_-$} [tl] at 238 407
\pinlabel {$\gamma_-$} [tl] at 478 404

\small
\pinlabel {(a)} at 80 245
\pinlabel {(b)} at 240 245
\pinlabel {(c)} at 480 245

\endlabellist
\centerline{\includegraphics[scale=.7]{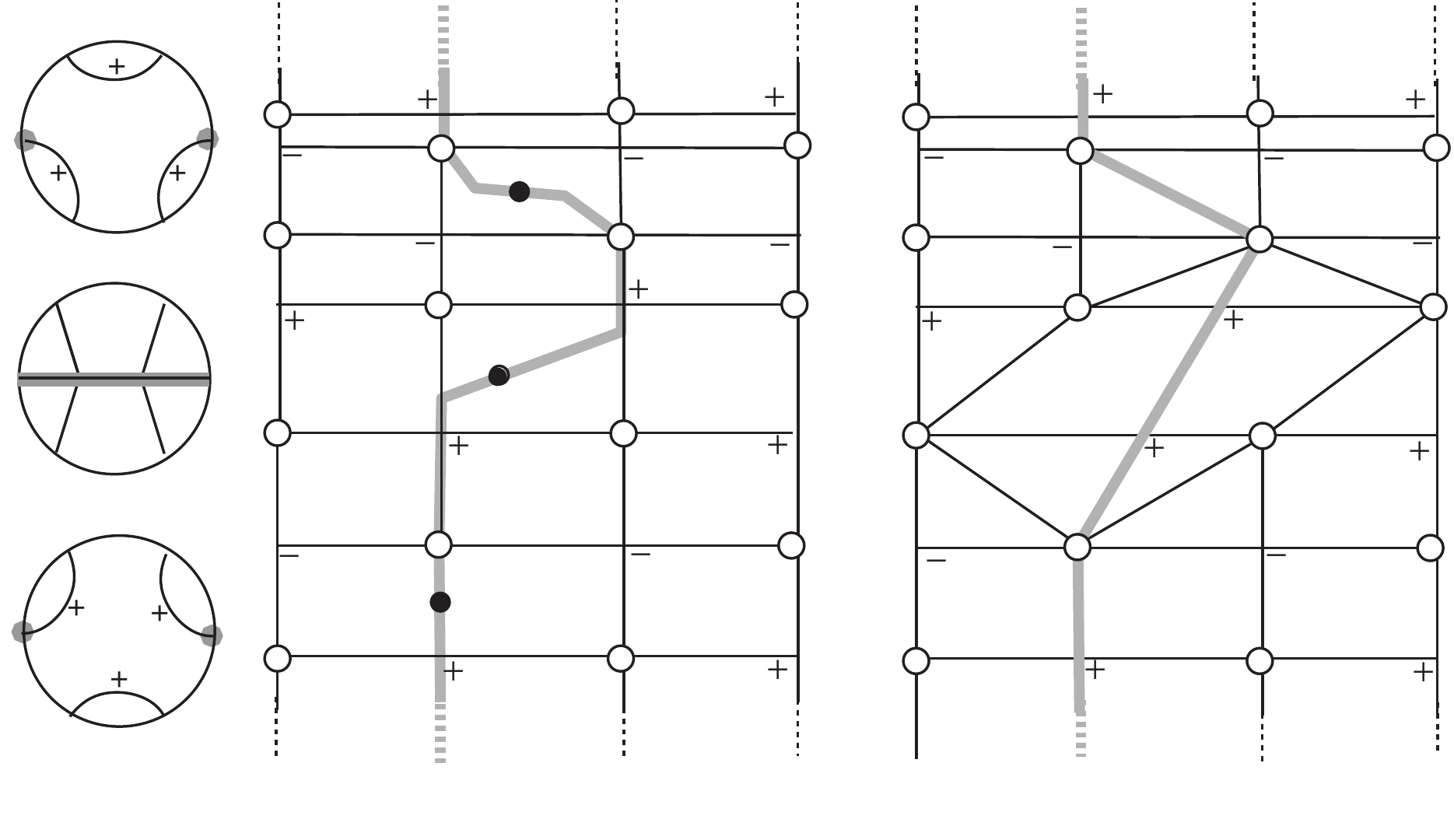}}
\caption{Positioning the clasp arc $\glow$ into a finite number of leaves (singular and/or non-singular)}
\label{figure:normal3}
\end{figure}  

We now claim that our procedure for creating $(\Nhigh, \Nlow)$ also
creates the conditions necessary for pushing the
associated clasp intersection arc into a union of leaves.  To see this let us review our construction of the
normal neighborhood pair.  See \fullref{figure:normal1}.  Recall the notation: 
$e^+ e^\prime , e^\prime e^{\prime\prime}, \dots , e^{\prime \dots \prime} e_-
\subset \glow \subset \Tlow$ is a subdivision of
the clasp arc into $k$ subarcs such that each subarc crosses one singular leaf.  Observe that each of the
subarcs
$e^+ e^\prime , e^\prime e^{\prime \prime}, \dots , e^{\prime \dots \prime} e^-
\subset \glow \subset \Tlow$   has a neighborhood which has the foliation of one of the tiles in
\fullref{figure:annulus2}.  So except for a neighborhood
around each $e^\prime , e^{\prime \prime}, \dots , e^{\prime \dots \prime}$ in
$\ge$, we can push $\ge$ into a finite number of leaves. The
parity of each of these singular leaves is the same as the parity of the pierce end-point of
$\ge$.  (See the end of \fullref{section:getting started} for the definition of the sign).  \fullref{figure:normal3}(b) illustrates an
example of a  clasp arc having a positive pierce end-point that has been partially
pushed into singular leaves. In the example which is illustrated the pierce point is positive, and since there is a
negative singularity on the `vertical' chain of singular leaves the clasp arc is forced to bend in the
manner that is illustrated in sketch (b).   \fullref{figure:normal3}(c) illustrates that, after a change in fibration which introduces singular leaves that
have two singular points (\fullref{figure:normal3}(a)) it is possible to
push our clasp arc into a finite union of leaves by pushing the remaining
arc neighborhoods of the points $ e^\prime , e^\prime , \dots , e^{\prime \dots \prime}$ into leaves of the foliation.  These
two-singularity tiles  arise naturally during the change of
fibration of \fullref{figure:first change in foliation} in the following way: as we perform the change
of fibration in  \fullref{figure:first change in foliation} we must pass through a tile having six sides and
a two-singularity singular leaf.  (This is explained carefully in
Birman--Finkelstein \cite[Figure~2.2]{B-F}).
The two singular points will always have common parity.  In the example, after the change we can reposition the clasp
arc so that it only passes through positive singularities, as in sketch (c).  Once we allow the use of such tiles in the foliation of
$\pa$, the final assertion in the proof of \fullref{proposition:normal} follows.   \end{proof}

\begin{remark}
 \label{remark:clasps into union of leaves}
{\rm When we push the clasp arcs into a union of leaves,  we lose some of the symmetry of normal neighborhoods.  In this regard we note the following.  The symmetry will be needed in one place: for the construction of `thin annuli'  in  \fullref{subsubsection:forming thin annuli}. The building blocks for those annuli are the regions which are illustrated in \fullref{figure:normal7}.   To construct those regions, all we need is a normal neighborhood of that part of a positive clasp arc $\gamma^i_+$ which begins on $\Xhigh$, follows an $a_+$ arc and along $\gamma_+$  to a little bit after it crosses the first singular leaf, together with the corresponding segment on $\gamma^i_-$.  In such a subset of a normal neighborhood, the required symmetry will always be present, even when it does not extend to a more global symmetry.}
\end{remark}

\section{Pushing across $\ta$}
\label{section:pushing across ta}

The machinery has been set up.  In this section we learn how to use the foliation of $\ta$ and $\pa$ to push $\Xhigh$ across the clasp annuli $\ta$ to $\Xlow$.  That is, we show how we use the foliation  to look inside the black box of Markov towers.  Here is an overview of the work in this section:  

In \fullref{subsection:the complexity function c(Xhigh,Xlow,ta)}
we learn how to produce the sequences (2) and (3) of
\fullref{theorem:MTWS}.  Those sequences use the simplest moves of the
MTWS, exchange moves and destabilization,  to simplify the foliation
outside the normal neighborhoods of the clasp arcs, to find the modified
braids $\Xlow^\prime$ and $\Xhigh^\prime$.  We alert the reader to 
\fullref{remark:exchange moves and the winding operation}, where the connection between sequences (1) and  (2) of the MTWS and the `winding' phenomenon that was illustrated  in  \fullref{figure:exchange3} is explained.  

  In \fullref{subsection:microflypes} we introduce a new basic move, the `microflype', which generalizes the example of a microflype which we gave earlier in \fullref{subsection:a key example of the basic construction} and \fullref{figure:getstarted4}.  In particular, we identify microflype `regions' in the foliated annulus $\pa$.   In \fullref{subsection:pushing across a thin annulus} we construct more complicated regions which we call `thin annuli' and learn how to use sequences of microflypes, stabilizations and destabilizations, which can then be amalgamated into admissible flypes, to push $\Xhigh$ across the thin annuli. At the same time we show how to construct a template that represents the push across a thin annulus.  In \fullref{subsection:constructing the template from the total annulus} we generalize this construction, setting up one of the main tools that will be used in the proof of the MTWS.  
  
In \fullref{subsection:G-flype moves} we describe the
foliation that $\ta$ must have so that
the Markov tower corresponds to a G--flype, and work out the foliation of $\ta$ and $\pa$ that gave us 
the G--flype template of \fullref{figure:template2}.   In
\fullref{subsection:pushing across regions with a G-exchange foliation} we describe the foliation
$\ta$ and $\pa$ must have so that the Markov tower corresponds to a G--exchange move, working out in detail 
the explicit foliation that yielded the G--exchange template of \fullref{figure:exchange6}. 
Lastly, in \fullref{subsection:pushing across a standard annulus} we discuss the isotopy associated
with cyclic templates,  and in particular the example that we gave earlier in \fullref{figure:template0}.

\subsection{The complexity function $c(\Xhigh,\Xlow,\ta)$}
\label{subsection:the complexity function c(Xhigh,Xlow,ta)}
By \fullref{proposition:normal} we may assume that each clasp arc pair is contained in a normal
neighborhood pair, and can be pushed into the union of leaves.  Assume from now on that has been done.
Define the
{\it complexity} $c(\Xhigh,\Xlow,\ta)$  
to be the lexicographically ordered pair $(c_1, c_2)$, where $c_1$ is the number of singular points on the clasp
arcs and $c_2$ is the number of singular points that are outside normal neighborhoods.  Note that if the
clasp arcs are transverse to the fibers of $\fib$, then $c_1$ can be interpreted as the number of singular
leaves crossed by $\ghigh$ and $\glow$.

The reader may wonder why we do not include a count of vertices in the foliation?  The reason is
simple: the vertex count is determined by the singularity count,  using the Euler characteristic of the annulus.  
See the argument at the beginning of the proof of \fullref{lemma:Euler characteristic equality}.  

The reader may also wonder why our complexity function ignores all the singularities that are in normal
neighborhoods but not on the clasp arcs.  The reason is again simple: when we constructed normal neighborhoods
we created lots of inessential $b$--arcs.  If we omit them, we will be forced to give up
normal neighborhoods, and we don't want to do that now.  

Finally, the reader may wonder why we are not including braid index in our complexity function.  That is a
more subtle matter.  At this stage in the work it suffices to say that we will introduce it later (see
the  {\it augmented} complexity function of \fullref{section:the proof}.)

In the manuscript \cite{B-M_V} the authors proved that when $\cX$ is the
$\mu$--component unlink, exchange
moves and destabilization suffice to reduce any closed braid representative $\Xhigh$ to 
the identity braid in the $\mu$--strand braid group.  In this section we see how far we can go in the
simplification of our clasp annulus $\ta$ with the help of exchange moves and
destabilization.  The moves that we use here will create the two subsequences (1) and (2) of the MTWS.
The reason that we
had to separate the sequences (1) and (2) in the statement of the MTWS from the sequence (3) is that it may happen
that the given braid $\Xlow$ does not have minimum complexity with respect to exchange moves, for example it may be
wound up in the manner illustrated in \fullref{figure:exchange2}.  If so, then if we
simply tried to modify the given $\Xhigh$ to $\Xlow$, it might not be possible to do it without increasing complexity
at some point. 

We shall regard changes in foliation and braid isotopy to be `trivial moves'.
On the other hand, in general exchange
moves modify  the braid isotopy class and destabilizations change the braid index
too.  Our goal will be to minimize the complexity, using all four moves: braid isotopy, changes in
foliation, exchange moves and destabilizations. 

\begin{figure}[htpb]
\labellist\small
\pinlabel {$G_{++}$} at 150 590
\pinlabel {$G_{--}$} at 270 590
\pinlabel {$G_{+-}$} at 390 590
\pinlabel {$G_{-+}$} at 510 590
\endlabellist
\centerline{\includegraphics[scale=.65]{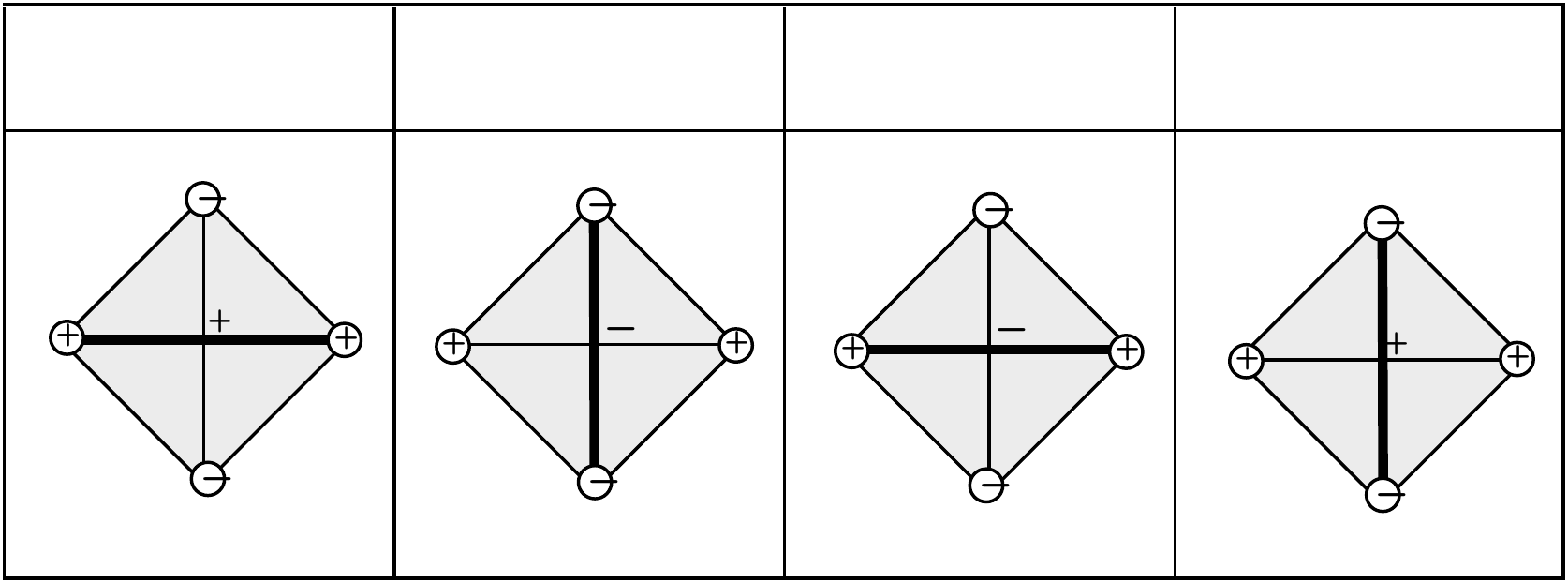}}
\caption{The graph $G_{\epsilon,\delta}$ includes all singular leaves through vertices of sign $\epsilon$
and singularities of sign $\delta$.  The thick edges in this example illustrate the passage of
$G_{+,+}, G_{-,-}, G_{+,-}, G_{-,+}$ through a $bb$--tile.}
\label{figure:graphs1}
\end{figure}
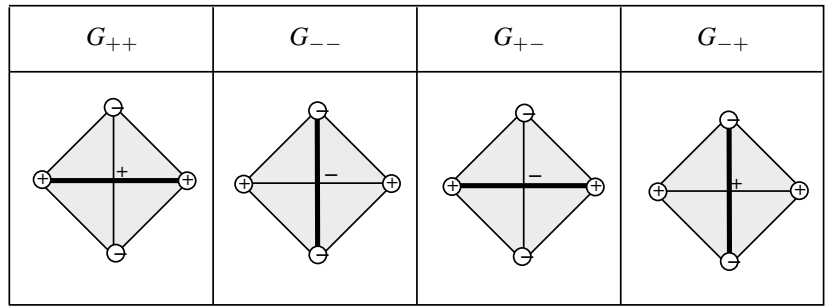

One expects the singular leaves to contain key information about a foliation.  In the case of braid foliations
more is true, because the singular leaves divide nicely into subsets which are characterized by the signs of
their vertices and singularities.  Let $G_{\epsilon,\delta}$, where $\epsilon$ and $\delta$ are $\pm$, be the set of
all singular leaves which pass through only vertices of sign $\epsilon$ and singularities of sign $\delta$.  We
consider the four subgraphs of the graph of the singular leaves $G_{+,+},  G_{-,-}, G_{+,-}, G_{-,+}$.
By definitions
$G_{\epsilon,\delta}\cap G_{-\e,-\delta} = \emptyset$.  See
\fullref{figure:graphs1}, which illustrates how the 4 graphs intersect a
$bb$--tile. 
 Similar graphs appeared in Bennequin \cite{Bennequin}, in connection with his studies of the characteristic
foliation of surfaces bounded by knots.

The intersection of a clasp arc $\ge$ with a subdisc $\Delta$ of $\pa$ is {\it good} if $(\Delta
\cap\ge)\subset(\Delta\cap G_{\epsilon,\delta})$ for some $(\epsilon,\delta) \in \{(+,+),(+,-),(-,+),(-,-)\}$
and if no component of $\Delta\cap\ge$ includes the puncture endpoint of $\ge$.
A subdisc $\Delta$ of $\pa$ is {\it good} if every component of intersection of $\Delta$ with a clasp arc $\ge$ is good.
A subdisc that is disjoint from the clasp arcs is, of course, good.  A disc that contains a puncture endpoint of a clasp arc is not good.

A vertex is said to be {\it near $X_\e$ }, where  $\e=\pm$, if it is the endpoint of a leaf
of type $a_\e$ in the foliation of $\pa$.  A vertex $v$ is said to be an {\it interior vertex} it is
not near either $\Xhigh$ or $\Xlow$.
In both cases we define {\it link$(v)$} to be  the
closure of the union of all $b$--arcs and $a_\e$--arcs which meet $v$.

To begin our work, we study $ab$ and $bb$--exchange moves and changes in foliation in the presence of clasp arcs,
under the hypothesis that the moves are supported in good discs: 

{\bf The $(ab)^\star$ exchange move}\qua
Let $v$ be a vertex of valence 2 in the foliation of $\pa$ which is near
$X_\e$ and is the endpoint of a $b$--arc.
Assume that
$\Delta =$ link$(v)$ is a good disc which has non-empty intersection with a clasp arc $\ge$.
Consult \fullref{figure:exchange5}(a), which illustrates the case when $v$ is a positive vertex.  The disc
$\Delta$ is the closed disc bounded by the singular leaves $wqp$ and $wst$ in \fullref{figure:exchange5}(a) and
by the subarc $pt$ of $\Xhigh$.  If $\e=-$ all vertex signs are reversed and the roles of
$\Xhigh$ and
$\Xlow$ are interchanged, but the underlying phenomena are unchanged.  We are interested in whether we can
push $\Xhigh$ across $\Delta$ in the presence of the clasp arc $\gamma$? 

Notice that all 4 graphs  $G_{+,+}, G_{+,-}, G_{-,+}$ and $G_{-,-}$ intersect
$\Delta$ in \fullref{figure:exchange5}(a).  By our definition of a good disc,
$\gamma\cap\Delta$ is a subset of some $G_{\epsilon,\delta}$, however the case
$\epsilon=+$ is impossible.  The reason is: $G_{+,\pm}$ ends at $v$ in the interior
of $\Delta$ (because the two singularities in $\Delta$ have opposite signs), also by hypothesis $\Delta$ is a 
good disc, so that $\gamma$ cannot have its puncture endpoint in $\Delta$.  Therefore $(\ge \cap \Delta)\subset
G_{-,+}$ or $G_{-,-}$.  Both cases can occur, and the two cases are essentially the
same.  In both cases
$\Delta$ is embedded, because $\ge$ lies in a normal neighborhood and so its partner $\g_{-\e}$ cannot
intersect $\Delta$.  In fact normal neighborhoods tell us more: $\ge$ is the only clasp arc that can intersect
$\Delta$.  So
$\gamma$ begins at, say, 
$p\in \Xhigh$ and passes through $q$ to $w$, always in $\partial\Delta$.  Just as in the embedded case, we can
push
$\Xhigh$ across a neighborhood of $\Delta$, removing two vertices and two singularities from the
foliation.  This is the `light bulb' move of \fullref{figure:exchange5}(b).  It shortens $\ge$ and so
also shortens its partner $\g_{-\e}$, reducing the entry $c_1$ in $c(\Xhigh,\Xlow,\ta)$.  The second sheet of $\pa$
is unchanged by the move, except for the fact that its clasp arc is shortened by pulling in its puncture
endpoint.  The link
$\Xhigh$ remains embedded throughout the move.  The argument is identical if $\Xhigh$ is replaced by $\Xlow$.
This is the move that we call the $(ab)^\star$ exchange move.  It is
exactly the same as the $ab$--exchange move,
but in the presence of clasp arcs and with the assumption of good discs.

{\bf The $(bb)^\star$ exchange move}\qua  Let $v$ be an interior vertex of valence 2 in the foliation of
$\pa$.  Assume that $\Delta =$ link$(v)$ is a good disc which has non-empty intersection with a clasp arc
$\gamma$.  Consult \fullref{figure:foliations7}, which illustrates the case when $v$ is a negative vertex.  The disc
$\Delta$ is the closed disc bounded by the singular leaves $w_1sw_2$ and $w_1qw_2$.  All
vertex signs could be reversed, it will not matter, so we assume they are as illustrated.  We are interested in
whether we can do the $bb$--exchange move and collapse of the pocket, as illustrated in the passage from the left
to right sketches in \fullref{figure:foliations7}, in the presence of the clasp arc $\gamma$, when $\Delta$
is a good disc?

As in the case of the $(ab)^\star$ move, the intersection of our clasp arc with $\Delta$ is assumed to be a subset
of the intersection of one of the 4 graphs with $\Delta$.
Since the puncture point of $\gamma$ cannot be in $\Delta$
the only possibilities are the graphs $G_{+,+}$ or
$G_{+,-}$, because the other two possibilities lead to a puncture point in the interior of $\Delta$.
Assume without
loss of generality that $\gamma \subset G_{+,+}$, for example $\gamma$ might pass through the points 
$pw_1sw_2p'$ in
the left sketch in \fullref{figure:foliations7}.  Using the hypothesis that
$\gamma$ is in a normal neighborhood, we can (by \fullref{proposition:normal}) push
$\gamma$ to a new position which we call $pw_1's'w_2'p'$ where it is transverse
to the leaves of the foliation.  But then the $bb$--exchange move goes through as
before, as described in full detail in sketches 2.15--2.19 and the
accompanying text of Birman--Finkelstein \cite{B-F}.  Since that
argument is long and technical, and since a complete reference is available, we do not repeat it here.  After the
move the `pocket may be collapsed', as illustrated in the passage from the
left to the middle sketch in \fullref{figure:exchange4} of this paper, and
also in the passage from the left to the right sketch in
\fullref{figure:foliations7}.  The clasp arc can then be pushed back into
a union of leaves, ie to the position 
$pw_2p'$ in the right sketch.  The complexity
$c(\Xhigh,\Xlow,\ta)$ has been reduced because the singularity $s$ is no longer on $\gamma$.  This is the move that
we call the $(bb)^\star$ exchange move.  It is exactly the same as the
$bb$--exchange move, but in the presence
of clasp arcs and with the assumption of good discs.

\textbf{The change of foliation of \fullref{lemma:first change in foliation} in the presence of clasp
arcs}\qua As in the case of the $ab$ and $bb$--exchange moves, we will refer
to the manuscript by Birman and Finkelstein \cite{B-F} for all
details, since the required changes in foliation were justified very carefully there, in the proof of Theorem
2.1 of that paper, under the hypothesis that the surface in question was embedded.  But the immersed case is
really no different from the embedded case because the changes in foliation that we need here can always be thought
of as being induced by changes in the choice of disc fibers for the fibration
$\fib$ of 3--space minus the braid axis, and not by any change in the
position of $\ta$ in 3--space.  See \fullref{figure:first change in foliation}, which illustrates the case when we are interested in changing the
foliation on two $bb$--tiles which meet in a common $b$--arc.  The sketches show the changes when there are no clasp
arcs.  The fact that the clasp arcs have been positioned in normal neighborhoods shows us that we can change the
foliation in one of the sheets without changing the other sheet.  Assume
that the union of our two $bb$--tiles is a
good disc.  There are three cases:

\begin{enumerate}

\item [(i)] Both $s_1$ and $s_2$ are contained in a single normal neighborhood as
consecutive singularities, on a
clasp arc.  After the
change of foliation, the path of the clasp arc contains one fewer
singularity, so $c(\Xhigh,\Xlow,\ta)$ is reduced.  In particular, if $ \epsilon = +$ and we have
a symmetric clasp arc then the resulting change of foliation is the one that is labeled (ii) in 
\fullref{figure:first change in foliation}. In this case the 
clasp arc contains the two right-most positive vertices and the
right-most singular point.  If $ \epsilon = -$, then the resulting change of foliation is the
one that is labeled (i) in  \fullref{figure:first change in foliation}, and the
clasp arc contains the two left-most negative vertices and the
left-most singular point.

\item [(ii)]  One of the two singularities, say $s_1$, is contained in a normal neighborhood of a
clasp arc, but not on the clasp arc.  After the change
in foliation, the path of the clasp arc is unchanged.  For example,
if $\epsilon = + $ and the clasp arc contains $s_1$ (resp.\  $s_2$)  then we have a change of
foliation that takes us to (i) (resp.\  (ii)).

\item [(iii)]
The two singularities are contained in disjoint normal neighborhoods of different clasp
arcs.  The path of the two clasp arcs is unchanged by the change in foliation.
For example, if $\epsilon = +$ and the bottom singularity was contained
in the clasp arc for the beginning tilings in the sequences
in \fullref{figure:first change in foliation}
then the resulting change of foliation is (ii) in all three cases.
\end{enumerate}

Armed with this knowledge, we are now able to adapt \cite[Theorem~3.1]{B-F}
to our clasp annulus $\ta$ and its foliated preimage $\pa$.  We will use the symbol 
$\Nboth$ for a normal neighborhood of the clasp arc $\ge$ and the symbol $\cN$ for the union of all normal
neighborhoods of all clasp arcs.

\begin{proposition}
\label{proposition:behavior of graphs} Each of the following holds for
all four graphs $G_{+,+}, G_{-,-}, G_{-,+}$ and $G_{+,-}$: 
\begin{enumerate}
\item  $G_{\epsilon,\delta} \cap
G_{-\epsilon,-\delta}=\emptyset$.
\item  Every singular point and every vertex in the 
foliation of $\pa$ is in
$G_{+,+}$ or $G_{-,-}$ (and also in $G_{+,-}$ or $G_{-,+}$).
\item  $G_{\epsilon,\delta}$ has no interior isolated 
vertex $v$.
\item  $G_{\epsilon,\delta} \setminus (G_{\epsilon,\delta} \cap \cN) $
has no interior endpoint vertex $v$.
\ee

Assume that all clasp arcs are positioned in normal neighborhoods, and that initially they have been pushed into
unions of leaves of the foliation.  Then after some number of exchange moves (combined with changes in foliation and
isotopies in the complement of the axis) the following holds for all four graphs at once:
\be
\item [{\rm (5)}] $G_{\epsilon,\delta}$ contains no closed loop 
$l$ which bounds a good
disk $\Delta\subset\pa$.
\ee
After some number of exchange moves and destabilizations, the following also holds for all four graphs at once:
\be
\item  [{\rm (6)}] There is no closed loop $l$ which is the union of an edgepath
$\cE_1\subset G_{\epsilon,+}$ and an edgepath $\cE_2\subset G_{\epsilon,-}$ which
bounds a good disc $\Delta\subset\pa\setminus\cN$.
\item  [{\rm (7)}] Let $\cE_1,\cE_2,\cE_3$ be connected arcs, with $\cE_1\subset G_{\e,\e},
\ \cE_2\subset G_{\e,-\e},\  \cE_3\subset X_\e.$ Then there is no
closed loop $l=\cE_1\sqcup \cE_2\sqcup \cE_3$ or $\cE_2\sqcup \cE_3$ or $\cE_1\sqcup \cE_3$ which
bounds a good disc $\Delta$ on $\pa$.
\ee
\end{proposition}

\begin{proof}  The proof is almost identical with the proof of
\cite[Theorem~3.1]{B-F}.
Outside normal neighborhoods of the clasp arcs it is identical.  Inside normal neighborhoods the key
concepts which makes it possible to carry over arguments used in
\cite{B-F} are good discs.  In the presence of clasp arcs we simply use the
exchange moves $(ab)^\star$ and $(bb)^\star$ instead of $ab$ and $bb$ and the complexity function
$c(\Xhigh.\Xlow,\ta)$.
An application of ${(ab)}^\star$ can shorten the length of a clasp
arc $\ge$, reducing the number of singular leaves it intersects.
It will necessarily shorten the length of $\gamma_{-\epsilon}$ and
decrease the number of tiles of $\pa$ that are in ${\bf N}_{-\e}$.
Thus, there may be some inessential $b$--arcs that are now away from
normal neighborhoods, and they can be eliminated.

Also, an application of ${(bb)}^\star$ can shorten the length of a
clasp arc $\ge$, reducing the number of singular leaves
it intersects.  But, it does not immediately reduce the number of
singular leaves $\gamma_{-\e}$ intersects.  Moreover, neither
$\ge$ nor $\gamma_{-\e}$ may be in normal neighborhoods
anymore.  To re-establish the symmetry between the the lengths of
the two clasp arcs in $\pa$ we take a tab neighborhood around
$\ge$ (since a tab neighborhood will have the shortest
clasp length) and re-apply the normal neighborhood construction of
\fullref{proposition:normal}.\end{proof}  

Our first corollary produces the modified braids $\Xhigh^\prime$ and $\Xlow^\prime$ in the statement of the MTWS.  In the course of changing $X_\pm$ to $X_\pm^\prime$, it also produces the sequences (1) and (2) of the MTWS, which modify $X_\pm \to X_\pm^\prime$:

\begin{corollary} 
\label{corollary:minimum complexity up to exchange moves and destabilization}  
Assume that all changes are made in the complement of normal neighborhoods of the clasp arcs.   Then with the help of:
\bi
\item complexity-reducing exchange moves on $\Xhigh$ and $\Xlow$,
\item complexity-reducing destabilizations of $\Xhigh$,  
\item changes in foliation (which modify $\ta$ but do not change $\Xhigh$ or $\Xlow$), 
\ei
the triplet $(\Xhigh,\Xlow,\ta)$ may be changed to a new triplet $(\Xhigh',\Xlow',\ta')$ which has minimum complexity, and for which the foliation of $\pa$ has the following properties:  
\be

\item [{\rm (1)}]  There are no $a_\epsilon s$ or $a_\epsilon a_\epsilon$
singularities in $\pa$ for $\epsilon = +$ or $-$.

\item [{\rm (2)}]  If $\alpha$
is a path in  $G_{\epsilon,\delta}$ which begins and ends on
$X_{\epsilon}$, and  $\Delta$ is the disc on $\pa$ which is split off by
$\alpha$ and a subarc of $X_\epsilon$, then $\intr(\Delta)$ contains a puncture end-point of a $\gamma_{-\epsilon}$ clasp arc.

\item [{\rm (3)}]   If $\alpha \subset G_{\epsilon,\delta}$
is a loop which bounds a disc $\Delta$ in $\pa$, then $\intr(\Delta)$ contains at least two puncture end-points of clasp arcs, one
from some $\gamma^i_{\epsilon}$ and the other from $\gamma^i_{-\epsilon}$.  

\item [{\rm (4)}]  If $v \in G_{\epsilon,\delta}$
is an endpoint vertex, then either $link(v)$ contains a puncture
endpoint of a clasp arc or $v$ lies in a  normal neighborhood of a clasp arc.
\ee
\end{corollary}

\begin{proof}   In the case of an embedded surface, details on how to pass
from the analogue of  \fullref{proposition:behavior of graphs} to that of
the Corollary  are given in the proof of \cite[Lemma~3.8]{B-F}.
We have modified the statement of  \cite[Lemma~3.8]{B-F}, by saying
that if a disc $\Delta$ exists with the stated undesirable properties,
then $\Delta$  is not a good disk.  In the case of assertion (3) of
the Corollary we need the puncture endpoints of both $\gamma^i_+$ and
$\gamma^i_-$ because if both are not in $\Delta$ then $\Delta$ will be
embedded, but no such embedded disc exists.  \end{proof}

\begin{remark} \label{remark:exchange moves and the winding operation}  {\rm  Armed with \fullref{corollary:minimum complexity up to exchange moves and destabilization}, the reader is advised to go back and look at \fullref{figure:exchange3} again.  It was intuitively clear, when we first presented this sketch,  that the unwinding process which is illustrated in the passage from left to right ought to be `complexity-reducing'.  Now that we have finally defined an appropriate complexity function we can understand exactly what that means.  

We need the modifications of \fullref{corollary:minimum complexity up to exchange moves and destabilization} for the following reason: If it should happen that the given braid $\Xlow$ is wound up some number of times as on the left, and if we were to modify $\Xhigh$ but leave $\Xlow$ unaltered, then we might be able to use exchange moves and destabilizations to modify $\Xhigh$ to a minimum braid index representative, using only complexity-decreasing moves,  but if we did not modify $\Xlow$ too then we would have to increase complexity at the end.  That is the reason that we need the initial complexity-reducing modifications of {\it both} $\Xhigh$ and $\Xlow$.  That is the content of \fullref{corollary:minimum complexity up to exchange moves and destabilization}.}
\end{remark}

An annulus component of $\ta$ is {\it standard}
 if every component $E \subset G_{\epsilon,\delta}$
on the annulus satisfies the following conditions:

\be
\item $E$ is homeomorphic to either $S^1$ or $[0,1]$, ie either a circle or a line.
\item If $E$ is a circle then it is a core circle of the annulus component.
\item If $E$ is a line with $\partial E = p \cup p^{\prime}$ then $p \in X_{\epsilon}$
and $p^\prime$ is near $X_{-\epsilon}$.
\ee

\begin{figure}[htpb!]
\labellist\small
\pinlabel {standard annulus} at 78 474
\pinlabel {trivially foliated annulus} at 79 318
\pinlabel {type $1_d$} at 267 471
\pinlabel {type $0_{d,1}$} at 494 471
\pinlabel {type $-1_d$} at 267 337
\pinlabel {type $0_{d,2}$} at 494 337

\tiny
\pinlabel {$\delta$} [tr] at 43 549
\pinlabel {$\delta$} [tr] at 86 607
\pinlabel {$\delta$} [tr] at 86 505
\pinlabel {$-\delta$} [r] at 76 561
\pinlabel {$-\delta$} [r] at 76 537
\pinlabel {$-\delta$} [b] at 95 561
\pinlabel {$-\delta$} [t] at 95 537

\pinlabel {$s$} [r] at 180 570
\pinlabel {$-\delta$} [br] at 206 552
\pinlabel {$\delta$} [tr] at 227 516
\pinlabel {$-\delta$} [br] at 249 552
\pinlabel {$\delta$} [tr] at 272 516
\pinlabel {$-\delta$} [br] at 294 552
\pinlabel {$\delta$} [tr] at 315 514
\pinlabel {$-\delta$} [br] at 339 552
\pinlabel {$s'$} [l] at 370 570

\pinlabel {$s$} [r] at 407 570
\pinlabel {$-\delta$} [br] at 434 550
\pinlabel {$\delta$} [tr] at 454 513
\pinlabel {$-\delta$} [br] at 476 550
\pinlabel {$\delta$} [tr] at 500 513
\pinlabel {$-\delta$} [br] at 521 550
\pinlabel {$\delta$} [tr] at 544 513
\pinlabel {$s'$} [l] at 571 495

\pinlabel {$s$} [r] at 177 367
\pinlabel {$\delta$} [tr] at 204 387
\pinlabel {$-\delta$} [br] at 225 423
\pinlabel {$\delta$} [tr] at 246 387
\pinlabel {$-\delta$} [br] at 268 423
\pinlabel {$\delta$} [tr] at 292 387
\pinlabel {$-\delta$} [br] at 313 423
\pinlabel {$\delta$} [tr] at 337 387
\pinlabel {$s'$} [l] at 367 367

\pinlabel {$s$} [r] at 407 367
\pinlabel {$\delta$} [tr] at 434 384
\pinlabel {$-\delta$} [br] at 457 421
\pinlabel {$\delta$} [tr] at 478 384
\pinlabel {$-\delta$} [br] at 500 421
\pinlabel {$\delta$} [tr] at 524 384
\pinlabel {$-\delta$} [br] at 545 421
\pinlabel {$s'$} [l] at 571 439
\endlabellist
\centerline{\includegraphics[scale=.65]{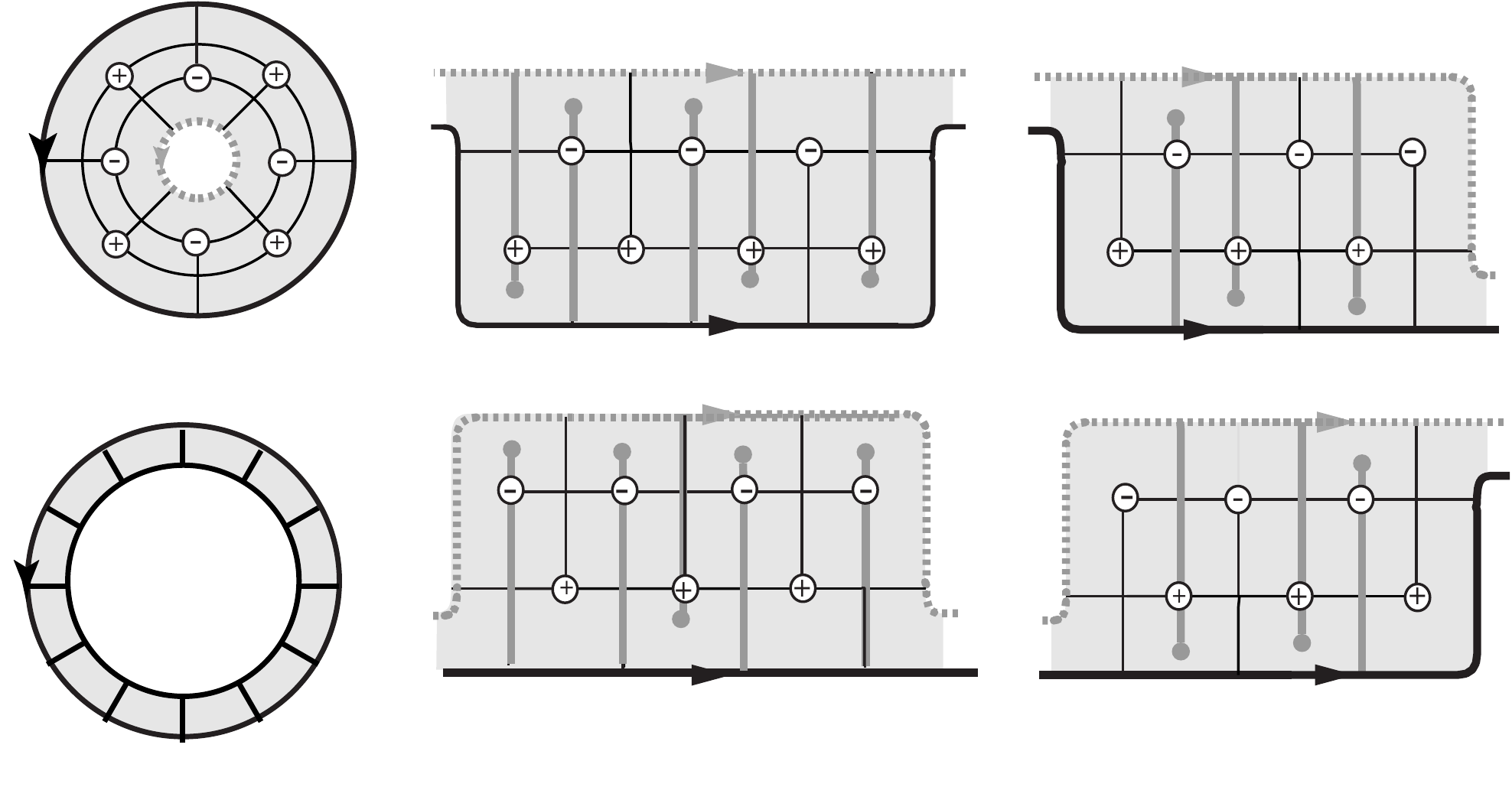}}
\caption{Examples of a standard annulus, a trivially foliated annulus, and regions in a thin annulus. Thin annuli will be defined in \fullref {subsection:pushing across a thin annulus}.}
\label{figure:normal7}
\end{figure}

An annulus component of $\ta$ is {\it trivially foliated} if it is
foliated entirely by $s$--arcs.
Examples of a standard annulus, and of a trivially foliated annulus, are given in \fullref{figure:normal7}.  

Our second Corollary relates to modifications which may be needed after all the clasp arcs have been shortened, that is $c_1=0$.  

\begin{corollary}
\label{corollary:standard annulus}
Let $(\Xlow,\Xhigh,\ta)$ be of minimal
complexity with $c_1 = 0$.  Then each component
of $\pa$ is either an annulus which is foliated entirely with $s$--arcs or a
standard annulus.
\end{corollary}

\begin{proof}   The statement follows directly from \fullref{corollary:minimum complexity up to exchange moves and destabilization}, when we add the assumption that $c_1=0$.  \end{proof} 

\begin{remark}
The problem of pushing $\Xhigh$ across such a standard annulus will be treated in \fullref{subsection:pushing across a standard annulus}.
\end{remark}

\begin{remark}
\label{remark:example,ta-tiled}
We know exactly how the foliation of $\pa$ looks inside a normal neighborhood of a clasp
arc.  We now ask what can happen outside the union of normal neighborhoods, when the puncture points are in the
interior of
$\pa$, as in the left and right sketches in \fullref{figure:normal2}.  An example is illustrated in \fullref{figure:template4}. The shaded regions are unions of overlapping normal
\begin{figure}[htpb!]
\labellist\small
\pinlabel {$U$} [b] at 159 666
\pinlabel {$V$} [b] at 181 666
\pinlabel {$Y$} [b] at 203 666
\pinlabel {$Z$} [b] at 226 666
\pinlabel {$J$} [b] at 249 666
\pinlabel {$K$} [b] at 290 666
\pinlabel {$W$} [t] at 191 625
\pinlabel {$W'$} [t] at 269 625
\pinlabel {$W''$} [t] at 355 625
\endlabellist
\centerline{\includegraphics[scale=.8]{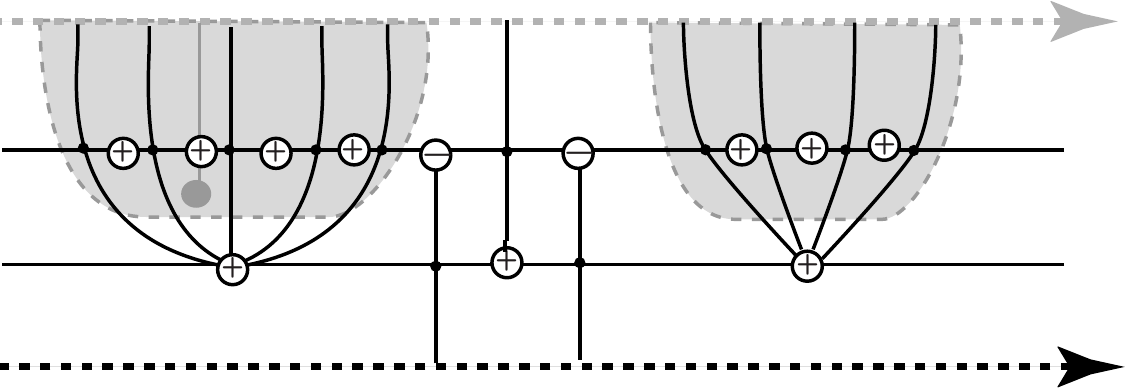}}
\caption{The foliation of $\pa$, near the puncture point of a clasp arc $\gamma_-$ but outside its normal
neighborhood.  It is assumed here that the puncture point of $\gamma_-$ is in the interior of $\pa$.}
\label{figure:template4}
\end{figure}
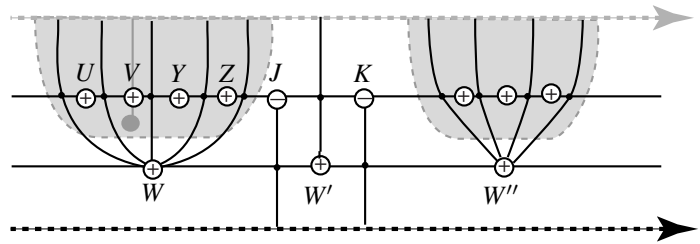
neighborhoods of clasp arcs.   The clasp arcs in question are negative.  We show only one of them, to avoid clutter and
enable us to focus on the features which are of interest now.  It is to be compared with the right sketch in \fullref{figure:normal2}.
The boundaries of the region we are studying are indicated with dotted black and grey arcs transverse to the
foliation because in general they are in the interior of $\pa$ (although in special cases there is no reason why
they could not be $\Xhigh$ and $\Xlow$.   This sketch will be used in \fullref{section:the proof} in the proof of
the finiteness of $\cT(m)$.\hfill $\diamondsuit$
\end{remark}

\subsection{Pushing across a microflype region}
\label{subsection:microflypes}

Having simplified the foliation of $\ta$  as much as possible with the use of exchange moves and
destabilizations, we begin to investigate new features of the foliation of $\ta$.

 A clasp arc $\gamma$ in $\ta$ has two preimages in $\pa$, namely $\ghigh$ and $\glow$.  Recall that 
$\ge, \ \epsilon = \pm$, begins on $X_\epsilon$ and ends at an interior point of $\pa$.  Recall also that $\ge$ was
defined to be short if it does not intersect any singular leaves.  We now define $\ge$ to be 
{\it long} if the puncture endpoint of $\ge$ is on an $a_{-\e}$--arc.  We
say that $\gamma$ is {\it doubly long} if both $\ghigh$ and $\glow$ are long.  An example of a doubly long
clasp arc was given in \fullref{figure:getstarted4}.  On the other hand, in \fullref{figure:template3} there are 4 clasp arcs, and among the eight preimages  $\glow^1,
\glow^2, \ghigh^3$ and $\ghigh^4$ are long, but their partners $\ghigh^1,
\ghigh^2, \glow^3$ and $\glow^4$ are not, so that there are no doubly long clasp arcs.  (Remark: the foliated
annulus in \fullref{figure:template3} gave rise to the 6--braid G--flype template which we described in the
introduction to this paper, in \fullref{figure:template2}).

 A clasp arc $\ge$ in $\pa$ is {\it intermediate} if it is neither short nor long.  This implies that its
puncture endpoint is on a $b$--arc.  Look ahead to \fullref{figure:exchange7} for examples.  Eight clasp arcs
are depicted there, and the 16 preimages are all intermediate.  (Remark:
the foliated annulus \fullref{figure:exchange7} gave rise to the 6--braid
G--exchange template which we described in the introduction to this
paper, in \fullref{figure:exchange6}).   The
{\it length} of the clasp arc $\ghigh$ is the number of singular leaves which $\ghigh$
crosses.  We will not be concerned with the length of $\glow$.  Note that if long clasp arcs occur, then $c_1 > 0$. 

In this subsection we consider the case when there is exactly
one clasp arc $\gamma$ and it is doubly long, so that $c(\Xhigh,\Xlow,\ta)=(1,0)$.  Since $c_2=0$, we may
assume that outside normal neighborhoods of $\ghigh$ and $\glow$, the foliation of $\pa$ consists
entirely of bands which are foliated by $s$--arcs.  But then, since there are no other clasp arcs in $\pa$, we may
just as well simplify  the normal neighborhoods by the deletion of
inessential $b$--arcs.  Thus $\pa$ is a union of an
$a_+s$--tile and an $a_-s$--tile, joined up by bands of $s$--arcs.  This very simple foliated annulus was  introduced long ago, as the `key example' at the end of \fullref{section:getting started}, in \fullref{figure:getstarted4}.
Understanding microflypes will allow us, later, to consider very much more complicated annuli $\pa$.
Recall  that microflypes are braid-index preserving moves which replace very simple Markov towers, as can be
seen from \fullref{figure:exchange2}, so that they might be expected to be basic to our work.

Of course, when we first encountered the microflype region in $\pa$ in \fullref{section:getting started}, as
an example of the basic construction, we did not have available to us the machinery of
braid foliations.  By our work in \fullref{subsection:normal} we
know that the signs of the singularities in the two tiles are opposite.  

The leaves in the braid foliation of
$\pa$ are level sets for the embedding of the two tiles which make up
$\pa$.  Referring back to \fullref{figure:getstarted4}, it should now be completely clear that the bottom row of
sketches represents the  immersion in 3--space of the two foliated tiles in the top row of sketches, and that the
motion of $\Xhigh$ across the two illustrated discs is indeed realized by a flype which takes $\Xhigh$ to a new
position which is separated from $\Xlow$ by a band of $s$--arcs.  The associated block strand diagram is clearly that
for a flype, which in the case which is illustrated is a negative flype.  The braid inside the braid block is a single
full twist of two strands.  The sign of the full twist that is illustrated is negative.

Summarizing: A {\it microflype region} 
in $\pa$ is a subset of $\pa$ which is a union of an $a_+s$--tile and
an $a_-s$--tile, intersected by paired clasp arcs $\ghigh$ and $\glow$, each of which intersects a single
singular leaf.    A {\it microflype} 
on a block-strand diagram is a flype in which the
braid in the braid block $R$ consists of exactly one full twist of either sign
on two strands.  Such a  braid block $R$ will be called a {\it microblock}.
All strands have weight 1.  The {\it sign of a microflype} 
is the pair $(\alpha,\beta)$ where $\alpha$ is sign of the half-twist which is outside the braid block and
$\beta$ is the sign of the full twist which is inside the braid block.  In 
\fullref{subsection:pushing across a thin annulus}, below,  we will show that microflypes with their
associated microblocks are the building blocks of the most general flypes, with arbitrary braids in the braid block and
arbitrary weights on the strands.

By the construction in \fullref{proposition:normal}, it is possible to push the clasp arcs into the
associated singular leaves, giving a more symmetric embedding in 3--space.

\subsection{Pushing across thin annuli}
\label{subsection:pushing across a thin annulus}

While
we have been able, up to now,  to make the tacit assumption that we are working with knots, we now
return to the general case of links. 

Our  task in this section is to learn how to use flypes to shorten 
the length of long (but not necessarily doubly long) clasp arcs in the foliation of $\pa$,
thereby reducing the integer $c_1$ in the complexity pair $(c_1,c_2)$.
We have already seen that in the situation where there is exactly one clasp arc, we may use a
microflype to push $\Xhigh$ across the clasp annulus $\ta$.  The situation
which we face now has two factors which make it significantly more complicated.  The first is
that we must allow for the possibility that there are $k\geq 1$ clasp arcs.  The
second is that if $\Xhigh$ has $\mu$ components, then $\ta$ will be the image of $\mu$ annuli
under an immersion, and we must allow for the possibility of clasp intersections between
distinct annuli.

\subsubsection{Constructing the thin annuli} 
\label{subsubsection:forming thin annuli} 
Preparing for the shortening of long clasp arcs, we will
construct a family $\cS$ of `thin annuli' which is a subset of $\pa$.  Normal neighborhoods will play a
key role in the construction.  We will prove:

\begin{lemma}
\label{lemma:construction of cS}  Assume
 that $\cX$ has $\mu$ components.  Suppose that $\pa$ contains  long clasp arcs $\ghigh^i , \ghigh^j, i
\not= j$.  Then there exists a family $\cS$ of $\mu$ annuli, each a subannulus of $\pa$, and each with a component
of
$\Xhigh$ as one of its boundaries, such that every annulus in
$\cS$ is either trivially foliated by $s$--arcs, or is a standard annulus, or has a foliation
satisfying the following (see \fullref{figure:normal7} for examples):
\be
\item [(i)] Each non-trivially foliated annulus contains at least one long clasp arc.  Moreover, all of its clasp arcs
are doubly long and of length 1, with respect to the induced foliation of $\cS$.
\item [(ii)] After an isotopy of $\pa$ which leaves $\Xhigh$ and the other components
of $\partial\cS$ fixed, each $b$--arc in $\cS$ may be
assumed to have at least one of its endpoints on a clasp arc.
\ee
\end{lemma}

\begin{proof}   The case when there is exactly one long clasp arc was just discussed in \fullref{subsection:microflypes}, so we assume that there are at least two long clasp arcs. 
Let $\ghigh^i$ be a long clasp arc, so that
$\glow^i$ begins on
$\Xlow$ and ends near $\Xhigh$.  Assume that $\ghigh^i$ and $\glow^i$ have been pushed into a union of 
leaves in their normal neighborhoods.  This is possible, by the construction in \fullref{proposition:normal}.  Let
$\Nlow^i$ be the normal neighborhood of $\glow^i$.  We focus our attention on a rectangle which we call
$\thintilelow^i$.  See sketches (1) and (2) in \fullref{figure:normal4}.  It is a subset of  $\Nlow^i$, and it
coincides with
$\Nlow^i$ when the clasp arc has length 2.  It has three edges which are in $\partial \Nlow^i$.  We choose
an arc which is  transverse to the foliation as its other `horizontal' boundary.

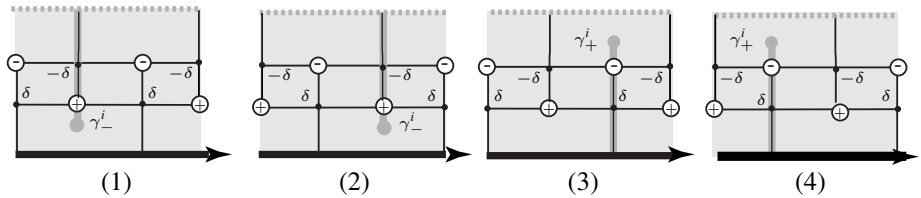
\begin{figure}[htpb!]
\labellist\small
\pinlabel {(1)} at 129 452
\pinlabel {(2)} at 259 452
\pinlabel {(3)} at 384 452
\pinlabel {(4)} at 508 452

\tiny
\pinlabel {$\delta$} [bl] at 74 496
\pinlabel {$-\delta$} [tr] at 107 519
\pinlabel {$\delta$} [bl] at 143 496
\pinlabel {$-\delta$} [tr] at 174 519

\pinlabel {$-\delta$} [tl] at 208 518
\pinlabel {$\delta$} [br] at 238 494
\pinlabel {$-\delta$} [tl] at 274 518
\pinlabel {$\delta$} [br] at 307 494

\pinlabel {$-\delta$} [tr] at 364 516
\pinlabel {$\delta$} [bl] at 332 494
\pinlabel {$-\delta$} [tr] at 431 516
\pinlabel {$\delta$} [bl] at 400 494

\pinlabel {$-\delta$} [tl] at 455 516
\pinlabel {$\delta$} [br] at 487 492
\pinlabel {$-\delta$} [tl] at 521 516
\pinlabel {$\delta$} [br] at 555 494

\pinlabel {$\gamma_-^i$} [l] at 110 485
\pinlabel {$\gamma_-^i$} [l] at 279 484
\pinlabel {$\gamma_+^i$} [r] at 397 530
\pinlabel {$\gamma_+^i$} [r] at 482 530

\endlabellist
\centerline{\includegraphics[scale=.7]{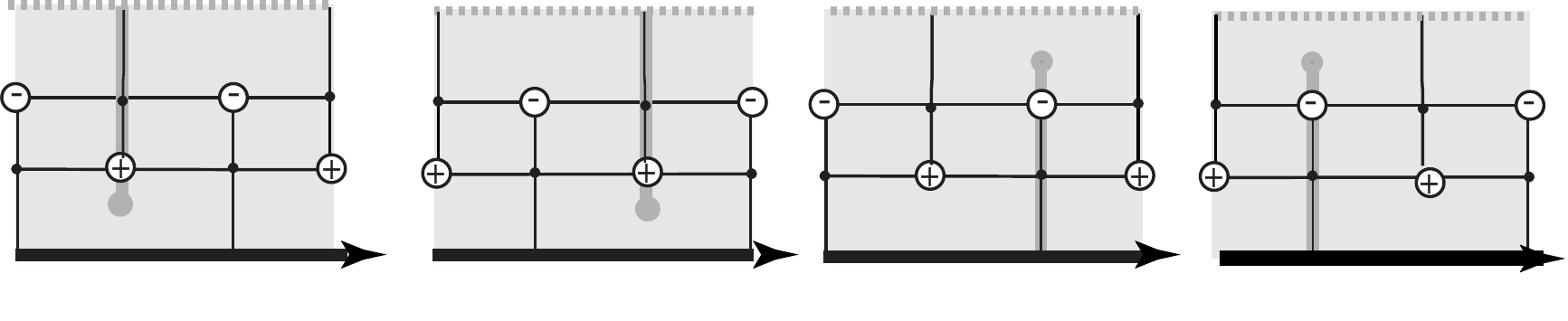}}
\caption{Sketches (1) and (2) illustrate the two possibilities for $\thintilelow^i$.
Sketches (3) and (4) illustrate the two possibilities for $\thintilehigh^i$.}
\label{figure:normal4}
\end{figure}

There are corresponding subrectangles $\thintilehigh^i$ and we show the two
possible arrangements in sketches  (3) and (4) of \fullref{figure:normal4}.  As in the
case of $\thintilelow^i$, the rectangle $\thintilehigh^i$ has 3 boundary edges which are in $\partial
\Nhigh^i$.  We choose its
$4^{th}$ boundary edge to be an arc which is transverse to the foliation, so that the rectangle
has 4 vertices and 4 singularities.

In the special case when $\ghigh^i$ has length 2 the puncture endpoint of the arc
$\ghigh^i$ will be an interior point of $\thintilehigh^i$, but if $\ghigh^i$ has length $\geq 3$
then  $\ghigh^i\cap \thintilehigh^i$ will have both of its endpoints on the boundary.  We correct this by modifying
$\ghigh^i\cap \thintilehigh^i$ to a subarc of
$\ghigh^i$ which has the same image in $\ta$ as the intersection of $\glow^i$ with its induced normal neighborhood.
This will give us a shortened induced arc which (by an abuse of notation) we continue to refer to as $\ghigh$.
It begins on
$\Xhigh$ and ends at a point in the interior of $\thintilehigh^i$.  By construction, $\glow^i$ and the new $\ghigh^i$ have
the same image in $\ta$, and determine the clasp intersection between the rectangles $\thintilelow^i$ and
$\thintilehigh^i$ which
is induced by the clasp intersection in $\ta$ corresponding to $\g^i$.  The four cases which are illustrated in
\fullref{figure:normal4} will be referred to as types (1),(2),(3),(4).

Let $\cN_\cS$ be the union of all of the $\thintilelow^i$ and $\thintilehigh^i$.  Note that, while
$\thintileboth^i \cap (\glow^j \cup \ghigh^j) = \emptyset$ for all $ j\not= i$, there is no
reason why 
$\thintileboth^i \cap \thintileboth^j$ should be empty.  This leads us to the following preliminary
definition of a connected collection of normal neighborhoods:
Choose $\thintile, \thintile^\prime \in \cN_\cS$.  We say that $\thintile$ and $\thintile^\prime$
are {\it connected}, and write $\thintile \connect \thintile^\prime$,
if $\thintile \cap \thintile^\prime \not= \emptyset$ in $\pa$.
Two examples are given in \fullref{figure:normal5}.
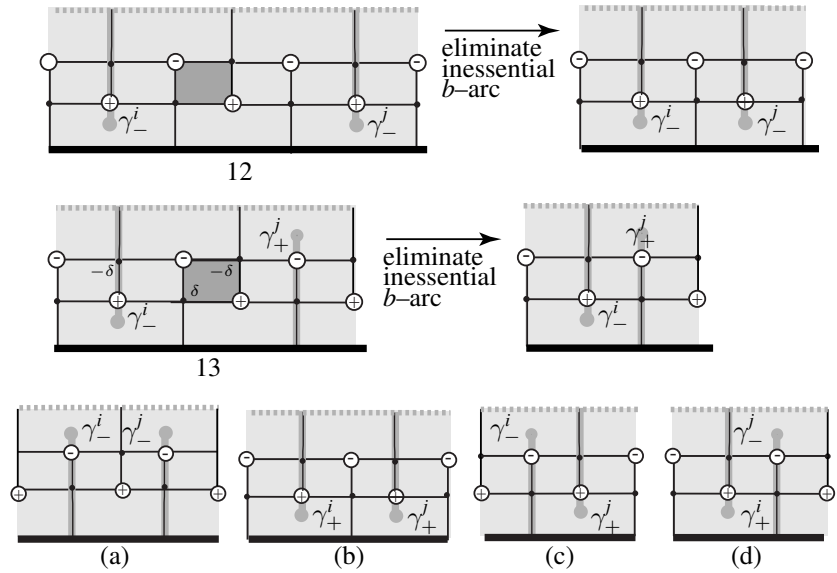
\begin{figure}[htpb!]
\labellist\small
\pinlabel {12} at 229 442
\pinlabel {13} at 212 335
\pinlabel {\parbox{45pt}{eliminate inessential $b$--arc}} [t] at 372 519
\pinlabel {\parbox{45pt}{eliminate inessential $b$--arc}} [t] at 340 405
\pinlabel {(a)} at 161 231
\pinlabel {(b)} at 289 231
\pinlabel {(c)} at 404 231
\pinlabel {(d)} at 506 231

\hair=2pt
\pinlabel {$\gamma_-^i$} [l] at 160 468
\pinlabel {$\gamma_-^j$} [l] at 296 468
\pinlabel {$\gamma_-^i$} [l] at 451 469
\pinlabel {$\gamma_-^j$} [l] at 509 469

\pinlabel {$\gamma_-^i$} [l] at 165 361
\pinlabel {$\gamma_+^j$} [r] at 260 407
\pinlabel {$\gamma_-^i$} [l] at 422 362
\pinlabel {$\gamma_+^j$} [r] at 460 409

\pinlabel {$\gamma_-^i$} [l] at 140 300
\pinlabel {$\gamma_-^j$} [r] at 185 300
\pinlabel {$\gamma_+^i$} [l] at 266 253
\pinlabel {$\gamma_+^j$} [l] at 318 253
\pinlabel {$\gamma_-^i$} [r] at 386 300
\pinlabel {$\gamma_+^j$} [l] at 418 253
\pinlabel {$\gamma_+^i$} [l] at 500 253
\pinlabel {$\gamma_-^j$} [r] at 519 300

\tiny
\pinlabel {$-\delta$} [tr] at 163 393
\pinlabel {$\delta$} [bl] at 200 372
\pinlabel {$-\delta$} [tr] at 228 393
\endlabellist
\centerline{\includegraphics[scale=.7]{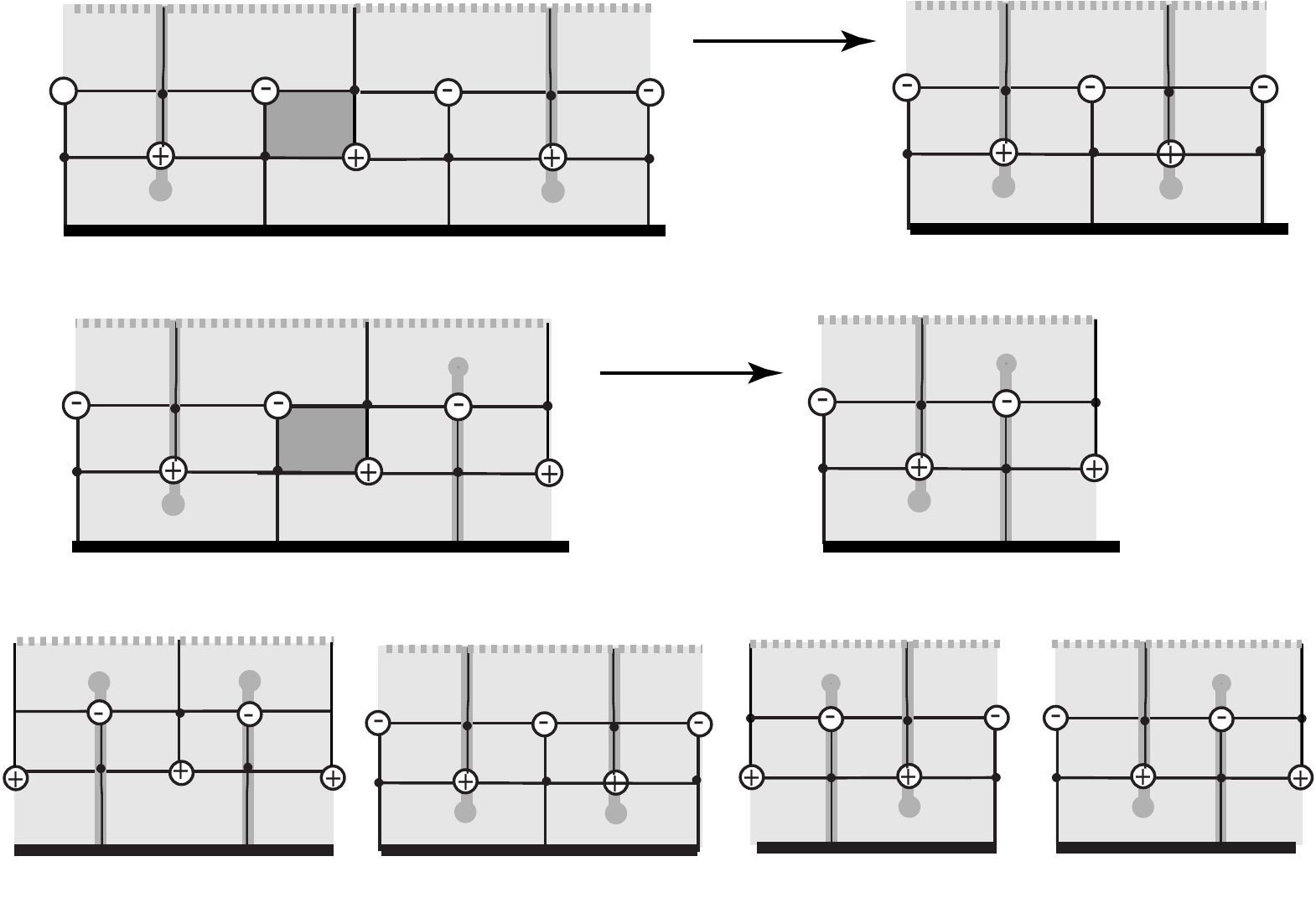}}
\caption{Adjacent pairs of foliated connected normal neighborhoods.  Sketches a,b,c,d in the bottom
row show all possibilities, after the elimination of inessential $b$--arcs.}
\label{figure:normal5}
\end{figure}
In the top row, $\thintile$ has
type 1 and $\thintile^\prime$ has type 2 and $\thintile \cap \thintile^\prime =
\partial \thintile \cap \partial \thintile^\prime$.
In the middle row
$\thintile$ has type 1, $\thintile^\prime$ has type 3 and they intersect along two singular leaves and the disc
between them.  Observe that this is the maximal possible intersection, because the normal
neighborhood of a clasp arc never intersects another clasp arc.  It follows
that the possible sequences in a connected set are 12, 13, 14, 23, 24, 34, 21, 31, 41, 32, 42, 43.  A collection of normal neighborhoods $\{\thintile^1,\ldots, \thintile^p\}\subset
\cN_\cS$ is connected if there is a connecting path between any two neighborhoods in the
collection, ie if
$\thintile^\prime ,\thintile^{\prime\prime} \in
\{\thintile^1, \ldots \thintile^p\}$  then there exists a subcollection
$\{ \thintile^{i_1}, \ldots , \thintile^{i_q} \} \subset \{\thintile^1, \ldots, \thintile^p \}$ such that
$ \thintile^\prime \connect \thintile^{i_1} \connect \thintile^{i_2} \connect \cdots
\connect \thintile^{i_q} \connect \thintile^{\prime\prime} .$  
A connected component of $\cN_\cS$ is called a {\it region} and is denoted by
the symbol $\cR$.

We now observe that each
$b$--arc in each $\cR$ is in a normal neighborhood of some clasp arc.  We distinguish between
two types of $b$--arcs: those whose endpoints are vertices which do not meet a clasp arc, and
those which have at least one vertex endpoint which is on a clasp arc.  Let's look first at the
former.  Examples can be seen in the darkly shaded subrectangles in the left sketches in rows 1 and 2 of \fullref{figure:normal5}.  Recall that our normal neighborhoods were created by the repeated
use of finger moves, which necessarily created some inessential $b$--arcs.  But
all of the $b$--arcs which do not intersect clasp arcs are inessential, and may be deleted
by an isotopy of $\pa$ which is supported on a disc in the interior of the connected region, as in
the passage from the left to the right in \fullref{figure:normal5}.  Therefore we may assume
that $\cR$ contains no such $b$--arcs.  Notice that we have chosen our definition of
complexity so that this modification does not alter the complexity.
With the modification, it is easy to see that there are precisely 4 possible sequences of two
modified normal neighborhoods, as illustrated in the bottom row in \fullref{figure:normal5}, sketches (a),(b),(c),(d).

A region $\cR$ is either an annulus or a rectangle, as illustrated in
\fullref{figure:normal7}.  If it is an annulus, then it satisfies
properties (i)--(ii) of \fullref{lemma:construction of cS}.  Assume it's a rectangle.
The lower horizontal boundary of $\cR$ is a subarc of a component of
$\Xhigh$, and so the connected components $\cR_1,\dots,\cR_q$ associated to any
given component of $\Xhigh$ have a natural cyclic order on
$\Xhigh$.  We would like to use this natural order to join them by bands
of $s$--arcs to
obtain annuli.  However the vertical edges of the rectangles are not $s$--arcs.  The
following observation saves the day: the grey dotted horizontal boundary of each rectangle was
chosen in a rather arbitrary way as an arc which is transverse to the foliation and in the
interior of $\pa$, and if we now modify our choices by replacing the vertical edges of the rectangles in
(a),(b),(c),(d) by the sketches in (e),(f),(g),(h) of \fullref{figure:normal6} we will be in business.

\begin{figure}[htpb!]
\labellist\small
\pinlabel {(e)} at 197 450
\pinlabel {(f)} at 234 450
\pinlabel {(g)} at 312 450
\pinlabel {(h)} at 383 450
\pinlabel {$(\mathrm e')$} at 200 342
\pinlabel {$(\mathrm f')$} at 239 342
\pinlabel {$(\mathrm g')$} at 322 342
\pinlabel {$(\mathrm h')$} at 381 342

\hair=2pt
\pinlabel {$\gamma_+^i$} [l] at 195 514
\pinlabel {$\gamma_+^j$} [r] at 235 514
\pinlabel {$\gamma_-^i$} [r] at 325 470
\pinlabel {$\gamma_-^j$} [l] at 363 470

\pinlabel {$\gamma_-^i$} [l] at 200 363
\pinlabel {$\gamma_-^j$} [r] at 240 363
\pinlabel {$\gamma_+^i$} [r] at 330 410
\pinlabel {$\gamma_+^j$} [l] at 375 410

\endlabellist
\centerline{\includegraphics[scale=.85]{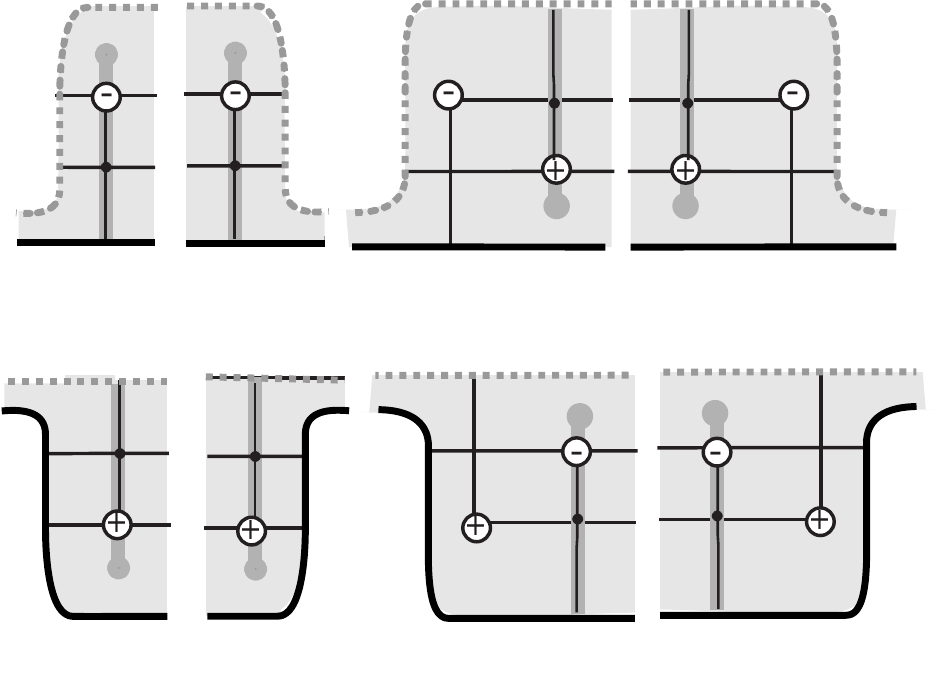}}
\caption{Possible arrangements of clasp arcs in adjacent intersecting
normal neighborhoods in $\pa$, near the right and left
boundaries of a connected collection of normal neighborhoods}
\label{figure:normal6}
\end{figure}
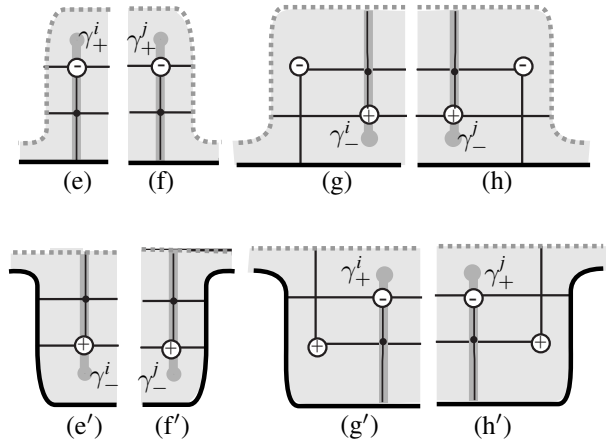
Sketches (e) and (f) are obtained from sketch (a) of \fullref{figure:normal6} by modifying the
grey boundary on the left and right respectively.  Note that the modified grey boundary is
everywhere transverse to the leaves of the foliation.  On the other hand, if we attempt to do the
same thing in the situation of (b), choosing the grey boundary to be close
to the left (resp.\ 
right) clasp arc, a point of tangency with leaves of the foliation will be introduced, so it is
necessary to include the singular leaf which is on the left (resp.\  right) in the modified
connected region, as illustrated in (g) (resp.\  (h)).  We leave it to the reader to check that (e)
and (h) are modifications of (c), and that (f) and (g) are modifications of (d).
In (e)--(h), which are based directly upon the sketches in \fullref{figure:normal4}, all have their rightmost singular leaf in
$G_{-,\delta}$. But there is a second possibility: we are given the foliation near $\Xhigh$, and we do not know
much about it outside the normal neighborhoods.  In fact there are 4 other possible patterns for
the right and left boundaries, also restricted by the fact that $\Xhigh$ must be transverse to
fibers of $\fib$.  They are illustrated in sketches (e$'$),(f$'$),(g$'$),(h$'$) of \fullref{figure:normal6}, and
were obtained from (e),(f),(g),(h) by interchanging the roles of grey and black.

Let's examine the possibilities for the regions.  \fullref{figure:normal7} shows the foliated annulus and the four possible foliated disc regions in
$\cS$, up to the number of pairs of components of $G_{-,\delta}$ and $G_{+,-\delta}$.
In each, we give examples of how the clasp arcs might be placed.  To construct
$\cS$ from the various non-trivially foliated regions we first join disc regions which are
consecutive as one travels along a component of
$\Xhigh$ with bands of $s$--arcs which run between $\Xhigh$ and $\Xlow$.  We may also need some
number of standard annuli which have long clasp arcs.  (See \fullref{figure:normal7} again.)
Finally we may need some number of annuli
which are foliated entirely by $s$--arcs.  This completes the proof of \fullref{lemma:construction of cS}.\end{proof}  

Let
$\cS$  be an annular subset of
$\pa$ which contains $\Xhigh$ (or a component of $\Xhigh$ in the case when $\cX$ is a link) as one of its boundary
components.  We consider the foliation of $\cS$ which is induced by the foliation on $\pa$.  A
(possibly empty) family of $s$--arcs ${\bf S} =
\{ s_1, \dots , s_l: s_i\subset\cS\}$ is a {\it complete collection of
$s$--arcs} in $\cS$ if (i) no two $s$--arcs in the collection split off a sub-band of
$\cS$ that is foliated entirely by $s$--arcs, and (ii) for any other $s$--arc $s
\subset\pa$ there exists an $s_i \in {\bf S}$ such that $s\cup  s_i$ splits off a
sub-band of $\cS$ that is foliated entirely by $s$--arcs.
It is immediate that cutting $\cS$ open along a complete collection ${\bf S}$ of
$s$--arcs decomposes $\cS$ into a disjoint union of {\it thin regions} and
{\it bands of $s$--arcs}.  This construction will be used in what follows.

The grey boundary of the thin annuli will in general be in the interior of $\pa$,
although in special cases it will coincide with $\Xlow$.  
Observe that the difference in braid index $b(\Xhigh)-b(\Xlow)$ across the union of the annuli which make up $\cS$ is
 the total number of positive vertices in $\cS$ minus the total
number of negative vertices in $\cS$.  From \fullref{figure:normal7} we see that this
difference is always either $0,1$ or $-1$ in a single non-trivially foliated region.  Thus the
motion of  $\Xhigh$ across $\cS$ to  $\Xlow$  increases braid index if and only if there are more regions  of type $-1$ than of type $+1$ in  $\cS$.  The regions of type 0 do not affect the count.  We call
our regions types $1_d, -1_d, 0_{d,1}$ and $0_{d,2}$, the subscript indicating that the
region is a disc.  (The two cases of type $0_d$ are distinguished by the placement of
their clasp arcs.)  There is also the special case of the standard annulus.  The standard annulus first appeared 
in \fullref{corollary:standard annulus} as an embedded annulus.  We are now allowing for the
occurrence of clasp arcs.

The previous construction of a thin annular subset $\cS$ in $\pa$ has $\Xhigh$ as one of its boundary components.
But, we can also use the same construction to produce a thin annular subset which has $\Xlow$ as one of its
boundary components.  Let $\cShigh$ (resp.\  $\cSlow$) be the thin annular subset of $\pa$ having $\Xhigh$
(resp.\  $\Xlow$) as one of its boundary components.  Notice that $\cSlow$, like $\cShigh$, will have 
type $0_d, 1_d, -1_d$ regions that are connected by $s$--bands.  However, the dotted grey boundary will be
replaced by a solid grey boundary since it will now be $\Xlow$; and the solid black boundary will be replace
by a dotted black boundary since it may or may not be $\Xhigh$.

Referring to the regions in \fullref{figure:normal7} we consider the cyclic ordering of the singularities
in $\fib$.  We say a type $0_{d,1}$ or $1_d$ (resp.\  $0_{d,2}$ or $-1_d$) region is a {\it fan} if
all of its singularities of parity $-\delta$ (resp.\  $\delta$) occur in sequence in the fibration, followed by all of
its singularities of parity $\delta$ (resp.\  $-\delta$).

We now have a proposition that allows
us to use $b(\Xhigh)$ and $b(\Xlow)$ to limit the occurrences of these regions in $\cShigh$ and 
$\cSlow$.  It will play a key role, later, when we prove the finiteness of the set of templates. (Look ahead to part (6) of \fullref{lemma:growth of V's and E(s)}  in \fullref{section:the proof}.) 

Let $\rho\co \pa \to \ta $ be the natural immersion.

\begin{proposition}
\label{proposition:index count on thin annuli}
The annular regions $\cShigh, \cSlow \subset \pa$ contribute to $b(\Xhigh)$ and $b(\Xlow)$
in the following ways.
\be
\item[\rm{(1)}] The number of type $1_d$ regions in $\cShigh$ is bounded by $b(\Xhigh)$.
\item[\rm{(1$^\prime$)}] The number of type $-1_d$ regions in $\cSlow$ is bounded by $b(\Xlow)$.
\item[\rm{(2)}] If $\cR_1 \subset \cShigh$ is a type $0_{d,1}$ region and $\cR_2 \subset \cShigh$ is a type $0_{d,2}$
region with
$\rho(\cR_1) \cap \rho(\cR_2) \not= \emptyset$, then the pair $(\cR_1, \cR_2)$ contributes at least $+1$
to $b(\Xhigh)$.
\item[\rm{(2$^\prime$)}] If $\cR_1 \subset \cSlow$ is a type $0_{d,1}$ region and $\cR_2 \subset \cSlow$ is a type $0_{d,2}$
region with
$\rho(\cR_1) \cap \rho(\cR_2) \not= \emptyset$, then the pair $(\cR_1, \cR_2)$ contributes at least $+1$
to $b(\Xlow)$.
\item[\rm{(3)}] Suppose $\cR_1, \cR_2 \subset \cShigh$ are regions of type $0_{d.1}$ and $0_{d,2}$ and
that $\cR_3 \subset \cShigh$ is a region of any type that is a fan. Assume that
\be
\item $\rho(\cR_1) \cap \rho(\cR_3) \not= \emptyset$,
\item $\rho(\cR_2) \cap \rho(\cR_3) \not= \emptyset$,
\item $\cR_1$ and $\cR_2$ are adjacent to a common $s$--band $\cR_s$.
\ee
Then the triple $(\cR_1, \cR_2, \cR_s)$ contributes at least $+1$ to $b(\Xhigh)$.

\item[\rm{(3$^\prime$)}]  Suppose $\cR_1, \cR_2 \subset \cSlow$ are regions of type $0_{d.1}$ and $0_{d,2}$ and that
$\cR_3 \subset \cSlow$ is a region of any type that is a fan.  Assume that
\be
\item $\rho(\cR_1) \cap \rho(\cR_3) \not= \emptyset$,
\item $\rho(\cR_2) \cap \rho(\cR_3) \not= \emptyset$,
\item $\cR_1$ and $\cR_2$ are adjacent to a common $s$--band $\cR_s$.
\ee
Then the triple $(\cR_1, \cR_2, \cR_s)$ contributes at least $+1$ to $b(\Xlow)$.
\ee
\end{proposition}

\textbf{Proof of (1)}\qua Consider a type $1_d$ region.  Its boundary
is a union of two $s$--arcs, labeled $s$ and $s^\prime$, and subarcs
$x_+,  x_-$ of $\Xhigh, \Xlow$.  Let $H_\theta, H_{\theta^\prime}$
be the fibers which contain $s, s^\prime$. Then $x_+$ (resp.\  $x_-$)
begins on $H_\theta$, travels around the braid axis some number of
times, staying transverse to fibers, and ends at $H_{\theta^\prime}$.
Since there are more positive than negative vertices in a type $1_d$
region, it follows that $x_+$ makes at least one more circuit about
the braid axis than $x_-$. In particular, $x_+$ makes at least one full
circuit about $\axis$ and so contributes at least 1 to $b(\Xhigh)$.

\textbf{Proof of (2)}\qua   Since $\rho(\cR_1)$ and $\rho(\cR_2)$
intersect along a common clasp arc, we can cut open their preimages in
the thin annulus subset $S$ of $\pa$ along the two components of the clasp
arc.  Paste them together in the unique new way to obtain two new regions
which resemble type $1_d$ and $-1_d$ regions in $S$.  Observe that,
while the new regions look like they are type $\pm 1_d$, their images
in $\ta$ will have the property: the $X_\epsilon$ boundary, $\epsilon =
\pm$, say  $x_\epsilon$ will be a union of two arcs, $x_{\epsilon,1}\cup
x_{\epsilon,2}$, where the final endpoint of $x_{\epsilon,1}$ and the
initial endpoint of $x_{\epsilon,2}$ are on the same fiber, but are at
different points of that fiber.  We call this a `surgery'.  After the
surgery, we may use the same argument that we used to prove (1), to show
that the angular support of $x_{+,1}\cup x_{+,2}$ is $2\pi$ greater than
the angular support of $x_{-,1}\cup x_{-,2}$. Therefore $x_{+,1}\cup
x_{+,2}$ contributes at least 1 to $b(\Xhigh)$.

\textbf{Proof of (3)}\qua We  have the following list of possibilities:  

\textbf{Case (i)}\qua  $\cR_1$ and $\cR_2$ have different types, say $\cR_1$
is type $0_{d,1}$ and $\cR_2$ is type $0_{d,2}$. 
\bi 
\item[(ia)] As we traverse $\Xhigh$ on $\cR_3$  we encounter an
$a_+$--arc that is adjacent to a clasp arc that $\cR_3$ shares with
$\cR_1$ before we encounter an $a_+$--arc that is adjacent to a clasp arc that $\cR_3$ shares with $\cR_2$.  

We surger our three regions $\cR_1$, $\cR_2$ and $\cR_3$, in $\pa$, along the pre-images of the clasp arcs in $\rho(\cR_1) \cap \rho(\cR_3)$ and $\rho(\cR_2) \cap \rho(\cR_3)$.   Since the angular support for $\cR_s$ is either
contained within the angular support of $\Xhigh \cap \cR_3$, or is outside the angular support of $\Xhigh \cap \cR_3$,  the angular length of $\Xhigh \cap (\cR_1 \cup \cR_s \cup \cR_2)$ must be greater than
$2\pi$.  (Conceptually, $\cR_3$ creates an `adequate amount of spacing' between $\cR_1$ and $\cR_2$ so  that the triple $(\cR_1, \cR_2 ,\cR_s)$ contributes $+1$ to $b(\Xhigh)$.)

\item[(ib)] As we traverse $\Xhigh$ on $\cR_3$
we encounter an $a_+$--arc that is adjacent to a clasp arc that $\cR_3$ shares with $\cR_2$ before
we encounter an $a_+$--arc that is adjacent to a clasp arc that $\cR_3$ shares with $\cR_1$. 

We need to consider the positioning of the $s$--band, $\cR_s$, that runs
from $\cR_1$ to $\cR_2$.  If the angular length of $\cR_s$ is greater than $2\pi$ then we are done.  So assume
that the angular length of $\cR_s$ is less than $2\pi$.  Next, we can assume that the angular length of
$\Xhigh \cap \cR_3$ is also less than $2\pi$.

Now, a convenient way of accounting for braid index is to line up the foliations of $\cR_1$, $\cR_2$ and
$\cR_3$ (which are still thought of as being in $\pa$) so that corresponding clasp arcs
line up.  (For examples the reader should jump ahead to Figures \ref{figure:amalgamate2} and
\ref{figure:amalgamate4}). Since our argument is concerned with the measurement of
angular length we need only focus on where we start and end this measurement.  
It is convenient to use the endpoint of a clasp arc in $\cR_1$ (or $\cR_2$)
as our starting point and the endpoint a clasp arc
in $\cR_2$ (or, reversing the interplay, $\cR_1$) as our ending point.
Thus, as a simplifying measure in drawing such a illustration we
can focus in on the portion of the foliation of $\cR_1$ and $\cR_2$ where these
clasp arcs are placed, for each such portion of the foliation has
just two vertices, two singularities and a single clasp arc.
Consider the position of $\cR_s$ in this superimposed picture.   If $\cR_s$
is not next to $\cR_3$ then we can again do a cut-and-paste to produce a
topological annulus that contributes $+1$ to the braid index of $b(\Xhigh)$.  If $\cR_s$ is next to
the foliation of $\cR_3$ then it is easy to see that $\Xhigh \cap (\cR_1 \cup \cR_s \cup \cR_2)$ contributes
$+1$ to $b(\Xhigh)$. In the superimposed image of $\cR_1 \cup \cR_s \cup \cR_2$ in
$\cR_3$ it can be seen that $\Xhigh \cap (\cR_1 \cup \cR_s \cup \cR_2)$  bounds a type $1_d$ region.
Again, the assumption that $\cR_3$ is a fan is not needed.
(In the right sketch of \fullref{figure:amalgamate2}(ii), the short $s$--band that runs from dot $g$ to dot $c'$ is
an example of an $s$--band that would be superimposed in the foliation of the right sketch in
\fullref{figure:amalgamate2}(ii).

\ei

\textbf{Case (ii)}\qua $\cR_1$ and $\cR_2$ are both type $0_{d,1}$ (or both type $0_{d,2}$). Suppose regions $\cR_1,\cR_2$ are of type $0_{d,1}$ and $\cR_3$
is a region that is a fan.  We now employ this procedure of
superimposing the foliations of
$\cR_1$ and $\cR_2$ onto the foliation of $\cR_3$.  We again consider the positioning of $\cR_s$.
As before we assume that the angular support of $\cR_s$ and $\Xhigh \cap \cR_3$ is less than $2\pi$.
When we superimpose $\cR_1 \cup \cR_s \cup \cR_2$ onto $\cR_3$ we will see two copies of a type $0_{d,1}$
region side-by-side.  (Please refer to \fullref{figure:normal7} for a understanding of the labels
$s$ and $s'$.)  We have two sub-possibilities: 

\bi
\item[(iia)] $\cR_s$ runs from the $s'$
label of the left region (which will be $\cR_1$ in our side-by-side image) to the $s$ label of the right region,
which is $\cR_2$. 

Notice that the cyclic ordering of the singularities in $\cR_3$ imposes an order on the
singularities of $\cR_1$ and $\cR_2$ where they intersect $\cR_3$.  Thus, since $\cR_3$ is a fan the
$-\delta$ singularity in $\cR_1$ must occur after the $\delta$ singularity in $\cR_2$.  (Recall we are
making the simplifying assumption that $\cR_1$ and $\cR_2$ only have two singularities.)  So the angular
support of $\cR_s$ must overlap with the angular support of $\partial \cR_1$ forcing
$\Xhigh \cap (\cR_1 \cup \cR_s)$ to have angular
length greater than $2\pi$.  Thus, the triple $(\cR_1,\cR_2,\cR_s)$
contributes $+1$ to $b(\Xhigh)$.

\item[(iib)] $\cR_s$ runs from the $s$
label of the left region, which is $\cR_1$, to the $s'$ label of the right region,
which is $\cR_2$.

We cut-and-paste $\cR_1 ,\cR_2 ,\cR_3$ to produce a type $0_{d,1}$ region
having a portion of $\Xhigh$ in its boundary but sheared in two places.  If we adjoin $\cR_s$ to this
new $0_{d,1}$ region we see that the resulting $\Xhigh$ boundary contributes $+1$ to $b(\Xhigh)$.
\ei  
The proofs of (1$^\prime$), (2$^\prime$) and (3$^\prime$) are almost
identical to the proofs of (1),(2),(3), and we leave them to the
reader.\hfill\qed

\subsubsection{Using flypes to push $\Xhigh$ across $\cS$}
\label{subsubsection:using flypes to push Xhigh across cS}
In this subsection we assume that $\cS$ is a collection of thin (but not standard) annuli.  The case of standard
annuli will be considered separately, in  \fullref{subsection:pushing across a standard annulus}.

To simplify the notation, we use
`{\bf bold-faced type}' for the black boundary and all its auxiliary structures; and `{\rm Roman type}' for the grey
boundary and all its auxiliary structure.  In particular, $\cS$ is
foliated by $s$--arcs, $\ba$--arcs (adjacent
to $\bX$), $b$--arcs and $\ra$--arcs (adjacent to $\rX$).
(We will specify further auxiliary structure in due course.)
Thus, if $\cS \subset \ta$ is thin, then $\partial \cS = \bX \cup \rX$ where $\bX \subset
\ta$ is on the $\Xhigh$ side of $\cS$ in $\ta$, and $\rX \subset \ta$ is on the $\Xlow$ side of $\cS$ in $\ta$.

The main goal of this section is to establish that the isotopy of $\bX$ across $\cS$ to $\rX$ is the result
of a collection of flypes, not necessarily all admissible.
Specifically, given $\cS$ we will produce a recipe
for constructing a template $(\bD,\rD)$ where, via a common braiding assignment to the blocks of $\bD$ and $\rD$,
we have that $\bD$ carries $\bX$ and $\rD$ carries $\rX$.  Our recipe will tell us how to designate
blocks of both flavors (fixed and moving), and how to designate the strands connecting the blocks.
Moreover, we will show how the thin structure
of $\cS$ gives rise to a collection of flypes that carries $\bD$ to $\rD$.

To motivate our work, we begin with examples.  The key to understanding the
moving blocks will be a `block amalgamation' process.  Here we give two
examples which illustrate: (i) how the clasp
arcs of
$\cS$ are used to designate the microblocks first mentioned in
\fullref{subsection:microflypes}; (ii) how
the foliation of $\cS$ is used to amalgamate the microblocks into larger moving blocks; and
how the foliation of $\cS$ is used to flype these moving blocks.

\begin{figure}[htpb!]
\centerline{\includegraphics[scale=.65]{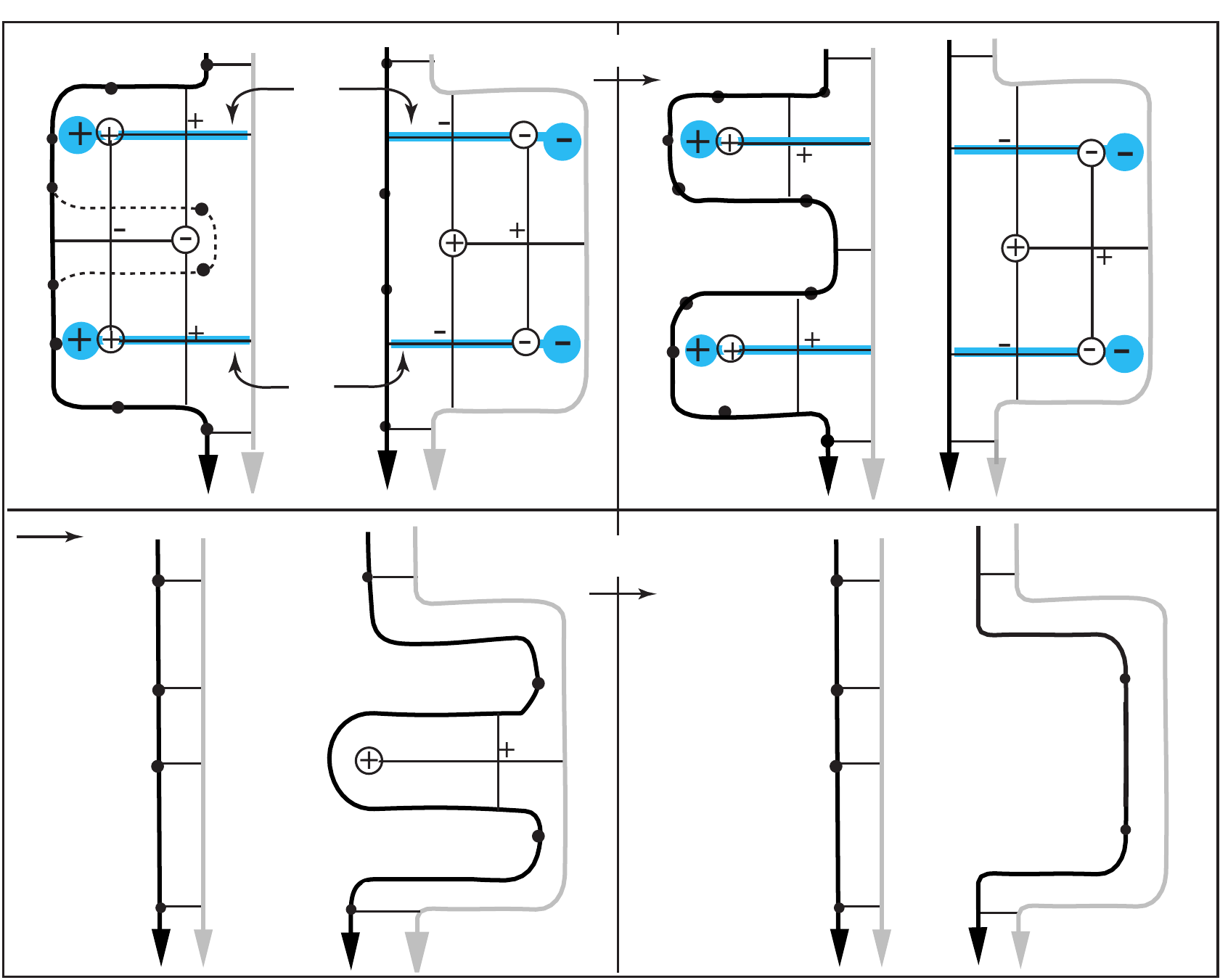}}
\caption{The isotopy of $\Xhigh$ across $\cS$ in a simple example}
\label{figure:amalgamate2}
\end{figure}

\begin{figure}[htpb!]
\labellist\tiny
\hair=2pt
\pinlabel {$a$} [r] at 62 539
\pinlabel {$e$} [r] at 69 516
\pinlabel {$b$} [r] at 68 495
\pinlabel {$d$} [b] at 97 537
\pinlabel {$f$} [bl] at 82 494
\pinlabel {$g'$} [l] at 167 543
\pinlabel {$c'$} [bl] at 170 496
\pinlabel {$d'$} [tl] at 80 484
\pinlabel {$f'$} [t] at 98 444
\pinlabel {$g$} [tl] at 170 481
\pinlabel {$c$} [l] at 165 438
\pinlabel {$a'$} [r] at 68 484
\pinlabel {$b'$} [r] at 66 441
\pinlabel {$e'$} [r] at 69 465
\pinlabel {axis} [r] at 131 488

\pinlabel {$a$} [r] at 192 540
\pinlabel {$e$} [r] at 198 516
\pinlabel {$b$} [r] at 198 495
\pinlabel {$f$} [l] at 212 494
\pinlabel {$d$} [bl] at 226 537
\pinlabel {$g'$} [l] at 296 543
\pinlabel {$a'$} [r] at 198 484
\pinlabel {$e'$} [r] at 199 466
\pinlabel {$b'$} [r] at 193 442
\pinlabel {$d'$} [l] at 209 484
\pinlabel {$f'$} [tl] at 227 444
\pinlabel {$c$} [l] at 294 438

\pinlabel {$a$} [r] at 343 539
\pinlabel {$b$} [r] at 349 494
\pinlabel {$a'$} [r] at 349 481
\pinlabel {$b'$} [r] at 345 440
\pinlabel {$f$} [l] at 359 494
\pinlabel {$d'$} [l] at 357 482
\pinlabel {$g'$} [b] at 445 538
\pinlabel {$c$} [t] at 441 448

\pinlabel {flype} [b] at 451 488

\pinlabel {$a$} [b] at 459 538
\pinlabel {$b'$} [l] at 460 441
\pinlabel {$g'$} [l] at 560 538
\pinlabel {$c$} [l] at 561 442
\endlabellist
\centerline{\includegraphics[scale=.7]{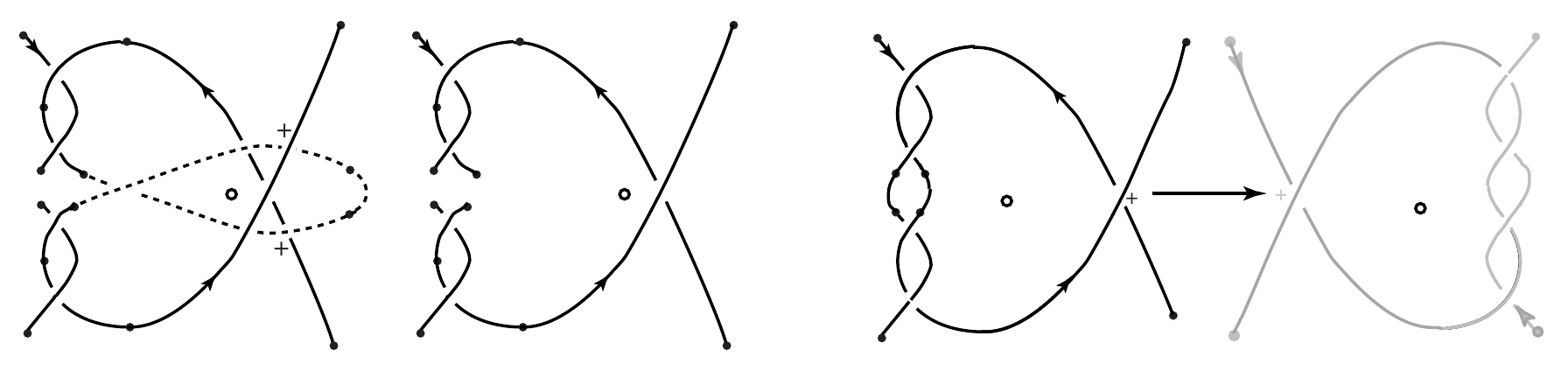}}
\caption{The 3--space embedding of $\Xhigh$ and $\Xlow$ in the Example in \fullref{figure:amalgamate2}}
\label{figure:amalgamate3}
\end{figure}
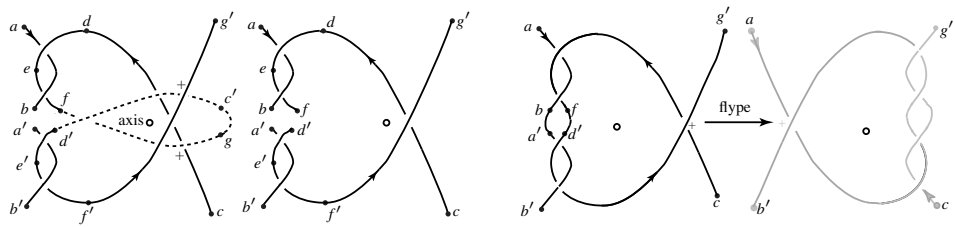

Example 1 is illustrated in Figures \ref{figure:amalgamate2} and \ref{figure:amalgamate3}.
\fullref{figure:amalgamate2}(i) depicts two regions of the preimage of $\cS$ 
sharing two clasp intersections, $\g^1$ and $\g^2$.  For each clasp intersection in $\cS$
there is a readily identifiable pair of tab neighborhoods in the preimage
(just look at what remains after stabilizing
$\bX$ (resp.\  $\rX$) along the singular leaf $\theta_4$
(resp.\  $\theta_3$)).  So for each microflype we have an
associated microblock that contains a positive full twist.  \fullref{figure:amalgamate2}(i)
coupled with the left sketch in \fullref{figure:amalgamate3} illustrates
the strands in the two microblocks:  strands $ab$ and $df$ are associated with the microblock for
$\g^1$; and strands $a^\prime b^\prime$ and $d^\prime f^\prime$ are associated with the microblock for
$\g^2$.

Next, notice that the foliation of $\cS$ yields the following sequence of isotopies for moving $\bX$
to $\rX$: in \fullref{figure:amalgamate2}(i)--(ii) we stabilize $\bX$ along the $\theta_4$ singular
leaf; in (ii)--(iii) we perform two microflypes; and in (iii)--(iv)
we destabilize $\bX$.  The geometric realization of \fullref{figure:amalgamate2}(ii) is depicted
in the left sketch of \fullref{figure:amalgamate3}.  (The alphabetic labeling of points on $\bX$
and $\rX$ are meant to correspond between the two figures and the reader is encouraged to check
the details of this correspondence.)

Now observe that the right two sketches in \fullref{figure:amalgamate3} illustrates an amalgamation
between the two full twists of the microblocks, allowing us to consolidate this
stabilization, microflypes, destabilization sequence into a single flype.  Also, observe that the
two regions in \fullref{figure:amalgamate2}(i) will be fans. This can be verified by checking
the ordering of the angular support of the edgepaths $c d e f d' e' f' g'$ and $a b a' b'$ that are
in $\bX$, along with the corresponding angular support of the analogous edgepaths in $\rX$,
since all singular leaves have an endpoint on the boundary of these regions.

\begin{figure}[htpb]
\labellist\small
\pinlabel {(i)} [b] at 87 304
\pinlabel {(ii)} [b] at 189 304
\pinlabel {(iii)} [b] at 319 304
\pinlabel {(iv)} [b] at 504 304

\tiny
\pinlabel {$a$} [r] at 153 562
\pinlabel {$b$} [r] at 153 461
\pinlabel {$a'$} [r] at 153 430
\pinlabel {$b'$} [r] at 153 367
\pinlabel {$c$} [r] at 110 550
\pinlabel {$d$} [b] at 62 533
\pinlabel {$e$} [r] at 51 508
\pinlabel {$f$} [t] at 63 476
\pinlabel {$g$} [r] at 110 470
\pinlabel {$c'$} [r] at 110 421
\pinlabel {$d'$} [b] at 65 404
\pinlabel {$e'$} [r] at 50 379
\pinlabel {$f'$} [t] at 64 347
\pinlabel {$g'$} [r] at 110 340

\pinlabel {$a$} [r] at 264 524
\pinlabel {$c$} [t] at 411 340
\pinlabel {$b$} [r] at 263 451
\pinlabel {$d$} [b] at 312 519
\pinlabel {$e$} [r] at 272 489
\pinlabel {$f$} [b] at 319 357
\pinlabel {$g$} [b] at 423 525
\pinlabel {$a'$} [r] at 258 419
\pinlabel {$c'$} [t] at 384 340
\pinlabel {$b'$} [r] at 266 347
\pinlabel {$d'$} [b] at 347 512
\pinlabel {$e'$} [r] at 272 379
\pinlabel {$f'$} [t] at 312 348
\pinlabel {$g'$} [b] at 385 527

\pinlabel {$a$} [b] at 448 529
\pinlabel {$b'$} [t] at 447 329
\pinlabel {$c$} [t] at 570 325
\pinlabel {$c'$} [t] at 550 328
\pinlabel {$g$} [r] at 549 488
\pinlabel {$g'$} [l] at 568 488

\endlabellist
\centerline{\includegraphics[scale=.65]{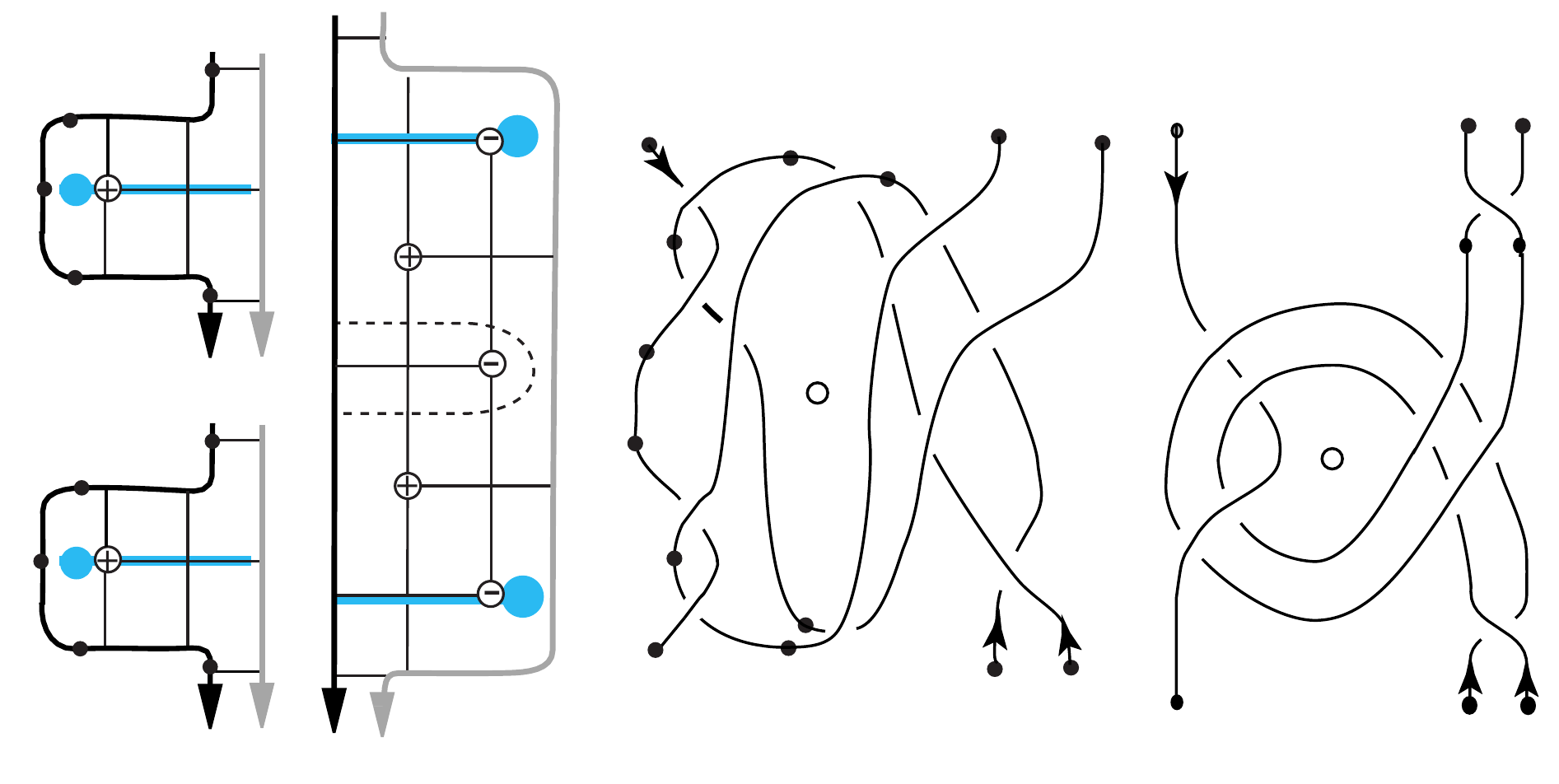}}
\caption{Amalgamating microflypes can result in flypes with weighted strands}
\label{figure:amalgamate4}
\end{figure}
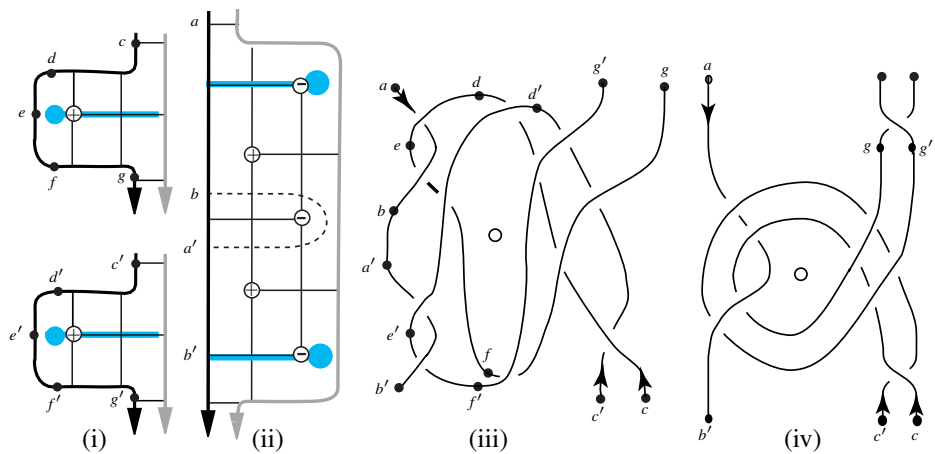

The three regions of $\cS$ for Example 2 are illustrated in \fullref{figure:amalgamate4}(i) and (ii).  Again, there are two clasp intersections
in $\cS$ along with two pairs of tab neighborhoods which clearly delineate the
microflypes and microblocks associated with each clasp intersection.  \fullref{figure:amalgamate4}(iii)
shows the geometric realization of $\bX$: strands $ab$ and $df$ properly contains a full twist in one
microblock; and strands $a^\prime b^\prime$ and $d^\prime f^\prime$ properly contains a full twist in
the another microblock.  And, again the foliation of $\cS$ supplies us with a sequence of stabilizations,
microflypes, destabilization for moving $\bX$ to $\rX$: the starting stabilization is along the
segment $b a^\prime \subset \bX$; perform the two microflypes; and then do the two remaining destabilizations.
The question is, can we amalgamate the two microblocks to achieve a three strand block that would contain the
braiding depicted in \fullref{figure:amalgamate4}(iv)?  The answer is that the information
in the foliation is ambiguous; we have not assigned any angular information to the 
singularities in \fullref{figure:amalgamate4}(i) and (ii).  In \fullref{figure:amalgamate4}(iii)
we have illustrated the point $f$ (resp.\  $d^\prime$) as occurring
after (resp.\  before)
the point $b^\prime$ (resp.\  $d$) in the fibration, however there is no information in \fullref{figure:amalgamate4}(i) and (ii) that forces this choice.  The occurrences might have been
reversed.  If they were reversed then the amalgamation of \fullref{figure:amalgamate4}(iv) would not have been
possible.  Thus, the foliation of \fullref{figure:amalgamate4}(i) and (ii) could depict either one braid
-index-decreasing flype with a three strand block, or two elementary flypes followed by a
destabilization. 

As in Example 1 of \fullref{figure:amalgamate2}, it is easy to see that the region in \fullref{figure:amalgamate4}(ii) will be a fan.  This will be independent of how the
amalgamation of blocks occurs. 

\begin{figure}[htpb!]
\labellist\tiny
\pinlabel {$W$} at 111 675
\pinlabel {$X$} at 108 656
\pinlabel {$Y$} at 108 594
\pinlabel {$Z$} at 105 579
\pinlabel {\rotatebox{180}{$W$}} at 447 567
\pinlabel {\rotatebox{180}{$X$}} at 454 622
\pinlabel {\rotatebox{180}{$Y$}} at 459 638
\pinlabel {\rotatebox{180}{$Z$}} at 457 691
\endlabellist
\centerline{\includegraphics[scale=.75]{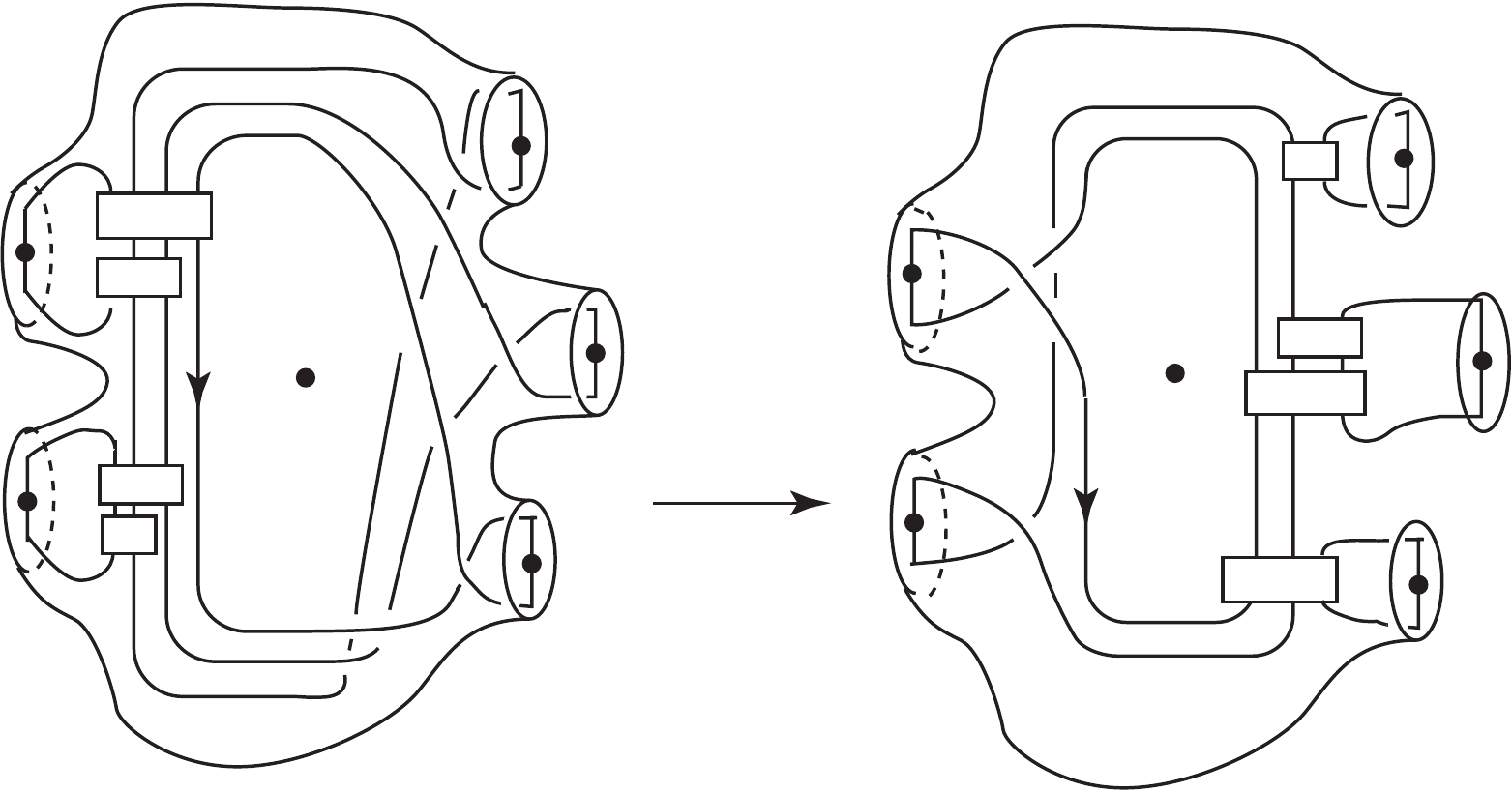}}
\caption{Four independent flypes which cannot be amalgamated}
\label{figure:amalgamate1}
\end{figure}
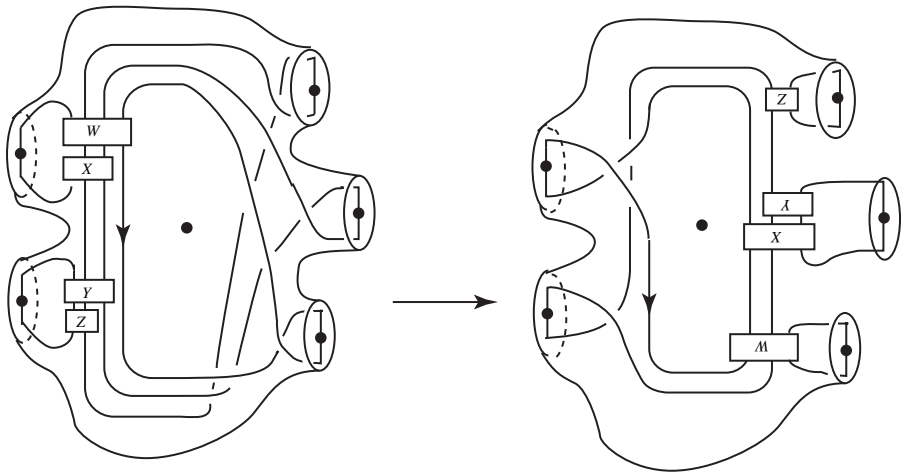

Example 3 (see \fullref{figure:amalgamate1}) shows a situation where successive flypes cannot be amalgamated.
We see four flypes, with associated braid blocks X,Y,Z,W.  In the initial diagram it looks
as if it might be possible to amalgamate $X$ and $Y$, but if we study the final diagram we see
that this is impossible.

\begin{remark}
\label{remark:properties of thin annuli}
Observations based on these examples lead us to the following
remarks about  thin annuli $\cS$:
\bi
\item[{\rm(i)}] The foliation of $\cS$ enables us to readily identify tab neighborhoods of clasp arcs.
\item[{\rm(ii)}] For each pair of tab neighborhoods (where the pairing is via the pairing of the clasp arcs)
we have an associated microflype and microblock.  By our parity conditions on the singularities of
a region (refer to \fullref{figure:normal7}) we know that all of the microflypes associated with
a particular region are of the same parity, ie either all positive flypes or all negative flypes.
In particular, any type $0_{d,1}, 0_{d,2},  1_d ,  -1_d$ region used in a flype will be a fan.
\item[{\rm(iii)}] We can move $\bX$ across the regions of $\cS$ to $\rX$ through a sequence of stabilizations
of $\bX$, microflypes and destabilizations of the resulting new
$\bX$. \hfill$\diamondsuit$
\ei
\end{remark}

Recall that our main goal in this subsection is to establish that the isotopy of $\bX$ across $\cS$ to $\rX$ is the
result of a collection of flypes, not necessarily all admissible.  For that we need to understand when it is
possible to amalgamate a collection of microblocks into a larger block that is moved by a flype.  So let
us finally formalize the definition of `amalgamation'.  To do this we first need to re-characterize microblocks in
term of the foliation of $\cS$.

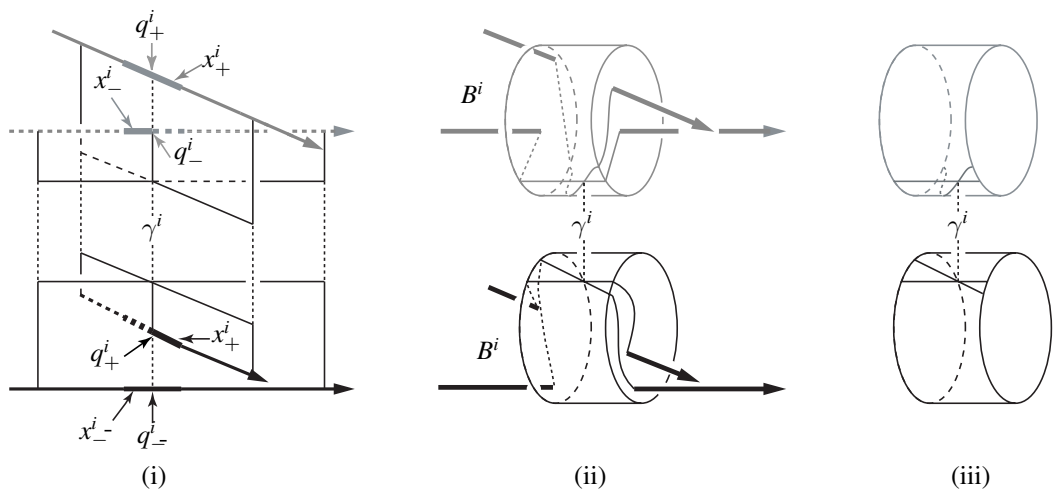
\begin{figure}[htpb]
\labellist\small
\hair=1pt
\pinlabel {$x_-^i$} [b] at 43 493
\pinlabel {$q_+^i$} [b] at 69 528
\pinlabel {$q_-^i$} [l] at 82 458
\pinlabel {$x_+^i$} [l] at 100 516
\pinlabel {$\gamma^i$} at 70 411
\pinlabel {$q_+^i$} [r] at 52 331
\pinlabel {$x_+^i$} [l] at 106 343
\pinlabel {$x_-^i$} [tr] at 44 292
\pinlabel {$q_-^i$} [t] at 70 290
\pinlabel {$B^i$} at 271 496
\pinlabel {$B^i$} at 281 335
\pinlabel {$\gamma^i$} at 340 413
\pinlabel {$\gamma^i$} at 576 413
\pinlabel {(i)} at 71 255
\pinlabel {(ii)} at 345 255
\pinlabel {(iii)} at 584 255
\endlabellist
\centerline{\includegraphics[scale=.6]{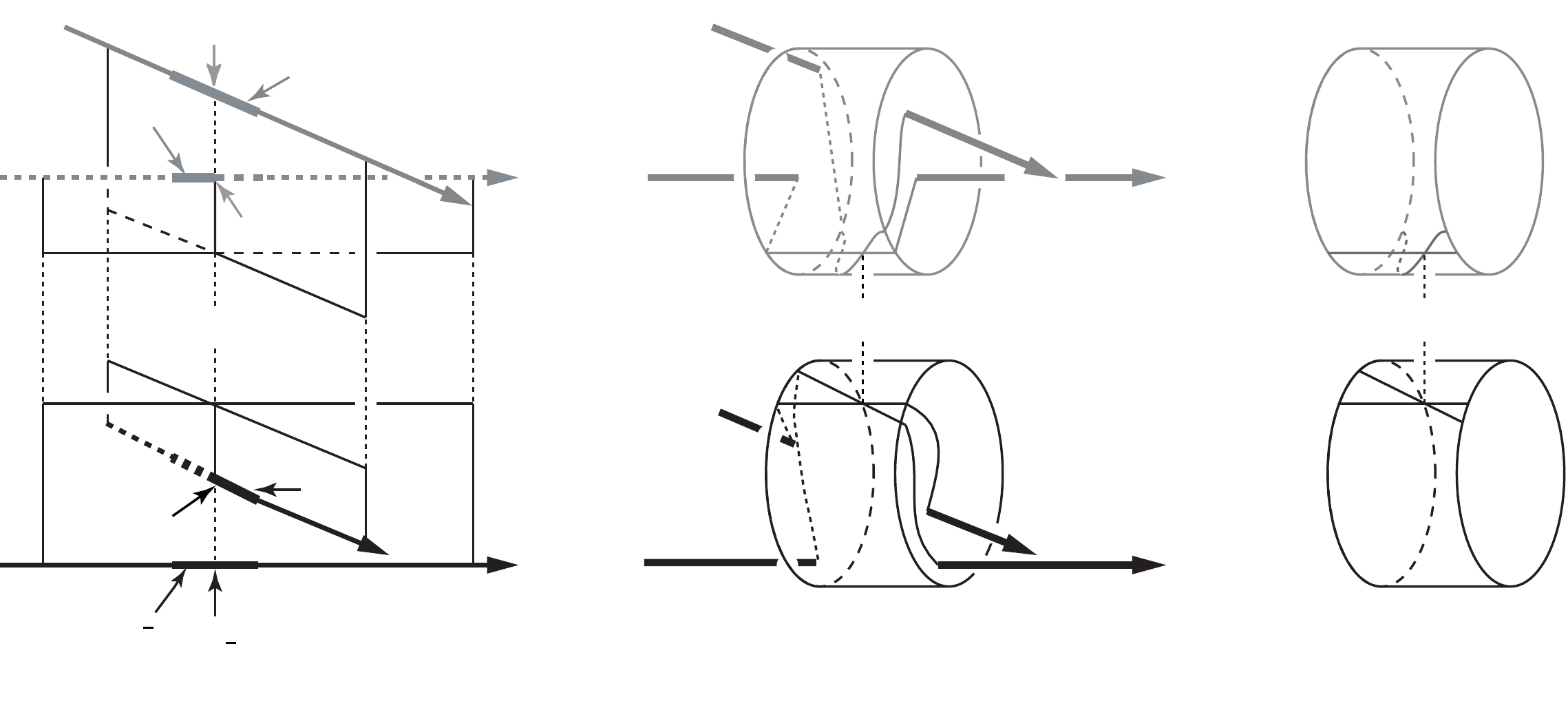}}
\caption{(i) The four microstrands associated to a doubly long clasp arc.  (ii) The
associated microblocks.  (iii) The braid projection of the microstrands on
the cylinder walls.}
\label{figure:amalgamate5}
\end{figure}

By hypothesis, all clasp arcs in $\cS$ are doubly long.  In particular, there are no short clasp arcs.  For
any pair $\ghigh^i, \glow^i$ we have pushed both arcs into a chains of
leaves.  (For thin annuli these chains will always have exactly two leaves.)
Since $\ghigh^i$ begins on $\bX$, one of the endpoints of $\ghigh^i$ is a point
$\qhigh^i\in \bX$.  The fact that $\gamma^i$ is doubly long shows that one of the endpoints
of $\glow^i$ is also near  $\bX$.  The singular leaf which contains it ends at
a point $\qlow^i\in \bX$.  Since our clasp arc is doubly long, there is an analogous
picture near $\rX$.  Both are illustrated in \fullref{figure:amalgamate5}.  The part of $\ta$ that
belongs to the two normal neighborhoods runs all the way from the black boundary to the grey boundary.  The two normal
neighborhoods intersect along the clasp arc.

We fix subarcs $\xhigh_+^i$ and  $\xhigh_-^i$ of $\bX$,
where  $\xhigh^i_+$ (resp.\  $\xhigh^i_-)$) is a closed neighborhood
of $\qhigh^i$ (resp.\  
$\qlow^i$) on $\bX$.  The subarc $\xhigh^i$ is chosen so that it does not
intersect the subarc associated to adjacent singular leaves on $\bX$, also so that the arcs
$\xhigh_+^i$ and $\xhigh_-^i$ have the same angular support
$[\theta_{i_1},\theta_{i_2}]\subset\fib$.  We call these arcs the {\it black microstrands}.
There are, of course, similar thick {\it grey microstrands} 
$\xlow_+^i$ and $\xlow_-^i$ of X.  Notice that the isotopy of $\bX$ across $\cS$
pushes $\bX$ across $\Nhigh^i$ (resp.\  $\Nlow^i$) to $\rX$, so the two black microstrands are mapped to the two grey microstrands.

Construct disjoint solid cylinders $\bB^i, \rB^i$ in 3--space which have the structure of blocks, as defined in
\fullref{section:introduction} of this paper.  These cylinders are, of course, microblocks, as defined in \fullref{subsection:microflypes}.  The microblock $\bB^i$ is foliated by discs and contains the braided arcs $\xhigh_+^i$ and
$\xhigh_-^i$, which meet the disc fibers transversally inside $\bB^i$.  It intersects the fibers of
$\fib$ in the interval $[\theta_{i_1},\theta_{i_2}]$.  We also have, without further work, a grey
microblock.  (See \fullref{figure:amalgamate5}(ii).)

There is an additional and very important feature of the geometry: The fact that our clasp arcs are doubly long
tells us, immediately, that the two tab neighborhoods of the clasp arcs reach all the way from $\bX$ to $\rX$,
joining the black and grey microblocks.  In fact, they intersect the
side boundaries of the black (resp.\  grey) microblocks, and the intersection is a pair of intersecting arcs.  We
label the double points of the black and grey projection with the index of $\g^i$.  (See \fullref{figure:amalgamate5}(iii).)  We can think of these crossed arcs as `local projections' of the microblocks.
The very interesting feature of this projection is that the black and grey projections are joined to one-another
by  the tab neighborhoods of $\ghigh$ and $\glow$.  (This is easy to see in \fullref{figure:amalgamate5}(ii), even though we omitted it to keep the picture as simple as possible).

Construct the black and grey microblocks, $\{\bB^1, \ldots, \bB^r \}$ and  $\{\rB^1, \ldots , \rB^r \}$, one for each clasp arc in $S$.  Let
$\bB
\subset S^3
\setminus \axis$ be a $3$--ball having the structure of a $2$--disc cross an interval, $\D \times [0,1]$.
Decompose  $\partial \bB$  as $\bc \cup \bt \cup \bb$, where  $\bc = \partial\D \times [0,1]$; $\bt = \D \times
\{0\}$; and $\bb = \D \times \{1\}$.  (Our notation was chosen to suggest $\bt$ for `top'; $\bc$ for `cylinder';
and $\bb$ for `bottom'.)  Then $\bB$ is an {\it amalgamating block} for the microblocks $\{\bB^1, \ldots , \bB^r \}$
associated to the microstrands $\{{\bf x}^1_+,{\bf x}^1_- , \ldots {\bf x}^r_+,{\bf x}^r_-\}$ of $\cS$ if the following hold:

\be

\item  Each ${\bB}^i \subset \intr(\bB)$ for $i \leq i \leq r$.

\item Its support, ie $\{ \theta \ | \ \bB \cap H_\theta \not= \emptyset \}$, is a closed interval having
length less than $2\pi$.  Every component of $\bX \cap \bB$ contains a microstrand.

\item  The top $\bt$ and bottom  $\bb$ of $\partial\bB$ are contained
in generic disc fibers of $\fib$.  Moreover, each component of $\bt \cap \cS$
(resp.\  $\bb \cap \cS$) is a subarc of an $\ba$--arc having one endpoint on 
$\partial \bt$
(resp.\  $\partial \bb$) and one endpoint on $\bX$.  The leaves of the
induced foliation on $\bc \subset \partial\bB$ are circles.  The only
non-singular leaves that $\bc$ intersect are $\ba$--arcs and $\bc \cap \bX 
= \emptyset$. 

\item If ${\bf x}^i_+$ (resp.\  ${\bf x}^i_-) \subset  \ \intr (\bB)$, then
its partner ${\bf x}^i_-$ (resp.\  ${\bf x}^i_+)$ is also a subset of \ $\intr(\bB)$. Also,  the dotted segment between each ${\bf q}^i_+$ and ${\bf 
q}^i_-$ is a subset of $\intr(\bB)$.

  \item  Replacing 
$$\{\bB,\bt, \bc, \bb, \ba, \{\bB^1, \ldots, \bB^r \},\{\bx^1_+,\bx^1_-,
\ldots,\bx^r_+,\bx^r_-\} \}$$
by 
$$\{ \rB, \rt, \rc, \rb, \ra, \{\rB^1, \ldots , \rB^r
\},\{\rx^1_+,\rx^1_-, \ldots, \rx^r_+, \rx^r_-\} \}$$ 
we can also define the amalgamation of microblocks $\{\rB^1, \ldots , \rB^r \}$.

 \ee

Our next question is: what conditions need to be met
so that we can say the isotopy across $\cS$ takes $\bB$ to $\rB$?  A hint comes from \fullref{figure:amalgamate5}(iii).  Consider the projection of the microstrands in $\bB^i$ onto its cylinder boundary and
the projection of the microstrands in $\rB^i$ onto its cylinder boundary.  Notice that these two projections
are reflections of each other.  This observation leads us to the following two definitions.

First, let $\bB$ be an amalgamating block. The {\it braid projection} $\pi(\bB)$ of $\Xhigh$ onto $\bc \subset \partial \bB$  is the graph $\bc \cap \cS$, regarded as a subset of $\bc$ together with the clasp arc index labeling of this graph's double
points.  We define $\pi(\rB)$ in a similar fashion.

\begin{figure}[htpb]
\labellist\small
\hair=1pt
\pinlabel {$\gamma^1$} [r] at 154 408
\pinlabel {$\gamma^2$} [l] at 283 408
\pinlabel {$\gamma^1$} at 410 417
\pinlabel {$\gamma^2$} at 465 417
\endlabellist
\centerline{\includegraphics[scale=.6]{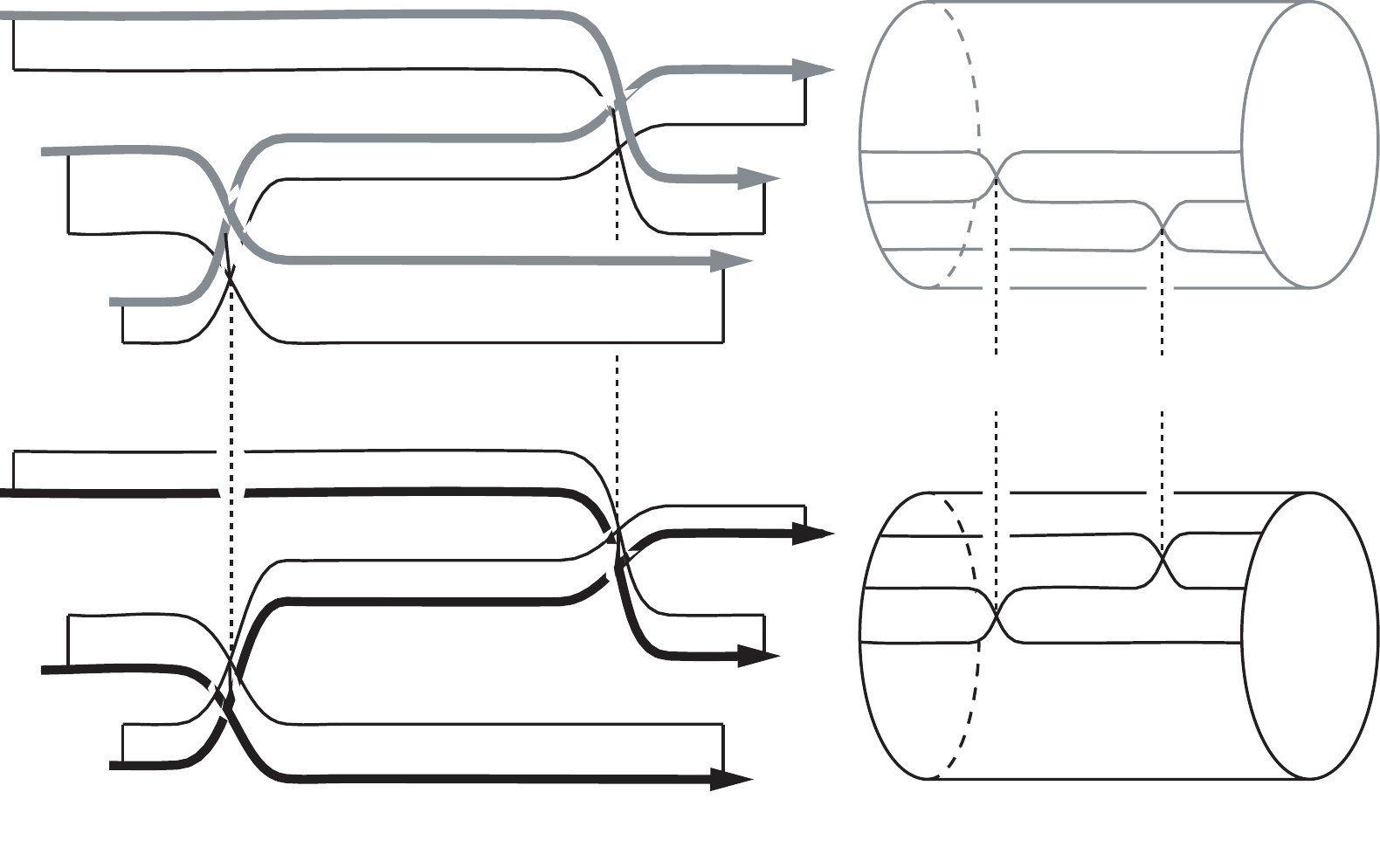}}
\caption{The projection of related microstrands onto the cylindrical walls in the boundaries of amalgamating blocks}
\label{figure:amalgamate6}
\end{figure}
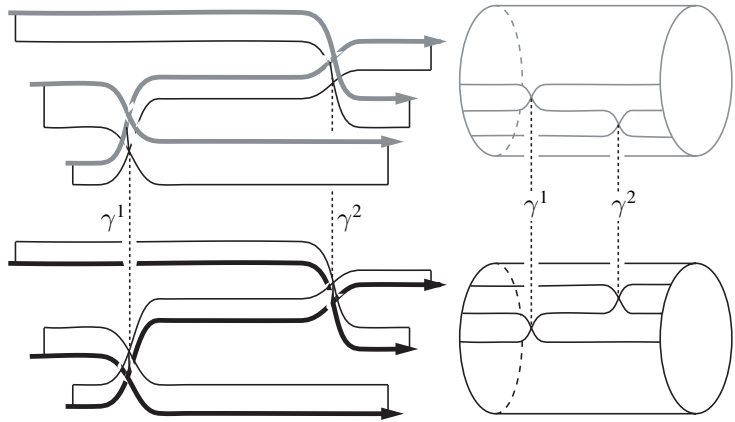

Second, amalgamating blocks $\bB$
and $\rB$ are a pair of {\it related} amalgamating blocks if the following hold:
\bi
\item[(a)] $\{\bB^1,\ldots,\bB^r\}$ and $\{\rB^1,\ldots,\rB^r\}$ are microblocks associated with
clasp intersections $\{\g^1,\ldots,\g^r\}$.
\item[(b)] $\bB$ is a block amalgamation of microblocks $\{\bB^1,\ldots,\bB^r\}$.
\item[(c)] $\rB$ is a block amalgamation of microblocks $\{\rB^1,\ldots,\rB^r\}$.
\item[(d)] $\pi(\bB)$ and $\pi(\rB)$ are reflections of each other.
\ei
Recall \fullref{figure:amalgamate1} and see \fullref{figure:amalgamate6} for an example.

\begin{proposition}
\label{proposition:sequence of flypes}
Let $\cS$ be a collection of thin annuli with boundary components $\bX$ and $\rX$.
Then the motion of $\bX$ to $\rX$ may be realized by a collection of flypes and destabilizations.
\end{proposition}

\begin{proof}  We will construct a pair of candidates for block-strand
diagrams $(\bD,\rD)$ such that $\bD$ (resp.\  $\rD$) carries $\bX$
(resp.\  $\rX$) and the motion $\bD\to\rD$ is a sequence of flypes.  We say that our diagrams  $(\bD,\rD)$ are merely  {\it candidates} for block-strand diagrams is that condition (4) in the definition of a block-strand diagram may not hold.  Indeed, it may be difficult to verify whether  (4)  does or does not hold at this stage of the argument.  For example, it's possible that we could amalgamate blocks, and so produce new blocks which have full braid index. That is a matter which we are not ready to discuss at this time.  (Look ahead to \fullref{subsection:the braid index of blocks in a template} for more on this. See, in particular, \fullref{corollary:no spanning blocks}).  However, if we can construct $\bD$ and $\rD$ so that they carry $\bX$ and $\rX$, and prove that the motion $\bD\to\rD$ is a sequence of flypes and destabilizations, we will have accomplished our task. 

To construct  $\bD$ and $\rD$ we need to understand several aspects of their structure.  Observe that the strands
of $\bX$ which change position during the passage to $\rX$ are precisely those which bound the `tiled' part
of $\ta$, that is the part which is away from the bands of $s$--arcs. On the other hand, the part that stays fixed is the part
that intersects the bands of $s$--arcs.  So we can consider, separately, the `moving' and `fixed' parts of 
$(\bD,\rD)$: 

\textbf{Moving blocks}\qua An amalgamating block $\bB$ will be moved to an
amalgamating block $\rB$ if and only if they are related amalgamating blocks, as defined above.  Note that there
may be some choices involved when we select the amalgamating blocks.  We make those choices in such a way
that the set of all moving blocks has minimal cardinality.

\textbf{Moving Strands}\qua The moving strands in $\bD$ are all subarcs of $\bX$ which
are not amalgamated into moving blocks and are not in the bands of $s$--arcs.

\textbf{Fixed blocks}\qua Here is the intuitive idea: The fixed blocks are associated to braiding between the
strands of $\bX$ (and so also of $\rX$) in the part of $\cS$ which is
foliated by bands of $s$--arcs.  The
braiding between bands of $s$--arcs is not part of our geometry.  $\cS$ only detects the places where
$\bX$ and
$\rX$ differ in a non-trivial way, but in regions where there are bands of
$s$--arcs they essentially coincide.
However, what might happen is that there is braiding between the bands of
$s$--arcs, but that the braiding is
interrupted by, for example, strands of
$\bX$ which separate two potential blocks.  In order to determine the angular regions where the blocks occur,
we therefore look for singularities in the tiled part of $\cS$.  The block subdivision so-obtained will be too
extensive, and we will then need to amalgamate fixed blocks.

We consider the $H_\theta$--sequence for $\cS$.
Let $\cS_{\text{tiled}} \subset \cS$ be the portion of
$\cS$ that is {\em not} foliated by $s$--arcs.  Let $\{\theta_0,\theta_1,\ldots,\theta_\zeta\} \subset [0,2\pi)$ be
a cyclic listing of all the angles at which the corresponding $H_{\theta_i} \in \fib$ contains a singularity,
ordered according to their natural cyclic order in $\fib$.  Since  the bands
of $s$--arcs are foliated without singularities, they must be located in
$$ H' = \bigcup^{i=\zeta}_{\theta\in (\theta_i,\theta_{i+1}), \ i=0 } \{ \
H_\theta \ \setminus \  H_\theta \ \cap \cS_{\text{tiled}} \} .$$ 
Each component of $H'$ is a disc, and it contains only $s$--arcs.  Let $C$ be
a connected component of
$H'$.  Then
$C$ has a $\D \times [\theta_i , \theta_{i+1}]$ structure, and it may
contain some number of bands of $s$--arcs which
braid with one-another inside $C$.  If so, we amalgamate this braiding of
$s$--arc bands in $C$ into a single block
$B(C)$, which will be a fixed block in the template $(\bD,\rD)$.  If $C$
has no $s$--arcs then
$B(C)$ is vacuous.  If $C$ has a single $s$--arc then $B(c)$ is a single fixed strand.

Suppose that there is another connected component $C'$ in $H'$, with its fixed block $B(C')$, and suppose further 
that $\{s-{\rm arcs}\} \cap C' \cap H_{\theta_{i+1}} \subset \{s-{\rm arcs}\} \cap C\cap H_{\theta_{i+1}}$.
If this happens, the singularity at
$\theta_{i+1}$ could have been ignored.  Another way to say this is that we can amalgamate the blocks
$B(C)$ and
$B(C^\prime)$ into a larger fixed block of the template $B(C \cup C^\prime)$ of $(\bD,\rD)$.  There may be some
choices involved, and we make them so that the set of all fixed blocks has minimal  cardinality.

\textbf{Fixed strands}\qua It may happen that we have two
blocks,
$B(C)$ and $B(C')$, for which $(C\cap H_{\theta_{i+1}})\cap (C'\cap H_{\theta_{i+1}}) \not=\emptyset$, but 
$(C\cap H_{\theta_{i+1}})\not= (C'\cap H_{\theta_{i+1}})$, so that the amalgamating condition fails.  In this
situation there are strands that run between $B(C)$ and $B(C')$, however
not all of the $s$--bands
coming out of $B(C)$ go into $B(C')$, and/or not all of the strands going into $B(C')$ are strands
that emerge from $B(C)$.  Any braiding between the $s$ bands that run between $B(C)$ and $B(C')$ can obviously be
pushed into either $B(C)$ or $B(C')$ (eliminating $H_{\theta_{i+1}}$ as a singular fiber) and  they will then
become fixed strands.

It is clear from the construction of $(\bD,\rD)$ that the template $(\bD,\rD)$ carries $(\bX,\rX)$.

We still need to show that for each moving block pair $(\bB , \rB)$
we have a flype (possibly with weighted strands) taking $\bB$ to $\rB$ as in \fullref{figure:flype2}(a).
To do this we first need to isolate each moving black block. 
Specifically, we wish to subdivide the foliation of $\cS$ so that there is a subannular region
$\cS^\prime \subset \cS$ such that if ${\bf y} \subset \bX$ is a strand of $\bD$
then ${\bf y}$ intersects $s$--arcs in the foliation of $\cS^\prime$.  We refer to
\fullref{figure:normal8} to see how $\cS^\prime$ can be obtained through the
stabilization of $\rX$.  Suppose we have two clasp intersections $\g^i$ and $\g^j$ where
the black microblock associated with $\g^i$ is in an amalgamated black block $\bB(i)$ and
the black microblock associated with $\g^j$ is in an amalgamated black block $\bB(j)$.
\begin{figure}[htpb!]
\labellist\small
\pinlabel {(a)} [r] at 152 523
\pinlabel {(b)} [r] at 152 444
\pinlabel {(c)} [r] at 152 374
\pinlabel {(d)} [r] at 152 306
\tiny
\pinlabel {$g_+^i$} [l] at 210 544
\pinlabel {$g_+^j$} [l] at 257 544
\pinlabel {$g_+^i$} [r] at 341 543
\pinlabel {$g_+^j$} [l] at 388 543
\pinlabel {$g_-^i$} [l] at 198 430
\pinlabel {$g_-^j$} [l] at 244 430
\pinlabel {$g_-^i$} [l] at 340 430
\pinlabel {$g_-^j$} [l] at 461 430
\pinlabel {$-\delta$} [tr] at 340 458
\pinlabel {$g_+^i$} [r] at 213 395
\pinlabel {$g_-^j$} [l] at 237 357
\pinlabel {$g_+^i$} [r] at 327 400
\pinlabel {$g_-^j$} [l] at 422 359
\pinlabel {$g_-^i$} [l] at 218 284
\pinlabel {$g_+^j$} [l] at 244 327
\pinlabel {$g_-^i$} [l] at 349 283
\pinlabel {$g_+^j$} [l] at 436 323
\pinlabel {$-\delta$} [tr] at 350 313
\endlabellist
\centerline{\includegraphics[scale=.7]{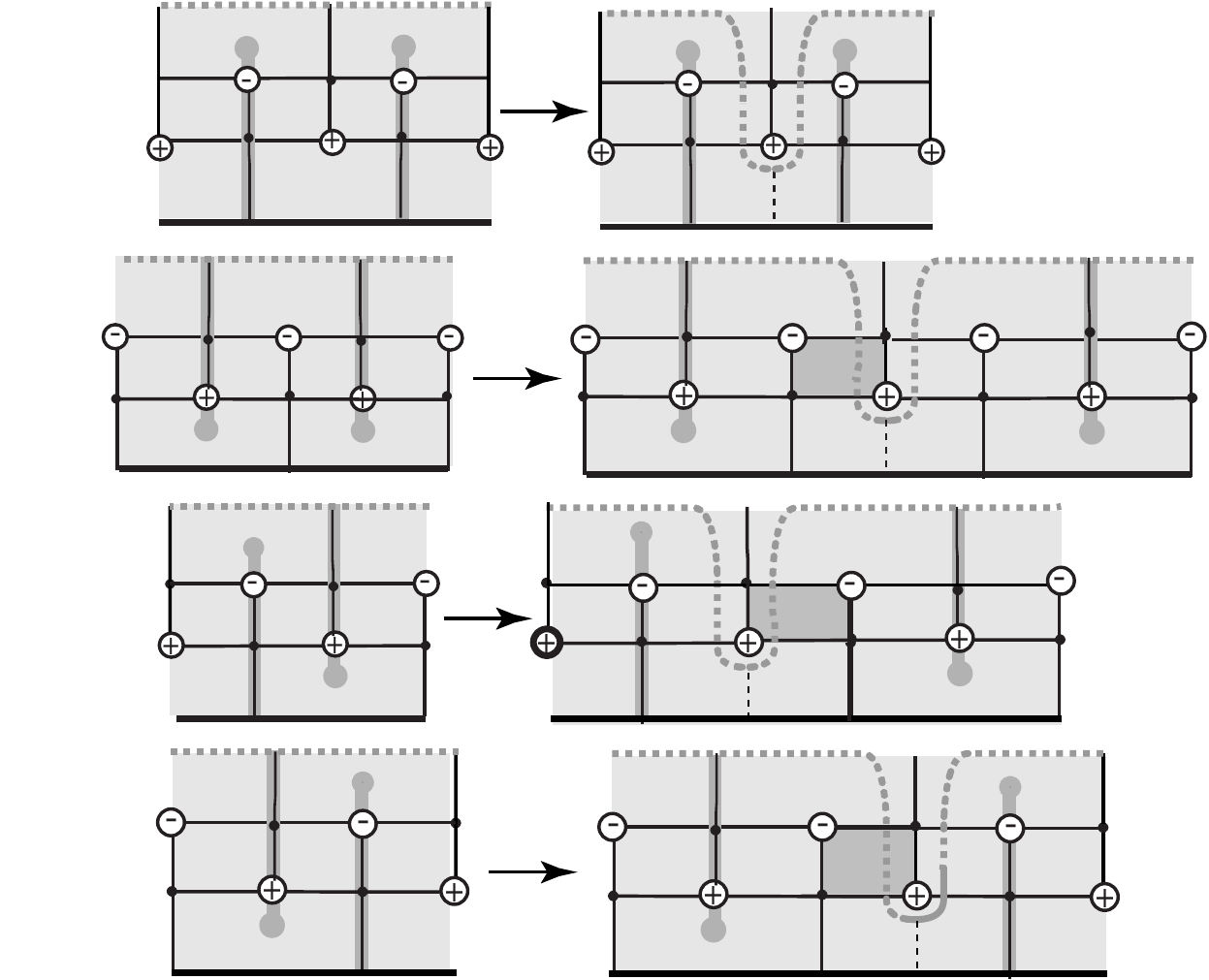}}
\caption{Subdividing non-trivially foliated regions in $\cS$. The black
(resp.\  dotted grey) horizontal arcs are the arcs we refer to in
the text as `${\bf y}$' (resp.\  `y'). The intersections of the clasp arcs with the regions of interest are labeled $g_\pm^i$ and $g_\pm^j$.}
\label{figure:normal8}
\end{figure}
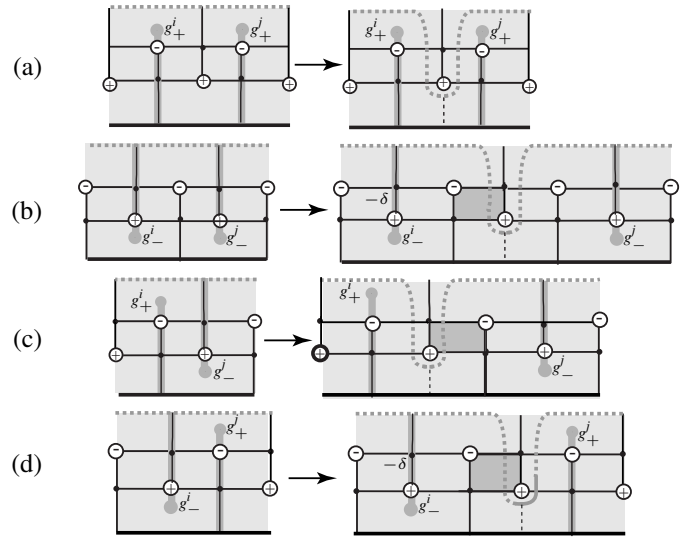
Let ${\bf y} \subset \bX$ be a strand of $\bD$ that has its endpoints on $\bB(i)$ and $\bB(j)$.
Since ${\bf y}$ is a strand it cannot intersect any microstrands of $\bX$.  Thus, it either
intersects $s$--arcs in the foliation of $\cS$ (in which case no subdivision is necessary), or
the only non-singular leaves it intersects are $\ba$--arcs.  To describe the needed
subdivision in the latter situation, we
look at the corresponding grey microblocks associated with our two clasp intersections;
the corresponding grey amalgamating blocks $\rB(i)$ and $\rB(j)$; and the corresponding grey strand
$ {\rm y} \subset \rX$ which has endpoints on $\rB(i)$ and $\rB(j)$.  By assumption, the only
non-singular leaves intersected by the strand ${\rm y}$ are $\ra$--arcs in the foliation of $\cS$.  If
${\rm y}$ intersects an $ab$--singular leaf (as illustrated in \fullref{figure:normal8}(a)) then we
can stabilize $\rX$ along this singular leaf to produce a subannulus of $\cS$ that has ${\bf y}$
intersecting $s$--arcs.  If the strand y does not intersect an $\ra b$--singular leaf then we can
enlarge the foliation of $\cS$ thru the addition of two vertices and singularities, as illustrated in
the three remaining sequences, \fullref{figure:normal8}(b)--(d).  (The new
$b$--arc will necessarily be
inessential.)  The strand y will now intersect the endpoint of an $\ra b$--singular leaf and a stabilization
of $\rX$ is possible.  The corresponding black strand ${\bf y}$ will then
intersect $s$--arcs in the resulting
$\cS^\prime$, isolating the block $\bB(i)$ from the block $\bB(j)$.

\begin{figure}[htpb!]
\labellist\tiny
\pinlabel {\rotatebox{90}{stabilization}} [r] at 120 534
\pinlabel {flype 1} [t] at 206 555
\pinlabel {braid isotopy} [t] at 302 555
\pinlabel {flype 2} [t] at 403 558
\pinlabel {\rotatebox{270}{destabilization}} [l] at 490 545
\pinlabel {\small flype} [t] at 310 481
\pinlabel {$P$} at 132 485
\pinlabel {$Q$} at 132 472
\pinlabel {$P$} at 134 606
\pinlabel {$Q$} at 132 583
\pinlabel {\rotatebox{180}{$P$}} at 274 611
\pinlabel {$Q$} at 233 589
\pinlabel {\rotatebox{180}{$P$}} at 381 589
\pinlabel {$Q$} at 335 588
\pinlabel {\rotatebox{180}{$P$}} at 478 599
\pinlabel {\rotatebox{180}{$Q$}} at 477 612
\pinlabel {\rotatebox{180}{$P$}} at 478 475
\pinlabel {\rotatebox{180}{$Q$}} at 479 489
\endlabellist
\centerline{\includegraphics[scale=.8]{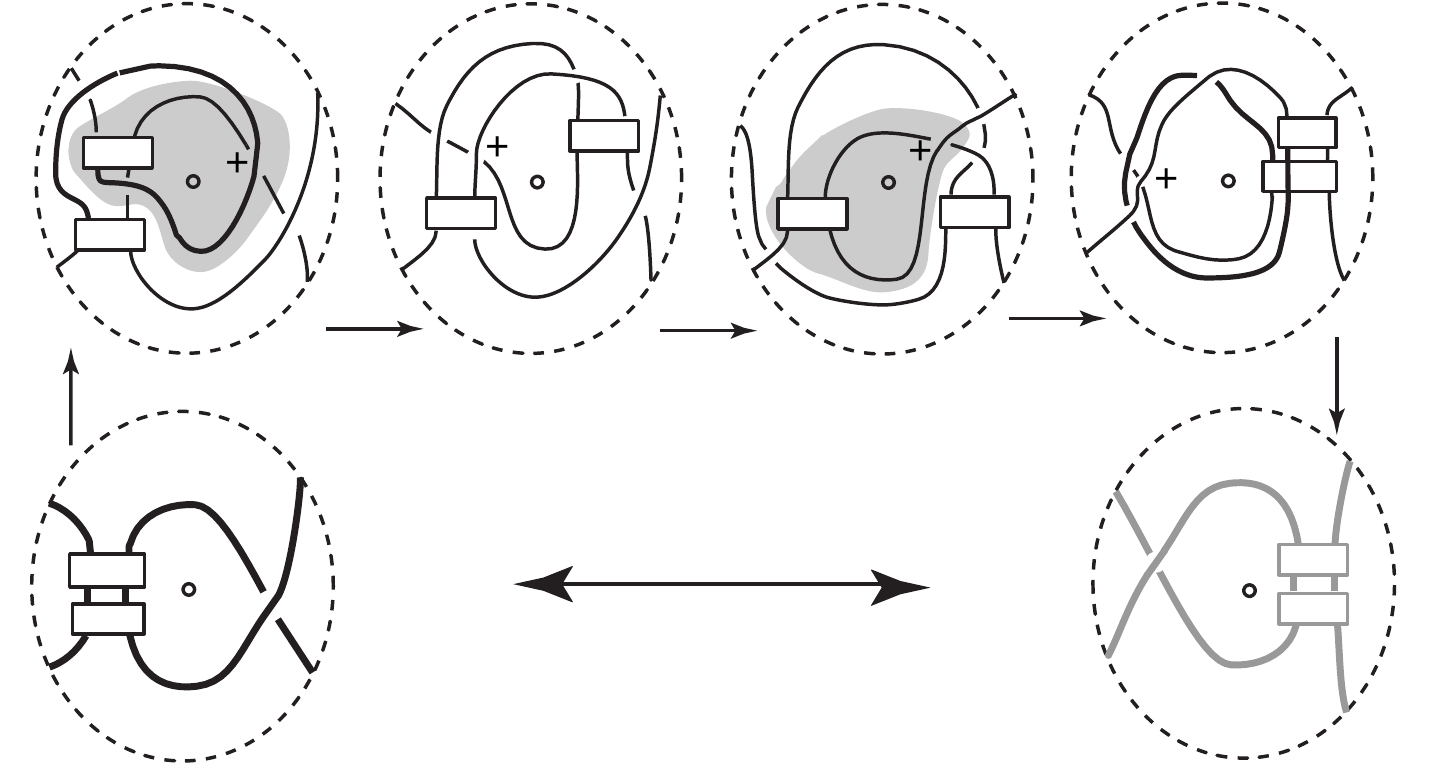}}
\caption{An example that shows how a stabilization, two microflypes and a destabilization are
amalgamated into a flype. The bottom left and right sketches show the corresponding subsets of 
${\bf  D}$ and $\rD$. The top right sketch illustrates the part of the diagram that we refer to as $\rD^\prime$ in the text.}
\label{figure:amalgamate0}
\end{figure}
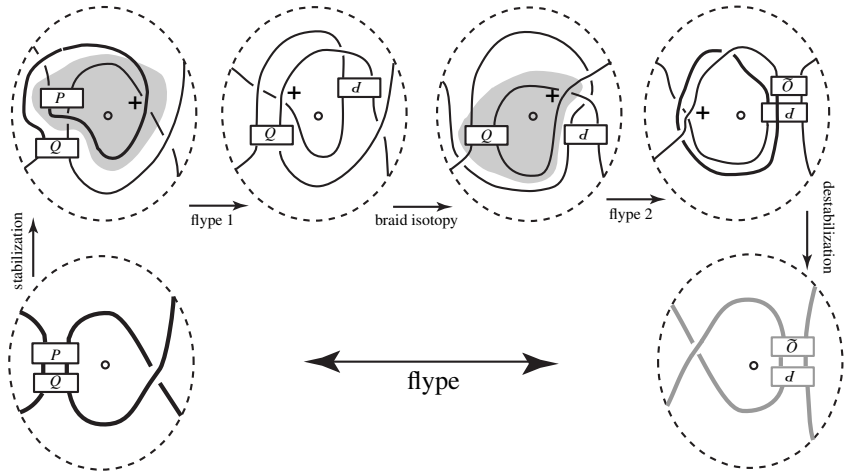

Let $\cS^\prime \subset \cS$ be thin subannuli that isolate all the moving blocks.  Then there is a pair $(\bD,\rD^\prime)$ where $\rD^\prime$ is obtained from $\rD$ by the
stabilizations in the subdivision of \fullref{figure:normal8}. Part of it is illustrated in \fullref{figure:amalgamate0}.  
Saying the moving blocks of
$(\bD,\rD^\prime)$ are isolated is saying that every strand of $\bD$ or
$\rD$ intersects $s$--arcs
in the foliation of $\cS^\prime$.  Now, if we can show that moving across $\cS^\prime$ is a collection of flypes,
then it will follow that the destabilizations come from going from $\rD^\prime$ to $\rD$.

To show that the isotopy across $\cS^\prime$ takes $\bD$ to $\rD^\prime$, using a collection of flypes,
we proceed as follows: since the blocks are isolated, we need only consider what the isotopy is for one block.  This
means we are back in a situation similar to that of Example 1 and Example 2 at the beginning of this section.
We have a subcollection of regions $\{\cR_1, \ldots , \cR_h\} \subset \cS^\prime$
in the foliation of $\cS^\prime$ (as defined in 
\fullref{subsubsection:forming thin annuli}) that are grouped together by clasp intersections.
These regions are all either type $0_{d,1}, 0_{d,2}, 1_d$ or ${-1}_d$.
All of the associated microblocks of the clasp intersections of this collection of regions can
be amalgamated into a block pair $(\bB, \rB)$.  Our claim is that the isotopy of $\bB$ across
the regions $\{\cR_1, \ldots , \cR_h\}$ to $\rB$ is a flype.  By statement (ii) of \fullref{remark:properties of
thin annuli} we know that we can subdivide the microblocks in $\bB$ and flype them across to the
microblocks in $\rB$.  Also, the fact $\pi(\bB)$ and $\pi(\rB)$ are reflections tells us that we can re-amalgamate
them into $\rB$.  In the simplest case of a subdivision into two
blocks (\fullref{figure:amalgamate0})
we see that the entire isotopy of statement (ii) of \fullref{remark:properties of thin annuli}
can be consolidated into a single flype.  This situation can obviously be generalized into a larger number of
blocks in the subdivision by an iterated nesting, the first iteration being:
before doing `flype 1' we subdivide block P into two new
blocks and, since block P in the upper-left sketch could have been flyped,
perform the entire sequence on these two new blocks first.
Thus, $\bB$ is flyped to $\rB$ and the proof of \fullref{proposition:sequence of flypes} is complete.\end{proof} 
\begin{example}
The projection criterion is essential to choosing block amalgamations that maintain their
integrity during a flype isotopy.  This is clear from the example that we gave earlier,  in \fullref{figure:amalgamate1}.  Before the flype it looks as if $X$ and $W$ (and also $Y$ and $Z$) can be amalgamated, but after the flype that is impossible.  Indeed,  the requirement that the projections be reflections of each other will be violated if we
attempt to amalgamate $X$ and $W$, or $Y$ and $Z$. \hfill $\diamondsuit$
\end{example}

\subsection{Constructing a candidate for a template $(\Dhigh,\Dlow)$ from $\ta$}
\label{subsection:constructing the template from the total annulus}

We interrupt the flow of the argument, briefly.  The procedure that we gave in \fullref{subsubsection:using flypes to push Xhigh across cS} 
for constructing a candidate for a template from thin annuli is almost identical to the
procedure for constructing a candidate for a template $(\Dhigh,\Dlow)$ from the foliation of $\ta$.  We will need to know how to do it when we reach the end of \fullref{section:the proof}, and since the alert reader  will have the essential ideas in mind right now, this seems like a good moment to go through the steps in the construction.  The argument will be repetitious, because it is modeled on that in \fullref{subsubsection:using flypes to push Xhigh across cS}, however it seems better to be repetitious than to keep asking the reader to turn back.   

Assume that we are given $\ta$, with boundary $ \Xhigh - \Xlow$. Assume that the complexity $(c_1,c_2)$ has been minimized by the use of exchange moves.  Assume further that we have removed all short clasp arcs, and that all non-short clasp arcs are in normal neighborhoods.  
We will show that we can construct a pair of diagrams $(\cD_+, \cD_-)$, and an isotopy $\cD_+ \to \cD_-$ such that $(\cD_+,\cD_-)$ carries $(\Xhigh,\Xlow)$,  and  the motion $\cD_+ \to \cD_-$ sends $\Xhigh\to\Xlow$.   As in \fullref{proposition:sequence of flypes} we do not know whether $\cD_+$ and $\cD_-$ are block strand diagrams because we do not know whether condition (4) for a block-strand diagram holds.  That matter will not concern us now. (Look ahead to \fullref{subsection:the braid index of blocks in a template} for more on this. See, in particular, \fullref{corollary:no spanning blocks}).  

As in \fullref{subsubsection:using flypes to push Xhigh across cS} there are four
ingredients in the structure of $(\Dhigh,\Dlow)$: the moving blocks, the fixed blocks, the
moving strands,  and the fixed strands.  The precise description of how we construct  each is very close to that in \fullref{subsubsection:using flypes to push Xhigh across cS}, except that: 
\bi
\item $\ta$ replaces $\cS$,
\item $\ta_{\text{tiled}}$ (the portion of $\ta$ not foliated by
$s$--arcs) replaces $\cS_{\text{tiled}}$. As a consequence, we are constructing the pair $(\cD_+,\cD_-)$ (instead of the pair $(\bD,\rD)$). Also the motion that takes $\cD_+\to\cD_-$  now takes $\Xhigh\to\Xlow$ (instead of the motion that takes  $\bD\to\rD$ taking  $\bX\to\rX)$. 
\ei

\textbf{(i)\qua The Moving Blocks}\qua If there are no doubly long clasp
arcs, then there also are no moving blocks, and this part of the
construction ends. So we assume that there are doubly long clasp arcs.
Therefore the two tab neighborhoods of the clasp arcs reach all the
way from $\Xhigh$ to $\Xlow$. We are in the situation of
\fullref{figure:amalgamate5},
and we may construct black and grey microstrands, exactly as we did in
\fullref{subsubsection:using flypes to push Xhigh across cS}, only now
the black and grey microstrands are subarcs of $\Xhigh$ and $\Xlow$.
Notice that the isotopy of  $\Xhigh$ across
$\ta$ pushes $\Xhigh$ across $\Nhigh^i$ to $\Xlow$, so the two $\Xhigh$
microstrands are always mapped to the two $\Xlow$ microstrands.  Exactly as in \fullref{subsubsection:using flypes to push
Xhigh across cS} we have cylinders $\bB^i, \rB^i$ for our microblocks.  The local projection
of the microstrands in $\bB^i$ (resp.\  $\rB^i$) are determined by
the intersections of $\cShigh$ (resp.\  $\cSlow$) with the walls of
$\bB^i$ (resp.\  $\rB^i$).  As in \fullref{subsubsection:using flypes to push Xhigh across cS}, the double points of the local projections are labeled with the index of $g^i$.  

Recall the amalgamation of microblocks, and the example of a template  in \fullref{figure:amalgamate1}.  Before the two flypes it looks as if $X$ and $W$ (and also $Y$ and $Z$) can be amalgamated, but after the flype we see that is impossible.  The requirement that we introduced, that the projections be reflections of each other, is violated if we attempt to amalgamate $X$ and $W$, or $Y$ and $Z$.
We developed a `projection criterion' for when this is possible.  The essential idea was that the block amalgamations that we used maintained their integrity during the passage from the black side to the grey side. 

Let's begin the construction of the black microblocks.  The blocks $\bB$ and $\rB$ have the same structure as in 
\fullref{subsubsection:using flypes to push Xhigh across cS}.  We continue to use the terms `top',
`cylinder', `bottom', with the same meaning as in \fullref{subsubsection:using flypes to push Xhigh across cS}.  Then $\bB$ is an amalgamation of microblocks
$\{\bB^1, \ldots , \bB^r \}$ for microstrands $\{{\bf x}^1_+,{\bf x}^1_- , \ldots
{\bf x}^r_+,{\bf x}^r_-\}$ of $\ta$ if the following hold:
\be

\item As in \fullref{subsubsection:using flypes to push Xhigh across cS}:
Each ${\bB}^i \subset \intr(\bB)$ for $i \leq i \leq r$.
\item As in \fullref{subsubsection:using flypes to push Xhigh across cS}:
Its support, ie $\{ \theta \ | \ \bB \cap H_\theta \not= \emptyset \}$, is a closed interval having
length less than $2\pi$.  Every component of $\bX \cap \bB$ contains a microstrand.

\item A subtle change from (3) in \fullref{subsubsection:using flypes to push Xhigh across cS}: note that  $\bt \cap \ta$ and $\bb\cap\ta$ are each still a
collection of subarcs of $\ahigh$--arcs and, possibly, subarcs in $b$--arcs.  
This necessarily implies that $\bt$ and $\bb$ are contained
in a generic disc fiber of $\fib$. We can assume that neither $\bt$ nor
$\bb$ intersect any
microstrands.   Continuing, we have the immersion $\rho \co \pa 
\longrightarrow \ta$.
Each component of ${\bf R}\subset \rho^{-1}(\bB \cap \ta)\subset\pa$
is a rectangular region that does not contain any singular points.
Specifically:
\bi
\item[(a)] If ${\bf R} \cap {\Xhigh} \not= \emptyset$ then ${\bf R}$ is trivially foliated 
with $\partial {\bf R} = \beta_{+1} \cup \beta_{+2} \cup
\beta_{+3} \cup \beta_{+4}$,
where (i) $\beta_{+1} \subset \Xhigh$, and $\beta_{+1}$ intersects some microstrands in our
specified set. (ii) $\beta_{+2}$ and
$\beta_{+4}$ are subarcs of $\ahigh$--arcs. (iii) $\beta_{+3}$ is transverse
to the
foliation of $\pa$. Moreover the only non-singular leaves it intersects
are $\ahigh$--arcs .
\item[(b)] If ${\bf R} \cap {\bf X} = \emptyset$ then ${\bf R}$ is  trivially foliated, 
 with $\partial {\bf R} = \beta_{+1} \cup \beta_{+2}
\cup \beta_{+3}
\cup \beta_{+4}$ where $\beta_{+1}$ and $\beta_{+3}$ are in leaves; also
$\beta_{+2}$
and $\beta_{+4}$ are transverse to the foliation of $\pa$.
\ei
\item As in \fullref{subsubsection:using flypes to push Xhigh across cS}:
If ${\bf x}^i_+$ (resp.\  ${\bf x}^i_-) \subset  \ \intr (\bB)$, then its
partner ${\bf x}^i_-$ (resp.\  ${\bf x}^i_+)$ is also a subset of \ $\intr(\bB)$. Also,  the dotted segment between each ${\bf q}^i_+$ and ${\bf 
q}^i_-$ is a subset of $\intr(\bB)$.

\item As in \fullref{subsubsection:using flypes to push Xhigh across cS}:  Replacing 
$$\{\bB,\bt, \bc, \bb, \ba, \{\bB^1, \ldots, \bB^r \},\{\bx^1_+,\bx^1_-,
\ldots,\bx^r_+,\bx^r_-\} \}$$
by 
$$\{ \rB, \rt, \rc, \rb, \ra, \{\rB^1, \ldots , \rB^r
\},\{\rx^1_+,\rx^1_-, \ldots, \rx^r_+, \rx^r_-\} \},$$ 
we can also define the amalgamation of microblocks $\{\rB^1, \ldots , \rB^r \}$.
\ee
As before, we also have related concepts for
$\rB$, the amalgamation of microblocks $\{\rB^1, \ldots , \rB^r \}$.

Our definition of the braid projection of $\Xhigh$ onto $\cB \subset \partial \bB$ has changed a little bit
because of the new conditions in (3) above.  For a given amalgamating block $\bB$ we consider
$\{{\bf R}^1, \ldots , {\bf R}^l\} \subset \rho^{-1}(\bB \cap \ta)$.  These are the regions referred to in condition (3)(a).  Let
$\{\beta_{+3}^1, \ldots , \beta_{+3}^l \}$ be the corresponding $\beta_{+3}$ boundary sides of these regions.  Then the graph $\rho(\cup_{1 \leq i \leq l} \beta_{+3}^i) \subset \bc$, along with the clasp arc index labeling of the double points, are in $\pi(\bB)$.  Similarly, $\pi(\rB)$ is just the labeled graph
$\rho(\cup_{1 \leq i \leq l} \beta_{-3}^i) \subset \rc$.  Using this definition of the braid projection of $\Xhigh$ onto $\cB \subset \partial \bB$,
our definition of `related amalgamating blocks' is the same as before, viz: 

Let $\bB$ be an amalgamating block. The braid projection $\pi(\bB)$ of $\Xhigh$ onto $\bc \subset \partial \bB$  is the graph $\bc \cap \cS$. regarded as a subset of $\bc$ together with the clasp arc index labeling of this graph's double
points.  We define $\pi(\rB)$ in a similar fashion.
Amalgamating blocks $\bB$
and $\rB$ are a pair of related amalgamating blocks if the following hold:
\bi
\item[(a)] $\{\bB^1,\ldots,\bB^r\}$ and $\{\rB^1,\ldots,\rB^r\}$ are microblocks associated with
clasp intersections $\{\g^1,\ldots,\g^r\}$.
\item[(b)] $\bB$ is a block amalgamation of microblocks $\{\bB^1,\ldots,\bB^r\}$.
\item[(c)] $\rB$ is a block amalgamation of microblocks $\{\rB^1,\ldots,\rB^r\}$.
\item[(d)] $\pi(\bB)$ and $\pi(\rB)$ are reflections of each other.
\ei
Look back to \fullref{figure:amalgamate6} for an example.

\textbf{(ii)\qua Fixed Blocks}\qua The fixed blocks are associated to braiding between the
strands of $\bX$ (and so also of $\rX$) in the part of $\cS$ which is
foliated by bands of $s$--arcs.  The
braiding between bands of $s$--arcs is not part of our geometry.  $\cS$ only detects the places where
$\bX$ and
$\rX$ differ in a non-trivial way, but in regions where there are bands of
$s$--arcs they essentially coincide.
However, what might happen is that there is braiding between the bands of
$s$--arcs, but that the braiding is
interrupted by, for example, strands of
$\bX$ which separate two potential blocks.  In order to determine the angular regions where the blocks occur,
we therefore look for singularities in the tiled part of $\cS$.  The block subdivision so-obtained will be too
extensive, and we will then need to amalgamate fixed blocks.  For further details, see the proof of \fullref{proposition:sequence of flypes}, which applies without any changes at all.

\textbf{(iii)\qua Moving Strands}\qua The moving strands in $\bD$ are all subarcs of $\bX$ which
are not amalgamated into moving blocks and are not in the bands of $s$--arcs.   We will give details on how they are placed in in the proof of \fullref{proposition:finiteness of cT(m)}.

\textbf{(iv)\qua Fixed Strands}\qua It may happen that we have two
blocks,
$B(C)$ and $B(C')$, for which $(C\cap H_{\theta_{i+1}})\cap (C'\cap H_{\theta_{i+1}}) \not=\emptyset$, but 
$(C\cap H_{\theta_{i+1}})\not= (C'\cap H_{\theta_{i+1}})$, so that the amalgamating condition fails.  In this
situation there are strands that run between $B(C)$ and $B(C')$, however
not all of the $s$--bands
coming out of $B(C)$ go into $B(C')$, and/or not all of the strands going into $B(C')$ are strands
that emerge from $B(C)$.  Any braiding between the $s$ bands that run between $B(C)$ and $B(C')$ can obviously be
pushed into either $B(C)$ or $B(C')$ (eliminating $H_{\theta_{i+1}}$ as a singular fiber) and  they will then
become fixed strands.

We are now in a position to take $\ta$ and construct an associated candidate for a template $(\Dhigh,\Dlow)$. It is clear from the construction that $(\Dhigh$ and $\Dlow)$ will carry
$\Xhigh$ and $\Xlow$.

\begin{remark}
\label{remark:constructing candidates for templates}
Go back to \fullref{figure:template2}. The bottom two sketches show the actual pair of block-strand diagrams $(\cD_+, \cD_-)$ that make up the template in \fullref{figure:template2}, whereas the other 5 sketches give it meaning by showing the actual combination of flypes (not necessarily admissible), exchange moves   and destabilizations that explain the isotopy from $\cD_+$ to $\cD_-$. 

The construction that we just gave falls short of doing what we did in \fullref{figure:template2} in two ways. The first has been noted before this, and we stress it again, but the second is new: 
\be
\item As noted earlier, the blocks may not either be `consolidated' or `optimal', two  matters that  will be discussed in \fullref{section:the proof}.  That's why we have said, repeatedly, that our $\cD_+$ and $\cD_-$ are {\it candidates} for block-strand diagrams, and the pair $(\cD_+,\cD_-)$ is a {\it candidate} for a template.
\item We have not constructed the intermediate stages that explain the
isotopy. That is clear, because all we dealt with is the long clasp arcs,
but to explain the intermediate isotopy we would have had to also look at
clasp arcs which are not, initially, doubly long but become doubly long
after part of $\ta$ has been crossed. That construction remains for future
work, for readers who may be interested in, for example, in working on our
`open problem (3)' in
\fullref{section:open problems} of this paper.  \hfill$\diamondsuit$
\ee
\end{remark}

\subsection{Pushing across regions with a G--flype foliation}  
\label{subsection:G-flype moves}

We have seen that if $\pa$ contains a long
clasp arc, then there is a family of thin annuli $\cS^1\subset\pa$ having $\partial\cS^1 = \bX^0 \cup \bX^1$, where $\bX^0
= \Xhigh$, also $\bX^0$ can be pushed across $\cS^1$ to $\bX^1$ using a sequence of flypes.  The complement of
$\cS^1$ in
$\pa$  will be a new family of annuli which has boundary
$\bX^1
\cup
\Xlow$.  It may happen that we can then apply \fullref{proposition:behavior of graphs} to simplify $\cS^1$
via ${(ab)}^\star$ exchange moves; or it may happen that $\pa
\setminus
\cS^1$ has long clasp arcs, in which case we can iterate the construction.  Putting the two cases together we can
construct a second annulus
$\cS^2$ with $\partial\cS^2 = \bX^1 \cup \bX^2 $ such that the movement across $\cS^2$ is either an ${(ab)}^\star$
exchange move or a sequence of flypes.  This procedure can be iterated until either we produce thin annuli that
have $\Xlow$ as a boundary component, or there are no more $(ab)^\star$ exchange moves, or there are no more long clasp
arcs.  If we do have thin annuli with $\Xlow$ in the boundary then we will have decomposed $\ta$ into a
sequence of thin annuli.  And, we have moved across $\ta$ using flypes and $ab$ or $ab^\star$ exchange moves.  The
difficulty in moving across
$\pa$ this way is that it may happen that we needed to use inadmissible flypes.

With that difficulty in mind, we say that a family of thin annuli $\cS^k\cup\cS^{k+1}\cup\cdots\cup
\cS^{k+l}$ supports a
{\it G--flype foliation} if (i) it is possible to cross the region using only flypes and $ab$ or $(ab)^\star$
exchange moves, and (ii) there exist integers $k,q$ such that $b(\bX^{k+i}) \geq b(\bX^{k})$ for $1\leq i \leq
(q-1)$ and
$b(\bX^{k+q}) \leq b(\bX^{k})$.  In this situation the combined flyping and $ab$ or $(ab)^\star$ motion across the
region 
$\cS^k\cup\cS^{k+1}\cup\cdots\cup \cS^{k+q}$ is a {\it G--flype}.
It is readily seen that if $(\bB , \rB)$ is a pair of related amalgamating blocks for the pair $(\Xhigh,\Xlow)$
then $\bB$ maintains its integrity through successive flypes and is isotopied to $\rB$.

\begin{example}
The template in \fullref{figure:template3} is an
illustration of a G--flype.  
\begin{figure}[htpb] 
\labellist\small
\pinlabel {\parbox{105pt}{(a)\qua from strand 3 of block $Y$ to strand 1
of block $Z$}} [l] at 10 479
\pinlabel {\parbox{105pt}{(b)\qua from strand 2 of block $Z$ to strand 2
of block $Y$}} [l] at 373 355
\pinlabel {\parbox{105pt}{(c)\qua from strand 3 of block $Z$ to strand 3
of block $Y$}} [l] at 75 207
\pinlabel {\parbox{105pt}{(d)\qua from strand 2 of block $Y$ to strand 2
of block $Z$}} [l] at 355 195
\tiny
\pinlabel {$1$} [tr] at 100 582
\pinlabel {$3$} [tl] at 140 541
\pinlabel {$4$} [tr] at 180 582
\pinlabel {$5$} [bl] at 140 583
\pinlabel {$13$} [bl] at 180 542
\pinlabel {$14$} [bl] at 100 542
\pinlabel {$\gamma^3$} [l] at 121 514
\pinlabel {$\gamma^4$} [l] at 197 514
\pinlabel {$+$} at 107 595
\pinlabel {$+$} at 187 595
\pinlabel {$+$} [tl] at 140 583
\pinlabel {$+$} at 140 542

\pinlabel {$1$} [br] at 399 408
\pinlabel {$2$} [br] at 399 518
\pinlabel {$4$} [tr] at 399 558
\pinlabel {$7$} [br] at 467 506
\pinlabel {$7$} [r] at 483 483
\pinlabel {$8$} [tr] at 467 460
\pinlabel {$10$} [r] at 357 483
\pinlabel {$12$} [l] at 440 484
\pinlabel {$14$} [tr] at 400 449
\pinlabel {$15$} [tl] at 329 460
\pinlabel {$17$} [r] at 313 484
\pinlabel {$18$} [l] at 329 506
\pinlabel {$\gamma^3$} [l] at 411 570
\pinlabel {$\gamma^1$} [l] at 420 396
\pinlabel {$+$} [tr] at 399 408
\pinlabel {$+$} [br] at 399 558
\pinlabel {$+$} [l] at 313 484
\pinlabel {$+$} [l] at 484 483
\pinlabel {$+$} at 399 517
\pinlabel {$+$} at 330 460
\pinlabel {$+$} at 468 460

\pinlabel {$3$} [bl] at 194 261
\pinlabel {$5$} [br] at 194 371
\pinlabel {$6$} [tr] at 193 411
\pinlabel {$6$} [br] at 263 359
\pinlabel {$8$} [r] at 279 336
\pinlabel {$9$} [tr] at 263 313
\pinlabel {$9$} [r] at 153 336
\pinlabel {$11$} [l] at 237 336
\pinlabel {$11$} [tr] at 194 302
\pinlabel {$16$} [tl] at 125 313
\pinlabel {$17$} [l] at 125 359
\pinlabel {$18$} [r] at 109 337
\pinlabel {$\gamma^4$} [l] at 207 422
\pinlabel {$\gamma^2$} [l] at 214 249
\pinlabel {$+$} [br] at 193 411
\pinlabel {$+$} [tr] at 194 261
\pinlabel {$+$} [l] at 279 336
\pinlabel {$+$} [l] at 109 337
\pinlabel {$+$} at 195 370
\pinlabel {$+$} at 264 312
\pinlabel {$+$} at 126 312

\pinlabel {$2$} [bl] at 457 243
\pinlabel {$10$} [tr] at 497 243
\pinlabel {$12$} [bl] at 457 283
\pinlabel {$13$} [tr] at 417 243
\pinlabel {$15$} [tr] at 417 283
\pinlabel {$16$} [tr] at 497 283
\pinlabel {$\gamma^1$} [l] at 429 303
\pinlabel {$\gamma^2$} [l] at 511 303
\pinlabel {$+$} at 498 242
\pinlabel {$+$} at 418 242
\pinlabel {$+$} at 429 229
\pinlabel {$+$} at 511 229
\pinlabel {$+$} [tr] at 457 243

\endlabellist
\centerline{\includegraphics[scale=.50]{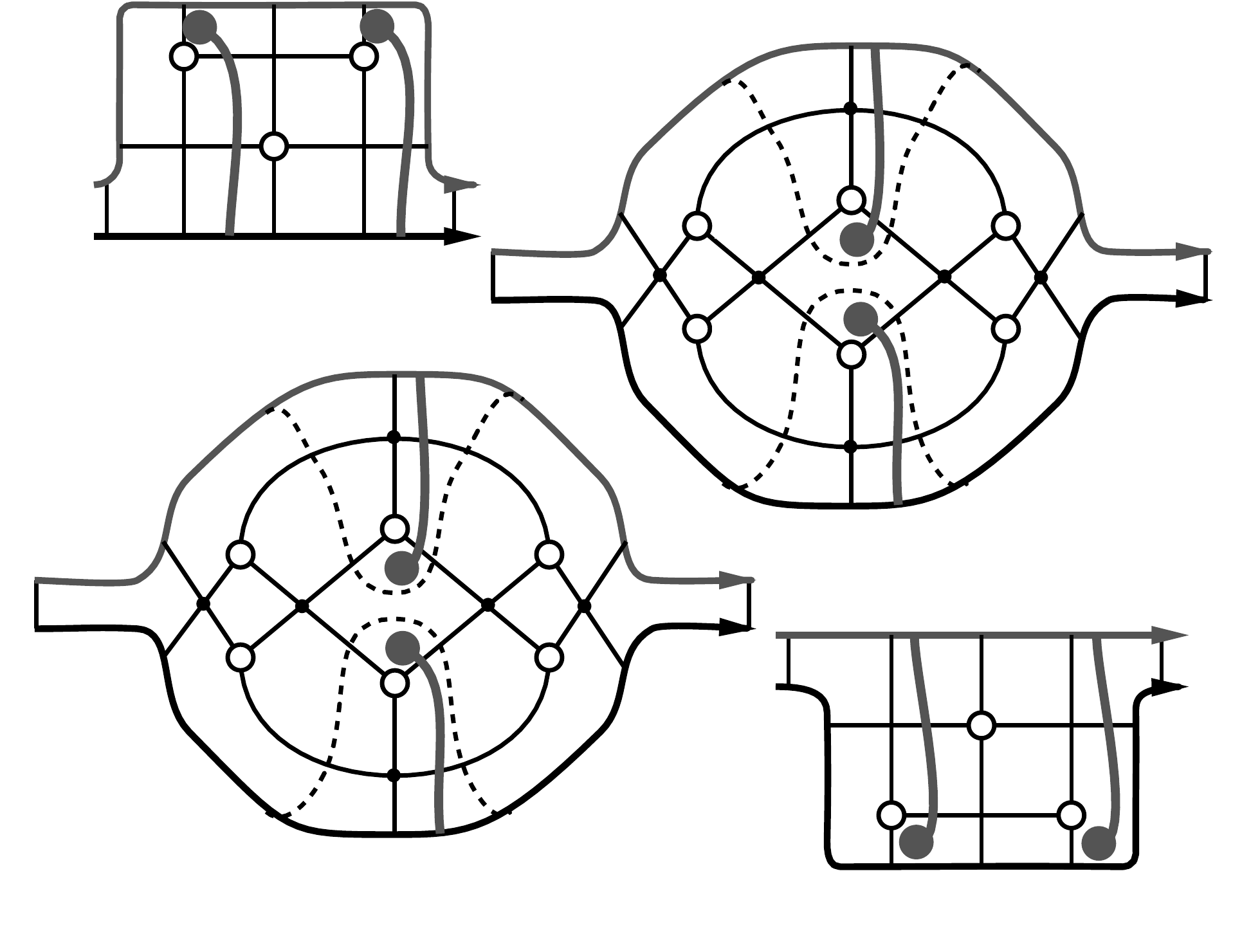}}
\caption{The foliation of $\pa$, in the situation of the
6--braid template of \fullref{figure:template2}}  
\label{figure:template3} 
\end{figure}

It shows the
foliated annulus $\pa$ which was used to construct the template which we saw in 
\fullref{section:introduction} in the boxed sketches at the bottom of \fullref{figure:template2}.  The labels `block 
$Y$' and `block $Z$' refer to the blocks in \fullref{figure:template2}.  In the interest 
of keeping the foliation as simple and understandable as possible, we 
decided to position the clasp arcs only in tab neighborhoods, and not 
in normal neighborhoods.  Thus, the clasp arcs are {\em not} contained in 
the finite union of leaves.  (The addition of 12 vertices and 
12 singularities would be needed in order to position the clasps in normal 
neighborhoods.)  In order to remove distracting details from the picture, we have eliminated inessential
$b$--arcs, collapsing the normal neighborhoods of the clasp arcs $\gamma^1_+$ and $\gamma^2_+$.  This is always
possible, and we did it to save space.  The arcs
$\g^1$ and
$\g^2$ in sketch (d) are the clasp arcs used to construct the first thinly foliated annuli that are
split off by the dotted curves in sketches (b) and (c).  After we have flyped across this initial thin
annulus (the first negative flype of \fullref{figure:template2}) we will have four resulting  vertices
where we can perform $(ab)$ exchange move.  These are the two double-strand exchange moves in the isotopy sequence of
\fullref{figure:template2}.  Finally, the remaining portion of $\ta$ will just be the region in sketch (a) and tab
neighborhoods around clasp arcs $\g^3$ and $\g^4$ in sketches (b) and (c).    All of
these individual motions were across thin annuli, giving us a decomposition of $\ta$ into thin annuli.
Since the braid index starts at $b(\Xhigh) = 6$, goes to $7$ after the first flype, remains at $7$ for all of the
exchange moves and only returns to $6=b(\Xlow)$ after the last flype, all of the thin annuli in this
decomposition go to make up the motion across $\pa$.  The combination of
all of these isotopies is the 6--braid
G--flype which we illustrated in \fullref{figure:template2} and \fullref{figure:template3}.

Notice that there are no moving blocks in \fullref{figure:template2} because there are no doubly
long clasp arcs in \fullref{figure:template3}.  The following features of
the template in \fullref{figure:template2} are outside the support of the
isotopy: the four blocks $X,Y,Z,W$ and the following strands
that join them: all strands that enter and leave $W$ and $X$, also strand
1 entering $Y$, strand 1 leaving $Y$, strands 3
and 4 entering $Z$ and strands 1 and 4 leaving $Z$.  We need to account for the changes in strands 2 and 3
entering $Y$, strands 2 and 3 leaving $Y$, strands 1 and 2 entering $Z$
and strands 2 and 3 leaving $Z$.  These are
all described completely by the data in the foliated annulus.
\end{example}

\subsection{Pushing across regions with a G--exchange foliation}  
\label{subsection:pushing across regions with a G-exchange foliation}  
In this subsection all clasp arcs have intermediate length, that is the pierce points of both $\ghigh$ and $\glow$
are on $b$--arcs.

 We begin
with an example which illustrates how G--exchange moves arise.  The example is the foliated subsurface of $\pa$ that
supports the G--exchange move of \fullref{figure:exchange6}.  It will not be difficult to understand this figure, now that the main tools in this paper,
ie the foliated immersed annulus $\ta$ and its foliated preimage $\pa$, are in place.

\begin{figure}[htpb!]
\labellist\small
\pinlabel {$a_-$} [l] at 304 499
\pinlabel {$a_+$} [l] at 305 474
\pinlabel {$s_1'$} [l] at 272 487
\pinlabel {$b_-$} [l] at 518 498
\pinlabel {$b_+$} [l] at 519 474
\pinlabel {$s_2'$} [l] at 487 484
\pinlabel {$c_-$} [l] at 298 354
\pinlabel {$c_+$} [l] at 297 330
\pinlabel {$s_3'$} [l] at 263 342
\pinlabel {$d_-$} [l] at 512 357
\pinlabel {$d_+$} [l] at 513 333
\pinlabel {$s_4'$} [l] at 482 345

\pinlabel {$s_1$} [r] at 141 399
\pinlabel {$s_2$} [r] at 357 395
\pinlabel {$s_3$} [r] at 133 253
\pinlabel {$s_4$} [r] at 355 255

\pinlabel {$1$} [r] at 151 441
\pinlabel {$2$} [b] at 202 499
\pinlabel {$3$} [t] at 215 385
\pinlabel {$4$} [l] at 266 439

\pinlabel {$5$} [r] at 365 441
\pinlabel {$6$} [b] at 416 498
\pinlabel {$7$} [t] at 428 384
\pinlabel {$8$} [l] at 479 444

\pinlabel {$3$} [r] at 143 294
\pinlabel {$4$} [b] at 194 354
\pinlabel {$5$} [t] at 208 240
\pinlabel {$6$} [l] at 257 295

\pinlabel {$7$} [r] at 361 301
\pinlabel {$8$} [b] at 410 358
\pinlabel {$1$} [t] at 424 244
\pinlabel {$2$} [l] at 476 302

\endlabellist
\centerline{\includegraphics[scale=.7]{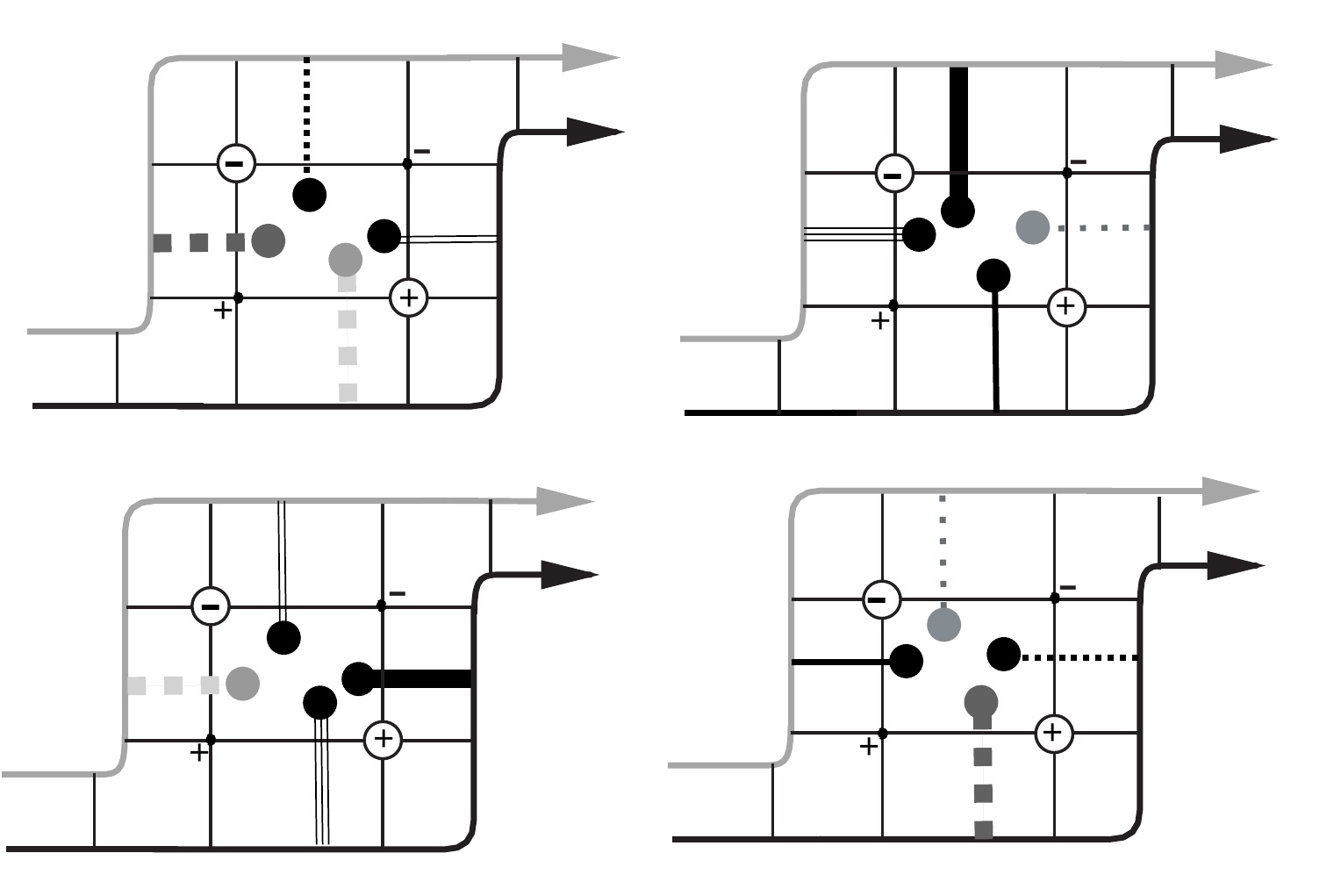}}
\caption{Foliated subsurface of $\pa$ which are the support of the G--exchange move of \fullref{figure:exchange6}}
\label{figure:exchange7}
\end{figure}

The foliated subsurface of $\pa$ that supports the G--exchange move of \fullref{figure:exchange6} is illustrated in \fullref{figure:exchange7}, which shows four discs on $\pa$.  Label the discs $\cR_1$ (top  left), $\cR_2$ (top right), $\cR_3$ (bottom left) and $\cR_4$
(bottom right).   No $\cR_i$ is good,
because each contains the puncture endpoints of  clasp arcs.  On the other hand, each is an
embedded subset of $\pa$ (because no $\cR_i$  contains both preimages of a clasp arc pair).  The clasp arcs are all doubly-intermediate
(because all of the puncture endpoints are on $b$--arcs).

Ignoring the clasp arcs momentarily, we see that each disc is topologically equivalent to the shaded disc in
the upper left sketch in \fullref{figure:exchange5}, ie it  contains a positive vertex of valence 2 and type $ab$, and
singularities of opposite sign.  Each of our 4 regions contains in its boundary a subarc of
$\Xlow$ and and a subarc of $\Xhigh$, colored dark and light respectively.  These arcs are labeled 
$a_\pm, b_\pm, c_\pm$ and $d_\pm$.
If the clasp arcs were not there we could use four $ab$--exchange moves to
push $a_+$ (resp.\  $b_+, c_+, d_+$) across its disc to $a_-$
(resp.\  $ b_-, c_-, d_-$). 

Since $\cR_1$ intersects $\cR_4$ along clasp arcs 1 and 2, and
$\cR_3$ along clasp arcs 3 and 4, we cannot complete the move on strand $a_+$ until we begin the
moves on strands $c_+$ and $d_+$.  But then, the motions of strands $c_+$ and $d_+$
across discs $\cR_3$ and $\cR_4$ cannot be completed because those discs have clasp
intersections with $\cR_2$, and in fact no one of the motions can be
completed until all of the others are completed too.  That is, we have a
G--exchange move.  This particular very simple G--exchange move was
discussed in \fullref{section:introduction} to this paper, and was
depicted in the 6--braid example in \fullref{figure:exchange6}.  
 
\begin{remark}
While we have 
gone to some pains to insure that all clasp arcs are in normal neighborhoods, in the example just given of a
G--exchange move the normal neighborhoods are ignored for reasons of space, as they would enlarge the
pictures to the point where they would obscure the features that are of interest.  Since normal neighborhoods
were created by adding many many inessential $b$--arcs; going the other way, they can also be deleted by an
isotopy of the embedded part of $\ta$. \hfill$\diamondsuit$
\end{remark}

To generalize this example we first need to understand the foliated subregions of $\pa$  which lead to sequences of
exchange moves that carry a subarc of $\Xhigh$ over a `rooted block and strand tree'.  It will be easiest to study
them first without the clasp arcs, and then add the clasp arcs later.  Then we will need to understand the associated
block and strand diagrams.  Since the definitions are somewhat detailed, it may be helpful to see them worked out in a
few special cases first.

\begin{figure}[htpb!]
\labellist\small
\pinlabel {$s$} [r] at 95 692
\pinlabel {$s'$} [l] at 285 692
\pinlabel {$w$} [b] at 152 718
\pinlabel {$w'$} [b] at 222 718
\pinlabel {$v$} [t] at 156 666
\pinlabel {$v'$} [t] at 224 670
\pinlabel {$r'$} [bl] at 186 714
\pinlabel {$r$} [bl] at 189 672

\pinlabel {$s$} [r] at 219 569
\pinlabel {$r'$} [bl] at 275 591
\pinlabel {$r$} [tr] at 329 549
\pinlabel {$s'$} [l] at 385 568
\pinlabel {$w$} [b] at 266 594
\pinlabel {$w'$} [b] at 335 594
\pinlabel {$v$} [t] at 267 547
\pinlabel {$v'$} [t] at 338 547

\pinlabel {(a)} at 173 618
\pinlabel {(b)} at 439 618
\pinlabel {(c)} at 144 565

\tiny
\pinlabel {$B_i$} at 376 651
\pinlabel {$B_{i\!+\!1}$} at 435 652
\pinlabel {$B_{i\!+\!2}$} at 494 652
\pinlabel {$t(i,i{+}1)$} [b] at 385 720
\pinlabel {$t(i{+}1,i{+}2)$} [b] at 499 720
\pinlabel {$s(i,i{+}1)$} [t] at 406 644
\pinlabel {$s(i{+}1,i{+}2)$} [t] at 466 644

\pinlabel {$B_i$} at 251 569
\pinlabel {$B_{i\!+\!1}$} at 304 569
\pinlabel {$B_{i\!+\!2}$} at 357 569
\endlabellist
\centerline{\includegraphics[scale=.7]{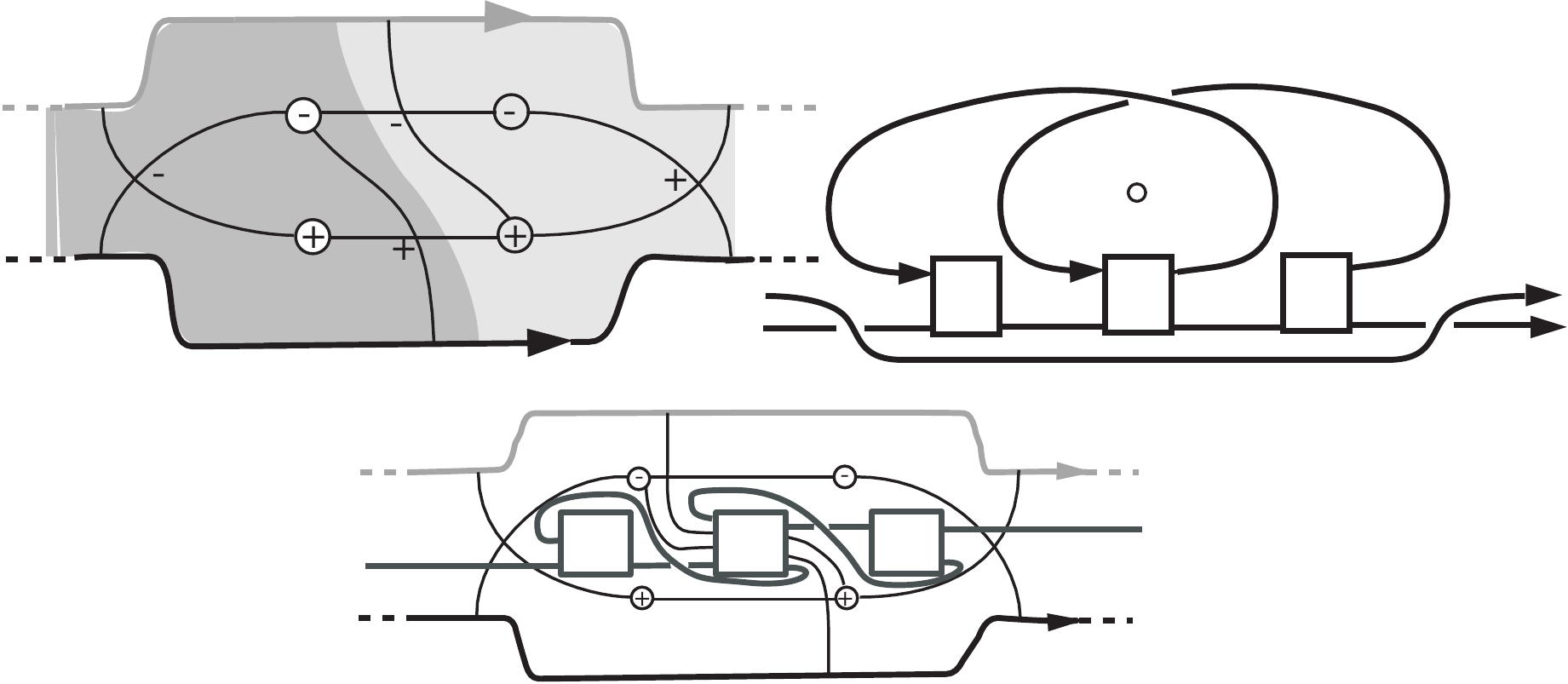}}
\caption{(a) The root foliation on $\pa$.  (b) The root diagram. (c) The two viewpoints are
put together.}
\label{figure:handle2}
\end{figure}
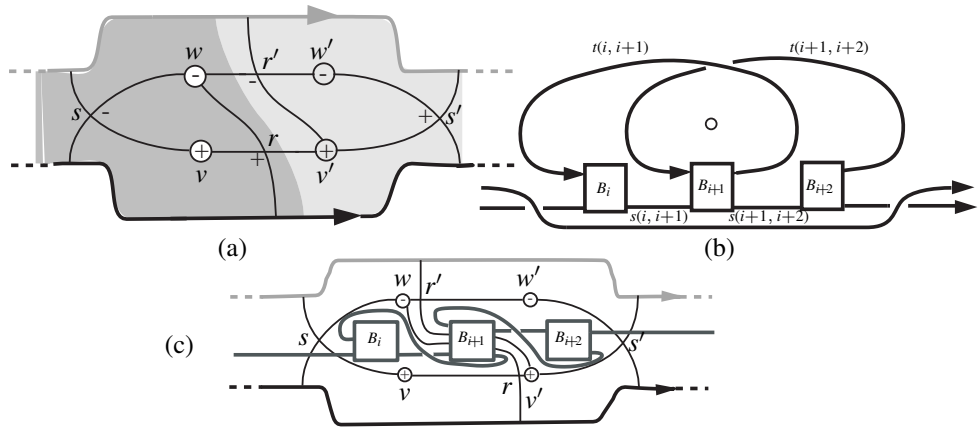

In \fullref{figure:handle2}(a)  a subarc of $\Xhigh$  can be pushed across the shaded regions by two
exchange moves: the first across the darkly shaded disc (containing the vertices $v$ and $w$ and
the singularities $s$ and $r$), the second across the lighter shaded disc, which contains the vertices $v',w'$ and the
singularities $s',r'$.
The support of the first exchange move, ie the dark shaded region,  is a `pouch' which is pierced twice by
the axis, at $v$ and $w$.  If the
foliation of $\pa$  contains only essential $b$--arcs, then the closed braid must wrap
around the braid axis (perhaps with many strands traveling along together, say
$t$ strands in all) in between the two pierce-points
$v$ and $w$.  The first exchange move is a push of a subarc of $\Xhigh$ across the  pouch, crossing the axis twice as
it does so at $v$ and $w$.  Then there is a second exchange move across a second pouch, crossing the axis twice at
$v'$ and $w'$.  There will also be
$t'$ braid strands wrapping about the axis in between $v'$ and
$w'$, and perhaps a braid in between the weighted strands $t$ and $t'$.

Keeping all this  in mind we turn to \fullref{figure:handle2}(b).
It shows a picture of a root diagram which is
part of a closed braid diagram.  There are 3 blocks, labeled $B_i, B_{i+1}, B_{i+2}$ and two
weighted strands, labeled $t(i,i+1)$ and
$t(i+1,i+2)$, traveling around the braid axis.  These blocks and strands
are there because if not the $b$--arcs which
foliate the pouch would not be essential.  The pouch is not shown in this picture.

In \fullref{figure:handle2}(c) we put together the information in \fullref{figure:handle2}(a) and (b).  The region in (a) is to be thought of as a very
flexible disc with two pouches.  We are looking  through the pouches to
the block and strand tree which is visible inside them.  The  braid axis $\axis$ pierces the
pouches in axis pieces $vw$ and $v'w'$.  The darker pouch (we called it $C$) is the support of the
first $ab$--exchange move in sketch (a).  It covers the  darker
$\Delta$--disc.  The lighter shaded pouch
$C'$ covers the lighter $\Delta$--disc.  The motion of our subarc of $\Xhigh$ is the sequence of two exchange moves
over the two pouches.  The subarc is like a `handle' which moves  over the braid blocks $B_i, B_{i+1}, B_{i+2}$ and the
weighted strands. 

\begin{figure}[htpb!]
\labellist\tiny
\pinlabel {$t^1(i,i{+}1)$} [br] at 178 555
\pinlabel {$t^1(i{-}1,i)$} [b] at 210 550
\pinlabel {$t^1(i{+}1,i{+}2)$} [b] at 270 555
\pinlabel {$t^2(i,i{+}1)$} [bl] at 312 557

\pinlabel {$s(i{-}2,i{-}1)$} [r] at 134 497
\pinlabel {$s(i{-}1,i)$} [r] at 134 440
\pinlabel {$s(i,i{+}1)$} [r] at 134 383
\pinlabel {$s(i{+}1,i{+}2)$} [r] at 134 324
\pinlabel {$s(i{+}2,i{+}3)$} [r] at 134 270

\pinlabel {$B_{i-1}$} at 144 468
\pinlabel {$B_{i}$} at 144 411
\pinlabel {$B_{i+1}$} at 144 352
\pinlabel {$B_{i+2}$} at 144 295

\pinlabel {$B'$} at 398 396
\pinlabel {\rotatebox{25}{$B_1''$}} at 472 425
\pinlabel {\rotatebox{330}{$B_2''$}} at 475 368
\endlabellist
\centerline{\includegraphics[scale=.7]{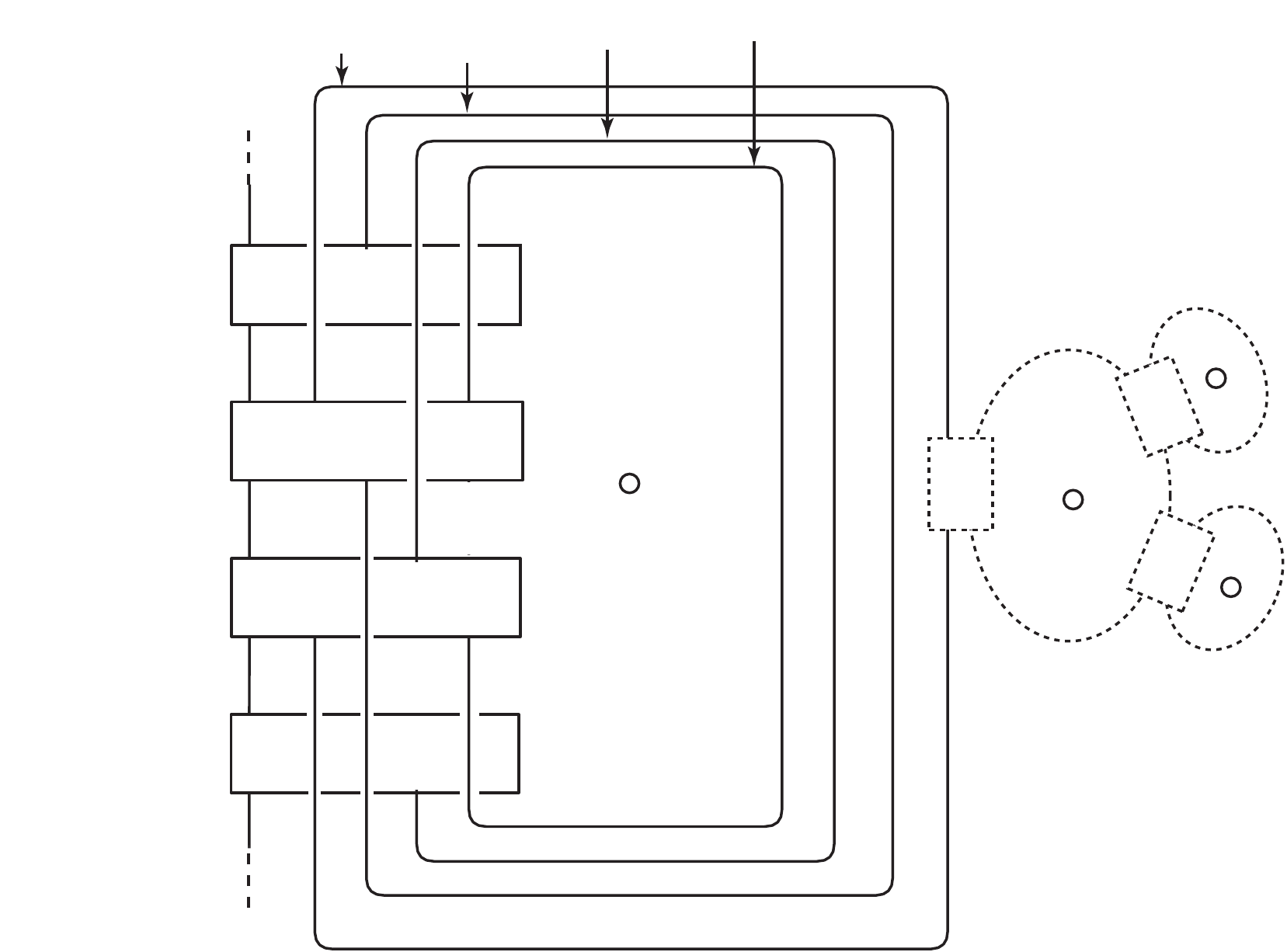}}
\caption{A block and strand tree. The dotted parts will be discussed later.}
\label{figure:handle1}
\end{figure}
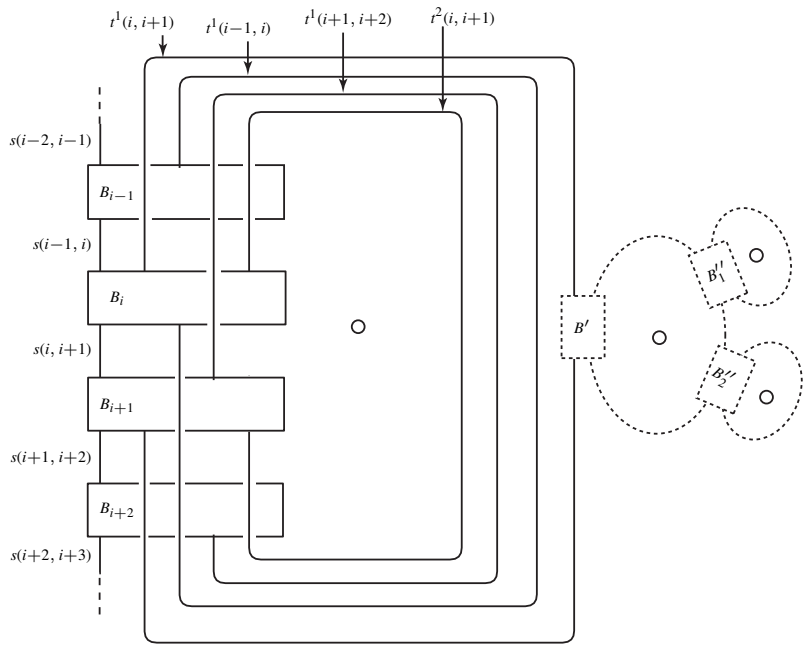

Now for a formal definition of roots, branches and a block and strand tree.  \fullref{figure:handle1} should be helpful.
A collection of braid blocks
$\{ B_1, \ldots , B_l \}$ and weighted strands
$\{ s(1,2),\ldots,s(i,i+1), \ldots \}, \ 1\leq i\leq l$, and additional weighted strands
$\{t^1(1,2), \ldots ,t^{k_1}(1,2),t^1(2,3), \ldots , t^{k_2}(2,3), \ldots, 
t^1(i,i+1), \ldots , t^{k_i}(i,i+1), \ldots \}$ is a \emph{root}  
if $s(i,i+1), \ 1\leq i\leq l$ has
endpoints at the bottom of $B_i$ and at the top of
$B_{i+1}$; $t^j(i,i+1), \ 1\leq j \leq k_i$ has endpoints at the top of $B_i$ and at the
bottom of $B_{i+1}$, and if there exist embedded discs $\Delta (i,i+1,j) \subset S^3$ satisfying 
the following further conditions:

\bi
\item The $\Delta$--discs have disjoint interiors.  Also, for each $t^j(i,i+1), \
1\leq j\leq k_i$  there is one
associated disc.

\item The braid axis $\axis$ intersects each $\Delta$--disc transversally in a single point.

\item $\partial \Delta(i,i+1,j) \subset t^j(i,i+1)\cup B_i \cup s(i,i+1) \cup B_{i+1}.$
Also $\Delta(i,i+1,j) \cap s(i,i+1) = s(i,i+1)$ and 
$\Delta(i,i+1,j) \cap t^j(i,i+1) = t^j(i,i+1)$.

\item If $\Delta(i,i+1,j)$ and $\Delta(i,i+1,j^\prime)$ are distinct 
$\Delta$--discs which
are intersected in succession by some meridian loop of $s(i,i+1)$ then 
there exist $\Delta(i-1,i,m)$
and $\Delta(i+1,i+2,q)$ such that the unoriented $\axis$ 
intersects first
$\Delta(i,i+1,j)$,  then $\Delta(i-1,i,m)$ and $\Delta(i+1,i+2,q)$ 
(in either order); then
$\Delta(i,i+1,j^\prime)$; then all other $\Delta$--discs.
\item For each $\Delta(i,i+1,j)$ there exist a $\Delta(x,y,z)$ with either 
$(x,y,z) = (i-1,i,m)$
or $(x,y,z) = (i+1,i+2,q)$.  Moreover, the unoriented $\axis$ intersects 
in succession
$\Delta(i,i+1,j), \Delta(x,y,z)$; then all remaining $\Delta$--discs.
\ei

Given a braid structure
$(\fib,\axis)$, a {\it radial sphere} 
is a 2--sphere that is transversally intersected by
$\axis$ twice and is transverse to all of the disc fibers of
$\fib$.  An {\it axis piece} $\alpha$ in a radial sphere $S$ is a closed arc whose
interior is transverse to the disc fibers of $\fib$ and which has empty
intersection with at least one disc fiber.  Axis pieces $\alpha,
\beta \subset S$, are loop equivalent if $\partial\alpha =
\partial\beta = \alpha \cap \beta$ and if $\alpha \cup \beta$ bounds a 2--disc
$\delta \subset S$ such that $\delta \cap \axis = \emptyset$.  

A {\it branch} 
is a block $B$ with associated weighted strands $t$, along with a 2--disc $\Delta$ such that:
\bi
\item $\Delta$ is transversally intersected by $\axis$ at one point.

\item $\partial \Delta = t \cup a$ where $a \subset \partial B$.
Specifically, $a$ is an arc made up of three segments, $a = a_T \cup a_S
\cup a_B$ where: $a_T$ is on the top of $B$; $a_S$ is on the side of
$B$; and $a_B$ is on the bottom of $B$.

\item $\intr(\Delta) \cap {B} = \emptyset$
\ei

The braid block $B_1^{\prime\prime}$ and the weighted strand which emerges from it and loops around the axis
is an example.

Iterating the construction of attaching a branch to a root, we obtain a 
{\it rooted block-strand tree}.

\begin{figure}[htpb!]
\labellist\small
\pinlabel {(a)} at 140 332
\pinlabel {(b)} at 388 332

\pinlabel {$w$} [br] at 108 475
\pinlabel {$v$} [tl] at 192 409
\pinlabel {$x$} [br] at 137 451
\pinlabel {$y$} [tl] at 165 430

\pinlabel {$s$} [r] at 303 433
\pinlabel {$r'$} [br] at 487 493
\pinlabel {$r$} [tl] at 440 384
\pinlabel {$v$} [t] at 341 379
\pinlabel {$v'$} [t] at 487 381
\pinlabel {$w$} [b] at 393 500
\pinlabel {$x$} [b] at 374 463
\pinlabel {$x'$} [br] at 418 463
\pinlabel {$y$} [t] at 357 416
\pinlabel {$y'$} [tl] at 440 415
\endlabellist
\centerline{\includegraphics[scale =.7]{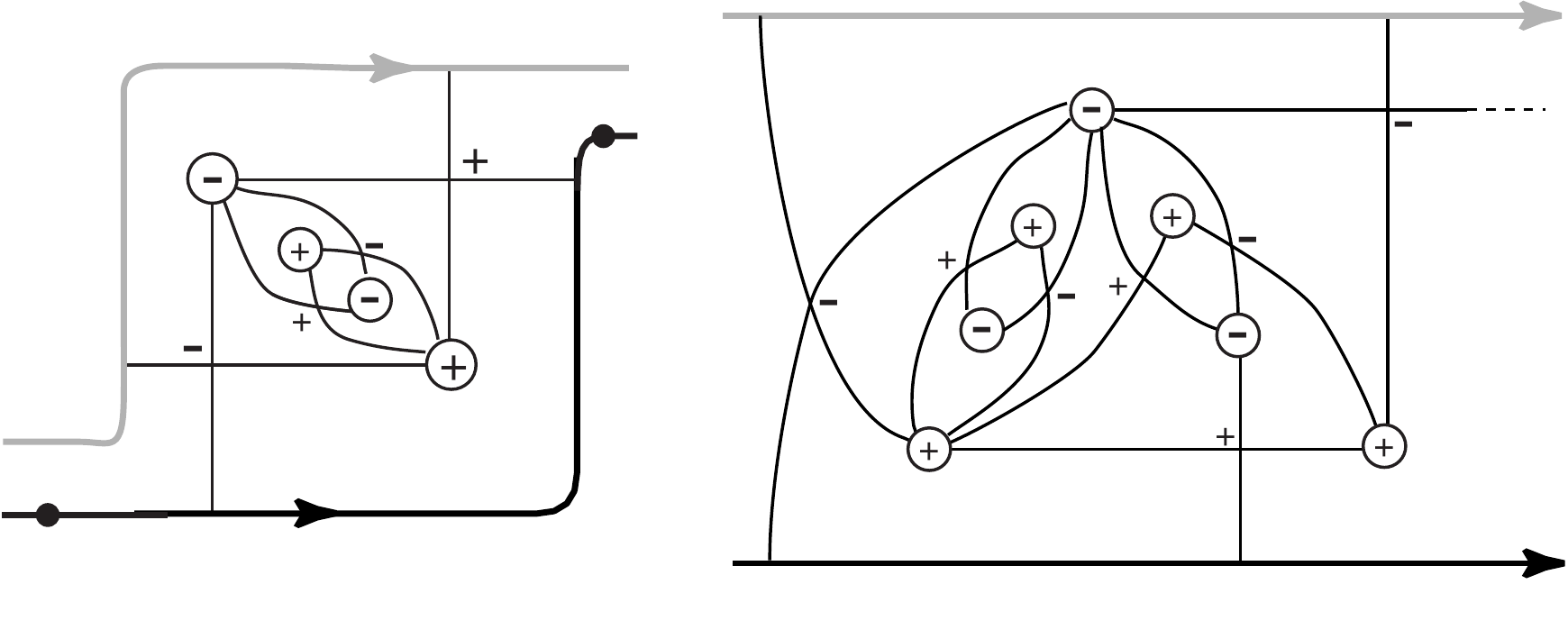}}
\caption{Expanding root foliations to tree foliations}
\label{figure:handle3}
\end{figure}

Next we need to understand the foliated subregions of $\pa$  which lead to sequences of
exchange moves that carry a handle subarc  of $\Xhigh$ over a rooted block and strand tree.  Again it will be helpful
to see examples before we give the general definition.  We first illustrate how the very simple pouch $\cR_1$ in \fullref{figure:exchange7} might itself develop a pouch.  As before, we ignore the clasp arcs.  \fullref{figure:handle3}(a)
shows the double pouch foliation of the region $\cR$ of \fullref{figure:handle3}(a). 
In the expanded foliation the handle arc (the arc $\alpha$ of Figures \ref{figure:exchange5}) cannot be
pushed across $\cR$ by $ab$ exchanges.  Fortunately $bb$--exchange moves come to the rescue.  (See Figures
\ref{figure:exchange4} and \ref{figure:foliations7}.)  The new vertex $x$ has valence 2 and type
$(b,b)$.  \fullref{lemma:ab and bb exchange}, part (2), applies.  We can do a $bb$ exchange, and then
remove the resulting inessential $b$--arcs.  Now an $ab$--exchange is possible.  As for the corresponding braid
picture, the presence of the  new vertices $x,y$ means that the dotted root block $B'$ of sketch
\fullref{figure:handle1} has grown new branches $B_1^{\prime\prime}$ and $B_2^{\prime\prime}$. The move over the
tree in \fullref{figure:handle3}(a) will be a $bb$--exchange followed by
an $ab$ exchange. In the presence of clasp arcs, the G--exchange move of Figures \ref{figure:exchange6} and \ref{figure:exchange7} will of course become much more complicated when new branches are added. In this regard we note: an important feature is that every time new branches are added the braid index increases.  Thus very complicated block and strand trees will only be encountered at very high braid index. 

A slightly more complicated example, in \fullref{figure:handle3}(b), shows the changes in foliation when we
grow new branches in two different ways in the foliation of the region $\cR$ of \fullref{figure:handle2}(a).  The changes are supported inside the region $wsvv'r'w$.  The branch
associated to the new vertices $x$ and $y$ (resp.\  $x'$ and $y'$) is attached to the strand
joining blocks
$B_{i}$ and $B_{i+1}$ (resp.\  blocks $B_{i+1}$ and $B_{i+2}$).  The root diagram of \fullref{figure:handle2}(b) has changed to a block-strand tree.  In the foliation, the growth has
all been `inward'.  This time two $bb$ exchanges and two $ab$ exchanges are needed to realize
the G--exchange move over the block-strand tree.  The reason G--exchange moves can be hard to visualize is
because the part of the surface that undergoes the change in foliation is always far away from
the block-strand tree in the closed braid.  Putting this is another way, the foliated surface
points out the way to organize very big sequences of exchange moves, some of which can be quite
difficult to see in the closed braid diagram.

Finally, we come to the general definition.  Let ${\bf S} \subset \pa$ be
a complete collection of $s$--arcs.
For present purposes a region
$\cR \subset \pa \setminus {\bf S}$
is either a rectangular shaped subdisc or a subannular region of $\pa$.
Thus, as before, if $\cR$ is a subdisc then 
$\partial \cR =  s \cup  Y_+ \cup s^\prime \cup Y_- $,
where: $s$ and $s^\prime$ are
subarcs of leaves in the foliation of $\pa$;  and $Y_+$ and $Y_-$ are
oriented arcs transverse to the foliation in the positive direction.
If $\cR$ is a subannulus then $\partial \cR = Y_+\cup Y_- $
where $Y_\pm$ are
oriented curves transverse to the foliation in the positive direction.
Inside these regions we will have the induced foliation.  While our regions 
may be intersected by clasp arcs, we are not concerned with them at this time.

Let $\cR \subset \pa
\setminus{\bf S}$ be as above.  Let $C$ be a component of $\cR\cap G_{\epsilon,\delta}$, where the graph $G_{\epsilon,\delta}$ was defined pictorially in \fullref{figure:graphs1}.  We say that
$\cR$ has a {\it root foliation} if: 
\be
\item $C$ is homeomorphic to either $S^1$ or $[0,1]$.
\item If $C$ is homeomorphic to $S^1$ then 
$\cR$ is an annulus and $C$ is homotopically equivalent to a core circle of $\cR$.
\item If $C$ is homeomorphic to $[0,1]$ then $C$ has an endpoint on $Y_\pm$ and
an endpoint near $Y_\mp$.
\ee

\begin{proposition}
\label {proposition:root foliations}
Let $\cR \subset \pa \setminus{\bf S}$ be a component which is embedded.
Assume that all $b$--arcs are essential.  Assume that  $\cR$ has a root foliation.  Then the isotopy
which corresponds to pushing a component of ${\bf Y}_+$ across $\cR$ is a
sequence of $ab$--exchange moves over a root.
\end{proposition}

\begin{proof}  The isotopy of the braid across a region which has a root foliation can be realized by a
sequence of $ab$ exchanges.  To see this, notice that at least one  pouch $\cP_i$ is associated to each such
region, and since the regions are crossed in a definite order the pouches can be joined in the
same order.  The assumption that each $b$--arc is essential implies that $\cP_i$ cannot be removed
by isotopy, and is associated to a new braid block or blocks.  The union
of all of the $\cP_i's$ gives a disc region $\cR$ with a subarc of ${\bf Y}_+$ in its boundary.  The
union of all of the blocks is a {\it root}.  The isotopy of the subarc of ${\bf Y}_+$ across $\cP$ is a
sequence of exchange moves across this root.  This completes the proof.
\end{proof}

\begin{remark}
It might happen, in the situation of \fullref{proposition:root foliations}, that ${\bf S}$ is empty on some component of $\cR$.  If that occurs, the component in question will be a standard annulus. Pushing across a standard annulus will be treated in  \fullref{subsection:pushing across a standard annulus}.
\hfill$\diamondsuit$
\end{remark}

 A region $\cR$ has a {\it tree foliation} if,
after a sequence of $bb$--exchanges the foliation is reduced to a root foliation.

The next proposition shows how a tree
foliation imposes a block decomposition on our two braids:
\begin{proposition}
\label{proposition:tree foliations}
Let $\cR \subset \pa \setminus{\bf S}$ be a component whose image is embedded.  Assume that ${\bf
S}\not=\emptyset$.  Assume that all $b$--arcs are essential, and that $\cR$ has a tree foliation.  Then the isotopy
which corresponds to pushing ${\bf Y}_+$ across $\cR$ is a sequence of exchange moves over a block-strand tree.
\end{proposition}

\begin{proof}  The  isotopy of the braid through the
sequence of $bb$--exchanges corresponds to collapsing a tree to its root.  The collapsing is realized by
$bb$--exchange moves followed by braid isotopy to remove inessential $b$--arcs.  After that, an isotopy of
${\bf Y}_+$ across the new $\cR$  corresponds to a sequence of $aa$ and
$ab$--exchange moves over this root. \end{proof} 

\begin{remark} 
\label{remark:tabs and flypes}
Two remarks are in order.  They relate to Propositions \ref
{proposition:root foliations} and 
\ref{proposition:tree foliations}.
\begin{enumerate}
\item In both propositions the basic assumption is that the region $\cR$ is
embedded.  Of course this will be the case if no clasp arcs intersect $\cR$.  But it will also be
the case if (a) clasp arcs intersect $\cR$ but no clasp arc pair intersects $\cR$, or (b) there are no long clasp arcs which intersect  $\cR$.
\item After a G--exchange move across $\cR$, valence 1 vertices may
be revealed.  See \fullref{figure:foliations4}.   Later, when we get to the final steps in the proof, we will remove them
by destabilization.\hfill$\diamondsuit$
\end{enumerate} 
\end{remark}

We are finally ready to introduce clasp arcs into the picture.  In the most general case there will be
several related regions
$\cR_1,\dots,\cR_k$ which are intersected by paired clasp arcs.  Thus the exchange moves across one
$\cR_i$ will have to be interrupted midway to do part of an exchange move along an associated
$\cR_j$.  Let $\{\cR_1,\dots,\cR_k \} \subset \pa$ be a collections of  regions
such that each $\intr(\cR_i), \ 1\leq i \leq k,$ is embedded in $S^3 \setminus (\Xhigh \cup 
\Xlow)$ and each $\cR_i$ has a tree foliation.  Assume each $\cR_i$ has at least one clasp 
arc with an endpoint on $\Xhigh$ and at least one clasp arc with an endpoint 
on $\Xlow$, and with the puncture endpoints on $b$--arcs in $\cR_i$, so that in particular no clasp
arc is long.  Moreover, assume that the image of
$\{\cR_1,
\ldots , \cR_k\}$  in $\ta$ is connected.  Then the collection of regions is said to have a 
{\it G--exchange foliation}.  The motion across a region with a
G--exchange foliation is a {\it G--exchange
move}.  As we have shown, in the absence of clasp arcs it would be a sequence of exchange moves which carries a
 subarc of $\Xhigh$ over a rooted block and strand tree.

\begin{proposition}
\label{proposition:G-exchange move}
Let $\{\cR_1,\dots, \cR_k\} \subset \pa$ be a collection of
regions which, taken together, have a G--exchange foliation.  Then an isotopy
of $\Xhigh$ to  $\Xlow$ across the regions $\{\cR_1, \dots ,\cR_k\}$ 
is realized by a G--exchange move.
\end{proposition}

\begin{proof}   Since each region $\cR_i$ is embedded, the isotopy of any one arc $\cR_i \cap \Xhigh$
across
$\cR_i$ to $\cR_i \cap \Xlow$ corresponds to a G--exchange move.
Moreover, since each region contains clasp arcs having endpoints on both $\Xhigh$ and $\Xlow$ 
and since their image in $\ta$ is connected, we have an interdependence of the isotopies 
across all the regions.  Thus, an isotopy across the collection $\{\cR_1, \ldots , \cR_k\}$ 
is a G--exchange move. \end{proof} 

\subsection{Pushing across a standard annulus}
\label{subsection:pushing across a standard annulus}  

\begin{figure}[htpb]
\labellist\small
\pinlabel {stabilize} [t] at 214 410
\pinlabel {\parbox{45pt}{3 exchange moves}} [t] at 380 410
\pinlabel {destabilize} [t] at 479 407

\tiny
\pinlabel {$\delta$} [br] at 146 511
\pinlabel {$-\delta$} [r] at 120 487
\pinlabel {$-\delta$} [r] at 173 488
\pinlabel {$\delta$} [tr] at 96 462
\pinlabel {$\delta$} [bl] at 198 459
\pinlabel {$-\delta$} [r] at 120 434
\pinlabel {$-\delta$} [r] at 174 434
\pinlabel {$\delta$} [tl] at 146 410

\pinlabel {$\delta$} [br] at 287 508
\pinlabel {$-\delta$} [r] at 261 484
\pinlabel {$-\delta$} [r] at 315 485
\pinlabel {$\delta$} [tr] at 237 458
\pinlabel {$\delta$} [tl] at 339 457
\pinlabel {$-\delta$} [r] at 261 431
\pinlabel {$-\delta$} [r] at 315 431

\pinlabel {$-\delta$} [r] at 434 429
\endlabellist
\centerline{\includegraphics[scale=.7]{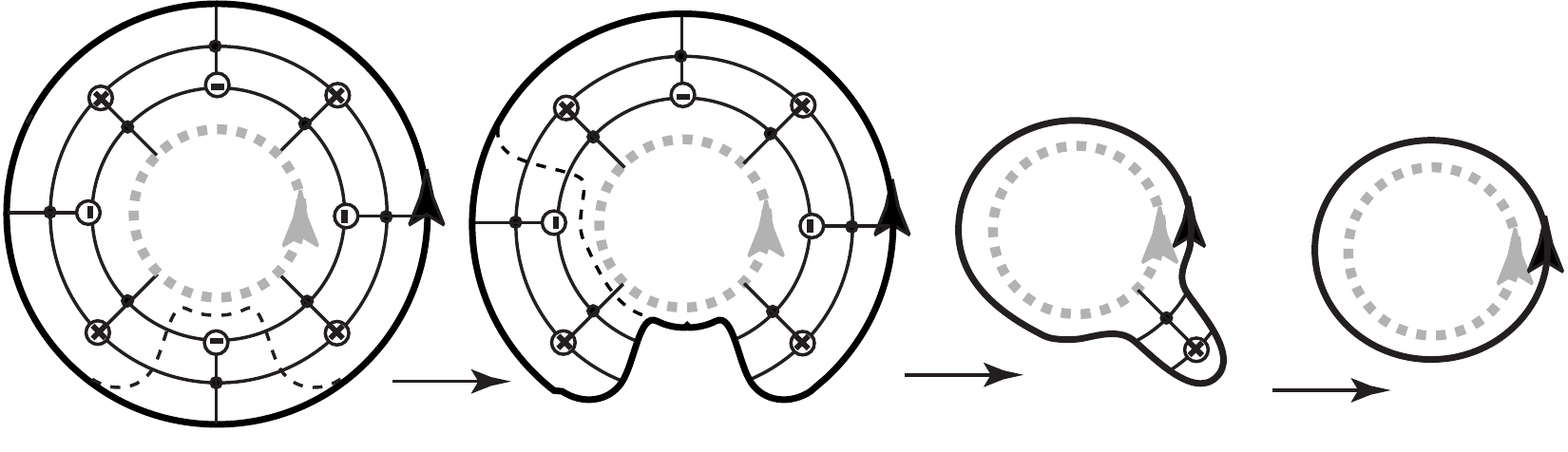}}
\caption{Pushing across a standard annulus}
\label{figure:annulus3}
\end{figure}
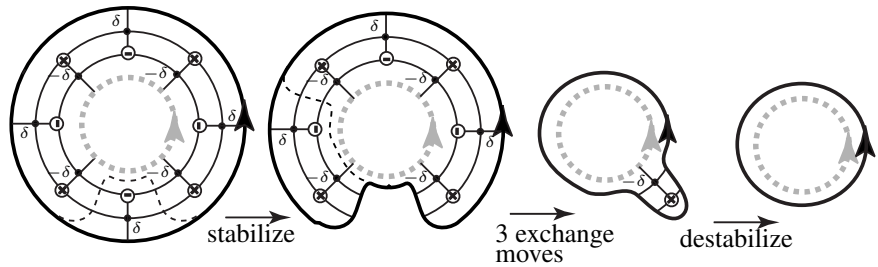

In constructing the thin annuli, we omitted  the case of a standard annulus (with or without clasp
arcs). Motions across a standard annulus will in general occur at the very end of the isotopy from
$\Xhigh\to\Xlow$, after all singularities that are on clasp arcs have been eliminated, and all short clasp arcs
have been removed.   To give an example of such a motion, let us suppose that $\cX$ is a link of
$\mu\geq 2$ components and that the first component, $\Xhigh^1$ cobounds with $\Xlow^1$ a standard
annulus with $2k$ vertices, with or without clasp arcs.  \fullref{figure:annulus3} shows an example without clasp arcs, when $k=4$.

\fullref{figure:annulus3} shows that a stabilization
along one of the singular leaves in $G_{-,\e}$, followed by a sequence of $k-1$ exchange moves,
followed by a destabilization, suffices to move $\Xhigh^1$ to $\Xlow^1$. The total
change in braid index is zero.   We have a Markov tower similar to the
one in \fullref{figure:exchange2}. If there are clasp arcs, we may assume that  they have been pushed into singular leaves.  Keeping in mind the stabilization move in \fullref{figure:annulus3}, where $\Xhigh$ is pushed across a singular leaf, we note that when there is a pair of clasp arcs a microflype pushes both the black boundary arc and the dotted grey boundary arc, simultaneously,  across a pair of singular leaves.  After all of the clasp arc pairs have been removed with the help of microflypes, the rest of the standard annulus  can be crossed with the help of complexity-reducing destabilizations and  exchange moves.  
In \fullref{section:introduction} of this paper, in  \fullref{figure:template0}, we constructed the Markov tower and the associated template for the crossing of a standard annulus without clasp arcs.  If there are clasp arcs, then as was just shown the annulus can be crossed using microflypes, exchange moves and destabilizations, and no other template is needed.  The blocks shown  in  \fullref{figure:template0}
allow for the possibility that other components of $\cX$ braid with $X^1$. If there
were no other components, then the entire move in \fullref{figure:template0} can be realized by braid isotopy.  If there are other components, or if $k>4$, this may not be the case.
The cyclic move exists for every $k$. The associated templates are the 
cyclic templates of the MTWS.

\section{The proof of the MTWS}
\label{section:the proof}
The machinery which we need to prove \fullref{theorem:MTWS}, the MTWS, has been set up.  In this section we
give our proof.  The reader is referred to \fullref{subsection:statement of results} for the statement of
the theorem.

\subsection[]
{Constructing the sequences \eqref{sequence:unwinding Xlow},
\eqref{sequence:exchange-reducibility on Xhigh} and \eqref{sequence:flype
sequence on Xhigh'} and the templates in $\cT(m)$}

\label{subsection:the sequences and the templates}

We are given a $\mu$--component oriented link type $\cX$ in oriented $S^3$ or $\reals^3$.  We are also given closed
braid representatives $\Xhigh \in \cB(\cX)$ and $\Xlow\in\cB_{\min}(\cX)$.  Indexing the components of $\cX$ as
$\cX^1,\dots,\cX^\mu$, we choose corresponding indices for the components of $\ta$ and $\Xhigh$ and $\Xlow$, so
that each annulus $\cA^j$ in $\ta$ has $\partial\cA^j = \Xhigh^j-\Xlow^j$.

By \fullref{proposition:get started for links}, which describes the `basic construction' for links, 
we know how to construct the clasp annulus $\ta$.  By the results in \fullref{section:braid foliations of
the immersed annulus} we may assume that $\ta$ supports a braid foliation.  In particular,
by \fullref{proposition:normal} we may assume that each clasp arc has a normal
neighborhood, and that it has been pushed into a union of leaves.  Let
$c(\Xhigh,\Xlow,\ta) = (c_1,c_2)$ be the complexity of the triple $(\Xhigh,\Xlow,\ta)$, as defined in
\fullref{subsection:the complexity function c(Xhigh,Xlow,ta)}.  Thus $c_1$ is the number of singularities that
are on clasp arcs and $c_2$ is the number of singularities that are outside normal neighborhoods of the clasp arcs.
The first step is:

\textbf{Construction of the sequences \eqref{sequence:unwinding Xlow} and 
\eqref{sequence:exchange-reducibility on Xhigh}}\qua We make as many modifications as are
possible, using only exchange moves  on
$\Xlow$ and only exchange moves and destabilizations on $\Xhigh$.  With these restrictions, Corollary 
\ref {corollary:minimum complexity up to exchange moves and destabilization} tells us that we may find  sequences  
$\Xlow = \Xlow^1 \to \cdots\to \Xlow^p
= \Xlow'$ as in \eqref{sequence:unwinding Xlow} and  $\Xhigh = 
\Xhigh^1\to\dots \to \Xhigh^q = \Xhigh' $ as in
\eqref{sequence:exchange-reducibility on Xhigh} of the MTWS such that 
$c(\Xhigh',\Xlow',\ta')$ is minimal up to exchange moves and destabilizations.  Exchange moves preserve
braid index and destabilizations reduce it.  Exchange moves and destabilizations are  both 
strictly reducing on
$c_2$.  Exchange moves in the presence of clasp arcs (ie the moves $ab^\star$ and $bb^\star$) preserve or reduce
$c_1$.  So our sequences are strictly complexity reducing with respect to $c(\Xhigh,\Xlow,\ta)$.

We next turn our attention to the construction of the sequence \eqref{sequence:flype sequence on Xhigh'} of
the MTWS.  We begin with a weak version.

\textbf{Claim}\qua
For all triplets $(\Xhigh',\Xlow',\ta')$ as above, we may find a sequence 
\begin{equation}
\label{sequence:using inadmissible flypes} 
\Xhigh' = X^1 \to X^2 \to \cdots \to X^r = \Xlow'
\end{equation}
which is strictly complexity reducing with respect to $c(\Xhigh',\Xlow',\ta')$ such
that every adjacent pair
$X^i \to X^{i+1}$ is related by a single destabilization, exchange move,
cyclic move, G--exchange move or flype.  Notice
that we do {\it not} require that the flypes be admissible, so that the sequence may not be non-increasing on braid
index.  

\begin{remark}
\label{remark:decomposition of ta}
\textbf{(a)}\qua  We draw the reader's attention to the similarities between this claim and the discussion at the beginning of
\fullref{subsection:G-flype moves} where we dealt with the special case when $\ta$
admits a decomposition into thin annuli.
Conceptually, subannuli of $\ta$ that are thin annuli can be thought of as clasp annuli in their own right, because
the boundary components represent $\cX$ and satisfy the conclusions of
\fullref{lemma:construction of the two annuli} and \fullref{proposition:normal}.  What this
claim is asserting is that $\ta^\prime$ admits a decomposition into clasp subannuli.  That is, for
each pair $(X^i , X^{i+1})$, the link $X^i - X^{i+1}$ bounds a clasp annulus that has its foliation
satisfying \fullref{proposition:normal}.  In the same fashion as in 
\fullref{subsection:G-flype moves} we use the notation $\cS^i$ for the clasp subannulus associated with $(X^i , X^{i+1})$.
The claim then asserts that $\cS^i$ corresponds to one of the following isotopies:
destabilization, exchange move, cyclic move, G--exchange move or flype.
\end{remark}  

We interrupt the remark to prove the claim. 

\begin{proof}[Proof of the claim]
The proof is by induction on $c(\Xhigh',\Xlow',\ta') = (c_1,c_2)$, using
lexicographical ordering.

To begin the induction, assume that $(c_1,c_2) = (0,0)$.  Since $c_1=0$, it follows that if there are clasp
arcs they must all be short, in which case \fullref{lemma:eliminating short clasp arcs} says that we may eliminate them by braid isotopy.  Therefore $\ta'$ is
embedded and foliated without singularities.  The foliation then consists
entirely of $s$--arcs.  But then 
$\Xhigh^\prime$ and $\Xlow^\prime$ represent the same braid isotopy class, and the MTWS is trivially true. 

\textbf{Inductive hypothesis}\qua Assume that the complexity of the
family of foliated annuli $\pa$ is $(c_1,c_2)$. Assume that the
claim is true whenever the
complexity is less than $(c_1,c_2)$.  There are several cases:

\textbf{Case (a)\qua $c_1>0$, and there are long clasp arcs}\qua By the construction
which was given in \fullref{lemma:construction of  cS}, we find a family $\cS$ of subannuli of $\pa$ such
that each component supports a thin foliation or is a standard annulus.
Some of these components may be trivially foliated, but since there is at least
one long clasp arc they are not all trivially foliated.  If $\cS$ is the union of standard annuli then we
can push across $\cS$ by a cyclic move as shown in \fullref{subsection:pushing across a standard annulus}.
Otherwise, $\cS$ is not trivially foliated but does contain $s$--arcs.
Then, by \fullref{proposition:sequence of flypes} we may use a
sequence of flypes to push $\Xhigh^\prime$ across $\cS$.  (Note: the flypes may not be admissible.)  This reduces
complexity.  By the induction hypothesis, the claim is true.

\textbf{Case (b)\qua $c_1>0$, but there are no long clasp arcs}\qua
In the presence of clasp arcs, but the absence
of long clasp arcs, \fullref{corollary:minimum complexity up to
exchange moves and destabilization} tells us that there must be at least
one G--exchange region.  Pushing across
it reduces complexity.  

\textbf{Case (c)\qua $c_1=0$, but $c_2$ is arbitrary}\qua
We have already removed
all short clasp arcs, so there are none.  But then there also are no normal neighborhoods.  The annulus $\ta$
is embedded.  By \fullref{corollary:standard annulus},  each component of $\pa$ will either be either an
annulus which is trivially foliated (in which case we are done as before) or there is a component which is a
standard annulus without clasp arcs.  In the latter case the motion is realized by a cyclic template. Again, the claim is true.   

While we have found a strictly monotonic complexity-reducing sequence of closed braids (sequence
\eqref{sequence:using inadmissible flypes}) joining $\Xhigh'$ to $\Xlow'$, we have not established sequence
\eqref{sequence:flype sequence on Xhigh'} of the MTWS because we have not established
that the braid index is non-increasing.

To deal with this deficiency we  introduce the {\it augmented complexity function}
$$C(\Xhigh',\Xlow',\ta) = (c_0,c_1,c_2),$$
where $c_0 = b(\Xhigh')$ and $c_1$ and $c_2$ are as before.  This is the complexity function
which is referred to in the statement of \fullref{theorem:MTWS}.
Notice that: 
\bi
\item Since the sequence \eqref{sequence:unwinding Xlow} is complexity-reducing with respect to
$c(\Xhigh,\Xlow,\ta)$ and since exchange moves preserve braid index, it is
complexity-reducing with respect to
$C(\Xhigh,\Xlow,\ta)$.
\item  Since the sequence \eqref{sequence:exchange-reducibility on Xhigh} is complexity-reducing with respect
to
$c(\Xhigh,\Xlow,\ta)$ and since destabilizations reduce braid index, it is 
complexity-reducing with respect to
$C(\Xhigh,\Xlow,\ta)$.
\ei  

Now notice that we can use the braid index entry of our augmented complexity function to pick out a subsequence of \eqref{sequence:using inadmissible flypes}:
\begin{equation}
\label{equation: final subsequence}
\Xhigh = X^{i_1} \to X^{i_2} \to \cdots \to X^{i_r} = \Xlow
\end{equation}
such that $0 < i_1 < i_2 < \cdots < i_r$ and $b(X^{i_{j}}) \geq b(X^{i_{j+1}})$, 
and
for any $X^k$ of \eqref{sequence:using inadmissible flypes} with $i_j < k < i_{j+1}$ we have
$ b(X^k) > b(X^{i_j})$.  Since $X^{i_j} - X^{i_{j+1}}$ cobound a clasp annulus $\cS^{i_j}$ which is the union of the
clasp subannuli mentioned in \fullref{remark:decomposition of ta} our
augmented complexity function is non-increasing on the triple $(X^{i_j}, X^{i_{j+1}}, \cS^{i_j})$.
The properties of sequence \eqref{sequence:flype sequence on Xhigh'} follow. 
\end{proof}

\addtocounter{remark}{-1}
\begin{remark}
\textbf{(b)}\qua Our decomposition of $\ta$ into
clasp subannuli $\cS^{i_1} \cup \cS^{i_2} \cup \cdots \cup \cS^{i_r}$ also
shows us the origin of templates in $\cT(m)$.  Specifically, an $\cS^{i_j}$ may
be a clasp subannulus of $\ta$ that is the union of any of the types of foliations of 
\fullref{subsubsection:forming thin annuli}, \ref{subsection:G-flype moves}, 
\ref{subsection:pushing across regions with a G-exchange foliation}
or  \ref{subsection:pushing across a standard annulus} along with destabilization.
This gives us a natural decomposition for $\cT(m)$.
Let $\cT(m,n)$ be the subset of all templates in $\cT(m)$ whose
initial braid has braid index $m$ and whose final braid has braid index $n \leq m$.  The union of the subsets
$\cT(m,n), \ n\leq m$ determines $\cT(m)$ because:
\begin{equation}
\label{equation:cT(m) and cT(m,n)}
\cT(m) = \cT(m,m) \  \sqcup \ \cT(m,m-1) \ \sqcup \ \cdots \  \sqcup \ \cT(m,1)
\end{equation}
Thus, employing the procedure given in \fullref{subsection:constructing the template from the total annulus}
for constructing a template, we can use each
triple $(X^{i_j}, X^{i_{j+1}}, \cS^{i_j})$ to construct a template in $\cT(b(X^{i_j}),b(X^{i_{j+1}}))$.
Moreover, using the triple $(\Xhigh, X^{i_{j+1}}, \cS^{i_1} \cup \cdots \cup \cS^{i_j})$ we can construct
a template in $\cT(m,b(X^{i_{j+1}}))$. \hfill$\diamondsuit$  
\end{remark}

Our proof of the MTWS is almost complete, except for two missing facts: 
the proof that no block in any template in $\cT(m)$  has full braid index, and the proof that
$\cT(m)$ is finite.  
We will deal with the former in \fullref{corollary:no spanning blocks}  and with the latter in \fullref{proposition:finiteness of cT(m)}.

We emphasize that the destabilization template is {\it not}  contained in $\cT(m)$.  On the other hand, for the remainder of the paper, it will be convenient to think of exchange moves, 
admissible flypes and G--exchange moves as templates in $\cT(m)$, rather than as separate moves.

\subsection{Cleaning up the templates in $\cT(m)$}  
\label{subsection:the braid index of blocks in a template}  

Our goal is to prove that no block in any template in $\cT(m)$  has full braid index, however we will do that in a setting that yields additional information for later use. The concept of a block-strand diagram is, at this moment, rather loose.  One can imagine that blocks could be slid around and amalgamated with one-another, and also subdivided.  We would like them to be more canonical.   But even if we succeed to make them canonical, a direct attack on the problem seems difficult.  We use a more roundabout approach.  When we organized the moves of the MTWS into three separate subsequences, we separated  the subsequence  \eqref{sequence:unwinding Xlow}, which requires only exchange moves, and the subsequence \eqref{sequence:exchange-reducibility on Xhigh}, which requires only exchange moves and destabilizations, from subsequence \eqref{sequence:flype sequence on Xhigh'}, which requires admissible flypes as well as the more general moves in $\cT(m)$.  We will see that, as a consequence, the templates in $\cT(m)$  have no blocks whose braid index is equal to the braid index of the diagram. In fact, we will be able to show more: that the templates of $\cT(m)$ are `consolidated' in a sense that will be made precise below, with the braid index condition on the blocks being  just one of the nice properties of consolidated templates.  

Let $\cD$ be a block-strand diagram having blocks $B$ and $B^\prime$.
We use the notation $t,b \subset \partial B$ (resp.\  $t^\prime, b^\prime \subset \partial B^\prime$) for the
top and bottom of $B$ (resp.\  top and bottom of $B^\prime$).
Recall our use of the term `amalgamating blocks' in
\fullref{subsubsection:using flypes to push Xhigh across cS}.  We now say
that $B^\prime$ (resp.\  $B$) can be 
amalgamated into $B$ (resp.\  $B^\prime$) if via a braid isotopy we can move $t^\prime \subset b$
(resp.\  $b \subset t^\prime$). (These are two ways of looking at the same phenomenon.)  We say the blocks of $\cD$ are {\it consolidated} if no two blocks of $\cD$
can be amalgamated.  Next, let $s$ and $s^\prime$ be strands in $\cD$ having endpoints on common blocks (still
called) $B$ and
$B^\prime$. We say that strands $s$ and $s^\prime$ are {\it topologically parallel} if there exists a rectangular disc
$R$ such that $R \cap \cD = \partial R = s \cup \beta_b \cup s^\prime \cup \beta_{t^\prime}$ where 
$\beta_b \subset b$ and $\beta_{t^\prime} \subset t^\prime$.  Note that if strands $s,s'$ are topologically parallel, it could happen, for example, that $s$ travels from $B$ to $B'$ without winding around the axis, but $s'$ winds about the axis as it does so.  To rule out this sort of complication, we say that $s$ and $s^\prime$ are
{\it braid parallel} if they are topologically parallel and $R \cap \axis = \emptyset$.
We say that $\cD$ is 
{\it consolidated} if the blocks of $\cD$ are consolidated and all of its topologically
parallel strands are braid parallel.

\begin{proposition}
\label{proposition:consolidated templates}
Given a template $(\Dhigh, \Dlow)$ with $\Dlow$ having minimal braid index
there exists a new template $(\Dhigh^\prime, \Dlow^\prime)$, where $\Dhigh^\prime$ is obtained from
$\Dhigh$ via a sequence of exchange moves and destabilizations; also $\Dlow^\prime$ is obtained from $\Dlow$ via a
sequence of exchange moves; also $\Dhigh^\prime$ and $\Dlow^\prime$ are both consolidated block-strand diagrams.
In particular, we can assume that the block-strand diagrams in every template in $\cT(m)$ are consolidated.
\end{proposition}

\begin{proof} 
Let $R$ be a rectangular disc that demonstrates that strands
$s$ and $s^\prime$ are topologically parallel for, say, $\Dhigh$.  We look at the induced foliation
of $\fib$ on $R$.  We can make the standard assumptions about $R$ being transverse to $\axis$ and
all but finitely many disc fibers of $\fib$ being transverse to $R$.  For the finitely many non-transverse
disc fibers, we can assume that each one contains a single saddle singularity.  We can then argue that
there are no leaves in the foliation that are circles.  The foliation of $R$ can
then be seen as a union of
$aa$--, $ab$--, $bb$-- and $sb$--tiles.  In particular, since the two sides $\beta_b$ and $\beta_{t^\prime}$ can be
assumed to be
$s$--arcs, $R$ will either be trivially foliated by $s$--arcs, or will contain
$sb$--singularities.

In the case where $R$ is not trivially foliated we assign an arbitrary orientation to $R$.
We can then talk about the $G_{\epsilon, \delta}$ graphs of $R$.  We wish to invoke the statements of
\fullref{proposition:behavior of graphs} to simplify the graph of $R$ through a sequence of
exchange moves and destabilizations.  Since $R$ is embedded all discs corresponding to those satisfying
statements (5), (6) and (7) of \fullref{proposition:behavior of graphs} will automatically be good.  It is
clear that after applying \fullref{proposition:behavior of graphs} repeatedly, $R$ will be trivially foliated, thus our resulting
strands will be braid parallel. 

We have a similar argument in the case where our strands are in $\Dlow$.  The only change is that our
application of \fullref{proposition:behavior of graphs} cannot yield any destabilizations because
$\Dlow$ has minimal braid index.  This establishes that we can replace a template $(\Dhigh,\Dlow)$ with
$(\Dhigh^\prime, \Dlow^\prime)$ via the sequence mentioned in statements 1 and 2, so that any topologically
parallel strands in $\Dhigh^\prime$ or $\Dlow^\prime$ are braid parallel.

Next, suppose there are two blocks, $B(1)$ and $B(2)$, in the template $(\Dhigh^\prime, \Dlow^\prime)$ such
that we can amalgamate $B(2)$ into $B(1)$ in, say $\Dlow^\prime$.  (These blocks may be moving or fixed.)
Such an amalgamation occurs because all of the strands entering the top of $B(2)$ are braid parallel and start
at the bottom of $B(1)$.
If we could also amalgamate them in $\Dhigh^\prime$, then there would be nothing to
stop us from doing the amalgamation
in both block-strand diagrams---our template would still be a template and would be simplified. 
The obstruction to amalgamating $B(2)$ into $B(1)$ in $\Dlow^\prime$ is that there is no similar amalgamation
of $B(2)$ into $B(1)$ in $\Dhigh^\prime$.  But, since the strands entering the top of $B(2)$ are all
braid parallel in $\Dlow^\prime$, they must be topologically parallel in $\Dhigh^\prime$.  But, we are assuming
that all topologically parallel strand in $\Dhigh^\prime$ are braid parallel.  So there is no obstruction.
Thus, the block amalgamation is possible.   After an iteration of amalgamations we may assume
that $(\Dhigh^\prime, \Dlow^\prime)$is consolidated. 
\end{proof}  

\begin{corollary}
\label{corollary:no spanning blocks}
Let $\cT = (\Dhigh,\Dlow)$ be a template in $\cT(m)$.   Then there is no block $B \subset \Dlow$ such that $b(B) = b(\Dlow)$. 
\end{corollary}

\begin{proof} 
Suppose that $\Dlow$ does contain a block $B$ such that the number of strands entering
the top of $B$ or leaving the bottom of
$B$ is equal to the braid index of $\Dlow$.  
Then we can amalgamate all other blocks in $\Dlow$ with $B$.   That is,
we can push all braiding of strands immediately before and after $B$ into $B$, thus
making all such braiding parallel.    This allows us to amalgamate any block immediately before
or after $B$ in $\Dlow$ with $B$.  We then iterate this procedure to amalgamate all blocks
in $\Dlow$ with $B$.  After this comprehensive amalgamation the consolidated
block-strand diagram $\Dlow^\prime$ will have a single block and will still carry $\Xlow$.
By \fullref{proposition:consolidated templates} we will thus have gone from
$\Xhigh$ (carried by $\Dhigh$) to $\Xlow$ (carried by the consolidated $\Dlow^\prime$)
through a sequence of braid isotopies, exchange moves and destabilizations. 
But note that these
are the moves that made up the sequences  \eqref{sequence:exchange-reducibility on Xhigh}
and \eqref{sequence:unwinding Xlow}.  This means that our $\Xhigh^\prime$
of sequence \eqref{sequence:exchange-reducibility on Xhigh} is braid isotopic to
our $\Xlow^\prime$ of sequence \eqref{sequence:unwinding Xlow}.  Thus,
sequence \eqref{sequence:flype sequence on Xhigh'} is vacuous and there is no
template to construct.  We conclude that we can assume that the templates
of $\cT(m)$, which are all constructed using sequence
\eqref{sequence:flype sequence on Xhigh'}, have no blocks of full braid index. \end{proof}

Our main goal for this section, the proof that no block in a template in $\cT(m)$ has full braid index, has been achieved. In so-doing, we learned a little bit more:  we may assume that the templates of $\cT(M)$ are consolidated.  In fact,  there is still is one more step we can take to make  them more canonical.  
Referring back to the discussion of {\em Fixed Blocks} in \fullref{subsection:pushing across a thin annulus}, in the first two  
paragraphs we used the occurrence of singularities in the fibration $\fib$
to designate fixed blocks.  In the third paragraph we then specified how some of these just-constructed
fixed blocks can be amalgamated.  Essentially, we were observing the existence of ``missing blocks''.  We make this concept and its implications rigorous below.  While this
is not necessary for pushing our finiteness argument forward, it does tie up
this loose end coming from \fullref{subsection:pushing across a thin annulus}.  Our brief detour will end at the end of this subsection.

Let $\cD$ be a block-strand diagram having braid structure
$(\axis,\fib)$.  Suppose there exist $3$--balls $B$ in the braid structure having a $\Delta \times [0,1]$
structure with the discs $\Delta \times {p}$ contained in fibers of $\fib$, $0 \leq p \leq 1$.  Assume that
$\intr(B)$ contains at least one block of $\cD$ and intersects only the strands of $\cD$, but
that $B \setminus (B \cap \cD)$ is not homeomorphic to an interval
cross a $2k$--punctured $2$--sphere.  We call such a $3$--ball $B$ a {\it missing block} in $\cD$.

Let $(\Dhigh, \Dlow)$ be a template that is consolidated. Then  $(\Dhigh, \Dlow)$ is an {\it optimal template}
if for every missing block of $\Dlow$, the ambient isotopy of $S^3$ that takes $\Dlow$ to $\Dhigh$ does not  
result in an missing block of $\Dhigh$.  In fact, we claim that we may assume that all templates of $\cT(M)$ are optimal.

To see this, let $(\Dhigh,\Dlow)= \cT \in \cT(m)$ and let $\{B^1,\ldots,B^k\}$ be a complete listing of all the
blocks in $\cT$.  As before, $b(B^i)$ is the braid index of the block $B^i$.  We define
the complexity of $\cT$  to be a lexicographically ordered $2$--tuple $(k, \sum_1^k (n - b(B^i ))$.
Now, if $B_-$ is a missing block of $\Dlow$ that is taken to a missing block $B_+$ of $\Dhigh$ by the 
ambient isotopy of $S^3$ that relates $\Dlow$ and $\Dhigh$ then by assumption $B_+$ must contain at least
one of the $B^i$ blocks plus something else.  That something else could be either additional blocks or
strands of $\Dhigh$ that do not intersect $B^i$ in the set $B_+ \cap \Dhigh \subset B_+$.  We now
replace $\cT \in \cT(m)$ with a new template
$$\cT^\prime =
( [(\Dhigh \setminus \intr(B_+)) \cup B_+] ,  [(\Dlow \setminus \intr(B_-)) \cup B_-]).$$  
Notice that
the complexity of $\cT^\prime$ is less than the complexity of $\cT$: if $B_+$ contains more than one
block then $\cT^\prime$ has fewer blocks than $\cT$; and, if $B_+$ contains just one block, $B^i \subset \Dhigh$,
along with some number of extra strands then $ m-b(B_+) < m - b(B^i)$.  In both situations the complexity is
reduced.

It is not hard to see that employing the recipe given in
\fullref{subsection:constructing the template from the total annulus}
for producing a $(\Dhigh,\Dlow)$ from an arbitrary triplet $(\Xhigh,\Xlow,\ta)$ does not necessarily
produce an optimal template.  We emphasize this point by ending this subsection with a useful
definition which we will need in our finiteness argument.

Let $\{ l_+^1, \ldots , l_+^h\}$ be a listing of the components of $\Xhigh
\cap \ta_{\text{tiled}}$ and, similarly,
$\{ l_-^1, \ldots , l_-^h\}$ be a listing of the components of $\Xlow \cap
\ta_{\text{tiled}}$ such that
the isotopy across $\ta$ has $l_+^i$ being taken to $l_-^i$ for $1\leq i \leq h$.  Specifically,
there is a component $\cR^i \subset \ta_{\text{tiled}}$ that has $l_\pm^i \subset \partial \cR^i$, $1 \leq i \leq h$.
Now, if we
collapse each of the $s$--arcs of $\ta - \ta_{\text{tiled}}$ to a point we can conceptually think of $\Xlow$ as being
obtained from $\Xhigh$ by replacing every $l_+^i$ in $\Xhigh$ with $l_-^i$, or to abuse notation
$[(\Xhigh \setminus (l_+^1 \cup \cdots \cup l_+^h) ) \cup ( l_-^1 \cup \cdots \cup l_-^h )] = \Xlow$.  (We will
continue this abuse of notation below.)

Recall that by Morton \cite[Theorem~1]{Morton1978} the
conjugacy class of $\Xhigh$ is determined by the link type $\Xhigh
\sqcup \axis$.
We say that a triplet has {\it unnecessary motion}
if its complexity can be reduced in the following manner.
Suppose now that there exists a proper subset
$\{l_+^{i_1}, \ldots , l_+^{i_{\rm l}} \} \subset \{ l_-^1, \ldots , l_-^h\}$ such that
$[(\Xhigh \setminus (l_+^{i_1} \cup \cdots \cup l_+^{i_{\rm l}}) ) \cup (l_-^{i_1} \cup \cdots \cup l_-^{i_{\rm l}})] 
\sqcup \axis $ has the same link type as $\Xhigh \sqcup \axis$.  Then let $\ta^\prime$ be the clasp annulus
that is obtained from $\ta$ by replacing each component $\{ \cR^{i_1},
\ldots , \cR^{i_{\rm l}} \} \subset \ta_{\text{tiled}}$
with an $s$--band that is parallel to the strands $\{l_-^{i_1}, \ldots , l_-^{i_{\rm l}}\}$.  Our new
clasp annulus is still cobounded by $\Xhigh-\Xlow$ and $C(\Xhigh,\Xlow, \ta^\prime) < C(\Xhigh,\Xlow, \ta)$.

\subsection{The set $\cT(m)$ is finite}
\label{subsection:finiteness of ct(m)}

We return to the main thread of the argument. There is still one very big unanswered question: how do
we know that $\cT(m)$ is finite?  Proving that it is, in fact, finite, is the main result in this
final subsection. See \fullref{proposition:finiteness of cT(m)}.    The proof of the MTWS will follow
immediately.

We begin with an investigation of restrictions on the tiling of $\ta$ which follow from the fact that it
is a topological annulus.  Let $W$ be a vertex in the foliation of $\pa$.
The {\it valence} $v$ of $W$ is the number
of singular leaves which have an endpoint on $W$.  Traveling around $W$ in either direction one will encounter a
cyclically ordered sequence $(t_1,t_2,\dots,t_v)$ of non-singular leaf types, where each $t_i$ is
either $a_\pm$ or $b$.  If this sequence includes
$\alpha$ edges of type $a_\pm$ and $v-\alpha$ edges of type $b$, then we say that $W$ has {\it type} $(\alpha,
v-\alpha)$.  Let $V(\alpha, v-\alpha)$ be the number of vertices of valence $v$ and type
$(\alpha, v-\alpha)$ in the foliation of $\pa$.  Let
$E(s)$ be the number of $s$--arcs which are boundary edges of a band of
$s$--arcs in the foliation of $\pa$.  For example, in the foliated annulus of \fullref{figure:normal7} we have $E(s) = 4$.   

\begin{lemma}
\label{lemma:Euler characteristic equality}
The vertices in the foliation of $\ta$ satisfy the following restriction:
\begin{multline}
\label{equation:positive equality} 
V(1,1) + 2 V(1,0) + 2 V(0,2) + V(0,3) = \\
2E(s) + V(2,1) + 2V(3,0) + 
\sum_{v=4}^\infty \sum_{\alpha=0}^v (v + \alpha -4) V(\alpha, v - \alpha),
\end{multline} 
with every term on both the LHS and the RHS non-negative.  Notice that every vertex
type appears and is counted, with the following 3 exceptions: vertices of type
$(1,2)$, $(2,0)$ and $(0,4)$ do not appear in this equation because in all three cases the
coefficient $(v+\alpha -4) =  0.$  
\end{lemma}

\begin{proof}   On each annular component of
$\pa$ the foliation determines a cellular decomposition which goes over to a cellular decomposition of $S^2$ on
shrinking the 2 boundary components to points.  Letting $V,E$ and $F$ be the number of vertices, edges and
tiles, the fact that $\chi(S^2) = 2$ shows that on each component of the foliated surface $\ta$ we have
$V+2 - E + F = 2$.  Each tile has four edges and each edge is an edge of exactly 2 tiles, so that
$E = 2F$.  Combining this with the previous equation we learn that $2V = E$.  Let
$E(a)$, $E(b)$, and $E(s)$ be the number of $a$--, $b$--, and $s$--edges, where we count both
$a_+$ and $a_-$ edges as being type $a$.  Then
\begin{equation}
\label{initial formula}
2V = E(a) + E(b) + E(s).
\end{equation}
  Since
\begin{eqnarray}
\label{equation:sum of vertices}
V & = & \sum_{v=1}^\infty \sum_{\alpha=0}^v V(\alpha,v-\alpha) \\
\label{equation:a edges}
E(a) & = & \sum_{v=1}^\infty \sum_{\alpha=0}^v \alpha V(\alpha,v-\alpha) \\
\label{equation:b edges}
2E(b) & = & \sum_{v=1}^\infty \sum_{\alpha=0}^v (v - \alpha) V(\alpha,v-\alpha)
\end{eqnarray}
we may combine equations \eqref{initial formula} through \eqref{equation:b edges} to obtain
\begin{equation}
\label{equation: second formula}
\sum_{v=1}^\infty \sum_{\alpha=0}^v (4 - v - \alpha) V(\alpha,v-\alpha) = 2E(s).
\end{equation}
Rearranging terms, we have proved the lemma.\end{proof} 

Our next goal is to learn which of the terms in Equation \eqref{equation:positive equality} are
bounded, and which terms can grow without bound, when we fix the braid indices
$b(\Xhigh)$ and $b(\Xlow)$.  For a triple $(\Xhigh , \Xlow ,\ta)$  of minimal complexity with $(b(\Xhigh), 
b(\Xlow)) = (m,n)$, let $\cN\subset\pa$ be the union
of all normal neighborhoods of the clasp arcs and let $\cN'$ be its complement in $\pa$.
As before, let $V(\alpha, \beta)$ be the number of vertices of valence $v$ and type
$(\alpha, \beta)$ in the foliation of $\pa$.  Let $V'(\alpha,\beta)$ denote the number of vertices which are in
$\cN'$ and contribute to $V(\alpha,\beta)$.

Our primary method for accounting for the braid index of the braid $\Xhigh \sqcup \Xlow$ is by
counting the number of $b$--arcs in non-singular disc fibers of $\fib$.
This number is independent of the choice of the fiber.  Let ${\bf b}$ be a
$b$--arc in $\ta\cap H_\theta$.  Then ${\bf b}$  has 2 endpoints on
$\axis$ and divides $H_\theta$ into 2 subdiscs, $\Delta_1$ and $\Delta_2$.
In \fullref{subsection:normal} we created normal neighborhoods of the
clasp arcs, at the expense of introducing lots of $b$--arcs into the
foliation of $\ta$. Those $b$--arcs are ones that can always be removed,
if we wish to do so, by reversing the isotopy.  We need to make sure that
they do not mess up our counting process.  For that reason, we now sharpen
the concept of an essential $b$--arc.  (Up to this point, a $b$--arc was defined to be essential if it was not inessential. )
\bi

\item A $b$--arc ${\bf b}$ is {\it strongly essential} if the interiors of
$\Delta_1$ and $\Delta_2$ each contain either an $a_\pm$-- arc or an
$s$--arc.   Note that in the creation of normal neighborhoods we never
introduced new strongly essential $b$--arcs.

\item A $b$--arc  ${\bf b}$ is {\it weakly essential} if the interiors of  
$\Delta_1$ and $\Delta_2$ do not contain $a_\pm$ or $s$ arcs.   The
reason these are only `weakly essential' is because, if their vertex endpoints are not adjacent on $\axis$, then one of discs $\Delta$ or $\Delta^\prime$
contains a family of $b$--arcs and the innermost member of that family
will be inessential. Thus a weakly essential $b$--arc can always be
changed to an inessential $b$--arc by an appropriate sequence of isotopies.

\item A $b$--arc  ${\bf b}$ is {\it inessential} if $v$ and $w$ are adjacent on
$\axis$.  This implies that neither $\Delta_1$ nor $\Delta_2$ contains an
$a_\epsilon$--arcs, an $s$--arc or a $b$--arc.

\ei

The construction of normal neighborhoods introduces numerous weakly
essential $b$--arcs in $\cN$. There may also be weakly essential
$b$--arcs in $\cN^\prime$, e.g. ${\bf b}$ may be in $\cN'$ but $\Delta$
and/or $\Delta^\prime$ may contain  $b$--arcs that are in $\cN$.  For this reason, we 
refine our count.  We already defined $V(\alpha,\beta)$. Recall that $V'(\alpha,\beta)$ is the number of vertices which are in
$\cN'$ and contribute to $V(\alpha,\beta)$.  Let $V^e(\alpha,\beta^e)$
denote the number of vertices which are adjacent to $\alpha$ $a$--arcs and
$\beta^e$ strongly essential $b$--arcs.  We will obtain the bounds that we need from $V(\alpha,\beta),  V^\prime(\alpha,\beta)$ and $V^e(\alpha,\beta^e)$

\begin{lemma}
\label{lemma:growth of V's and E(s)}
The  following  hold for the individual terms in Equation
{\rm \eqref{equation:positive equality}}, when the braid indices $m$ and $n$ are fixed
and $(\Xhigh,\Xlow,\ta)$ has minimal complexity through the application of exchange moves
(applied to $\Xhigh$ and $\Xlow$) and destabilizations (applied to $\Xhigh$):
\begin{enumerate}
\item [{\rm (1)}] $V(1,0)$ is zero.

\item [{\rm (2)}] $V(\alpha , v-\alpha)$ is zero for $\alpha \geq 2$.

\item [{\rm (3)}] If the complexity is minimal, then $V'(0,2) = V'(0,3) = 0.$

\item [{\rm (4)}]  $V^e(1,2)$ and $V^e(0,4)$ are bounded.

\item [{\rm (5)}] $E(s)$ is bounded.

\item [{\rm (6)}] $V(1,1)$ is bounded.

\item [{\rm (7)}] $\sum_{\beta=3}^\infty V^e(1, \beta^e) $ is bounded.

\item [{\rm (8)}] $\sum_{\beta=5}^\infty V^e(0, \beta^e) $ is bounded.
\end{enumerate}
\end{lemma}

\begin{remark}
We do not obtain bounds for $V(0,2), V(0,3), V(0,4), V(1,2), V'(1,2)$ and $V'(0,4)$, and indeed they cannot be bounded.
With this observation, we have
accounted for every possible term $V(\alpha,\beta)$ and $V'(\alpha,\beta)$
for $\alpha, \beta \geq 0$. \hfill$\diamondsuit$
\end{remark}

\begin{proof}  We consider the various inequalities in order:\\
\textbf{Proof of (1) and (2)}\qua These two statements follow from \fullref{corollary:minimum complexity up to exchange moves and destabilization}.  Specifically,
statement 1. of \fullref{corollary:minimum complexity up to exchange moves and destabilization} forces $V(1,0)$ to be zero.  If there is a vertex $v_\e$ that is adjacent
to the boundary of $\pa$ that contributes to the count of $V(\alpha , v-\alpha)$ for $\alpha \leq 2$
then there are two $a_\e$--arcs, $a_\e^1 , a_\e^2$, adjacent to $v_\e$
that are not isotopic to each other in
the foliation.  The subdisc, $\Delta(v_\e)$,
that $a_\e^1 \cup a_\e^2$ splits off in $\pa$ can only
contain $\g_\e$ clasp arcs.  The absence of $\g_{-\e}$ makes $\Delta(v_\e)$ embedded when
considered in $\ta$.  
This implies that all components of $G_{\e,\d} \cap \Delta(v_\e) \subset \Delta(V_\e)$
are simply connected, ie by statement
3. of \fullref{corollary:minimum complexity up to exchange moves and destabilization}
there are no loops in the graphs $G_{\e,\d} \cap \Delta(v_\e)$.  This further implies that there
exists a path $\alpha$ in either $G_{-\e, +}$ or $G_{-\e,-}$ satisfying
the assumptions in statement 2. of \fullref{corollary:minimum complexity up to exchange moves and destabilization}.
Since our subdisc does not contain a $\g_{-\e}$ clasp arc, we get a
contradiction of statement 2.  We conclude that $V(\alpha,v-\alpha)$ is zero.

\textbf{Proof of (3)}\qua By \fullref{proposition:behavior of graphs}, if a vertex $W$ has type (0,2) or (0,3),
and if its link is a good disc, then that vertex can be removed by changes of foliation
followed by exchange moves, reducing complexity. However, we are assuming minimum
complexity.  Therefore no such $W$ exists, unless link$(W)$ is not good,
ie it intersects $\cN$. 

\textbf{Proof of (4)}\qua  
Suppose that $W$ is a vertex that contributes to the count of
$V(1,2)$.  It is near $\Xhigh$ or $\Xlow$, so assume that it
is adjacent to a strongly essential $b$--arc.
An example was given in \fullref{figure:handle2}(a).  In this situation we showed in \fullref{figure:handle2}(b) that the corresponding embedding is a root.  Each new block
in the root contributes at least 1 to the braid index.
This forces the braid index of $\Xhigh
\sqcup \Xlow$ to grow, contradicting our assumption that it is fixed at $m+n$.  A similar argument
applies to $V(0,4)$.

\textbf{Proof of (5)}\qua Suppose that $E(s)$ is unbounded.  Then there will be an unbounded number of 
singularities of type $a_\epsilon s$ and/or  $sb$.  Notice that there is a direct
correspondence between the number of $a_\epsilon s$--singularities and the number of vertices
contributing to $V(1,0)$ and $V(\alpha , v - \alpha)$ for $\alpha \geq 2$.
But by statements (1) and (2) we know that $V(1,0)$ and
$V(\alpha , v - \alpha), \ \alpha \geq 2$ are bounded.
Thus, the only way that $E(s)$ can grow is if there is a growth in the number of singularities of type $sb$.

Previously in \fullref{subsubsection:forming thin annuli} we defined the notion of a complete
collection of $s$--arcs for thin annuli.  We have a similar notion for $\ta$ and $\pa$.  A
(possibly empty) family of $s$--arcs ${\bf S} =
\{ s_1, \dots , s_l: s_i \subset \ta \}$ is a {\it complete collection of
$s$--arcs} in $\ta$ if: (i) no two $s$--arcs in the collection split off a sub-band of
$\ta$ that is foliated entirely by $s$--arcs; and, (ii) for any other $s$--arc $s
\subset \ta$ there exists an $s_i \in {\bf S}$ such that $s \cup  s_i$ splits off a
sub-band of $\ta$ that is foliated entirely by $s$--arcs.
It is immediate that cutting $\ta$ open along a complete collection ${\bf S}$ of
$s$--arcs decomposes $\ta$ into a disjoint union of components that contain in one-to-one
correspondence the components of $\ta_{\text{tiled}}$ and
bands of $s$--arc.

Now fixing $\ta$, let ${\bf S}$
to be a complete collection of $s$--arcs in the foliation.  Assume that
$E(s)$ can be arbitrarily large.  We know that the number of components of $E(s)$ which have angular length
$\geq 2\pi$ is bounded, because each time that a band of $s$--arcs travels completely around $\axis$ it contributes
1 to $b(\Xhigh)$ and 1 to $b(\Xlow).$   This means that all of the growth in the cardinality $| \ta \setminus {\bf S} |$
(which is the same as $| \ta_{\text{tiled}} |$)
comes from components $\cC$ with the angular length of $\Xhigh\cap\cC$ and $\Xlow\cap\cC$ being
strictly less than $2\pi$.  In particular, for such a component the set of fibers $H_\theta$ for
which $H_\theta \cap (\Xhigh \cap \cC) \not= \emptyset$ coincides with the set of fibers for which   
$H_\theta \cap (\Xlow \cap \cC) \not= \emptyset$, since $\Xhigh\cap\cC$ and $\Xlow\cap\cC$ have their endpoints
on the same two $bs$--singularities.  We are thus seeing a growth in the
components of $\ta_{\text{tiled}}$ that are
characterized by the condition: 
\bi
\item There exists an $H_\theta \in \fib$ with $H_\theta \cap (X_\e \cap
\cC) = \emptyset, \ \e=+ \text{ and } -.$
\ei
We conclude that the following holds:
Let $\{\cC_1 , \ldots , \cC_{\rm J} \} \subset \ta \setminus {\bf S}$ be the set of components for which
there exists an $H_\theta \in \fib$ such that $H_\theta \cap (X_\pm \cap \cC_j) = \emptyset$ for all
$1 \leq j \leq {\rm J}$.  (Note that there may be a different $H_\theta$ for each $\cC_j$.)
If $E(s)$ grows then the index ${\rm J}$ must also grow.

Now suppose we have a subcollection of such components $\{\cC_{i_1} , \ldots , \cC_{i_{\rm L}} \} \subset
\{\cC_1 , \ldots , \cC_{\rm J}\}$ such that for every $H_\theta \in \fib$ we have that
$H_\theta \cap [\bigcup_{1 \leq j \leq {\rm L}} \cC_{i_j}]$ contains a
strongly essential $b$--arc.  Then
this subcollection contributes to the braid index of $\Xhigh \sqcup \Xlow$.  Since our braid index is fixed
such subcollections cannot grow in number.  Thus, as the index $J$ grows we can only have growth in a subset
$\{\cC_{i_1}^\prime , \ldots , \cC_{i_{\rm L^\prime}}^\prime \} \subset
\{\cC_1 , \ldots , \cC_{\rm J}\}$ with the property that there exists a fixed disc fiber $H_\theta^\prime$ with
$H_\theta^\prime \cap [\bigcup_{1 \leq j \leq L^\prime} \cC_{i_j}^\prime]$ being a union of weakly
inessential $b$--arcs.  Pushing these weakly essential $b$--arcs off of $H_\theta^\prime$ we see that the
isotopy across the components $\{\cC_{i_1}^\prime , \ldots , \cC_{i_{\rm L^\prime}}^\prime \}$ can be achieved
in the complement of $\axis$.  Thus, our original triplet $(\Xhigh,\Xlow,\ta)$ has unnecessary
motion and was not of minimal complexity.  But, we are assuming that we started with minimal
complexity.  Thus, $E(s)$ cannot grow.

\textbf{Proof of (6)}\qua We refer back to \fullref{proposition:index count on thin annuli}
to establish the boundedness of $V(1,1)$.  If we consider the construction of $\cShigh$ and $\cSlow$
we know from (1) and (2) of \fullref{proposition:index count on thin annuli}
that there can be at most $m$ type $1_d$ regions in $\cShigh$ and $n$ type $-1_d$ regions in $\cSlow$.
Each type $1_d$ (resp.\  $-1_d$) region in $\cShigh$ (resp.\  $\cSlow$)
accounts for two vertices that contribute to the count of $V(1,1)$.  So if $V(1,1)$ is able to increase
it must be through either the occurrence of type $0_{d,1}$ and $0_{d,2}$ regions in $\cShigh$ and $\cSlow$,
or through the occurrence of a subregion like those in \fullref{figure:exchange7}.
(The distinguishing feature between the two cases is whether the clasp arc intersecting the region
is long or not.)  An increase in
the latter will not increase $b(\Xhigh) + b(\Xlow)$ since it will be
associated with a G--exchange move.  So we need only consider growth in
$\cShigh$ and $\cSlow$.

If $V(1,1)$ is allowed to grow arbitrarily large then there will be one annular component of $\pa$ that will
contribute an arbitrarily large  number of vertices to the count of $V(1,1)$.  Thus, we will
have a single component of $\pa$ which will contribute an arbitrarily
large number of type $0_{d,1}$ or $0_{d,2}$ regions to the construction of
either $\cShigh$ or $\cSlow$.  Since this growth occurs on a single
component of $\pa$, any two type $0_d$ regions on that $\pa$ component
will be adjacent to a common $s$--band, (an assumption which is needed to apply (3) and (3$^{\prime})$ of \fullref{proposition:index count on thin annuli}).

Focusing on $\cShigh$, we know from (2) of \fullref{proposition:index count on thin annuli}
that we cannot have growth in the number of pairs of type $0_{d,1}$ regions and type $0_{d,2}$
regions that intersect each other.  From (5) of \fullref{proposition:index count on thin annuli}
we know that we cannot have a single region (see $\cR_3$
in \fullref{proposition:index count on thin annuli})
which is intersected by a growing number of type $0_d$ regions.
(By \fullref{remark:properties of thin annuli} any $\cR_3$ region in $\cShigh$ or $\cSlow$ will
be a fan, which is an assumption needed for the application of (3) and (3$^{\prime}$) of
\fullref{proposition:index count on thin annuli}.)

Thus, we can only have an increase in $V(1,1)$ if it comes from a pair of intersecting regions.
Dealing with the growth of $V(1,1)$ in $\cShigh$, we list the
possibilities: (i) a type $0_{d,1}$ could intersect another type $0_{d,1}$; (ii) a type $0_{d,1}$ could intersect
a type $0_{d,2}$; (iii) a type $0_{d,1}$ (or $0_{d,2}$) could intersect a type $1_d$; or (iv) a type $0_{d,1}$ (or $0_{d,2}$)
could intersect a type $-1_d$.  If we have possibility (i), for one of the $(1,1)$ vertices the
link($v$) will be a
good disc and we could have eliminated it by (7) of \fullref{proposition:behavior of graphs}.
This violates our minimal complexity assumption, so possibility (i) does not occur.
If we have possibilities (ii) or (iii)
then by (2) and (2$^\prime$) of \fullref{proposition:index count on thin annuli} there will be
a contribution of $+1$ to $b(\Xhigh)$.  So these occurrences are bounded.

Finally, we consider the growth of $V(1,1)$ in $\cShigh$ through an unbounded number of occurrences of
possibility (iv). Suppose there is growth in pairs of regions 
$(\cR_0,\cR_{-1_d})$ such that: $\cR_0 \subset
\cShigh$ is a type $0_{d,2}$ region; $\cR_{-1_d} \subset \cShigh$ is a type $-1_d$ region; and
$\rho(\cR_0) \cap \rho(\cR_{-1}) \not= \emptyset$.  We refer the reader back to Figures \ref{figure:normal5} and
\ref{figure:normal7}, and adapt them to our purpose at hand.  In \fullref{figure:normal7} we see an illustration
of a type $-1_d$ region.  Given such a region we can perform the inverse of the operation illustrated
in the top sketch of \fullref{figure:normal5} to introduce an inessential
$b$--arc that is positioned 
as the `left most' $b$--arc in $\cR_{-1_d}$.  (Referring to the type $-1_d$ region in \fullref{figure:normal7},
as we traverse the black boundary in the direction of its orientation,
with this introduction of an inessential $b$--arc, the first first singularity of parity $\delta$ will no longer
be associated with a clasp intersection.)  Now, if we stabilize $\Xhigh$ along this first singularity
we will eat into $\cR_{-1_d}$ and the remaining portion of $\cR_{-1_d}$ will be a type $0_{d,1}$ region which
we call $\cR_{-1_d}^\prime$.  The new stabilized $\Xhigh$ we call $\Xhigh^\prime$, and we will have
$b(\Xhigh^\prime) = b(\Xhigh) + 1$.  But, since $\cR_0$ is type $0_{d,2}$ and $\rho(\cR_0) \cap
\rho(\cR_{-1_d}^\prime)
\not= \emptyset$, by (3) of \fullref{proposition:index count on thin annuli} we know that
this intersection contributes $+1$ to $b(\Xhigh^\prime)$.  Since $\Xhigh^\prime$ came from $\Xhigh$ by a
single stabilization we know $2 \leq b(\Xhigh^\prime) = b(\Xhigh) + 1 $.  Thus, if we had $x$ such
$(\cR_0,\cR_{-1_d})$ pairs, for each pair we could have performed a similar stabilization on $\Xhigh$ to produce a
braid
$\Xhigh^\prime$; and we would know that $2x \leq b(\Xhigh^\prime) = b(\Xhigh) + x $.  This implies $x \leq b(\Xhigh)$.
So we have bounded $V(1,1)$.

\textbf{Proof of (7) and (8)}\qua  
We study Equation \eqref{equation:positive equality} and ask which
terms can grow without bound on both sides?  By statements (1) and (6) of this lemma we know that
the terms $V(1,0)$ and $V(1,1)$ on the LHS cannot grow without bound for fixed $m$ and $n$.  By
statement (3) of this Lemma we know that if $V(0,2)$ and/or
$V(0,3)$ on the LHS grow without bound, then the growth must occur inside the union $\cN$ of
all normal neighborhoods of clasp arcs.  By equation (7) we know that a growth in $V(1,1)$
will force a growth in $\sum_{\beta=3}^\infty V(1,\beta)$ or $E(s)$.

Passing to the RHS of Equation \eqref{equation:positive equality},
we know from statement (5) that $E(s)$ cannot grow without bound, for fixed
$m,n$.  By statement (2) we know that $V(2,1)$ and $V(3,0)$ are bounded.  But then, the only  terms
which might not be bounded, on the RHS of Equation \eqref{equation:positive equality}, are those in
the double sum.  However, of the terms in the double sum we know from statement (2) that
$V(\alpha,v-\alpha)$ is bounded if $\alpha\geq 2$.

Suppose that $\beta$ is bounded but that $V(1, \beta)$ increases without bound.  This means
that there is some fixed value of $\beta$ for which there are arbitrarily many vertices of
type $(1,\beta)$.  An example was illustrated earlier, in \fullref{figure:template4}.
In this illustration vertices $U,V,Y,Z$ are vertices `at the bottom' of a region that is composed of
normal neighborhoods (see the right sketch of \fullref{figure:normal2}),
and the vertices $W, W^{\prime\prime}$ can be thought of as `coning' these vertices and their associated
singular leaves.  (The vertex $W''$ should be thought of in a similar manner.)
Then $W$ and $W^{\prime\prime}$ contribute to the count of $V(0,\beta)$ or $V(1,\beta)$ and we are able to
see the interplay between these terms and $V(0,2)$ and $V(0,3)$ in equation (7).  The vertices $U,V,Y,Z$
are in $\cN$, but are adjacent to a vertex outside of $\cN$. They have valence 2.

Such a coning idea can be iterated.  Referring back to \fullref{figure:template5}, if we imagine
an additional vertex $W'''$ lying below the black dotted arc we could conceivably cone $W,W',W''$ and
all of the associated singular leaves to such a $W'''$.  Since the number of vertices in the shaded
region (vertices like $U,V,Y,Z$) can possibly grow, the number of vertices coning the 
bottom (or top) of a region composed of normal neighborhoods can also grow.  And, the number of
vertices coning the $W$--flavored vertices can then also grow, etc.  So we see that there is
no inherent reason why $\sum_{\beta=5}^\infty V(0, \beta) $ or $\sum_{\beta=3}^\infty V(1, \beta)$
should be bounded.  However, we need only establish that
$\sum_{\beta=5}^\infty V^e(0, \beta^e) $ or $\sum_{\beta=3}^\infty V^e(1, \beta^e)$
are bounded.

By construction the vertices that contribute to  $V^e(0, \beta^e)$ and $V^e(1, \beta^e)$ count
are outside normal neighborhoods.  We suppose that a growth in them is balanced in Equation
\eqref{equation:positive equality} by a growth in $V(0,2)$ and $V(0,3)$ that are associated with
the normal neighborhoods of $\ghigh$ arcs.

We need to go back to our original construction of $\ta$
and extract an embedded annulus from $\ta$ that contains the vertices that contribute to the
$\sum_{\beta=5}^\infty V^e(0, \beta^e) $ or $\sum_{\beta=3}^\infty V^e(1, \beta^e)$.
We do this by taking a tab neighborhood for each $\glow \subset \pa$ (see the left sketch
of \fullref{figure:normal1}) and removing it from $\pa$ and its image from $\ta$.
(This is equivalent to stabilizing $\Xlow$ along all of the singular leaves that are in the tab neighborhoods
of the preimages $\glow$ of the clasp arcs.)  Through an abuse of notation (and in keeping with 
\fullref{section:getting started}) we call this embedded annulus $\annhigh$.  

Recall the notation $\Xlow,\Xmed,\Xhigh,\annlow,\annhigh$ from the basic construction in 
\fullref{section:getting started}. Choose an annular neighborhood $\cA_0$ of $\Xmed$ in $\annhigh\cup\annlow$ which
does not intersect the clasp arcs. Then $\cA_0$ is embedded and has $\Xmed$ as its core circle. Since $\annhigh$
and $\annlow$ are both embedded, we may extend them to embedded annuli $\annhigh'=\annhigh\cup\cA_0$ and
$\annlow'=\annlow\cup\cA_0$ which have a common framing, also both are embedded and both have $\Xmed$ as a core
circle. From the construction in \fullref{section:getting started} we know that the algebraic linking number 
$Lk(\Xlow,\Xmed) = 0$. It then follows that the linking number $Lk(\Xhigh,\Xmed)$ is also $0$, and so $\annhigh'$
can be extended to a minimal genus Seifert surface $\seif_+$ (different from the one which we used in \fullref{section:getting started}) having $\Xhigh$ as its boundary. Observe that all of the vertices that contribute
to $\sum_{\beta=5}^\infty V(0, \beta) $ or $\sum_{\beta=3}^\infty V(1, \beta)$ are in $\annhigh$, and so also
are in $\seif_+$.

Now consider the count $\sum_{\beta=5}^\infty V^e(0, \beta^e)$ on $\annhigh \subset \seifhigh$, and
suppose we have a type $(0,2)$ vertex, $v \in \annhigh$.  The possibilities are:  that both
of the $b$--arcs that are adjacent to $v$ are strongly essential; or one is strongly essential and the
other is weakly essential; or both are weakly essential.  (Our notion of weakly essential
and strongly essential are now with respect to the surface $\seif_+$.)  We observe that if we
have the last case, when we eliminate either of the weakly essential $b$--arcs using an
exchange move and the surgery in \fullref{figure:foliations7}, the sum
$\sum_{\beta=5}^\infty V^e(0, \beta^e)$ will remain constant.  (Referring to the labeling in
\fullref{figure:foliations7}, any strongly essential $b$--arcs that were adjacent to the
vertex $w_1$ will be adjacent to the vertex $w_2$ after the surgery.  Thus, the essential
valence of $w_2$ will increase by exactly the essential valence of $w_1$.)

Next, we observe that the count $\sum_{\beta=5}^\infty V^e(0, \beta^e)$ is invariant under our change
in foliation illustrated in \fullref{figure:first change in foliation}.  
This is proved by examining the changes in the $H_\theta$--sequence under the change in foliation (which reverses
the order of the associated singularities).

So we allow in $\annhigh$ the conversion of type $(0,3)$ vertices to type $(0,2)$
vertices, and the elimination of type $(0,2)$ vertices if both of the
adjacent $b$--arcs
are weakly essential.  Such alterations to $\annhigh$ keep $\sum_{\beta=5}^\infty V^e(0, \beta^e)$
constant.  If the number of type $(0,2)$ vertices which we cannot eliminate in this manner
is arbitrarily large then, since for each such vertex there is a strongly essential (with respect
to $\seifhigh$) $b$--arc, the braid index of $\Xhigh$ will be unbounded.
(Basically, $\annhigh$ is forcing the existence
of a block and strand tree of arbitrarily high index as described in \fullref{subsection:G-flype moves}.)  Since we cannot have
an arbitrarily large number of such type $(0,2)$ vertices adjacent to
strongly essential $b$--arcs,
and they balance out $\sum_{\beta=5}^\infty V^e(0, \beta^e)$, this sum must be bounded.

A similar argument applies to $\sum_{\beta=3}^\infty V^e(1, \beta^e)$. Also, we can interchange the
role $\Xhigh$ with $\Xlow$, using $\annlow$ instead of $\annhigh$.  (A subtle point is that we will
have to change the orientation of $\annlow$ to match that of $\Xlow$.) This completes the proof of \fullref{lemma:growth of V's and E(s)}. \end{proof} 

\begin{remark}
\label{remark: euler characteristic formula}
It is interesting to note the similarities between \eqref{equation:positive equality}
and \cite[Equation~(7)]{B-M_VI}.  Given any Seifert surface, $\seif$ with Euler characteristic $\chi(\seif)$,
assume that $\seif$ is tiled by $aa$--, $ab$-- and $bb$--tiles. Using the notation introduced earlier for
$V(\alpha,\beta)$ , we have
\begin{multline}
\label{equation:Euler characteristic formula} 
V(1,1) + 2 V(0,2) + V(0,3) = \\
4\chi(\seif) + V(2,1) + 2V(3,0) + 
\sum_{v=4}^\infty \sum_{\alpha=0}^v (v + \alpha -4) V(\alpha, v - \alpha).
\end{multline} 
In our proof of statements (7) and (8) of \fullref{lemma:growth of V's and E(s)} we used the fact that
$\annhigh$ can be extended to a Seifert surface $\seifhigh$ bounded by $\Xhigh$.  In Equation
\eqref{equation:positive equality} we noticed that when dealing with a bounded braid index, growth in
the values $V(0,2)$, $V(0,3)$ is balanced by growth in the values
$V(\alpha, v - \alpha)$ for $\alpha=0,1$ and $4 \leq v$.
So in Equation \eqref{equation:Euler characteristic formula}, if there is any additional growth in
the values $V((\alpha, v - \alpha)$ for $0 \leq \alpha$ and $4 \leq v$ it must be balanced out
by (positive) growth in the value $-\chi(\seifhigh)$.
(If there are any vertices in $\seifhigh \setminus \annhigh$ that contribute to the count of
$V(0,2)$ or $V(0,3)$ then by our argument in \cite{B-M_VI} they would have been eliminated through
changes in foliation and exchange moves.)
In other words, one can think of $\annhigh$ as being the largest subannuli in $\seifhigh$ such that
when Equation \eqref{equation:Euler characteristic formula} is specialized to the surface
$\seifhigh \setminus \annhigh$, growth in the sum 
$V(2,1) + 2V(3,0) + 
\sum_{v=4}^\infty \sum_{\alpha=0}^v (v + \alpha -4) V(\alpha, v - \alpha)$ is balanced
precisely by growth in the genus of $\seifhigh$.  As the number of tiles of $\seifhigh$ grow the only way
$\annhigh \subset \seifhigh$ can intersect these additional tiles is by
$s$--bands going through
$aa$-- or $ab$--tiles parallel to the $\Xhigh$ boundary.  The braiding of
$s$--bands 
of $\annhigh$ that comes from them running through $aa$-- or $ab$--tiles of $\seifhigh$
can be seen as accounting for some of the braiding that occurs in the fixed blocks
of a block-strand diagram. \hfill$\diamondsuit$
\end{remark}

\begin{proposition} 
\label{proposition:finiteness of cT(m)}
Choose any positive integer $m$.  Then for each fixed positive integer $n\leq m$ the set of templates in
$\cT(m,n)$ is finite. Moreover, while the finitely many templates in $\cT(m,n)$ depend on the braid indices $m$ and $n$ of $\partial\ann$, they do not depend in any other way on the choice of $\ann$.
\end{proposition}

\begin{proof}   The idea behind the assertion that $\cT(m)$ is finite is that the parts of the foliated clasp annulus
$\ta$ which can grow without bound when we fix the braid index of the boundary are all inside the blocks.  In this
regard observe that a block of braid index $k < m$ can contain an
unbounded number of distinct $k$--braids, and of
course in any one example the $k$--braid assignment to the block contributes to the foliation of $\ta$.  The hard
part of the proof is to show that in {\it all} cases where aspects of the foliation of
$\ta$ grow without bound, the growth in a template $\cT = (\Dhigh,\Dlow)$ can be understood as occurring inside the
blocks of $\Dhigh$ (which are also the blocks of $\Dlow$).

We begin by defining a subset of $\ta$ which contains precisely the information that we need  to construct
a template in $\cT(m)$.  In this regard we remark that one of the beautiful features of block-strand
diagrams is that most of the detailed information about the links that they support is concealed in the blocks,
however we do not need to know details of what is in the blocks to construct the templates.  Therefore we really
need a rather limited amount of information from the foliation of $\ta$ to construct the templates in
$\cT(m)$.

A $b$--arc in $\ta$ is said to be {\it near} $X_\e, \ \e = \pm$,
if it has a vertex endpoint that meets an $a_\e$--arc.  The subset $\BS$ of $\ta$ which is of
interest to us now is the union of all $s$--arcs and $a_\pm$ arcs in $\ta$, together with
all $b$--arcs which are near $\Xhigh$ or $\Xlow$, enlarged to a closed neighborhood in $\ta$ which is chosen
so that its boundary (which include both $\Xhigh$ and $\Xlow$) is a union of simple
closed curves which are transverse to the foliation of $\ta$.  We call it the {\it boundary support} of
$\pa$, because it is the subset of $\pa$ which determines the embeddings
of $\Dhigh$ and $\Dlow$ in 3--space, by
\fullref{proposition:embedding of X} and the construction in \fullref{subsection:constructing the
template from the total annulus}.  In the example in \fullref{figure:template3} the boundary support includes
everything except the $b$--arcs which join vertices 5 and 11, and also vertices 2 and 14.  In generic examples
$\BS$ will be a very small subset of $\pa$. 

By definition $\BS$ includes the thin annuli, $\cShigh$ (resp $\cSlow$)
which are associated to $\Xhigh$ (resp.\  $\Xlow$), and all $b$--arcs in $\ta$ that are near $\Xhigh \sqcup \Xlow$.  A clasp arc pair
$(\ghigh,\glow)$ in $\pa$ induces a clasp arc pair $(\ghigh,\glow)$ in
$\BS$. It follows from the definition of $\BS$  that
clasp arcs in
$\ta$ are in $\BS$ if and only if they are doubly long. For, if it happens
that $\gamma^i_\epsilon$ is long but $\gamma^i_{-\epsilon}$ is not, then
$\BS$ will be embedded near $X_\epsilon$, even though $\ta$ is not
embedded near $\gamma^i$, and we are studying the part of $\pa$ which
determines the embedding of the template in 3--space.  Notice further that
an $s$--arc is in $\ta$
if and only if it is in $\BS$.  From now on we 
will drop the parity subscripts for $a$--arcs when we talk about the induced foliation on $\BS$.

As before, $\cN\subset\pa$ be the union
of all normal neighborhoods of the clasp arcs and let $\cN'$ be its
complement in $\pa$.    Initially, we let $\cN_{\text{dl}}$ be the union of all
normal neighborhoods of all doubly long clasp arcs in $\pa$.  These are
the clasp arcs with the property that we can push $\Xhigh$ across it all
the way to $\Xlow$ by a sequence of microflypes that amalgamate to a
flype.   (Other clasp arcs are not doubly long initially, but may become
doubly long after G--exchange moves.) 
We enlarge $\cN_{\text{dl}}$ by adjoining to it any $ab$--tile that has only
weakly essential $b$--arcs intersecting the
boundary of our initial $\cN_{\text{dl}}$.  We allow continual enlargement
of $\cN_{\text{dl}}$ in this fashion until
any $ab$--tile that intersects $\cN_{\text{dl}}$ and has only weakly essential
$b$--arcs is also in $\cN_{\text{dl}}$.  Next, let
$\cN_{\text{dl}}^\prime$ be the complement of $\cN_{\text{dl}}$ in $\cN$.  Let
$\cS$ be the union of all bands of $s$--arcs.  It will be convenient to divide the foliation of $\BS$  into
parts:
\bi
\item $\BS_1 =  \cS$, the union of all bands of $s$--arcs.
\item $\BS_2=\BS\cap\cN_{\text{dl}}$.
\item $\BS_3=\BS\cap\cN_{\text{dl}}^\prime$.
\item $\BS_4 = \BS\cap\cN$, the intersection of $\BS$ with the union of
  all normal neighborhoods of all clasp arcs.
\item $\BS_5 = \BS \ \setminus \ (\BS_1 \cup \BS_4)$, ie
  the part of $\BS$ that is non-trivially tiled and outside all of the
  normal neighborhoods.
\ei

If $A$ is a subset of $\pa$, let $|A|$ denote the number of connected components in $A$.

\textbf{Claim}\qua $|\BS_i|$ is bounded for $i=1,2,3,4,5$.

\begin{proof}[Proof of claim]
$|\BS_1| = |\cS|$ is bounded, because every band of
$s$--arcs has 2 $s$--edges, but by (5) of \fullref{lemma:growth of V's and
E(s)} we know that the number of $s$--edges
in the foliation of $\ta$ is bounded.
Suppose next that $|\BS_2| = |\BS\cap\cN_{\text{dl}}|$ is unbounded.  Since the clasp
arcs in $\cN_{\text{dl}}$ are doubly long, there must be some connected component
of $\pa$ split along the bounded set $\cS$ which has the property that as one travels along
$X_\pm$ in this component one passes from $\cN_{\text{dl}}$ to
$\cN_{\text{dl}}'$ an unbounded number of times.  However, studying
the regions in \fullref{figure:normal7} we see that this would violate (4), (5) and (6) of \fullref{lemma:growth of V's and E(s)}, so this cannot happen.  The identical argument shows that
$|\BS_3|=|\BS\cap\cN_{\text{dl}}^\prime|$ is also bounded.  (Remark: this does not say that there is a bound on the
number of clasp arcs, in fact no such bound exists.)  Since $\BS_4 = \BS_2 + \BS_3$ it follows that $|\BS_4|$ is
also bounded.  Since $\BS_5$ is the part of $\BS$ that is
non-trivially tiled and outside the union of all normal neighborhoods, it follows that as we travel along a
component of $X_\pm$  we will intersect components of  $\BS_1, \BS_5$ and $\BS_4$.  We never
pass from a component of $\BS_5$ to another component of $\BS_5$ without passing through a component of $\BS_1$
or $\BS_4$.  Since $|\BS_1|$ and $|\BS_4|$ are both bounded, it follows
that $|\BS_5|$ must be bounded too.  This proves the claim.
\end{proof}

Our next task is examine the contributions of the components of $\BS_i$ to
$\Dhigh$ and $\Dlow$.  For this we need to investigate in detail (in a more general setting) the construction in
\fullref{subsection:constructing the template from the total annulus}.  Recall that to construct a template
$(\Dhigh,\Dlow)$ we needed to understand four aspects of its structure: the moving blocks, the moving strands, the
fixed blocks and the fixed strands.  We analyze each separately.

\textbf{Moving blocks}\qua An amalgamating block $\bB$ will be moved to an
amalgamating block $\rB$ if and only if $\bB$ and $\rB$  are related amalgamating blocks as defined in
\fullref{subsection:pushing across regions with a G-exchange foliation}.  This
assumes that they are associated to clasp arcs which are doubly long in $\pa$.  Thus the moving blocks will be
associated to $\cN_{\text{dl}}$.  Note that there may be some choices involved when we select the amalgamating blocks.  We
make those choices in such a way that the set of all moving blocks has minimal cardinality.

The strands of $\Xhigh$ which are incorporated into an amalgamating block lie in the black boundary
of $\cN_{\text{dl}}$, ie in the subset $\BS_2$ of $\BS$.  We have already proved that $\BS_2$ is bounded.  We know that if
a subarc of
$\Xhigh\cap\cN_{\text{dl}}$ is related to a corresponding subarc of
$\Xlow\cap\cN_{\text{dl}}$, then their angular lengths
coincide.    So let $\{N_1, \ldots , N_r\} = \cN_{\text{dl}}$ be a listing of all of the components.
For each component $N_i$ we define its {\it angular support} $\angle N_i$
to be the interval $[\theta^0_i,\theta^1_i] \subset [0,2\pi)$ for which $\theta \in [\theta^0_i,\theta^1_i]$
iff $H_\theta \cap [(\Xhigh \cup \Xlow) \cap N_i] \not= \emptyset$.
Notice that if $\bB$ and $\rB$ are related amalgamating blocks then the angular support of every
$N_i$ component that intersects $\bB$ (or $\rB$) must overlap, ie if each component of the subcollection
$\{N_{j_1}, \ldots , N_{j_{\rm R}} \} \subset \{N_1, \ldots , N_r\}$
intersects $\bB$ then $ \angle N_{j_1} \cap \cdots \cap \angle N_{j_{\rm
R}} \not= \emptyset$.  Since $\cN_{\text{dl}}$
is a finite set there are a finite number of angular support intervals, and for those intervals
there are only a finite number of possible intersection subsets.  Thus, there are only a finite
number of possible moving blocks in any template of $\cT(m,n)$.

\textbf{Moving Strands}\qua Every subarc of $\Xhigh$ which is away from the
bands of $s$--arcs is potentially a moving
strand, however some of these potential moving strands have been amalgamated into moving
blocks.  We separate the surviving moving strands into two types:

\be
\item[(a)] Moving strands that are the subarcs of $\Xhigh$ which are in
$\partial\cN_{\text{dl}}$ but were not
amalgamated into moving blocks.  The finiteness of this set of strands follows from the argument we used to prove
that the number of moving blocks is finite. 

\item[(b)] Moving strands that are the subarcs of $\Xhigh$ which intersect $\BS_5$.  Since $|\BS_5|$ is
bounded, we can restrict ourselves to the strands which intersect a single component of $\BS_5$.
\ee

By \fullref{proposition:embedding of X},
the embedding of these strands is determined by the ordering and signs of the vertices and singularities which
belong to a tile in $\BS_5$ which intersects $\Xhigh$.  If we can show that the number of such vertices is bounded,
it will follow that the number of singularities is also bounded.  Since the number of distinct ways to assign orders
and signs to a finite set of vertices and singularities is finite, it will then follow that the number of
possible arrangements of the moving strands in set (b) is bounded.

The vertices in question contribute to the count of $V(1,\beta)$, and unfortunately could be unbounded.
For example, since we can have arbitrarily many clasp arcs, the number of vertices contributing to
the count of $V(1,2)$ can be arbitrarily high.  But, these vertices are adjacent to weakly essential
$b$--arcs which do not add information to the positioning of our type (a) or (b) strands.  Thus,
the only vertices that we need to be concerned with are the ones that contribute to the count of
$V^e(1,\beta^e)$.  By statement (7) of \fullref{lemma:growth of V's and E(s)}, we know that they
are bounded.  Thus, there are only a finite number of possibilities for the positioning of moving strands
in any template in $\cT(m,n)$.

\textbf{Fixed blocks}\qua The argument here is more subtle than the one for the moving blocks, because the fixed
blocks are associated to bands of $s$--arcs, and so there is no tiling to work with.
Nevertheless, we can relate the phenomenon of block amalgamation to the tiling, in the following way.  Recall (see the construction in 
\fullref{subsection:constructing the template from the total annulus}) how the fixed blocks are formed.
Let $\ta$ be the clasp annulus and let $\ta_{\text{tiled}}$ be the part that is
not foliated by $s$--arcs.  Let
$\{\theta_1,\dots,\theta_r\}$ be a listing of all of the singularities in
$\BS\cap\ta_{\text{tiled}}$.  For each $\theta$
in one of the intervals $[\theta_i,\theta_{i+1}]$ we know there are no singularities.  For each such $\theta$ we
also know that $H_\theta$ contains $a$--arcs and $b$--arcs and $s$--arcs.  So $(H_\theta \setminus
[H_\theta\cap\ta_{\text{tiled}}])$ is a collection of discs, each containing only
$s$--arcs (with
some discs possibly containing no $s$--arcs).  Each component of 
the union over all $\theta\in[\theta_i,\theta_{i+1}]$ of the set $\{ (H_\theta \setminus \
H_\theta\cap\ta_{\text{tiled}})\}$ has a
$D^2\times [\theta_i,\theta_{i+1}]$ structure.  Let $C$ be one such
component.  If $C$ contains bands of $s$--arcs.
amalgamate them into a single block $\bB(C)$.  If $C$ contained no bands or just a single band
then $\bB(C)$ is either vacuous or a single fixed strand.

Now it may happen that there is another connected component
$C'$ with its associated block $\bB(C')$, such that we have the
{\it amalgamation condition}
\begin{equation}
\label{equation:amalgamation condition}
\{s{\rm -arcs}\} \cap C' \cap H_{\theta_{i+1}} \subset \{s{\rm -arcs}\} \cap C\cap H_{\theta_{i+1}}.  
\end{equation}
If this happens, amalgamate $\bB(C')$ and $\bB(C)$ and delete the singularity at $\theta_{i+1}$, so that we have a single block
$\bB(C \cup C')$.  (If $\bB(C')$ was vacuous or a single strand then $\bB(C \cup C')$ is essentially
still $\bB(C)$.)  Continue this amalgamation process as long as possible.  Among all amalgamated 
blocks discovered
in this way, choose one such that the set of fixed blocks has minimum cardinality.  In this way we will have
eliminated some number of singularities, ie the ones which separated the new-amalgamated blocks.  Let
$\theta^1,\dots,\theta^p$ be the angles which remain.  We need to show that this listing is bounded, when $b(\Xhigh)$ and $b(\Xlow)$ are fixed.

To do this we consider the effect a singularity has in the
$H_\theta$--sequence when weakly essential $b$--arcs
are used.  We refer to \fullref{figure:foliations3}, where the
possibilities for the $b$--edges of $bb$-- and $ab$--tiles are illustrated. Each could be either strongly or weakly essential.
\begin{figure}[htpb]
\centerline{\includegraphics[scale=.8]{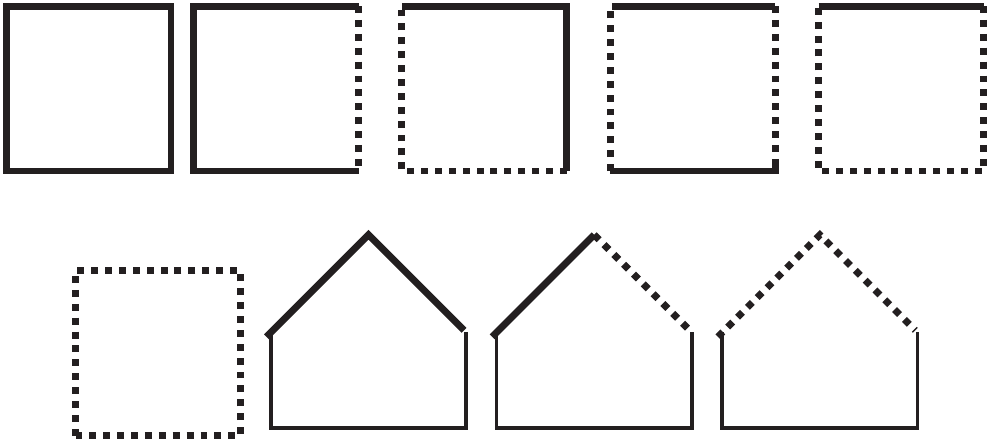}}
\caption{The possibilities for $b$--edges of $bb$-- and $ab$--tiles. A
black undotted $b$--edge means
strongly essential A black dotted $b$--edge means weakly essential.}
\label{figure:template5}
\end{figure}

First, note that occurrences of $aa$--singularities will effect the block amalgamation condition stated
previously, ie \eqref{equation:amalgamation condition} above. 
So the occurrence of $aa$--singularities will register in our listing of remaining angles
$\theta^1,\dots,\theta^p$.  But, by statements (1), (2) and (3)
of \fullref{lemma:growth of V's and E(s)} we will have a bounded
number of such singularities.

Second, notice that among the possibilities for $bb$--tiles, (1) and (2) will not effect our amalgamation
condition, because weakly essential $b$--arcs split off regions that
contain no $s$--arcs.
This is also true for (1) and (2) in our possibilities for $ab$ tiles.

Thus, the only possibilities that effect the amalgamation condition are
(4), (5) and (6) for $bb$--tiles,
and (3) for $ab$--tiles.

Now by statements (1), (2), (4), (6), (7) and (8) of \fullref{lemma:growth of V's and E(s)} there are a bounded number of
strongly essential $b$--arcs in a fixed fiber of $\fib$.  So there are a bounded number of singularities that can
effect our amalgamation condition.  Thus, the growth in our remaining angles $\theta^1,\dots,\theta^p$
is bounded.  (It is interesting to notice that since the tiles in the normal neighborhoods have all of
their sides labeled $wi$ they do not effect the amalgamation condition.)

We have in fact established more than just that there are a finite number of fixed blocks.  We have
established that there are only a finite number of possible positions for fixed blocks.  This
is because there are only a finite number of positions for vertices that contribute to
the count of $V^e(\alpha, \beta^e)$ and, thus, a finite number of singularities that
correspond to possibilities (3)--(5) for $bb$ tiles and possibilities (2) and (3) for $ab$ tiles.
Therefore any combinatorial information having to do with cyclic ordering of vertices on $\axis$ and
cyclic ordering of singularities in $\fib$ is also finite.

\textbf{Fixed strands}\qua  If the number of moving blocks, moving strands and fixed blocks is finite
then the number of fixed strands must be finite.
Using the observation we employed to establish the finiteness of positions
of fixed blocks, we can establish finiteness of the positions of fixed strands.

Thus everything is bounded, and so the number of block-strand diagram pairs is bounded.  The proof of \fullref{proposition:finiteness of cT(m)} is complete.  But then, so is the proof of \fullref{theorem:MTWS}.
\end{proof}

\section{Open problems}
\label{section:open problems}


(1)\qua\label{openproblem:Kirby Calculus}
As was pointed out in \fullref{subsection:the problem}, Markov's Theorem
is just one example of a class of theorems about 3--manifolds in which
some form of stabilization and destabilization play an important role.   A
very different and very important example is the Kirby Calculus
\cite{Kirby}, relating two surgery representations of a 3--manifold.  The stabilization move is the addition of an unknot with standard framing to the link that is to be surgered.  In this setting, what should be the analogue of the MTWS?  
Presumably, the first entry in the complexity function should be the number of components in the surgery link.    But unfortunately, in this regard, we do not know of any tools that could be put to work, to play the role that was played by the geometry of braid foliations and the group structure of the braid group in the proof of the MTWS.

(2)\qua\label{open problem:Reidemeister-Singer Theorem}
There are many analogies between the
study of knots via their closed braid representatives and the study of 3--manifolds
via their Heegaard diagrams (or equivalently via their `Heegaard gluing maps' in the mapping
class group of a closed orientable surface of genus $g$).  In the latter setting
equivalence classes of Heegaard splittings are in 1--1 correspondence with
double cosets in the mapping class group $\cM_g$ modulo the mapping class group $\cH_g$ of a
handlebody.  We pose as an open problem to find  moves which 
(like the moves in $\cT(m)$) change the equivalence class of a Heegaard
splitting of a 3--manifold without increasing its Heegaard genus.  A strategy for finding
such moves is given by the second author in \cite{Menasco2001a}, however (lacking an invariant) there is no proof
that this strategy actually produces inequivalent splittings.  In his PhD thesis
\cite{Zablow} Joel Zablow made a relevant contribution in his study of waves in Heegaard
diagrams.  Is there a tool which plays the role of braid foliations in the situation
of Heegaard splittings of 3--manifolds?  This seems to be a very interesting area for
future investigations.

(3)\qua\label{openproblem:construct the templates}
We pass to open questions about \fullref{theorem:MTWS} of this paper. In principle the templates in $\cT(m)$ can be enumerated,
but the actual enumeration is non-routine.  We pose this as an open
problem for $m=4,5$ and any other cases which prove to be computable.  A very interesting special case are to classify the templates in $\cT(m,m)$. In particular, these relate any two closed braid representatives when both have minimum braid index. There should be
applications.  For example, knowing that the only templates that we need
for braid index 3 are the flype and destabilization templates,  it is a
simple matter to classify links which are closed 3--braids (replacing
the complicated argument we used in \cite{B-M_III}), and it is to be expected that if $\cT(4)$ is
computed, then one would learn more about the classification of links of braid index 4.   One could expect many applications, if such an advance in knowledge could be achieved. 

(4)\qua Although actual enumeration of $\cT(m)$ may not be routine, if we restrict our attention to a
single type of isotopy such enumeration or characterization may be reasonably doable.  Specifically,
referring back to \fullref{subsection:pushing across a standard annulus},
the characterization of all knot complements that admit a cyclic move
predicated on the existence of an essentially embedded standard annulus in its complement would be
of interest.  Moreover, such knot complements could be divided into two classes: the first class would
have the components of $G_{++}$ and $G_{--}$ being homeomorphic to $[0,1]$; and, the second class would
have the components of $G_{++}$ and $G_{--}$ being homeomorphic to $S^1$.  The first class would use
a positive stabilization and destabilization to begin and end the cyclic move.  The second class would
use a negative stabilization and destabilization for the cyclic move.  The question is how to determine
when a knot is not in both classes, since the first class corresponds to a transversal isotopy and
the second class does not.  (See our article \cite{B-M_transverse}.)

(5)\qua In the manuscript \cite{B-M_IV} the authors proved that an arbitrary closed braid representative of a
composite knot or link may be modified by the use of exchange moves to a prime summand of the same braid
index.  We do not know whether the work in that paper can be incorporated into the proof of the
MTWS.  The difficulty is that we do not know whether the splitting 2--sphere which realizes the connect-sum
operation can be modified to one which intersects the clasp annulus $\ta$ in a `nice' way.

(6)\qua Some knots or links, for example the unlink \cite{B-M_V} and most iterated torus links \cite{Menasco2001b}, have
unique closed braid representatives of minimum braid index.  On the other hand, there are links of braid index
3 which have more than one conjugacy class of 3--braid representatives \cite{B-M_III}, and this pathology
clearly persists as braid index is increased.  We pose the open problem: find general conditions which suffice
for a knot or link type to have a {\it unique} conjugacy class of closed braid representative of minimum braid
index.

(7)\qua Referring back to \fullref{remark: euler characteristic formula},
it would be of interest to understand exactly how the structure of
a minimal genus Seifert surface can restrict the ability of a given link type having minimal braid index
$n$ to be carried by any template of $\cT(m,n)$.  In particular, if the Seifert surface has a foliation
composed of $aa$ tiles at minimal braid index is $\cT(m,n)$ necessarily empty for all values of $m$.

(8)\qua There are special knots and links, for example the unlink \cite{B-M_V} and most iterated
torus links \cite{Menasco2001b}, for which the MTWS is very simple: the moves are simply braid
isotopy and exchange moves.  We say that such links are {\it exchange-reducible}.  Are there other
examples of exchange-reducible links? Does the conclusion hold under weaker hypotheses? 

We remark that by the main result in Birman--Wrinkle \cite{B-W}, if a knot type $\cX$ is exchange-reducible, then
every transversal knot type ${\cal TX}$ associated to $\cX$ is transversally simple,
ie determined up to transversal isotopy by $\cX$ and the Thurston--Bennequin invariant.  Since
G--exchange moves and positive flypes are realized by transversal isotopy, it would be equally interesting if
the condition `exchange-reducible' was weakened to `exchange and positive flype-reducible'.

(9)\qua As noted in problem (4) above, the unlink is exchange-reducible.  This fact proves  that {\it
there exists} a monotonic and very rapid (perhaps even a quadratic) algorithm for recognizing the
unlink, through the use of exchange moves.  Unfortunately, however, the complexity function that would
translate this existence theorem into a working algorithm needs new techniques, as the complexity
function is concealed in the invisible family of discs  which the unlink bounds.  (One of these days
the first author will write a short note to show that it is also concealed in the auxiliary `extended
braid word' of Birman and Hirsch \cite{B-H}.)  We note that the unknot recognition algorithms  in
\cite{B-H} and  the finite unknot recognition algorithm in the very new paper by Ivan Dynnikov
\cite{Dynnikov2003}, which is based upon related foliation techniques,  are exponential.  A vague (but
we feel realistic) problem is to find an `energy functional' (AKA complexity function) which
uses the monotonic reduction process that is guaranteed to exist because of exchange-reducibility.  A
wild guess is that it is encoded in notions based upon Ricci curvature.

(10)\qua The MTWS begins with a choice of a closed braid representative $\Xlow\in\cX$ which has braid
index $b(\cX)$, however at this writing we do not know how to compute $b(\cX)$.  The most useful tool
that we know is the `Morton--Franks--Williams inequality', however  the MFW
inequality is doomed to fail in certain situations (see Kawamuro
\cite{Kawamuro}).  An important open problem is to develop new techniques for computing the braid index of a knot or link.   As can be seen from \cite{Kawamuro}, templates can be a very useful tool.

(11)\qua Conjecture: 
\bi
\item  Every block-strand diagram in a template in $\cT(m)$ has at least one block.  
\ei 
We remark that our attempts to find a counterexample have been unsuccessful, but we lack a proof
that it cannot happen.  Note that the number of counterexamples is necessarily finite.

(12)\qua\label{problem:arc presentations without stabilization}
In a standard annulus, the
graph $G_{-,\delta}$ (resp.\  $G_{+,\delta}$) is topologically equivalent to a circle  which cobounds with a
component of
$\Xhigh$ (resp.\  $\Xlow$) an embedded annulus. That annulus is foliated without singularities. From
this is follows that, if we regard $G_{-,\delta}$ and $G_{+,\delta}$ as defining knots in
$\reals^3$, then they will have the same knot type as the component in question of $\cX$. But in
fact more can be said. The graph $G_{-,\delta}$ (resp.\  $G_{+,\delta}$) is a union of arcs, each a branch
in a singular leaf, which join up a string of negative (resp.\  positive) vertices in a cycle.  Each
arc lies in a fiber of $\fib$ and has its endpoints on $\axis$, and so this representation of the component of
$\cX$ has an `arc presentation', in the sense defined by Cromwell \cite{C}
and Dynnikov \cite{Dynnikov2003}. 
Indeed, Ivan Dynnikov has been engaged in a project which begins with the introduction of braid
foliations, and goes on to study the foliations of the associated Seifert surface bounded by
$G_{-,\delta}$ and
$G_{+,\delta}$, adapting the braid foliation machinery in \cite{B-M_IV} and \cite{B-M_V} to arc
presentations.    One expects that there will be similar adaptations of the work in this paper to arc
presentations, although the adaptation is almost certainly non-trivial.

One reason why arc presentations are of interest is because they give a 
filtration of all knots and links, using the number of arcs as a measure of complexity,
and with that filtration there are always finitely many arc presentations which represent a
given knot type and have complexity at most the complexity of a given example. This fact is
important if one wishes to use the braid foliation machinery to construct algorithmic solutions to
the knot and link problem. 

\bibliographystyle{gtart}
\bibliography{link}

\end{document}